\documentclass{memo-l}
\usepackage{amsmath,amssymb,amscd,mathrsfs}

\newtheorem{theorem}{Theorem}[chapter]

\theoremstyle{plain}
\newtheorem{Lemma}[theorem]{Lemma}
\newtheorem{Theorem}[theorem]{Theorem}
\newtheorem{Corollary}[theorem]{Corollary}
\newtheorem{Conjecture}[theorem]{Conjecture}

\theoremstyle{remark}
\newtheorem{Remark}[theorem]{Remark}

\numberwithin{section}{chapter}
\numberwithin{equation}{chapter}

\makeatletter

\newif\ifmemoirs@
\memoirs@false

\def\op{\operatorname{op}}

\def\row{{\operatorname{row}}}
\def\ro{{\operatorname{row}}}

\def\col{{\operatorname{col}}}
\def\co{{\operatorname{col}}}
\def\ch{\operatorname{ch}}
\def\wt{\operatorname{wt}}
\def\res{\operatorname{res}}
\def\Col{\operatorname{Col}}
\def\Row{\operatorname{Row}}

\def\Tab{\operatorname{Tab}}
\def\Dom{\operatorname{Dom}}
\def\Std{\operatorname{Std}}
\def\rdet{\operatorname{rdet}}
\def\rad{\operatorname{rad}}

\def\cdet{\operatorname{cdet}}
\def\id{\operatorname{id}}

\def\height{\operatorname{ht}}

\def\V{\mathbb{V}}
\def\I{\mathbb{I}}
\def\C{{\mathbb F}}

\def\Z{{\mathbb Z}}
\def\N{{\mathbb N}}
\def\LL{{{\scriptscriptstyle\text{\rm{L}}}}}

\def\F{\mathrm{F}}

\def\pr{{\operatorname{pr}}}
\def\p{{\operatorname{p}}}
\def\diag{{\operatorname{diag}}}
\def\ad{\operatorname{ad}}
\def\gr{\operatorname{gr}}

\def\hom{{\operatorname{Hom}}}

\def\End{{\operatorname{End}}}

\def\sgn{{\operatorname{sgn}}}
\def\Mod{\operatorname{\text{-}mod}}

\def\Wh{\operatorname{Wh}}
\def\eps{{\varepsilon}}
\def\phi{{\varphi}}
\def\emptyset{{\varnothing}}

\def\bw{{{\textstyle\bigwedge}}}

\def\de{{\delta}}

\newcommand{\ba}{{\text{\boldmath{$a$}}}}
\newcommand{\bb}{{\text{\boldmath{$b$}}}}
\newcommand{\bc}{{\text{\boldmath{$c$}}}}

\newdimen\hoogte    \hoogte=8pt    
\newdimen\breedte   \breedte=8pt   
\newdimen\dikte     \dikte=0.5pt    

\newenvironment{young}{\begingroup
       \def\vr{\vrule height0.8\hoogte width\dikte depth 0.2\hoogte}
       \def\fbox##1{\vbox{\offinterlineskip
                    \hrule height\dikte
                    \hbox to \breedte{\vr\hfill##1\hfill\vr}
                    \hrule height\dikte}}
       \vbox\bgroup \offinterlineskip \tabskip=-\dikte \lineskip=-\dikte
            \halign\bgroup &\fbox{##\unskip}\unskip  \crcr }
       {\egroup\egroup\endgroup}
\def\diagram#1{\relax\ifmmode\vcenter{\,\begin{young}#1\end{young}\,}\else%
              $\vcenter{\,\begin{young}#1\end{young}\,}$\fi}

\newdimen\Hoogte    \Hoogte=12pt    
\newdimen\Breedte   \Breedte=12pt   
\newdimen\Dikte     \Dikte=0.5pt    

\newenvironment{Young}{\begingroup
       \def\vr{\vrule height0.8\Hoogte width\Dikte depth 0.2\Hoogte}
       \def\fbox##1{\vbox{\offinterlineskip
                    \hrule height\Dikte
                    \hbox to \Breedte{\vr\hfill##1\hfill\vr}
                    \hrule height\Dikte}}
       \vbox\bgroup \offinterlineskip \tabskip=-\Dikte \lineskip=-\Dikte
            \halign\bgroup &\fbox{##\unskip}\unskip  \crcr }
       {\egroup\egroup\endgroup}
\def\Diagram#1{\relax\ifmmode\vcenter{\,\begin{Young}#1\end{Young}\,}\else%
              $\vcenter{\,\begin{Young}#1\end{Young}\,}$\fi}

\begin{document}
\frontmatter
\title{Representations of Shifted Yangians\\ and\\ Finite {\boldmath $W$}-algebras}
\author{Jonathan Brundan}
\author{Alexander Kleshchev}
\address{Department of Mathematics, University of Oregon, Eugene, OR 97403, USA.}
\address{Department of Mathematics, University of Oregon, Eugene, OR 97403, USA.}
\email{brundan@uoregon.edu}
\email{klesh@uoregon.edu}
\urladdr{http://uoregon.edu/$\sim$brundan}
\urladdr{http://uoregon.edu/$\sim$klesh}
\thanks{Both authors supported in part by NSF grant no. DMS-0139019.}
\thanks{Part of this research was carried out
during a stay by the first author at the
(ex-) Institut\\\phantom{spaj} Girard Desargues, Universit\'e Lyon I in Spring 2004.
He would like to thank Meinolf Geck\\\phantom{skak} and the other members of the institute for their hospitality during this time.}
\date{25 August 2005}

\ifmemoirs@

\begin{abstract}
We study highest weight representations of
shifted Yangians over an algebraically closed field of characteristic $0$.
In particular, we classify the finite dimensional irreducible representations
and explain how to compute their Gelfand-Tsetlin characters in terms of
known characters of standard modules and
certain Kazhdan-Lusztig polynomials.
Our approach exploits the relationship between shifted Yangians and
the finite $W$-algebras associated to nilpotent orbits 
in general linear Lie algebras.
\end{abstract}

\maketitle
\tableofcontents
\mainmatter

\else

\chapter*{Representations of Shifted Yangians\\ and \\Finite \boldmath{$W$}-algebras}

\section*{Jonathan Brundan and Alexander Kleshchev}

\begin{center}
Department of Mathematics

University of Oregon

Eugene
OR 97403

USA
\end{center}

\vspace{8mm}

\begin{center}{\bf Abstract}\end{center}

\vspace{2mm}

\noindent
We study highest weight representations of
shifted Yangians over an algebraically closed field of characteristic $0$.
In particular, we classify the finite dimensional irreducible representations
and explain how to compute their Gelfand-Tsetlin characters in terms of
known characters of standard modules and
certain Kazhdan-Lusztig polynomials.
Our approach exploits the relationship between shifted Yangians and
the finite $W$-algebras associated to nilpotent orbits 
in general linear Lie algebras.

\vspace{6mm}

\begin{quote}
\begin{center}{\bf Contents}\end{center}

\vspace{2mm}

\noindent
\begin{tabular}{lr}
1. Introduction&\pageref{sintro}\\
2. Shifted Yangians&\pageref{syangians}\\
3. Finite $W$-algebras&\pageref{swalgebras}\\
4. Dual canonical bases&\pageref{sbases}\\
5. Highest weight theory&\pageref{shw}\\
6. Verma modules&\pageref{sverma}\\
7. Standard modules&\pageref{sstandard}\\
8. Character formulae\qquad\qquad\qquad\qquad\qquad\qquad\qquad\qquad&\pageref{sfdr}\\
References&\pageref{bib}\\
\end{tabular}
\end{quote}

\vspace{8mm}

\begin{center}{\bf Acknowledgements}\end{center}

\vspace{2mm}

\noindent
Part of this research was carried out
during a stay by the first author at the
(ex-) Institut Girard Desargues, Universit\'e Lyon I in Spring 2004.
He would like to thank Meinolf Geck 
and the other members of the institute for their hospitality during this time.

\noindent
This research was supported in part by NSF grant no. DMS-0139019.

\mainmatter

\fi

\chapter{Introduction}\label{sintro}

Following work of Premet,
there has been renewed interest recently in the representation theory
of certain algebras
that are associated to nilpotent orbits in complex semisimple Lie algebras. 
We refer to these algebras as {\em finite $W$-algebras}.
They should be viewed as analogues of universal 
enveloping algebras for the Slodowy slice
through the nilpotent orbit in question.
Actually, in the special cases considered in this article, the
definition of these algebras first appeared in 1979 in
the Ph.D. thesis of Lynch \cite{Ly}, extending the celebrated work of 
Kostant \cite{K} treating regular nilpotent orbits. However, despite 
quite a lot of attention by a number of authors since then,
see e.g. \cite{Kaw,Moe,Mat,BT,VD,GG, P, premet2, KS}, there is still
surprisingly little concrete information about the representation theory 
of these algebras to be found in the literature.
The goal in this article is to undertake a thorough 
study of finite dimensional representations of the finite 
$W$-algebras associated to nilpotent orbits in the 
Lie algebra $\mathfrak{gl}_N(\mathbb{C})$. We are able to make progress in this 
case thanks largely to the relationship 
between finite $W$-algebras and {\em shifted Yangians} 
first noticed in \cite{RS,BR} and developed in 
full generality in
\cite{BK}. 

Fix for the remainder of the introduction 
a partition $\lambda = (p_1 \leq \cdots \leq p_n)$ of $N$.
We draw the Young diagram of $\lambda$ in a slightly unconventional way, 
so that there are $p_i$ boxes in the $i$th row, numbering rows
$1,\dots,n$ from top to bottom in order of increasing length.
Also number the non-empty 
columns of this diagram 
by $1,\dots,l$ from left to right, and let $q_i$ denote the 
number of boxes in the $i$th column, 
so $\lambda' = (q_1 \geq \cdots \geq q_l)$ is the 
transpose partition to $\lambda$.
For example, if $(p_1,p_2,p_3) = (2,3,4)$ then the Young diagram of
$\lambda$ is
$$
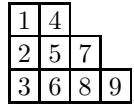
\Diagram{1&4\cr2&5&7\cr3&6&8&9\cr}
$$
and $(q_1,q_2,q_3,q_4) = (3,3,2,1)$.
We number the boxes of the diagram by 
$1,2,\dots,N$ down columns from left to right, and let $\row(i)$ and $\col(i)$ denote the row and column numbers 
of the $i$th box.

Writing $e_{i,j}$ for the $ij$-matrix unit in the Lie algebra
$\mathfrak{g} = \mathfrak{gl}_N(\mathbb{C})$, 
let $e$ denote the matrix $\sum_{i,j} 
e_{i,j}$
summing over all $1 \leq i,j \leq N$ such that $\row(i) = \row(j)$ and $\col(i) = \col(j)-1$. This is a nilpotent matrix of Jordan type $\lambda$.
For instance, if $\lambda$ is as above, then 
$e = e_{1,4}+e_{2,5}+e_{5,7}+e_{3,6}+e_{6,8}+e_{8,9}$.
Define a $\Z$-grading $\mathfrak{g} = \bigoplus_{j \in \Z} \mathfrak{g}_j$
of the Lie algebra $\mathfrak{g}$ by declaring that each $e_{i,j}$ is of degree $(\col(j)-\col(i))$.
This is a
{\em good grading} for $e\in \mathfrak{g}_1$ in the sense of \cite{KRW}
(see also \cite{EK} for the full classification). 
However, it is not the usual Dynkin grading arising from an $\mathfrak{sl}_2$-triple unless all the parts of $\lambda$ are equal. 
Actually, in the main body of the article, 
we work with more general good gradings than the one described here,
replacing the Young diagram of $\lambda$ with a more general diagram
called a {\em pyramid} and denoted by the symbol $\pi$; 
see $\S$\ref{sspyramids}. When the pyramid $\pi$ is left-justified, it coincides with the Young diagram of $\lambda$. We have chosen to focus just on this
case in the introduction,
since it plays a distinguished role in the theory.

Now we give a formal definition of
the finite $W$-algebra $W(\lambda)$ associated to this data.
Let $\mathfrak{p}$ denote the parabolic subalgebra 
$\bigoplus_{j \geq 0} \mathfrak{g}_j$ of $\mathfrak{g}$
with Levi factor $\mathfrak{h} =\mathfrak{g}_0$, 
and let $\mathfrak{m}$ denote the opposite nilradical $\bigoplus_{j < 0} \mathfrak{g}_j$.
Taking the trace form with $e$ defines a one dimensional representation
$\chi:\mathfrak{m} \rightarrow \mathbb{C}$.
Let $I_\chi$ be the two-sided ideal of the universal
enveloping algebra $U(\mathfrak{m})$ generated by $\ker \chi$.
Let $\eta:U(\mathfrak{p}) \rightarrow U(\mathfrak{p})$ be the
automorphism mapping
$e_{i,j} \mapsto e_{i,j}
+ \delta_{i,j}(n-q_{\col(j)}-q_{\col(j)+1}-\cdots - q_l)$
for each $e_{i,j} \in \mathfrak{p}$.
Then, by our definition,
$W(\lambda)$
is the following subalgebra of $U(\mathfrak{p})$:
$$
W(\lambda) = \{u \in U(\mathfrak{p})\:|\:
[x, \eta(u)] \in U(\mathfrak{g}) I_\chi
\text{ for all }x \in \mathfrak{m}\};
$$
see $\S$\ref{sswalgebras}.
The twist by the automorphism $\eta$ here is unconventional but quite convenient
later on; it is analogous to ``shifting by
$\rho$'' in the definition of Harish-Chandra homomorphism.
For examples, if the Young diagram of 
$\lambda$ consists of a single column and $e$ is the zero 
matrix, $W(\lambda)$ coincides with the entire universal enveloping algebra
$U(\mathfrak{g})$. At the other extreme, if the Young diagram of $\lambda$ 
consists of a single row
and $e$ is a regular nilpotent element, the work of Kostant \cite{K}
shows that $W(\lambda)$ is isomorphic to the 
center of $U(\mathfrak{g})$, in particular it is commutative.

For $u \in W(\lambda)$, right multiplication by $\eta(u)$ leaves
$U(\mathfrak{g}) I_\chi$ invariant, so induces a well-defined right
action of $u$ 
on the {\em generalized Gelfand-Graev representation}
$$
Q_\chi = U(\mathfrak{g}) / U(\mathfrak{g}) I_\chi
\cong U(\mathfrak{g}) \otimes_{U(\mathfrak{m})} \mathbb{C}_\chi.
$$
This makes $Q_\chi$ into a $(U(\mathfrak{g}), W(\lambda))$-bimodule.
The associated
algebra homomorphism $W(\lambda) \rightarrow \End_{U(\mathfrak{g})}(Q_\chi)^{\operatorname{op}}$ is actually an isomorphism, giving an
alternate definition of $W(\lambda)$ as an endomorphism algebra.

Another useful construction involves the homomorphism 
$\xi:U(\mathfrak{p}) \rightarrow U(\mathfrak{h})$
induced by the natural projection
$\mathfrak{p} \twoheadrightarrow \mathfrak{h}$.
The restriction of $\xi$ to
$W(\lambda)$ defines an {\em injective} algebra homomorphism
$W(\lambda) \hookrightarrow U(\mathfrak{h})$
which we call the {\em Miura transform}; see $\S$\ref{ssmiura}.
To explain its signigicance, we note that
$\mathfrak{h} = 
\mathfrak{gl}_{q_1}(\mathbb{C}) \oplus\cdots\oplus \mathfrak{gl}_{q_l}(\mathbb{C})$, so $U(\mathfrak{h})$ is naturally identified
with the tensor product
$U(\mathfrak{gl}_{q_1}(\mathbb{C})) \otimes\cdots\otimes U(\mathfrak{gl}_{q_l}(\mathbb{C}))$. Given $\mathfrak{gl}_{q_i}(\mathbb{C})$-modules $M_i$
for each $i=1,\dots,l$, the outer tensor product
$M_1 \boxtimes \cdots \boxtimes M_l$ is therefore a
$U(\mathfrak{h})$-module in the natural way. Hence, via the Miura transform,
$M_1 \boxtimes\cdots\boxtimes M_l$ is a $W(\lambda)$-module too.
This construction plays the role of tensor product in the
representation theory of $W(\lambda)$.

Next we want to recall the connection between
$W(\lambda)$ and shifted Yangians.
Let $\sigma$ be the upper triangular $n \times n$  
matrix with $ij$-entry $(p_j - p_i)$
for $i \leq j$. 
The {\em shifted Yangian} $Y_n(\sigma)$ associated to $\sigma$
is the associative algebra over $\mathbb{C}$ with generators
$D_i^{(r)}\:(1 \leq i \leq n, r > 0), 
E_i^{(r)}\:(1 \leq i < n, r > p_{i+1}-p_i)$ and
$F_i^{(r)}\:(1 \leq i <n, r > 0)$ subject 
to certain relations recorded explicitly in $\S$\ref{ssgenerators}.
In the case that $\sigma$ is the zero matrix, i.e. all parts of $\lambda$ are equal, $Y_n(\sigma)$ is precisely the usual Yangian $Y_n$
associated to the Lie algebra
$\mathfrak{gl}_n(\mathbb{C})$ 
and the defining 
relations are a variation on the Drinfeld presentation of \cite{D3}; see 
\cite{BKdrinfeld}. 
In general, the presentation of $Y_n(\sigma)$
is adapted to its natural {\em  triangular
decomposition}, allowing us to study representations 
in terms of highest weight theory. In particular, the subalgebra generated
by all the elements $D_i^{(r)}$ is a maximal commutative subalgebra
which we call the {\em Gelfand-Tsetlin subalgebra}.
We often work with the generating functions
$$
D_i(u) = 1 + D_i^{(1)} u^{-1} + D_i^{(2)} u^{-2} + \cdots
\in Y_n(\sigma)[[u^{-1}]].
$$
The main result of \cite{BK}
shows that 
the finite $W$-algebra $W(\lambda)$ is isomorphic to 
the quotient of $Y_n(\sigma)$
by the two-sided ideal generated by all
$D_1^{(r)}\:(r > p_1)$.
The precise identification of $W(\lambda)$ with this quotient
is described in $\S$\ref{ssquotients}.
Also in $\S$\ref{ssmiura}, we explain how the tensor product
construction outlined in the previous paragraph
is induced by the comultiplication of the Hopf algebra $Y_n$.

We are ready to describe the first results about
representation theory.
We call a vector $v$ in a $Y_n(\sigma)$-module $M$
a {\em highest weight vector} if 
it is annihilated by all
$E_i^{(r)}$ and each $D_i^{(r)}$ acts on $v$ by a scalar.
A critical point is that if $v$ is a highest weight vector in 
a $W(\lambda)$-module, viewed as a $Y_n(\sigma)$-module via the map
$Y_n(\sigma)\twoheadrightarrow W(\lambda)$, 
then in fact $D_i^{(r)} v = 0$
for all $r > p_i$. This is obvious for $i=1$, since the image of
$D_1^{(r)}$ in $W(\lambda)$ is zero by the definition of the map
for all $r > p_1$.
For $i > 1$, it follows from the following fundamental result
proved in $\S$\ref{svan}:
for any $i$ and $r > p_i$, the image of $D_i^{(r)}$
in $W(\lambda)$ is congruent to zero modulo the left ideal 
generated by all $E_j^{(s)}$.
Hence, if $v$ is a highest weight vector in
a $W(\lambda)$-module, then there exist scalars
$(a_{i,j})_{1 \leq i \leq n, 1 \leq j \leq p_i}$ such that
\begin{align*}
u^{p_1} D_1(u) v &= (u+a_{1,1}) (u+a_{1,2}) \cdots (u+a_{1,p_1}) v,\\
(u-1)^{p_2} D_2(u-1) v &= (u+a_{2,1}) (u+a_{2,2}) \cdots (u+a_{2,p_2}) 
v,\\\intertext{\vspace{-3mm}$$\vdots$$\vspace{-6mm}}
(u-n+1)^{p_n} D_n(u-n+1) v &= (u+a_{n,1}) (u+a_{n,2}) \cdots (u+a_{n,p_n}) 
v.
\end{align*}
Let $A$ be the $\lambda$-tableau obtained
by writing the scalars
$a_{i,1},\dots,a_{i,p_i}$ into the boxes on the $i$th row
of the Young diagram of $\lambda$. 
In this way, the highest weights
that can arise in $W(\lambda)$-modules
are parametrized by the set $\Row(\lambda)$ of 
{\em row symmetrized $\lambda$-tableaux}, i.e. tableaux of shape
$\lambda$ with entries 
from $\mathbb{C}$ viewed up to row equivalence.
Conversely, given any row symmetrized $\lambda$-tableau 
$A \in \Row(\lambda)$,
there exists a (non-zero) universal highest weight module
$M(A)$ generated by such a highest weight vector; see $\S$\ref{ssp}.
We call $M(A)$ the {\em generalized Verma module} of type $A$.
By familiar arguments, $M(A)$ has a unique irreducible quotient $L(A)$,
and then the modules $\{L(A)\:|\:A \in \Row(\lambda)\}$ give all irreducible
highest weight modules for $W(\lambda)$ up to isomorphism.

There is a natural abelian category $\mathcal M(\lambda)$
which is an analogue of the BGG category $\mathcal O$ for the
algebra $W(\lambda)$; see
$\S$\ref{ssgg}. 
(Actually, $\mathcal M(\lambda)$ is more like the category $\mathcal O^\infty$
obtained by weakening the hypothesis that a Cartan subalgebra acts semisimply in
the usual definition of $\mathcal O$.)
All objects in $\mathcal M(\lambda)$ are of finite length,
and the simple objects are precisely the irreducible highest weight modules,
hence the isomorphism classes $\{[L(A)]\:|\:A \in \Row(\lambda)\}$ give 
a canonical basis for the Grothendieck group $[\mathcal M(\lambda)]$ 
of the category $\mathcal M(\lambda)$. The generalized Verma modules belong to
$\mathcal M(\lambda)$ too, and it is natural to consider the
{\em composition multiplicities} $[M(A):L(B)]$ for $A, B \in \Row(\lambda)$.
We will formulate a precise combinatorial 
conjecture for
these, in the spirit of the Kazhdan-Lusztig conjecture, later
on in the introduction.
For now, we just record the following 
basic result about the structure of Verma modules; see $\S$\ref{sslink}.
For the statement, let $\leq$ denote the Bruhat ordering
on row symmetrized $\lambda$-tableaux; see $\S$\ref{sstableaux}.

\vspace{2mm}
{{\sc Theorem A} {\bf (Linkage principle)}. }{\em
For $A, B \in \Row(\lambda)$, the composition multiplicity 
$[M(A):L(A)]$ is equal to $1$, and
$[M(A):L(B)] \neq 0$ if and only if $B \leq A$ in the Bruhat ordering.
}
\vspace{2mm}

Hence,
$\{[M(A)]\:|\:A \in \Row(\lambda)\}$ is another natural basis for 
the Grothendieck group $[\mathcal M(\lambda)]$. 
We want to say a few words about the proof of Theorem A, since it 
involves an interesting technique.
Modules in the category $\mathcal M(\lambda)$
possess {\em Gelfand-Tsetlin characters}; see $\S$\ref{sscharblock}.
This is a formal notion that keeps track of the dimensions of the
generalized weight space decomposition of 
a module with respect to the Gelfand-Tsetlin subalgebra of $Y_n(\sigma)$,
similar in spirit to the $q$-characters of Frenkel and Reshetikhin \cite{FR}.
The map sending a module to its Gelfand-Tsetlin character
induces an embedding of the Grothendieck group $[\mathcal M(\lambda)]$ into 
a certain completion of the
ring of Laurent
polynomials
$\Z[y_{i,a}^{\pm 1}\:|\:i=1,\dots,n, a \in \mathbb{C}]$, for 
indeterminates $y_{i,a}$. The key step in our proof of Theorem A is the computation of the Gelfand-Tsetlin character of the Verma module $M(A)$ itself; see
$\S$\ref{sschars} for the precise statement.
In general, $\ch M(A)$ is an infinite sum of monomials in
the $y_{i,a}^{\pm 1}$'s involving 
both positive and negative powers, 
but the highest weight vector of $M(A)$ contributes
just the positive monomial
$$
y_{1,a_{1,1}} \dots y_{1,a_{1,p_1}} 
\times
y_{2,a_{2,1}} \dots y_{2,a_{2,p_2}} 
\times
\cdots
\times
y_{n,a_{n,1}} \dots y_{n,a_{n,p_n}},
$$
where $a_{i,1},\dots,a_{i,p_i}$ are the entries in the $i$th row of $A$ as above.
The highest weight vector of any composition factor 
contributes a similar such positive monomial.
So by analyzing the positive monomials
appearing in the formula for
$\ch M(A)$, we get information about the possible 
$L(B)$'s that can be composition factors of $M(A)$. The Bruhat ordering 
on tableaux emerges naturally out of these considerations.

Another important property of Verma modules has to do with tensor products.
Let $A\in\Row(\lambda)$ be a row symmetrized $\lambda$-tableau.
Pick any representative for it and let $A_i$ denote the $i$th column
of this representative 
with entries $a_{i,1},\dots,a_{i,q_i}$ read from top to
bottom.
Let $M(A_i)$ denote the 
usual Verma module for the Lie algebra $\mathfrak{gl}_{q_i}(\mathbb{C})$
generated by a highest weight vector
$v_+$ annihilated by all strictly upper triangular matrices
and on which $e_{j,j}$ acts as the scalar
$(a_{i,j}+n-q_i+j-1)$ for each $j=1,\dots,q_i$.
Via the Miura transform,  the tensor product 
$M(A_1) \boxtimes \cdots \boxtimes M(A_l)$ is then naturally a 
$W(\lambda)$-module 
as explained above, and the vector
$v_+ \otimes\cdots\otimes v_+$ is a highest weight vector
in this tensor product of type $A$.  
In fact, our formula for the Gelfand-Tsetlin character of $M(A)$ implies that
$$
[M(A)]= [M(A_1) \boxtimes \cdots \boxtimes M(A_l)],
$$
equality in the Grothendieck group $[\mathcal M(\lambda)]$. 
The first part of the 
next theorem, proved in $\S$\ref{sscenter}, is a
consequence of this 
equality;
 the second part is an application of
\cite{FO}. 

\vspace{2mm}
{{\sc Theorem B} {\boldmath\bf (Structure of center)}.
\em
Identifying $W(\lambda)$ with the endomorphism algebra 
of $Q_\chi$,
the natural multiplication map
$\psi:Z(U(\mathfrak{g})) \rightarrow \End_{U(\mathfrak{g})}(Q_\chi)$ defines an algebra isomorphism 
between the center of $U(\mathfrak{g})$
and the center of $W(\lambda)$. Moreover, $W(\lambda)$ 
is free as a module over its center.
}
\vspace{2mm}

Now we switch our attention to finite dimensional $W(\lambda)$-modules.
Let $\mathcal F(\lambda)$ denote the category of all finite dimensional
$W(\lambda)$-module, viewed as a 
subcategory of the category $\mathcal M(\lambda)$.
The problem of classifying all finite dimensional
irreducible $W(\lambda)$-modules reduces to determining precisely which
$A \in \Row(\lambda)$ have the property that 
$L(A)$ is finite dimensional.
To formulate the final result, we need one more definition.
Call a $\lambda$-tableau $A$ with entries in $\mathbb{C}$ 
{\em column strict} if in every column the entries
belong to the same coset of $\mathbb{C}$ modulo $\mathbb{Z}$ and 
are strictly increasing from bottom to top.
Let $\Col(\lambda)$ denote the set of all such column strict $\lambda$-tableaux.
There is an obvious map
$$
R:\Col(\lambda) \rightarrow \Row(\lambda)
$$
mapping a $\lambda$-tableau to its row equivalence class.
Let $\Dom(\lambda)$ denote the image of this map, the set of all
{\em dominant} row symmetrized $\lambda$-tableaux.

\vspace{2mm}
{{\sc Theorem C} {\bf (Finite dimensional irreducible representations)}. }
{\em
For $A \in \Row(\lambda)$, 
the irreducible highest weight module $L(A)$ is finite dimensional
if and only if $A$ is dominant.
Hence, 
the modules $\{L(A)\:|\:A \in \Dom(\lambda)\}$ form a complete set of pairwise non-isomorphic finite dimensional irreducible $W(\lambda)$-modules.
}
\vspace{2mm}

To prove this, there are two steps: one needs to show first that each $L(A)$
with $A \in \Dom(\lambda)$ is finite dimensional, and second that
all other $L(A)$'s are infinite dimensional.
Let us explain the argument for the first step.
Given $A \in \Col(\lambda)$, let $A_i$ be its $i$th column
and define $L(A_i)$ to be the unique irreducible quotient of 
the Verma module $M(A_i)$ introduced above.
Because $A$ is column strict, each $L(A_i)$ is a finite dimensional
irreducible $\mathfrak{gl}_{q_i}(\mathbb{C})$-module.
Hence we obtain a finite dimensional $W(\lambda)$-module
$$
V(A) = 
L(A_1) \boxtimes\cdots\boxtimes L(A_l),
$$ 
which we call the {\em standard module} corresponding to $A \in 
\Col(\lambda)$. It contains an obvious highest weight vector of type
equal to the row equivalence class of $A$.
This simple construction is enough to finish 
the first step of the proof.
The second step is actually much harder,
and is an extension of the proof due to Tarasov
\cite{Tarasov1, Tarasov2} and Drinfeld \cite{D3}
of the classification of finite dimensional irreducible
representations of the Yangian $Y_n$ by {\em Drinfeld polynomials}.
It is based on the remarkable fact that when $n=2$, i.e. the Young diagram 
of $\lambda$ 
has just two rows, {\em every} $L(A)\:(A \in \Row(\lambda))$
can be expressed as an irreducible tensor product; see $\S$\ref{sstr}.

Amongst all the standard modules, there are some special ones which are 
highest weight modules and whose isomorphism classes form a basis for the
Grothendieck group of the category $\mathcal F(\lambda)$. 
Let $A \in \Col(\lambda)$ be a column strict $\lambda$-tableau
with entries $a_{i,1},\dots,a_{i,p_i}$ in its $i$th row read from left to right.
We say that $A$ is {\em standard}
if $a_{i,j} \leq a_{i,k}$ for every $1 \leq i \leq n$ and
$1 \leq j < k \leq p_i$ such that $a_{i,j}$ and $a_{i,k}$ belong to the same
coset of $\mathbb{C}$ modulo $\mathbb{Z}$.
If all entries of $A$ are integers, this is the usual definition of a 
standard tableau: entries increase strictly up columns and weakly along rows.
Let $\Std(\lambda)$ denote the set of all standard $\lambda$-tableaux
$A \in \Col(\lambda)$.
Our proof of the next theorem is based on an argument due to Chari
\cite{Chari} in the context of quantum affine algebras; see $\S$\ref{sstp}.

\vspace{2mm}
{{\sc Theorem D} {\bf (Highest weight standard modules)}.}
{\em
For $A \in \Std(\lambda)$,
the standard module $V(A)$ is a highest weight module of highest weight 
equal to the row equivalence class of $A$.
}
\vspace{2mm}

Most of the results so far are analogous to well known results in the 
representation theory of Yangians and quantum affine algebras, and do not 
exploit the finite $W$-algebra side of the picture in any significant way. 
To remedy this, we need to apply {\em Skryabin's theorem} from \cite{Skry};
see $\S$\ref{ssskryabin}.
This asserts that the functor $Q_\chi \otimes_{W(\lambda)} ?$
gives an equivalence of categories
between the category of all 
$W(\lambda)$-modules and the category
$\mathcal W(\lambda)$ of all {\em generalized Whittaker modules},
namely, all
$\mathfrak{g}$-modules on which $(x-\chi(x))$ acts locally nilpotently
for all $x \in \mathfrak{m}$.
For any finite dimensional $\mathfrak{g}$-module $V$,
it is obvious that the functor $? \otimes V$ maps objects in
$\mathcal W(\lambda)$ to objects in $\mathcal W(\lambda)$.
Transporting through Skryabin's equivalence of categories, we obtain a functor $? \circledast V$ on $W(\lambda)\Mod$ itself; see $\S$\ref{ssti}. 
In this way, one can introduce {\em translation functors}
on the categories $\mathcal M(\lambda)$ and $\mathcal F(\lambda)$.
Actually, we just need some special translation functors, peculiar to 
the type $A$ theory and  denoted
$e_i, f_i$ for $i \in \mathbb{C}$, 
which arise from $\circledast$'ing
with the natural $\mathfrak{gl}_N(\mathbb{C})$-module 
and its dual. These functors fit into the axiomatic framework
developed recently by Chuang and Rouquier \cite{CR}; see $\S$\ref{sstf}.

Now recall the parabolic subalgebra $\mathfrak{p}$ of $\mathfrak{g}$ with
Levi factor $\mathfrak{h}$. We let $\mathcal O(\lambda)$ 
denote the corresponding
parabolic category $\mathcal O$, the category of all finitely generated
$\mathfrak{g}$-modules
on which $\mathfrak{p}$ acts locally finitely and $\mathfrak{h}$ acts 
semisimply. 
For $A \in \Col(\lambda)$
with entry $a_i$ in its $i$th box, we let $N(A) \in \mathcal O(\lambda)$ 
denote the {\em parabolic Verma 
module} generated by
a highest weight vector $v_+$ that is annihilated by all 
strictly upper triangular matrices in $\mathfrak{g}$
and on which $e_{i,i}$ acts as the scalar
$(a_i+i-1)$ for each $i=1,\dots,N$.
Let $K(A)$ denote the unique irreducible quotient of $N(A)$.
Both of the sets
$\{[N(A)]\:|\:A \in \Col(\lambda)\}$ and $\{[K(A)]\:|\:A \in \Col(\lambda)\}$
form natural bases for the Grothendieck group $[\mathcal O(\lambda)]$.
There is a remarkable functor
$$
\V:\mathcal O(\lambda) \rightarrow \mathcal F(\lambda)
$$
introduced originally (in a slightly different form)
by Kostant and Lynch.
We call it the {\em Whittaker functor};
see $\S$\ref{ssmain}.
It is an exact functor preserving central characters and
commuting with translation functors.
Moreover, it maps the
parabolic Verma module $N(A)$ to the standard module $V(A)$
for every $A \in \Col(\lambda)$. The culmination of this 
article is the following theorem.

\vspace{2mm}
{{\sc Theorem E} {\bf (Construction of irreducible modules)}.} 
{\em 
The Whittaker functor $\V:\mathcal O(\lambda) \rightarrow \mathcal F(\lambda)$
sends 
irreducible modules to irreducible modules or zero.
More precisely, take any $A \in \Col(\lambda)$ and 
let $B \in \Row(\lambda)$ be its row equivalence class.
Then
$$
\V(K(A)) \cong \left\{\begin{array}{ll}
L(B)&\text{if $A$ is standard},\\
0&\text{otherwise.}
\end{array}\right.
$$
Every finite dimensional 
irreducible $W(\lambda)$-module arises in this way.
}
\vspace{2mm}

There are three main ingredients to the proof of this
theorem. 
First, we need detailed information about
the translation functors $e_i, f_i$, much of which 
is provided by
\cite{CR} as an application of 
the representation theory of degenerate affine Hecke algebras.
Second, we need to know that the standard modules $V(A)$
have simple cosocle if $A \in \Std(\lambda)$, which follows from 
Theorem D. Finally, we need to apply 
the Kazhdan-Lusztig conjecture for the Lie algebra
$\mathfrak{gl}_N(\mathbb{C})$ in order to determine 
exactly when $\V(K(A))$ is non-zero.

Let us discuss some of 
the combinatorial consequences of Theorem E in more detail.
For this, we at last restrict our attention just 
to modules having integral central character. Let
$\Row_0(\lambda), \Col_0(\lambda), \Dom_0(\lambda)$ and $\Std_0(\lambda)$ denote the subsets
of $\Row(\lambda), \Col(\lambda), \Dom(\lambda)$ and $\Std(\lambda)$ consisting of the tableaux 
all of whose entries are integers.
The restriction of the 
map $R$ actually gives a bijection between the sets
$\Std_0(\lambda)$ and $\Dom_0(\lambda)$.
Let $\mathcal O_0(\lambda)$, $\mathcal F_0(\lambda)$ and $\mathcal M_0(\lambda)$
denote the full subcategories of $\mathcal O(\lambda)$, 
$\mathcal F(\lambda)$ and $\mathcal M(\lambda)$ consiting of objects
all of whose composition factors are of the form
$\{K(A)\:|\:A \in \Col_0(\lambda)\}$,
$\{L(A)\:|\:A \in \Dom_0(\lambda)\}$ and
$\{L(A)\:|\:A \in \Row_0(\lambda)\}$, respectively.
The isomorphism classes of these three sets of objects give
canonical bases for the Grothendieck groups
$[\mathcal O_0(\lambda)], [\mathcal F_0(\lambda)]$ and $[\mathcal M_0(\lambda)]$,
as do the isomorphism classes of the
parabolic Verma modules $\{N(A)\:|\:A \in \Col_0(\lambda)\}$,
the standard modules $\{V(A)\:|\:A \in \Std_0(\lambda)\}$, and the
generalized Verma modules $\{M(A)\:|\:A \in \Row_0(\lambda)\}$, respectively.

The functor $\V$ above restricts to an exact
functor $\V:\mathcal O_0(\lambda) \rightarrow \mathcal F_0(\lambda)$, 
and we also have the natural embedding $\I$ of 
the category $\mathcal F_0(\lambda)$ into
$\mathcal M_0(\lambda)$. At the level of Grothendieck groups, these functors
induce maps
$$
[\mathcal O_0(\lambda)] \stackrel{\V}{\twoheadrightarrow}
[\mathcal F_0(\lambda)] \stackrel{\I}{\hookrightarrow}
[\mathcal M_0(\lambda)].
$$
The translation functors $e_i, f_i$ for $i \in \Z$
(and more generally their divided powers 
$e_i^{(r)}, f_i^{(r)}$ defined as in \cite{CR})
induce maps also denoted $e_i, f_i$ on all these 
Grothendieck groups.
The resulting maps 
satisfy the relations of the Chevalley generators (and their divided powers) for the Kostant $\Z$-form $U_\Z$ of the universal enveloping 
algebra
of the Lie algebra $\mathfrak{gl}_\infty(\mathbb{C})$,
that is, the Lie algebra 
of matrices with rows and columns labelled by $\Z$ 
all but finitely many of which are zero.
The maps $\V$ and $\I$ are then $U_\Z$-module homomorphisms with respect to these
actions.

Now the point is that all of this categorifies a well known situation in linear algebra. Let $V_\Z$ denote the natural $U_\Z$-module, with basis
$v_i\:(i \in \Z)$.
We write $\bigwedge^{\lambda'}(V_\Z)$ for the tensor product
$\bigwedge^{q_1}(V_\Z) \otimes\cdots\otimes \bigwedge^{q_l}(V_\Z)$
and
$S^{\lambda}(V_\Z)$ for the tensor product
$S^{p_1}(V_\Z) \otimes\cdots\otimes S^{p_n}(V_\Z)$. 
These free $\Z$-modules have natural monomial bases 
denoted 
$\{N_A\:|\:A \in \Col_0(\lambda)\}$
and $\{M_A\:|\:A \in \Row_0(\lambda)\}$, 
respectively; see $\S$\ref{ssdcb}.
A well known consequence of the Littlewood-Richardson 
rule (observed already by Young long before) implies that the space
$$
\hom_{U_\Z}({\textstyle \bigwedge^{\lambda'}}(V_\Z), S^{\lambda}(V_\Z))
$$
is a free $\Z$-module of rank one; indeed, there is a 
canonical $U_\Z$-module homomorphism
$\V:\bigwedge^{\lambda'}(V_\Z) 
\rightarrow S^{\lambda}(V_\Z)$
that generates the space of all such maps.
The image of this map is $P^\lambda(V_\Z)$, a familiar $\Z$-form
for the 
{\em irreducible polynomial representation} of $\mathfrak{gl}_\infty(\mathbb{C})$ 
labelled by the partition $\lambda$.
So by definition $P^\lambda(V_\Z)$ is a subspace of $S^\lambda(V_\Z)$;
we denote the natural inclusion by
$\mathbb{I}$.
Recall $P^\lambda(V_\Z)$ also possesses a standard monomial basis
$\{V_A\:|\:A \in \Std_0(\lambda)\}$, defined 
from $V_A = \V(N_A)$.
Finally, we let $i:\bigwedge^{\lambda'}(V_\Z)
\rightarrow [\mathcal O_0(\lambda)]$,
$j:P^\lambda(V_\Z) \rightarrow [\mathcal F_0(\lambda)]$ and
$k:S^\lambda(V_\Z) \rightarrow [\mathcal M_0(\lambda)]$ 
be the $\Z$-module homomorphisms
sending $N_A \mapsto N(A),
V_A \mapsto [V(A)]$ and $M_A \mapsto [M(A)]$
for $A \in \Col_0(\lambda)$, $A \in \Std_0(\lambda)$ and $A \in \Row_0(\lambda)$,
respectively. 

\vspace{2mm}
{{\sc Theorem F} {\bf (Categorification of polynomial functors)}.}
{\em
The maps $i,j,k$ are all isomorphisms of $U_\Z$-modules, and 
the following diagram commutes:
$$
\begin{CD}
\bigwedge^{\lambda'}(V_\Z) &@>\V>>&P^\lambda(V_\Z) &@>\I>> & S^\lambda(V_\Z)\\
@ViVV&&@VVjV&&@VVkV\\
[\mathcal O_0(\lambda)]&@>>\V>&[\mathcal F_0(\lambda)]&@>>\I>&[\mathcal M_0(\lambda)].
\end{CD}
$$
Moreover, setting $L_A = j^{-1}([L(A)])$ for $A \in \Dom_0(\lambda)$, the basis
$\{L_A\:|\:A \in \Dom_0(\lambda)\}$ 
coincides with 
Lusztig's dual canonical basis/Kashiwara's upper global crystal basis
for the polynomial representation $P^\lambda(V_\Z)$.
}
\vspace{2mm}

Again, the Kazhdan-Lusztig conjecture plays the central role in the proof of this 
theorem.
Actually, we use the following 
increasingly well known reformulation of the Kazhdan-Lusztig conjecture in 
type $A$: setting $K_A = i^{-1}([K(A)])$, the basis
$\{K_A\:|\:A \in \Col_0(\lambda)\}$ 
coincides with the dual canonical basis
for the space $\bigwedge^{\lambda'}(V_\Z)$. 
In particular, this implies that the decomposition
numbers $[V(A):L(B)]$ for $A \in \Std_0(\lambda)$ and
$B \in \Dom_0(\lambda)$ 
can be computed in terms of certain Kazhdan-Lusztig polynomials
associated to the symmetric group $S_N$ evaluated at $q=1$.
From a special case, one can also recover the analogous result
for the Yangian $Y_n$ itself. We mention this, because it is 
interesting to compare the strategy followed here with that of
Arakawa \cite{A}, who also computes
the decomposition matrices of the Yangian in terms of Kazhdan-Lusztig polynomials starting from the Kazhdan-Lusztig
conjecture for the Lie algebra $\mathfrak{gl}_N(\mathbb{C})$, via \cite{AS}. There might also be a geometric approach 
to representation theory 
 of shifted Yangians 
in the spirit of \cite{var}.

As promised earlier in the introduction, let us now formulate a precise 
conjecture that explains how to compute
the decomposition numbers $[M(A):L(B)]$ for
all $A, B \in \Row_0(\lambda)$, 
also in terms of Kazhdan-Lusztig
polynomials associated to the symmetric group $S_N$.
Setting $L_A = k^{-1}([L(A)])$ for any $A \in \Row_0(\lambda)$, 
we conjecture that $\{L_A\:|\:A \in \Row_0(\lambda)\}$
coincides with the dual canonical basis for the
space $S^\lambda(V_\Z)$;
see $\S$\ref{ssgg}. 
This is a purely combinatorial reformulation in type $A$ of the
conjecture of de Vos and van Driel \cite{VD}
for arbitrary finite $W$-algebras, and is consistent with an idea of Premet
that there should be an equivalence of categories between the category
$\mathcal M(\lambda)$ here and a certain category
$\mathcal N(\lambda)$ considered by Mili\v c\'ic and Soergel \cite{MS}.
Our conjecture is known to be true in the special case that the Young diagram of
$\lambda$ consists of a single column: in that case it is 
precisely the Kazhdan-Lusztig conjecture for the 
Lie algebra $\mathfrak{gl}_N(\mathbb{C})$.
It 
is also true if the Young diagram of 
$\lambda$ has at most two rows, as can be verified 
by comparing 
the explicit construction of the simple highest weight
modules in the two row case from $\S$\ref{sstr}
with the explicit description of the dual canonical basis in this
case from \cite[Theorem 20]{qla}.
Finally, Theorem E would be an easy consequence
of this conjecture.

In a forthcoming article \cite{BKnew},
we will study the categories of
{\em polynomial} and {\em rational}
representations of $W(\pi)$ in more detail.
In particular, we will make precise the
relationship between polynomial representations of $W(\pi)$
and representations of
degenerate cyclotomic Hecke algebras, and we will 
relate the Whittaker functor $\V$
to work of Soergel \cite{Soergel} and Backelin \cite{Back}.
This should have applications to the representation theory
of {\em affine $W$-algebras} in the spirit of \cite{A2}.

\chapter{Shifted Yangians}\label{syangians}

 We will work from now on over 
an algebraically closed field $\C$ of characteristic $0$.
Let $\geq$ denote the partial order on $\C$ defined by
$x \geq y$ if $(x-y) \in \N$, where
$\N$ denotes $\{0,1,2,\dots\} \subset \C$.
We write simply $\mathfrak{gl}_n$ for the Lie algebra
$\mathfrak{gl}_n(\C)$.
In this preliminary chapter, we collect some basic definitions and
results about shifted Yangians, most of which are taken from \cite{BK}.
By a {\em shift matrix}\label{shift} we mean a matrix
$\sigma = (s_{i,j})_{1 \leq i,j \leq n}$ of
non-negative integers such that
\begin{equation}\label{shiftcon}
s_{i,j} + s_{j,k} = s_{i,k}
\end{equation}
whenever $|i-j|+|j-k|=|i-k|$.
Note this means that $s_{1,1}=\cdots=s_{n,n} = 0$, and the
matrix $\sigma$ is completely determined by the 
upper diagonal entries
$s_{1,2}, s_{2,3},\dots,s_{n-1,n}$
and the lower diagonal entries
$s_{2,1}, s_{3,2},\dots,s_{n,n-1}$.
We fix such a matrix $\sigma$ throughout the chapter.

\section{Generators and relations}\label{ssgenerators}

The {\em shifted Yangian} associated to the matrix $\sigma$
is the algebra $Y_n(\sigma)$\label{yangian} over $\C$
defined by generators \label{def}
\begin{align}
&\{D_i^{(r)}\:|\:1 \leq i \leq n, r > 0\},\label{gen1}\\
&\{E_i^{(r)}\:|\:1 \leq i < n, r > s_{i,i+1}\},\\
&\{F_i^{(r)}\:|\:1 \leq i < n, r > s_{i+1,i}\}\label{gen3}
\end{align}
subject to 
certain relations.
In order to write down these relations,
let 
\begin{equation}\label{didef}
D_i(u) := \sum_{r \geq 0} D_i^{(r)} u^{-r}
\in Y_n(\sigma)[[u^{-1}]]
\end{equation}
where $D_i^{(0)} := 1$, and then
define some  new elements $\widetilde{D}_i^{(r)}$ of $Y_n(\sigma)$
from the equation
\begin{equation}
\widetilde{D}_i(u) =
\sum_{r \geq 0} \widetilde{D}_i^{(r)}u^{-r}
:=-D_i(u)^{-1}.
\end{equation}
With this notation, the relations are as follows.
\begin{align}
[D_i^{(r)}, D_j^{(s)}] &=  0,\label{r2}\\
[E_i^{(r)}, F_j^{(s)}] &= \delta_{i,j} 
\sum_{t=0}^{r+s-1} \widetilde D_{i}^{(t)}D_{i+1}^{(r+s-1-t)} ,\label{r3}
\end{align}\begin{align}
[D_i^{(r)}, E_j^{(s)}] &= (\delta_{i,j}-\delta_{i,j+1})
\sum_{t=0}^{r-1} D_i^{(t)} E_j^{(r+s-1-t)},\label{r4}\\
[D_i^{(r)}, F_j^{(s)}] &= (\delta_{i,j+1}-\delta_{i,j})
\sum_{t=0}^{r-1} F_j^{(r+s-1-t)}D_i^{(t)} ,\label{r5}\end{align}
\begin{align}
[E_i^{(r)}, E_i^{(s+1)}] - [E_i^{(r+1)}, E_i^{(s)}] &=
E_i^{(r)} E_i^{(s)} + E_i^{(s)} E_i^{(r)},\label{r6b}\\
[F_i^{(r+1)}, F_i^{(s)}] - [F_i^{(r)}, F_i^{(s+1)}] &=
F_i^{(r)} F_i^{(s)} + F_i^{(s)} F_i^{(r)},\label{r7b}
\end{align}
\begin{align}
 [E_i^{(r)}, E_{i+1}^{(s+1)}]- [E_i^{(r+1)}, E_{i+1}^{(s)}] &=
-E_i^{(r)} E_{i+1}^{(s)},\label{r8}\\
[F_i^{(r+1)}, F_{i+1}^{(s)}] - [F_i^{(r)}, F_{i+1}^{(s+1)}] &=
 -F_{i+1}^{(s)} F_i^{(r)},\label{r9}\end{align}\begin{align}
[E_i^{(r)}, E_j^{(s)}] &= 0 \qquad\text{if }|i-j|> 1,\label{r10}\\
[F_i^{(r)}, F_j^{(s)}] &= 0 \qquad\text{if }|i-j|> 1,\label{r11}\\
[E_i^{(r)}, [E_i^{(s)}, E_j^{(t)}]] + 
[E_i^{(s)}, [E_i^{(r)}, E_j^{(t)}]] &= 0 \qquad\text{if }|i-j|=1,\label{r12}\\
[F_i^{(r)}, [F_i^{(s)}, F_j^{(t)}]] + 
[F_i^{(s)}, [F_i^{(r)}, F_j^{(t)}]] &= 0 \qquad\text{if }|i-j|=1,\label{r13}
\end{align}
for all meaningful $r,s,t,i,j$.
(For example, the relation (\ref{r8}) should be understood to hold for all
$i=1,\dots,n-2$, $r > s_{i,i+1}$ and $s > s_{i+1,i+2}$.)

It is often helpful to view $Y_n(\sigma)$ as an algebra graded by the root
lattice $Q_n$\label{qn} associated to the Lie algebra $\mathfrak{gl}_n$.
Let  $\mathfrak{c}$\label{c} be the (abelian) Lie subalgebra of
$Y_n(\sigma)$ spanned by the elements $D_1^{(1)},\dots,D_n^{(1)}$.
Let $\eps_1,\dots,\eps_n$\label{eps} be the basis for $\mathfrak{c}^*$ dual to the basis
$D_1^{(1)},\dots,D_n^{(1)}$.
We refer to elements of $\mathfrak{c}^*$ as {\em weights}\label{pn}
and elements of \begin{equation}
P_n := \bigoplus_{i=1}^n \Z \eps_i \subset \mathfrak{c}^*
\end{equation}
as {\em integral weights}.
The {\em root lattice} associated to the Lie algebra
$\mathfrak{gl}_n$ is then the $\Z$-submodule $Q_n$ of $P_n$
spanned by the {\em simple roots} 
$\eps_i - \eps_{i+1}$ for $i=1,\dots,n-1$.
We have the usual {\em dominance ordering} on $\mathfrak{c}^*$ defined 
by $\alpha \geq \beta$ if $(\alpha-\beta)$ is a sum of simple roots.
With this notation set up, the relations imply that 
we can define a $Q_n$-grading 
\begin{equation}\label{qgrading}
Y_n(\sigma) = \bigoplus_{\alpha \in Q_n} (Y_n(\sigma))_{\alpha}
\end{equation}
of the algebra $Y_n(\sigma)$
by declaring that the generators $D_i^{(r)}, E_i^{(r)}$ and $F_i^{(r)}$
are of degrees $0,\eps_i-\eps_{i+1}$ and $\eps_{i+1}-\eps_i$, respectively.

\section{PBW theorem}\label{sspbw}
For $1 \leq i < j \leq n$ and $r > s_{i,j}$ resp. $r > s_{j,i}$, 
we inductively define 
the {\em higher root elements} $E_{i,j}^{(r)}$ resp. $F_{i,j}^{(r)}$
of $Y_n(\sigma)$ from the formulae
\begin{align}\label{eij}
E_{i,i+1}^{(r)} := E_i^{(r)}, \qquad &
E_{i,j}^{(r)} := [E_{i,j-1}^{(r-s_{j-1,j})}, E_{j-1}^{(s_{j-1,j}+1)}],\\
F_{i,i+1}^{(r)} := F_i^{(r)}, \qquad &
F_{i,j}^{(r)} := [F_{j-1}^{(s_{j,j-1}+1)},F_{i,j-1}^{(r-s_{j,j-1})}].
\label{fij}\end{align}
Introduce the {\em canonical filtration}\label{filt}
$\F_0 Y_n(\sigma) \subseteq \F_1 Y_n(\sigma) \subseteq \cdots$
of $Y_n(\sigma)$ by declaring that
all $D_i^{(r)}, E_{i,j}^{(r)}$ and $F_{i,j}^{(r)}$ are of degree $r$, i.e.
$\F_d Y_n(\sigma)$ 
is the span of all monomials in these elements of total degree
$\leq d$.
Then \cite[Theorem 5.2]{BK} shows that 
the associated graded algebra $\gr Y_n(\sigma)$ is 
free commutative on generators
\begin{align}\label{gr1}
\{\gr_r D_i^{(r)}\:&|\:1 \leq i \leq n, s_{i,i} < r\},\\
\{\gr_r E_{i,j}^{(r)}\:&|\:1 \leq i < j \leq n, s_{i,j} < r 
\},\\
\{\gr_r F_{i,j}^{(r)}\:&|\:1 \leq i< j \leq n, s_{j,i} < r\}.
\label{gr3}
\end{align}
It follows immediately that the monomials in the elements
\begin{align}
\{D_i^{(r)}\:&|\:1 \leq i \leq n, s_{i,i} < r\},\label{bb1}\\
\{E_{i,j}^{(r)}\:&|\:1 \leq i < j \leq n, s_{i,j} < r \},\label{bb2}\\
\{F_{i,j}^{(r)}\:&|\:1 \leq i< j \leq n, s_{j,i} < r \}\label{bb3}
\end{align}
taken in some fixed order give a basis for the algebra $Y_n(\sigma)$.
Moreover, letting $Y_{(1^n)}$ resp.
$Y_{(1^n)}^+(\sigma)$ resp. 
$Y_{(1^n)}^-(\sigma)$ denote the subalgebra of $Y_n(\sigma)$
generated by the $D_i^{(r)}$'s resp. the
$E_i^{(r)}$'s resp. the $F_i^{(r)}$'s, the monomials just in the elements
(\ref{bb1}) resp. (\ref{bb2}) resp. (\ref{bb3}) taken in some fixed order
give bases for these subalgebras;
see \cite[Theorem 2.3]{BK}.
These basis theorems imply in particular that 
multiplication defines a vector space isomorphism
\begin{equation}\label{tridec}
Y_{(1^n)}^-(\sigma) \otimes Y_{(1^n)} \otimes Y_{(1^n)}^+(\sigma)
\stackrel{\sim}{\longrightarrow} Y_n(\sigma),
\end{equation}
giving us a {\em triangular decomposition} of the shifted Yangian.\label{triangular}
Also define the {\em positive} and {\em negative Borel subalgebras}
\begin{equation}
Y_{(1^n)}^\sharp(\sigma) 
:= Y_{(1^n)} Y_{(1^n)}^+(\sigma)
\text{,}
\qquad
Y_{(1^n)}^\flat(\sigma) 
:= Y_{(1^n)}^-(\sigma) Y_{(1^n)}.
\end{equation}
By the relations, these are indeed subalgebras of 
$Y_n(\sigma)$. Moreover,
there are obvious surjective homomorphisms
\begin{equation}\label{epis}
Y_{(1^n)}^\sharp(\sigma) \twoheadrightarrow 
Y_{(1^n)},
\qquad
Y_{(1^n)}^\flat(\sigma) \twoheadrightarrow Y_{(1^n)}
\end{equation} 
with kernels\label{kernels}
$K_{(1^n)}^\sharp(\sigma)$ and $K_{(1^n)}^\flat(\sigma)$
generated by all $E_{i,j}^{(r)}$ and all
$F_{i,j}^{(r)}$, respectively.

We now introduce a new basis for $Y_n(\sigma)$, which will play a central role in this article. First,
define the power series\label{defu}
\begin{align}
E_{i,j}(u) := \sum_{r > s_{i,j}} E_{i,j}^{(r)} u^{-r},\label{eiju}\qquad
&F_{i,j}(u) := \sum_{r > s_{j,i}} F_{i,j}^{(r)} u^{-r}
\end{align}
for $1 \leq i < j \leq n$, and set $E_{i,i}(u) = F_{i,i}(u) := 1$ by convention.
Recalling (\ref{didef}), let
$D(u)$ denote the $n \times n$ diagonal matrix
with $ii$-entry $D_i(u)$ for $1 \leq i \leq n$, let
$E(u)$ denote the $n \times n$
upper triangular matrix with
$ij$-entry $E_{i,j}(u)$ for $1 \leq i \leq j \leq n$, and let
$F(u)$
denote the $n \times n$
lower triangular matrix with
$ji$-entry $F_{i,j}(u)$ for $1 \leq i \leq j \leq 
n$. Consider the product
\begin{equation}\label{fde}
T(u) = F(u)D(u)E(u)
\end{equation}
of matrices with entries in $Y_n(\sigma)[[u^{-1}]]$.
The $ij$-entry of 
the matrix $T(u)$ defines a power series\label{tijind}
\begin{equation}\label{tij}
T_{i,j}(u) = \sum_{r \geq 0} T_{i,j}^{(r)} u^{-r}
:= \sum_{k=1}^{\min(i,j)} F_{k,i}(u) D_k(u) E_{k,j}(u)
\end{equation}
for some new elements $T_{i,j}^{(r)} \in \F_r Y_n(\sigma)$.
Note that $T_{i,j}^{(0)} = \delta_{i,j}$
and $T_{i,j}^{(r)} = 0$ for $0 < r \leq 
s_{i,j}$.

\begin{Lemma}\label{tbb}
The associated graded algebra $\gr Y_n(\sigma)$ is free commutative
on generators
$\{\gr_r T_{i,j}^{(r)}\:|\:1 \leq i,j \leq n, 
s_{i,j} < r\}$. Hence, the monomials
in the elements $\{T_{i,j}^{(r)}\:|\:1 \leq i,j \leq n, 
s_{i,j} < r\}$ taken in some 
fixed order form a basis for $Y_n(\sigma)$.
\end{Lemma}

\begin{proof}
Recall 
that $T_{i,j}^{(r)} = 0$ for $0 < r \leq s_{i,j}$.
Given this, it is easy to see, e.g. 
by solving the equation (\ref{fde}) in terms of quasi-determinants
as in \cite[(5.2)--(5.4)]{BKdrinfeld}, that
each of the elements $D_i^{(r)}, 
E_{i,j}^{(r)}$ and $F_{i,j}^{(r)}$ of $Y_n(\sigma)$
can be written as a 
linear combination of monomials of total 
degree $r$ 
in the elements $$
\{T_{i,j}^{(s)}\:|\:1 \leq i,j \leq n, 
s_{i,j} < s\}.
$$
Since we already know that $\gr Y_n(\sigma)$
is free commutative on the generators (\ref{gr1})--(\ref{gr3}), it follows
that the elements
$\{\gr_r T_{i,j}^{(r)}\:|\:1 \leq i, j \leq n,
s_{i,j} < r\}$ also generate $\gr Y_n(\sigma)$.
Now the lemma follows by dimension considerations.
\end{proof}

\section{Some automorphisms}\label{ssautomorphisms}
Let $\dot\sigma = (\dot s_{i,j})_{1 \leq i,j \leq n}$ 
be another shift matrix with $\dot s_{i,i+1} + \dot s_{i+1,i}
= s_{i,i+1}+s_{i+1,i}$ for all $i=1,\dots,n-1$.
Then the defining relations imply that there is a unique algebra isomorphism
\begin{equation}\label{iotadef}
\iota:Y_n(\sigma) \rightarrow Y_n(\dot\sigma)
\end{equation} 
defined on generators
by the equations
\begin{align}
\iota(D_i^{(r)})&= D_i^{(r)},\\
\iota(E_i^{(r)})&= (-1)^{s_{i,i+1}-\dot s_{i,i+1}}
E_i^{(r-s_{i,i+1}+\dot s_{i,i+1})},\\
\iota(F_i^{(r)}) 
&= (-1)^{s_{i+1,i}-\dot s_{i+1,i}}
F_i^{(r-s_{i+1,i}+
\dot s_{i+1,i})}.
\end{align}
This is not quite the same as the definition 
in \cite{BK}
(because of the extra signs), but the change causes no difficulties.

Another useful map is the anti-isomorphism
\begin{equation}\label{taudef}
\tau:Y_n(\sigma) \rightarrow Y_n(\sigma^t)
\end{equation}
where $\sigma^t$ denotes the {\em transpose} of the shift matrix $\sigma$,
defined on the generators by 
\begin{equation}
\tau(D_i^{(r)})= D_i^{(r)},
\qquad\tau(E_i^{(r)})= F_i^{(r)},
\qquad \tau(F_i^{(r)}) 
= E_i^{(r)}.
\end{equation}
Note that 
\begin{equation}\label{taudef2}
\tau(E_{i,j}^{(r)}) = 
F_{i,j}^{(r)},
\qquad
\tau(F_{i,j}^{(r)}) = 
E_{i,j}^{(r)},
\qquad
\tau(T_{i,j}^{(r)}) = T_{j,i}^{(r)}
\end{equation}
by (\ref{eij})--(\ref{fij}) and (\ref{tij}).

Finally for any power series
$f(u) \in 1 + u^{-1}\C[[u^{-1}]]$, it is easy to check from the relations that there is an automorphism
\begin{equation}\label{mufdef}
\mu_f:Y_n(\sigma)\rightarrow Y_n(\sigma)
\end{equation}
fixing each $E_i^{(r)}$ and $F_i^{(r)}$
and mapping
$D_i(u)$ to the product $f(u) D_i(u)$, i.e.
\begin{align}
\mu_f(D_i^{(r)})  &= \sum_{s=0}^r a_s D_i^{(r-s)}
\end{align}
if $f(u) = \sum_{s \geq 0} a_s u^{-s}$.

\section{Parabolic generators}\label{ssparabolic}
In this section, we recall some more complicated
{\em parabolic presentations} of $Y_n(\sigma)$ from \cite{BK}.
Actually the parabolic generators defined here will be needed
later on only in $\S$\ref{svan}.
By a {\em shape} we mean a tuple
$\nu = (\nu_1,\dots,\nu_m)$ of positive integers summing to $n$,
which we think of as the shape of 
the standard Levi subalgebra
$\mathfrak{gl}_{\nu_1}\oplus\cdots\oplus \mathfrak{gl}_{\nu_m}$
of $\mathfrak{gl}_n$.
We say that a shape 
$\nu = (\nu_1,\dots,\nu_m)$ is 
{\em admissible} (for $\sigma$) if 
$s_{i,j} = 0$ for all $\nu_1+\cdots+\nu_{a-1}+1 \leq i,j \leq \nu_1+\cdots+\nu_a$
and $a=1,\dots,m$, in which case we define
\begin{align}\label{spirit}
s_{a,b}(\nu) &:= s_{\nu_1+\cdots+\nu_a,\nu_1+\cdots+\nu_b}
\end{align}
for $1 \leq a,b \leq m$.
An important role is played by the
{\em minimal admissible shape} (for $\sigma$), namely, the admissible
shape whose length $m$ is as small as possible.

Suppose that we are given an admissible shape
$\nu = (\nu_1,\dots,\nu_m)$. 
Writing $e_{i,j}$ for the $ij$-matrix unit
in the space $M_{r,s}$ of $r \times s$ matrices over $\C$, define
\begin{multline}\label{matt}
{^\nu}T_{a,b}(u) :=\\ \sum_{\substack{1 \leq i \leq \nu_a\\ 1 \leq j \leq \nu_b}}
e_{i,j} \otimes T_{\nu_1+\cdots+\nu_{a-1}+i,\nu_1+\cdots+\nu_{b-1}+j}(u)
\in M_{\nu_a,\nu_b} \otimes Y_n(\sigma)[[u^{-1}]]
\end{multline}
for each $1 \leq a,b \leq m$.
Let
${^\nu}T(u)$ denote 
the $m \times m$ matrix 
with $ab$-entry ${^\nu}T_{a,b}(u)$.
Generalizing (\ref{fde}) (which is the special case
$\nu = (1^n)$ of the present definition), 
consider the Gauss factorization
\begin{equation}\label{gf}
{^\nu}T(u) = {^\nu}F(u) {^\nu}D(u) {^\nu}E(u)
\end{equation}
where ${^\nu}D(u)$ is an $m \times m$ diagonal matrix
with $aa$-entry denoted
${^\nu}D_a(u) \in M_{\nu_a,\nu_a}\otimes Y_n(\sigma)[[u^{-1}]]$,
${^\nu}E(u)$ is an $m \times m$ upper unitriangular matrix
with $ab$-entry denoted ${^\nu}E_{a,b}(u)
 \in M_{\nu_a,\nu_b}\otimes Y_n(\sigma)[[u^{-1}]]$ and
${^\nu}F(u)$ is an $m \times m$ lower unitriangular matrix
with $ba$-entry denoted ${^\nu}F_{a,b}(u)
 \in M_{\nu_b,\nu_a}\otimes Y_n(\sigma)[[u^{-1}]]$.
So, 
${^\nu}E_{a,a}(u)$ and 
${^\nu}F_{a,a}(u)$ 
are both the 
identity and
\begin{equation}\label{bowwow}
{^\nu}T_{a,b}(u) = \sum_{c = 1}^{\min(a,b)}
{^\nu}F_{c,a}(u) {^\nu}D_c(u) 
{^\nu}E_{c,b}(u).
\end{equation}
Also for $1 \leq a \leq m$ let 
\begin{equation}\label{rs}
{^\nu} \widetilde D_a(u) := -{^\nu}D_a(u)^{-1},
\end{equation}
inverse computed in the algebra $M_{\nu_a,\nu_a} 
\otimes Y_n(\sigma)[[u^{-1}]]$.
We expand
\begin{align}\label{mat1}
{^\nu}D_a(u) &= 
\sum_{1 \leq i,j \leq \nu_a}
e_{i,j} \otimes {^\nu} 
D_{a;i,j}(u) = 
\sum_{\substack{1 \leq i,j \leq \nu_a \\ r \geq 0}}
e_{i,j}\otimes{^\nu} D_{a;i,j}^{(r)} u^{-r},\end{align}\begin{align}\label{mat2}
{^\nu}\widetilde D_a(u) &= 
\sum_{1 \leq i,j \leq \nu_a}
e_{i,j} \otimes {^\nu} 
\widetilde D_{a;i,j}(u) = 
\sum_{\substack{1 \leq i,j \leq \nu_a \\ r \geq 0}}
e_{i,j}\otimes{^\nu}\widetilde D_{a;i,j}^{(r)} u^{-r},
\end{align}\begin{align}
\label{eg}
\qquad\:\;{^\nu}E_{a,b}(u) &= 
\sum_{\substack{1 \leq i \leq \nu_a\\ 1 \leq j \leq \nu_b}}
e_{i,j} \otimes {^\nu} 
E_{a,b;i,j}(u) = 
\sum_{\substack{1 \leq i \leq \nu_a\\ 1 \leq j \leq \nu_b\\r > s_{a,b}(\nu)}}
e_{i,j}\otimes{^\nu}E_{a,b;i,j}^{(r)} u^{-r},\end{align}\begin{align}\label{fg}
\qquad\;{^\nu}F_{a,b}(u) &= 
\sum_{\substack{1 \leq i \leq \nu_b\\ 1 \leq j \leq \nu_a}}
e_{i,j} \otimes {^\nu} 
F_{a,b;i,j}(u) = 
\sum_{\substack{1 \leq i \leq \nu_b\\ 1 \leq j \leq \nu_a\\r > s_{b,a}(\nu)}}
e_{i,j}\otimes{^\nu}F_{a,b;i,j}^{(r)} u^{-r},
\end{align}
where ${^\nu}D_{a;i,j}(u), {^\nu}\widetilde D_{a;i,j}(u),
{^\nu}E_{a,b;i,j}(u)$ and ${^\nu}F_{a,b;i,j}(u)$
are power series belonging to $Y_n(\sigma)[[u^{-1}]]$,
and
${^\nu}D_{a;i,j}^{(r)}$,
${^\nu}\widetilde D_{a;i,j}^{(r)}$,
${^\nu}E_{a,b;i,j}^{(r)}$
and
${^\nu}F_{a,b;i,j}^{(r)}$
are elements of $Y_n(\sigma)$.
We will usually omit the superscript $\nu$,
writing simply $D_{a;i,j}^{(r)}$, 
$\widetilde D_{a;i,j}^{(r)}$,
$E_{a,b;i,j}^{(r)}$ and $F_{a,b;i,j}^{(r)}$, and 
also abbreviate 
$E_{a,a+1;i,j}^{(r)}$ by $E_{a;i,j}^{(r)}$ 
and $F_{a,a+1;i,j}^{(r)}$
by $F_{a;i,j}^{(r)}$.
Note finally that the anti-isomorphism $\tau$ from (\ref{taudef})
satisfies
\begin{equation}
\tau(D_{a;i,j}^{(r)}) =
D_{a;j,i}^{(r)},\quad
\tau(E_{a,b;i,j}^{(r)}) =
F_{a,b;j,i}^{(r)},\quad
\tau(F_{a,b;i,j}^{(r)}) =
E_{a,b;j,i}^{(r)},
\end{equation}
as follows from (\ref{bowwow}) and (\ref{taudef2}).

In \cite[$\S$3]{BK}, we proved that
$Y_n(\sigma)$ is generated by the elements
\begin{align}\label{pgen1}
\{D_{a;i,j}^{(r)}\:&|\:a=1,\dots,m, 1 \leq i,j \leq \nu_a,
r > 0\},\\
\{E_{a;i,j}^{(r)}\:&|\:a=1,\dots,m-1, 1 \leq i 
\leq \nu_a, 1 \leq j \leq \nu_{a+1},
r > s_{a,a+1}(\nu)\},\\
\{F_{a;i,j}^{(r)}\:&|\:a=1,\dots,m-1, 1 \leq i 
\leq \nu_{a+1}, 1 \leq j \leq \nu_{a},
r > s_{a+1,a}(\nu)\}\label{pgen3}
\end{align}
subject to certain relations
recorded explicitly in \cite[(3.3)--(3.14)]{BK}.
Moreover, the monomials in the elements
\begin{align}
\{D_{a;i,j}^{(r)}\:&|\:1 \leq a \leq m, 1 \leq i,j \leq \nu_a, 
s_{a,a}(\nu) < r\},\label{pb1}\\
\{E_{a,b;i,j}^{(r)}\:&|\:1 \leq a < b \leq m, 
1 \leq i \leq \nu_a, 1 \leq j \leq \nu_b, s_{a,b}(\nu) < r\},\label{pb2}\\
\{F_{a,b;i,j}^{(r)}\:&|\:1 \leq a< b \leq m, 
1 \leq i \leq \nu_b, 1 \leq j \leq \nu_a, s_{b,a}(\nu) < r\}\label{pb3}
\end{align}
taken in some fixed order form a basis for $Y_n(\sigma)$.
Actually the definition of the higher root elements
$E_{a,b;i,j}^{(r)}$ and $F_{a,b;i,j}^{(r)}$ given here is 
different from the definition given in \cite{BK}. The equivalence of
the two definitions is verified by the following lemma.

\begin{Lemma}\label{thee}
For $1 \leq a < b-1 < m$, $1 \leq i \leq \nu_a, 1 \leq j \leq \nu_b$ 
and $r > s_{a,b}(\nu)$, we have that
$$
E_{a,b;i,j}^{(r)} = [E_{a,b-1;i,k}^{(r-s_{b-1,b}(\nu))}, 
E_{b-1;k,j}^{(s_{b-1,b}(\nu)+1)}]
$$
for any $1 \leq k \leq \nu_{b-1}$.
Similarly, 
for $1 \leq a < b-1 < m$, $1 \leq i \leq \nu_b, 
1 \leq j \leq \nu_a$ 
and $r > s_{b,a}(\nu)$, we have that
$$
F_{a,b;i,j}^{(r)} = [F_{b-1;i,k}^{(s_{b,b-1}(\nu)+1)}, 
F_{a,b-1;k,j}^{(r-s_{b,b-1}(\nu))}]
$$
for any $1 \leq k \leq \nu_{b-1}$.
\end{Lemma}

\begin{proof}
We just prove the statement about the $E$'s; the
statement about the $F$'s then follows
on applying the anti-isomorphism $\tau$.
Proceed by downward induction on the length
of the admissible shape $\nu = (\nu_1,\dots,\nu_m)$.
The base case $m = n$ is the definition
(\ref{eij}), so suppose $m < n$.
Pick $1 \leq p \leq m$ and
$x,y > 0$ such that $\nu_p = x+y$, then
let $\mu = (\nu_1,\dots,\nu_{p-1},x,y,\nu_{p+1},\dots,\nu_m)$,
an admissible shape of strictly longer length.
A matrix calculation from the definitions shows 
for each $1 \leq a < b \leq m$,
$1 \leq i \leq \nu_a$ and $1 \leq j \leq \nu_b$ that
\begin{align*}
{^\nu}E_{a,b;i,j}(u) &=
\left\{
\begin{array}{ll}
{^\mu}E_{a,b;i,j}(u)&\hbox{if $b < p$};\\
{^\mu}E_{a,b;i,j}(u)&\hbox{if $b=p, j \leq x$};\\
{^\mu}E_{a,b+1;i,j-x}(u)&\hbox{if $b=p, j > x$};\\
{^\mu}E_{a,b+1;i,j}(u)&\hbox{if $a < p, b > p$};\\
{^\mu}E_{a,b+1;i,j}(u)
\\
\qquad -\:\sum_{h=1}^{y}
{^\mu} E_{a,a+1;i,h}(u) {^\mu} E_{a+1,b+1;h,j}(u)\qquad
&\hbox{if $a=p, i \leq x$};\\
{^\mu}E_{a+1,b+1;i-x,j}(u)&\hbox{if $a=p, i > 
x$};\\
{^\mu}E_{a+1,b+1;i,j}(u)&\hbox{if $a > p$}.
\end{array}\right.
\end{align*}
Now suppose that $b > a+1$. We need to prove that
$$
{^\nu}E_{a,b;i,j}(u)
= [{^\nu}
E_{a,b-1;i,k}(u), {^\nu}
E_{b-1;k,j}^{(s_{b-1,b}(\nu)+1)}]u^{-s_{b-1,b}(\nu)} 
$$
for each $1 \leq k \leq \nu_{b-1}$.
The strategy is as follows: rewrite both sides of the
identity we are trying to prove in terms of the ${^\mu}E$'s and then use the 
induction hypothesis, which asserts that
$$
{^\mu}E_{a,b;i,j}(u)
= [{^\mu}
E_{a,b-1;i,k}(u), {^\mu}
E_{b-1;k,j}^{(s_{b-1,b}(\mu)+1)}]u^{-s_{b-1,b}(\mu)} 
$$
for each $1 \leq a < b-1 \leq m,
1 \leq i \leq \mu_a, 1 \leq j \leq 
\mu_b$ and 
$1 \leq k \leq \mu_{b-1}$.
Most of the cases follow at once on doing this; we just discuss the
more difficult ones in detail below.

\noindent
{\em Case one: $b < p$.} Easy.

\noindent
{\em Case two: $b = p, j \leq x$.} Easy.

\noindent
{\em Case three: $b=p, j > x$.}
We have by induction that
$$
{^\mu}E_{a,b+1;i,j-x}(u)
= [{^\mu}E_{a,b;i,h}(u),
{^\mu}E_{b;h,j-x}^{(s_{b,b+1}(\mu)+1)}]u^{-s_{b,b+1}(\mu)} 
$$
for $1 \leq h \leq x$.
Noting that $s_{b,b+1}(\mu) = 0$ and that
${^\mu}E_{b;h,j-x}^{(1)}
= {^\nu}D_{b;h,j}^{(1)}$, this shows that
${^\nu}E_{a,b;i,j}(u) = 
[{^\nu}E_{a,b;i,h}(u), {^\nu}D_{b;h,j}^{(1)}].$
Using the cases already considered and the relations, we get that
\begin{align*}
[{^\nu}E_{a,b;i,h}(u), {^\nu}D_{b;h,j}^{(1)}]
&=
[[{^\nu}
E_{a,b-1;i,k}(u), {^\nu}
E_{b-1;k,h}^{(s_{b-1,b}(\nu)+1)}], {^\nu}D_{b;h,j}^{(1)}]u^{-s_{b-1,b}(\nu)} \\
&=
[{^\nu}
E_{a,b-1;i,k}(u), {^\nu}
E_{b-1;k,j}^{(s_{b-1,b}(\nu)+1)}]u^{-s_{b-1,b}(\nu)} 
\end{align*}
for any $1 \leq k \leq \nu_{b-1}$.

\noindent
{\em Case four: $a < p, b > p$.} Easy if
$b > p+1$ or if $b = p+1$ and $k > x$.
Now suppose that $b = p+1$ and $k \leq x$.
We know already that
$$
{^\nu}E_{a,b;i,j}(u)
=  
[{^\nu}
E_{a,b-1;i,x+1}(u), {^\nu}
E_{b-1;x+1,j}^{(s_{b-1,b}(\nu)+1)}]u^{-s_{b-1,b}(\nu)}.
$$
Using the cases already considered to
express ${^\nu}E_{a,b-1;i,k}^{(r)}$ as a commutator then using 
the relation \cite[(3.11)]{BK}, 
we have that
$[{^\nu}E_{a,b-1;i,k}^{(r)}, 
{^\nu}E_{b-1;x+1,j}^{(s)}] = 0.$
Bracketing with ${^\nu}D_{b-1;k,x+1}^{(1)}$ and using the 
relations one deduces that
$$
[{^\nu}E_{a,b-1;i,x+1}^{(r)}, 
{^\nu}E_{b-1;x+1,j}^{(s)}]
=
[{^\nu}E_{a,b-1;i,k}^{(r)}, 
{^\nu}E_{b-1;k,j}^{(s)}].
$$
Hence, $$[{^\nu}
E_{a,b-1;i,x+1}(u), {^\nu}
E_{b-1;x+1,j}^{(s_{b-1,b}(\nu)+1)}]
=
[{^\nu}
E_{a,b-1;i,k}(u), {^\nu}
E_{b-1;k,j}^{(s_{b-1,b}(\nu)+1)}].$$ 
Using this we get that
${^\nu}E_{a,b;i,j}(u)
= 
[{^\nu}
E_{a,b-1;i,k}(u), {^\nu}
E_{b-1;k,j}^{(s_{b-1,b}(\nu)+1)}]u^{-s_{b-1,b}(\nu)} 
$ as required.

\noindent
{\em Case five: $a = p, i \leq x$.}
The left hand side of the identity we are trying to prove is equal to
$$
{^\mu}E_{a,b+1;i,j}(u) - \sum_{h=1}^{y}
{^\mu} E_{a,a+1;i,h}(u) {^\mu} E_{a+1,b+1;h,j}(u).
$$
The right hand side equals
$$
[{^\mu}E_{a,b;i,k}(u)- \sum_{h=1}^{y}
{^\mu}E_{a,a+1;i,h}(u) {^\mu} E_{a+1,b;h,k}(u)
,
{^\mu}E_{b;k,j}^{(s_{b,b+1}(\mu)+1)}]u^{-s_{b,b+1}(\mu)}.
$$
Now apply the induction hypothesis together with the fact from the
relations that
${^\mu}E_{a,a+1;i,h}(u)$
and 
${^\mu}E_{b;k,j}^{(s_{b,b+1}(\mu)+1)}$
commute.

\noindent
{\em Case six: $a = p, i > x$.} Easy.

\noindent
{\em Case seven: $a > p$.} Easy.
\end{proof}

We also introduce here one more family of elements of $Y_n(\sigma)$
needed in $\S$\ref{svan}.
Continue with $\nu = (\nu_1,\dots,\nu_m)$ being 
a fixed admissible shape for $\sigma$.
Recalling that
${^\nu}E_{a,a}(u)$
and
${^\nu}F_{a,a}(u)$ are both the identity,
we define 
\begin{align}\label{bare}
{^\nu}\bar E_{a,b}(u) &:= {^\nu}E_{a,b}(u) - 
\sum_{c=a}^{b-1}
{^\nu}E_{a,c}(u) {^\nu}E_{c,b}^{(s_{c,b}(\nu)+1)}
u^{-s_{c,b}(\nu)-1},\\\label{barf}
{^\nu}\bar F_{a,b}(u) &:= {^\nu}F_{a,b}(u) - 
\sum_{c=a}^{b-1}
{^\nu}F_{c,b}^{(s_{b,c}(\nu)+1)}
{^\nu}F_{a,c}(u) 
u^{-s_{b,c}(\nu)-1},
\end{align}
for $1 \leq a \leq b \leq m$.
As in (\ref{eg})--(\ref{fg}), we expand
\begin{align}\label{newgen}
{^\nu}\bar E_{a,b}(u) &= 
\sum_{\substack{1 \leq i \leq \nu_a\\ 1 \leq j \leq \nu_b}}
e_{i,j} \otimes {^\nu} 
\bar E_{a,b;i,j}(u) = 
\sum_{\substack{1 \leq i \leq \nu_a\\ 1 \leq j \leq \nu_b\\r > s_{a,b}(\nu)+1}}
e_{i,j}\otimes{^\nu}\bar E_{a,b;i,j}^{(r)} u^{-r},\\
{^\nu}\bar F_{a,b}(u) &= 
\sum_{\substack{1 \leq i \leq \nu_b\\ 1 \leq j \leq \nu_a}}
e_{i,j} \otimes {^\nu} 
\bar F_{a,b;i,j}(u) = 
\sum_{\substack{1 \leq i \leq \nu_b\\ 1 \leq j \leq \nu_a\\r > s_{b,a}(\nu)+1}}
e_{i,j}\otimes{^\nu}\bar F_{a,b;i,j}^{(r)} u^{-r},
\end{align}
where 
${^\nu}\bar E_{a,b;i,j}(u)$ and ${^\nu}\bar F_{a,b;i,j}(u)$
are power series in $Y_n(\sigma)[[u^{-1}]]$,
and
${^\nu}\bar E_{a,b;i,j}^{(r)}$
and
${^\nu}\bar F_{a,b;i,j}^{(r)}$
are elements of $Y_n(\sigma)$.
We usually 
drop the superscript $\nu$ 
from this notation.

\begin{Lemma}\label{barel}
For $1 \leq a < b-1 < m$, $1 \leq i \leq \nu_a, 1 \leq j \leq \nu_b$ 
and $r > s_{a,b}(\nu)+1$, we have that
$$
\bar E_{a,b;i,j}^{(r)} = [E_{a,b-1;i,k}^{(r-s_{b-1,b}(\nu)-1)}, 
E_{b-1;k,j}^{(s_{b-1,b}(\nu)+2)}]
$$
for any $1 \leq k \leq \nu_{b-1}$.
Similarly, 
for $1 \leq a < b-1 < m$, $1 \leq i \leq \nu_b, 
1 \leq j \leq \nu_a$ 
and $r > s_{b,a}(\nu)+1$, we have that
$$
\bar F_{a,b;i,j}^{(r)} = [F_{b-1;i,k}^{(s_{b,b-1}(\nu)+2)}, 
F_{a,b-1;k,j}^{(r-s_{b,b-1}(\nu)-1)}]
$$
for any $1 \leq k \leq \nu_{b-1}$.
\end{Lemma}

\begin{proof}
We just prove the statement about the $E$'s; the statement for the
$F$'s then follows on applying the anti-isomorphism $\tau$.
We need to prove that
$$
\bar E_{a,b;i,j}(u) =
[E_{a,b-1;i,k}(u), E_{b-1;k,j}^{(s_{b-1,b}(\nu)+2)}]
u^{-s_{b-1,b}(\nu)-1}.
$$
Proceed by induction on $b = a+2,\dots,
m$.
For the base case $b=a+2$, we have by the relation \cite[(3.9)]{BK} that
\begin{multline*}
[E_{a,b-1;i,k}(u), E_{b-1;k,j}^{(s_{b-1,b}(\nu)+2)}]
-
[E_{a,b-1;i,k}(u), E_{b-1;k,j}^{(s_{b-1,b}(\nu)+1)}]u
=\\
-[E_{a,b-1;i,k}^{(s_{a,b-1}(\nu)+1)},
E_{b-1;k,j}^{(s_{b-1,b}(\nu)+1)}]
u^{-s_{a,b-1}(\nu)}
- \sum_{h=1}^{\nu_{b-1}}
E_{a,b-1;i,h}(u) E_{b-1;h,j}^{(s_{b-1,b}(\nu)+1)}.
\end{multline*}
Multiplying by $u^{-s_{b-1,b}(\nu)-1}$ and using Lemma~\ref{thee}, 
this shows that
\begin{multline*}
[E_{a,b-1;i,k}(u), E_{b-1;k,j}^{(s_{b-1,b}(\nu)+2)}]
u^{-s_{b-1,b}(\nu)-1}
=E_{a,b;i,j}(u)
\\
-E_{a,b;i,j}^{(s_{a,b}(\nu)+1)} u^{-s_{a,b}(\nu)-1}
-
\sum_{h=1}^{\nu_{b-1}}
E_{a,b-1;i,h}(u) E_{b-1;h,j}^{(s_{b-1,b}(\nu)+1)}u^{-s_{b-1,b}(\nu)-1}.
\end{multline*}
The right hand side is exactly the definition
(\ref{bare}) of $\bar E_{a,b;i,j}(u)$
in this case.
Now assume that $b > a+2$ and
calculate using Lemma~\ref{thee}, 
relations \cite[(3.9)]{BK} and \cite[(3.11)]{BK}
and the induction hypothesis:
\begin{align*}
[E_{a,b-1;i,k}(u)&, E_{b-1;k,j}^{(s_{b-1,b}(\nu)+2)}]
u^{-1}\\
=&[[E_{a,b-2;i,1}(u), 
E_{b-2;1,k}^{(s_{b-2,b-1}(\nu)+1)}],
E_{b-1;k,j}^{(s_{b-1,b}(\nu)+2)}]u^{-s_{b-2,b-1}(\nu)-1}\\
=&[E_{a,b-2;i,1}(u), 
[E_{b-2;1,k}^{(s_{b-2,b-1}(\nu)+1)},
E_{b-1;k,j}^{(s_{b-1,b}(\nu)+2)}]]u^{-s_{b-2,b-1}(\nu)-1}\\
=&[E_{a,b-2;i,1}(u), 
[E_{b-2;1,k}^{(s_{b-2,b-1}(\nu)+2)},
E_{b-1;k,j}^{(s_{b-1,b}(\nu)+1)}]
]u^{-s_{b-2,b-1}(\nu)-1}\\
&-\sum_{h=1}^{\nu_{b-1}}
[E_{a,b-2;i,1}(u), E_{b-2;1,h}^{(s_{b-2,b-1}(\nu)+1)}
E_{b-1;h,j}^{(s_{b-1,b}(\nu)+1)}
]u^{-s_{b-2,b-1}(\nu)-1}
\end{align*}\begin{align*}
\hspace{19mm}=&[[E_{a,b-2;i,1}(u), 
E_{b-2;1,k}^{(s_{b-2,b-1}(\nu)+2)}],
E_{b-1;k,j}^{(s_{b-1,b}(\nu)+1)}]u^{-s_{b-2,b-1}(\nu)-1}\\
&-\sum_{h=1}^{\nu_{b-1}}
[E_{a,b-2;i,1}(u), E_{b-2;1,h}^{(s_{b-2,b-1}(\nu)+1)}]
E_{b-1;h,j}^{(s_{b-1,b}(\nu)+1)}u^{-s_{b-2,b-1}(\nu)-1}\\
=&[\bar E_{a,b-1;i,k}(u),
E_{b-1;k,j}^{(s_{b-1,b}(\nu)+1)}]
-\sum_{h=1}^{\nu_{b-1}}
E_{a,b-1;i,h}(u)
E_{b-1;h,j}^{(s_{b-1,b}(\nu)+1)}u^{-1}.
\end{align*}
Multiplying both sides by $u^{-s_{b-1,b}(\nu)}$ and
using the definition (\ref{bare}) together with 
Lemma~\ref{thee} once more gives the conclusion.
\end{proof}

\section{Hopf algebra structure}\label{ssgrownup}
In the special case that the shift matrix $\sigma$
is the zero matrix, we denote $Y_n(\sigma)$ simply by
$Y_n$. 
Observe that the parabolic generators 
$D_{1;i,j}^{(r)}$ of $Y_n$ defined from (\ref{gf}) relative to the
admissible shape 
$\nu = (n)$ are simply equal to the elements $T_{i,j}^{(r)}$
from (\ref{tij}).
Hence the parabolic presentation 
from \cite[(3.3)--(3.14)]{BK}
asserts in this case that the elements
$\{T_{i,j}^{(r)}\:|\:1 \leq i,j \leq n, r > 0\}$ generate $Y_n$
subject only to the relations
\begin{equation}\label{mr}
[T_{i,j}^{(r)}, T_{h,k}^{(s)}]
= \sum_{t=0}^{\min(r,s)-1}
\left( 
T_{i,k}^{(r+s-1-t)}T_{h,j}^{(t)} -
T_{i,k}^{(t)}T_{h,j}^{(r+s-1-t)}
\right)
\end{equation}
for every $1 \leq h,i,j,k \leq n$ and $r,s > 0$,
where $T_{i,j}^{(0)} = \delta_{i,j}$.
This is precisely the {RTT presentation} for
the {\em Yangian} associated to the Lie algebra $\mathfrak{gl}_n$
originating in the work of
Faddeev, Reshetikhin and Takhtadzhyan \cite{FRT};
see also \cite{D3} and \cite[$\S$1]{MNO}.
It is well known that the Yangian $Y_n$ is actually a Hopf algebra
with comultiplication
$\Delta:Y_n \rightarrow Y_n \otimes Y_n$ 
and counit $\eps:Y_n \rightarrow \C$ defined 
in terms of the generating function (\ref{tij}) by
\begin{align}\label{com2}
\Delta(T_{i,j}(u)) &= \sum_{k=1}^n T_{i,k}(u) \otimes T_{k,j}(u),\\
\eps(T_{i,j}(u)) &= \delta_{i,j}.
\end{align}
Note also that the algebra anti-automorphism $\tau:Y_n \rightarrow Y_n$
from (\ref{taudef2}) 
is a coalgebra anti-automorphism, i.e. we have
that
\begin{equation}\label{opp}
\Delta \circ \tau = P \circ (\tau \otimes \tau) \circ \Delta
\end{equation}
where $P$ denotes the permutation operator $x \otimes y \mapsto
y \otimes x$.

It is usually difficult to compute the comultiplication
$\Delta:Y_n \rightarrow Y_n \otimes Y_n$
in terms of the generators $D_i^{(r)}, E_i^{(r)}$
and $F_i^{(r)}$.
At least the case $n=2$ can be worked out explicitly
like in \cite[Definition 2.24]{Molev}: we have that
\begin{align}
\Delta(D_1(u)) &= D_1(u) \otimes D_1(u) + D_1(u) E_1(u)
\otimes F_1(u) D_1(u),\label{dmult1}\\
\Delta(D_2(u)) &= D_2(u) \otimes D_2(u) + \sum_{k \geq 1}
(-1)^k D_2(u) E_1(u)^k \otimes F_1(u)^k D_2(u),\label{dmult2}\\
\Delta(E_1(u)) &= 1 \otimes E_1(u) + \sum_{k \geq 1} (-1)^k
E_1(u)^k \otimes \widetilde D_1(u) F_1(u)^{k-1} D_2(u),\label{emult}\\
\Delta(F_1(u)) &= F_1(u) \otimes 1 + \sum_{k \geq 1} (-1)^k
D_2(u)E_1(u)^{k-1}\widetilde D_1(u) \otimes F_1(u)^k,\label{fmult}
\end{align}
as can be checked directly from (\ref{com2}) and (\ref{fde}).
The next lemma gives some further information about
$\Delta$ for $n > 2$; cf. \cite[Lemma 2.1]{CPmin}.
To formulate the lemma precisely, recall from (\ref{qgrading})
how $Y_n$ is viewed as a $Q_n$-graded algebra; the elements $T_{i,j}^{(r)}$
are of degree $(\eps_i-\eps_j)$ for this grading.
For any $s \geq 0$ and $m \geq 1$ with $m+s \leq n$ 
there is an algebra embedding \begin{equation}\label{psisdef}
\psi_s:Y_m \hookrightarrow Y_n,
\quad
D_i^{(r)} \mapsto D_{i+s}^{(r)}, 
E_i^{(r)} \mapsto E_{i+s}^{(r)}, 
F_i^{(r)} \mapsto F_{i+s}^{(r)}.
\end{equation}
A different description of this map in terms of the
generators $T_{i,j}^{(r)}$ of $Y_n$ is given in 
\cite[(4.2)]{BKdrinfeld}.
The map $\psi_s$ is {\em not} a Hopf algebra embedding:
 the maps
$\Delta \circ \psi_s$
and $(\psi_s \otimes \psi_s) \circ \Delta$ from $Y_m$
to $Y_n \otimes Y_n$ are definitely different
if $m < n$.

\begin{Lemma}\label{bettery}
For any $x \in Y_m$ such that $\psi_s(x) \in (Y_n)_{\alpha}$ 
for some $\alpha \in Q_n$, we have that
$$
\Delta (\psi_s(x))-
(\psi_s \otimes \psi_s) (\Delta(x))
\in
\sum_{0 \neq \beta \in Q_n^+} 
(Y_n)_{\beta} \otimes (Y_n)_{\alpha-\beta}
$$
where $Q_n^+$ here denotes the set of all
elements $\sum_{i=1}^{n-1} c_i (\eps_i-\eps_{i+1})$ of the root lattice $Q_n$
such that $c_i \geq 0$ for all $i \in \{1,\dots,s\}\cup\{m+s,\dots,n-1\}$.
\end{Lemma}

\begin{proof}
It suffices to prove the lemma in the two special cases
$s = 0$ and $m+s=n$.
Consider first the case that $s=0$.
Then $\psi_s:Y_m \hookrightarrow Y_n$ is just the 
map sending $T_{i,j}^{(r)} \in Y_m$ to $T_{i,j}^{(r)} \in Y_n$
for $1 \leq i,j \leq m$ and $r > 0$.
For these elements the statement of the lemma is clear
from the explicit formula for $\Delta$ from (\ref{com2}).
It follows in general
since $Y_m$ is generated by these elements
and $Q_n^+$ is closed under addition.

Instead suppose that $m+s=n$. 
Let $\widetilde T_{i,j}^{(r)} := - S(T_{i,j}^{(r)})$ where
$S$ is the antipode. Then by \cite[(4.2)]{BKdrinfeld}, 
$\psi_s:Y_m \hookrightarrow Y_n$
is the map sending $\widetilde T_{i,j}^{(r)} \in Y_m$ to 
$\widetilde T_{i+s,j+s}^{(r)} \in Y_n$
for $1 \leq i,j \leq m, r > 0$.
Since
(\ref{com2}) implies that
$$
\Delta(\widetilde T_{i,j}^{(r)})  =
-
\sum_{k=1}^n\sum_{t=0}^r 
\widetilde T_{k,j}^{(t)} \otimes \widetilde T_{i,k}^{(r-t)},
$$
the proof can now be completed as in the previous paragraph.
\end{proof}

Now we can formulate a very useful result describing 
the effect of $\Delta$ on the generators of $Y_n$ in general.
Recall from (\ref{epis}) that
$K_{(1^n)}^\sharp(\sigma)$ resp.
$K_{(1^n)}^\flat(\sigma)$ denotes the two-sided ideal of
the Borel subalgebra
$Y_{(1^n)}^\sharp(\sigma)$ resp. 
$Y_{(1^n)}^\flat(\sigma)$ generated by the
$E_i^{(r)}$ resp. the $F_i^{(r)}$;
in the case $\sigma$ is the zero matrix, we denote these
simply by
$K_{(1^n)}^\sharp$ and $K_{(1^n)}^\flat$.
Also define
\begin{equation}\label{hidef}
H_i(u) = \sum_{r \geq 0} H_i^{(r)} u^{-r} :=
\widetilde D_i(u) D_{i+1}(u)
\end{equation}
for each $i=1,\dots,n-1$. Since
$\widetilde D_i(u) = - D_i(u)^{-1}$, we have that $H_i^{(0)} = -1$.

\begin{Theorem}\label{utah}
The comultiplication $\Delta:Y_n\rightarrow Y_n\otimes Y_n$
has the following properties:
\begin{enumerate}
\item[\rm(i)]
$\Delta(D_i^{(r)}) \equiv \sum_{s=0}^r D_i^{(s)} \otimes D_i^{(r-s)}
\pmod{K_{(1^n)}^\sharp \otimes K_{(1^n)}^\flat}$;
\item[\rm(ii)]
$\Delta(E_i^{(r)}) \equiv 1 \otimes E_i^{(r)} -
\sum_{s=1}^r E_i^{(s)} \otimes H_i^{(r-s)}
\pmod{(K_{(1^n)}^\sharp)^2 \otimes K_{(1^n)}^\flat}$;
\item[\rm(iii)]
$\Delta(F_i^{(r)}) \equiv F_i^{(r)} \otimes 1 -
\sum_{s=1}^r H_i^{(r-s)} \otimes F_i^{(s)}
\pmod{K_{(1^n)}^\sharp \otimes (K_{(1^n)}^\flat)^2}$.
\end{enumerate}
\end{Theorem}

\begin{proof}
This follows from Lemma~\ref{bettery}, (\ref{dmult1})--(\ref{fmult}) 
and \cite[Corollary 11.11]{BK}.
\end{proof}

Returning to the general case, there is for any shift matrix
$\sigma = (s_{i,j})_{1 \leq i,j \leq n}$ a canonical embedding
$Y_n(\sigma) \hookrightarrow Y_n$ such that the 
generators $D_i^{(r)}, E_i^{(r)}$ and $F_i^{(r)}$
of 
$Y_n(\sigma)$ from
(\ref{gen1})--(\ref{gen3}) 
map to the elements of $Y_n$ with the same name.
However, the higher root elements
$E_{i,j}^{(r)}$ and $F_{i,j}^{(r)}$ of $Y_n(\sigma)$
do not in general map
to the elements of $Y_n$ with the same name under this embedding, 
and the elements $T_{i,j}^{(r)}$ of $Y_n(\sigma)$
do not in general map to the elements $T_{i,j}^{(r)}$ 
of $Y_n$. In particular, if $\sigma \neq 0$ 
we do {\em not} know 
a full set of relations for 
the generators $T_{i,j}^{(r)}$ of $Y_n(\sigma)$.

Write $\sigma = \sigma' + \sigma''$ where $\sigma'$
is strictly 
lower triangular and 
$\sigma''$ is strictly upper triangular.
Embedding the shifted Yangians
$Y_n(\sigma)$, $Y_n(\sigma')$ and $Y_n(\sigma'')$ into $Y_n$
in the canonical way, the first part of \cite[Theorem 11.9]{BK} 
asserts that the comultiplication
$\Delta:Y_n \rightarrow Y_n \otimes Y_n$
restricts to a map\label{hopfref}
\begin{equation}\label{dsub}
\Delta:Y_n(\sigma) \rightarrow Y_n(\sigma') \otimes Y_n(\sigma'').
\end{equation}
Also the restriction of the counit $\eps:Y_n \rightarrow \C$
gives us the {\em trivial representation}
\begin{equation}\label{triv}
\eps:Y_n(\sigma) \rightarrow \C
\end{equation}
of the shifted Yangian, with
$\eps(D_i(u)) = 1$ and
$\eps(E_i(u)) = \eps(F_i(u)) = 0$.

\section{The center of $Y_n(\sigma)$}\label{sscenter1}
Let us finally describe the center $Z(Y_n(\sigma))$ of $Y_n(\sigma)$.
Recalling the notation (\ref{didef}), let\label{centref}
\begin{equation}\label{centelts}
C_n(u) =
\sum_{r \geq 0} C_n^{(r)} u^{-r}
:= D_1(u)D_2(u-1)\cdots D_n(u-n+1) \in Y_n(\sigma)[[u^{-1}]].
\end{equation}
In the case of the Yangian $Y_n$ itself, there is 
a well known alternative description of the power series $C_n(u)$
in terms of quantum determinants due to Drinfeld \cite{D3}
(see also \cite[Theorem 8.6]{BKdrinfeld}).
To recall this, given
an $n \times n$ matrix
$A = (a_{i,j})_{1 \leq i,j \leq n}$ with entries in some (not necessarily commutative) ring, set\label{detref}
\begin{align}\label{rdetdef}
\rdet A := \sum_{w \in S_n} \sgn(w) a_{1,w1} a_{2,w2} \cdots a_{n,wn},\\
\cdet A := \sum_{w \in S_n} \sgn(w) a_{w1,1} a_{w2,2} \cdots a_{wn,n},\label{cdetdef}
\end{align}
where $S_n$ is the symmetric group.
Then, working in $Y_n[[u^{-1}]]$, we have that
\begin{align}\label{cnr}
C_n(u) &= 
\rdet
\left(
\begin{array}{ccccc}
T_{1,1}(u-n+1)&T_{1,2}(u-n+1)&
\cdots&T_{1,n}(u-n+1)\\
\vdots&\vdots&\ddots&\vdots\\
T_{n-1,1}(u-1)&T_{n-1,2}(u-1)&
\cdots&T_{n-1,n}(u-1)\\
T_{n,1}(u)&T_{n,2}(u)&\cdots&T_{n,n}(u)
\end{array}
\right)\\\label{cnc}
&=
\cdet
\left(
\begin{array}{ccccc}
T_{1,1}(u)&T_{1,2}(u-1)&\cdots&T_{1,n}(u-n+1)\\
\vdots&\vdots&\ddots&\vdots\\
T_{n-1,1}(u)&T_{n-1,2}(u-1)&\cdots&T_{n-1,n}(u-n+1)\\
T_{n,1}(u)&T_{n,2}(u-1)&\cdots&T_{n,n}(u-n+1)
\end{array}
\right).
\end{align}
In particular, in view of this alternative description, 
\cite[Proposition 2.19]{MNO} shows that
\begin{equation}\label{centercom}
\Delta(C_n(u)) = C_n(u) \otimes C_n(u).
\end{equation}

\begin{Theorem}
The elements $C_n^{(1)},C_n^{(2)},\dots$ are algebraically
independent and generate $Z(Y_n(\sigma))$.
\end{Theorem}

\begin{proof}
Exploiting the embedding $Y_n(\sigma)\hookrightarrow Y_n$,
it is known by \cite[Theorem 2.13]{MNO} that the elements
$C_n^{(1)},C_n^{(2)},\dots$ are algebraically independent
and generate $Z(Y_n)$ (see also
\cite[Theorem 7.2]{BKdrinfeld} 
for a slight variation on this argument). 
So they certainly belong to
$Z(Y_n(\sigma))$. The fact that $Z(Y_n(\sigma))$ is no larger
than $Z(Y_n)$ may be proved by 
passing to the associated graded algebra 
$\gr^\LL Y_n(\sigma)$ from \cite[Theorem 2.1]{BK}
and following the idea of
the proof of \cite[Theorem 2.13]{MNO}.
We omit the details since we give an alternative argument in 
Corollary \ref{centcor} below.
\end{proof}

Recall the automorphisms $\mu_f:Y_n(\sigma) \rightarrow Y_n(\sigma)$
from (\ref{mufdef}).
Define
\begin{equation}\label{synind}
SY_n(\sigma) := \{x \in Y_n(\sigma)\:|\:\mu_f(x) = x \hbox{ for all }f(u) \in 1 + u^{-1}\C[[u^{-1}]]\}.
\end{equation}
Like in \cite[Proposition 2.16]{MNO}, one 
can show 
 that multiplication defines
an algebra isomorphism
\begin{equation}\label{syn}
Z(Y_n(\sigma)) \otimes SY_n(\sigma) \stackrel{\sim}{\rightarrow} Y_n(\sigma).
\end{equation}
Recalling (\ref{hidef}), ordered monomials in the elements
$\{H_i^{(r)}\:|\:i=1,\dots,n-1,r > 0\}$,
$\{E_{i,j}^{(r)}\:|\:1 \leq i < j \leq n, r > s_{i,j}\}$
and
$\{F_{i,j}^{(r)}\:|\:1 \leq i < j \leq n, r > s_{j,i}\}$
form a basis for $SY_n(\sigma)$.

\chapter{Finite $W$-algebras}\label{swalgebras}

In this chapter we review the definition of the finite $W$-algebras associated
to nilpotent orbits in the Lie algebra $\mathfrak{gl}_N$,
then explain their connection to the shifted Yangians.
Again, much of this material is based closely on \cite{BK}, though 
there are some important new results too.
Throughout the chapter, we assume that $\pi$
is a fixed {\em pyramid} of {\em level} $l$, that is, a sequence
$\pi = (q_1,\dots,q_l)$ 
of integers such that
\begin{equation}\label{pyr1}
0 < q_1 \leq \cdots \leq 
q_k,\quad q_{k+1} \geq\cdots\geq q_l > 0
\end{equation}
for some fixed integer
$0 \leq k \leq l$.
We also choose an integer $n$
greater than or equal to the {\em height}
$\max(q_1,\dots,q_l)$ of the pyramid $\pi$.

\section{Pyramids}\label{sspyramids}
We visualize the pyramid $\pi$ by means of a diagram consisting of
$q_1$ bricks stacked in the first column,
$q_2$ bricks stacked in the second column, \dots,
$q_l$ bricks stacked in the $l$th column, where columns are numbered
$1,2,\dots,l$ from left to right.
For example, the diagram of the pyramid $\pi = (1,2,4,3,1)$ is
\begin{equation}\label{neeg}
{\begin{picture}(90, 46)%
\put(0,-15){\line(1,0){75}}
\put(0,0){\line(1,0){75}}
\put(15,15){\line(1,0){45}}
\put(30,30){\line(1,0){30}}
\put(30,45){\line(1,0){15}}
\put(0,-15){\line(0,1){15}}
\put(15,-15){\line(0,1){30}}
\put(30,-15){\line(0,1){60}}
\put(45,-15){\line(0,1){60}}
\put(60,-15){\line(0,1){45}}
\put(75,-15){\line(0,1){15}}
\put(7,-7){\makebox(0,0){1}}
\put(22,-7){\makebox(0,0){3}}
\put(37,-7){\makebox(0,0){7}}
\put(52,-7){\makebox(0,0){10}}
\put(67,-7){\makebox(0,0){11}}
\put(22,8){\makebox(0,0){2}}
\put(37,8){\makebox(0,0){6}}
\put(52,8){\makebox(0,0){9}}
\put(37,23){\makebox(0,0){5}}
\put(52,23){\makebox(0,0){8}}
\put(37,38){\makebox(0,0){4}}
\put(79,-8){\makebox(0,0){.}}
\end{picture}}
\vspace{4mm}
\end{equation}
\noindent
Also number the rows of the diagram of
$\pi$ by $1,2,\dots,n$ from top to bottom, so that the $n$th row
is the last row containing $l$ bricks, 
and let $p_i$
denote the number of bricks on the $i$th row.
This defines the tuple
$(p_1,\dots,p_n)$ of {\em row lengths}, with
\begin{equation}\label{rl}
0 \leq p_1 \leq \cdots \leq p_n = l.
\end{equation}
As in the above example,
we always number the bricks of the diagram $1,2,\dots,N$ down 
columns starting with the first column.
Let $\row(i)$ and 
$\col(i)$ denote the number of the row and
column containing the entry $\diagram{$\scriptstyle{i}$\cr}$
in the diagram.
We say that the pyramid is {\em left-justified} 
if $q_1 \geq \cdots \geq q_l$
and
{\em right-justified} if $q_1 \leq \cdots \leq q_l$.

Recalling the fixed choice of the integer $k$ from (\ref{pyr1}), we associate 
a shift matrix $\sigma = (s_{i,j})_{1 \leq i,j \leq n}$ to the pyramid $\pi$
by setting
\begin{equation}\label{pyr3}
s_{i,j} := \left\{
\begin{array}{ll}
\#\{c=1,\dots,k\:|\:
i >
n-q_c \geq j\}&\hbox{if $i \geq j$,}\\
\#\{c=k+1,\dots,l\:|\:
i \leq
n-q_c < j\}&\hbox{if $i \leq j$.}
\end{array}\right.
\end{equation}
To make sense of this formula, we just point out
that the pyramid $\pi$ 
can easily be recovered
given just this shift matrix $\sigma$ and the level $l$, since
its diagram consists of $p_i = l - s_{n,i}-s_{i,n}$ bricks on the $i$th
row indented $s_{n,i}$ columns from the left edge and
$s_{i,n}$ columns from the right edge. Finally, 
let
\begin{equation}\label{Sijdef}
S_{i,j} := s_{i,j} + p_{\min(i,j)}.
\end{equation}
We stress that almost all of the notation in this section and later on
depends implicitly
on the fixed choices of $n$ and $k$.

\section{Finite $W$-algebras}\label{sswalgebras}
Let $\mathfrak{g}$ denote the Lie algebra $\mathfrak{gl}_N$, 
equipped with the trace form
$(.,.)$.
Define a 
$\Z$-grading
$\mathfrak{g} = \bigoplus_{j \in \Z}\mathfrak{g}_j$ 
defined by declaring that the
$ij$-matrix unit $e_{i,j}$ is of degree
$(\col(j) - \col(i))$ for each $1 \leq i,j \leq N$.
Let $\mathfrak h := \mathfrak{g}_0,
\mathfrak{p} := \bigoplus_{j \geq 0} \mathfrak{g}_j$ and
$\mathfrak{m} := \bigoplus_{j < 0} \mathfrak{g}_j$.
Thus $\mathfrak{p}$ is a standard parabolic subalgebra of $\mathfrak{g}$
with Levi factor $\mathfrak{h} = \mathfrak{gl}_{q_1} \oplus\cdots\oplus
\mathfrak{gl}_{q_l}$, and $\mathfrak{m}$ is the opposite nilradical.
Let $e \in \mathfrak{p}$ denote the nilpotent matrix
\begin{equation}\label{edef}
e = \sum_{i,j} e_{i,j}
\end{equation}
summing over all pairs $\:\diagram{$\scriptstyle{i}$&$\scriptstyle{j}$\cr}\:$
of adjacent entries
in the diagram; for example if $\pi$ is as in (\ref{neeg})
then $e = e_{5,8}+e_{2,6}+e_{6,9}+e_{1,3}+e_{3,7}+e_{7,10}+e_{10,11}$.
The $\Z$-grading $\mathfrak{g}= \bigoplus_{j \in \Z}\mathfrak{g}_j$ 
is then a {\em good grading for $e \in \mathfrak{g}_1$}
in the sense of \cite{KRW,EK}.

The map
$\chi:\mathfrak{m} \rightarrow \C, x \mapsto (x,e)$ is a Lie algebra
homomorphism.
Let $I_\chi$ denote the kernel of the associated
homomorphism $U(\mathfrak{m})\rightarrow \C$.
Also 
let $\eta:U(\mathfrak{p}) \rightarrow U(\mathfrak{p})$ be the
algebra automorphism defined by
\begin{equation}\label{etadef}
\eta(e_{i,j}) = e_{i,j} + 
\delta_{i,j}(n-q_{\col(j)} - q_{\col(j)+1} - \cdots - q_l)
\end{equation}
for each $e_{i,j} \in \mathfrak{p}$.
Now we define the {\em finite $W$-algebra} corresponding to the pyramid 
$\pi$ to be the subalgebra
\begin{equation}\label{origdef}
W(\pi) := \{u \in U(\mathfrak{p})\:|\:
[x,\eta(u)] \in U(\mathfrak{g}) I_\chi\hbox{ for all }x \in \mathfrak{m}\}
\end{equation}
of $U(\mathfrak{p})$. 
Note this is slightly different from the
definition used in \cite[$\S$8]{BK}: there we did not include the
shift by the automorphism $\eta$ at this point. 

The definition of $W(\pi)$ originates in work of 
Kostant \cite{K} and Lynch \cite{Ly}, and is a special
case of the construction due to Premet \cite{P} and then
Gan and Ginzburg \cite{GG} of non-commutative filtered
deformations of the coordinate algebra of the Slodowy slice associated
to the nilpotent orbit containing $e$.
To make the last statement precise, we need to introduce
the {\em Kazhdan filtration} 
$$\F_0 U(\mathfrak{p}) \subseteq \F_1 U(\mathfrak{p}) \subseteq \cdots$$
of $U(\mathfrak{p})$.
This can be defined
simply by declaring
that each matrix unit $e_{i,j} \in \mathfrak{p}$ is of filtered degree 
$(\col(j)-\col(i)+1)$,
that is, $\F_d U(\mathfrak{p})$ is the span of all the monomials
$e_{i_1,j_1} \cdots e_{i_r,j_r}$ in $U(\mathfrak{p})$ such that
$$
\col(j_1)-\col(i_1)+\cdots+\col(j_r)-\col(i_r) +r \leq d.
$$
The associated graded algebra $\gr U(\mathfrak{p})$ is obviously identified
with the symmetric algebra $S(\mathfrak{p})$, viewed as a graded algebra
via the {\em Kazhdan grading} in which each $e_{i,j}$ is of graded degree $(\col(j)-\col(i)+1)$.
We get induced a filtration
$$
\F_0 W(\pi) \subseteq \F_1 W(\pi) \subseteq \cdots
$$ of
$W(\pi)$, also called the Kazhdan filtration,
by setting $\F_d W(\pi) := W(\pi) \cap \F_d U(\mathfrak{p})$;
so $\gr W(\pi)$ is naturally a graded subalgebra of $\gr U(\mathfrak{p})
= S(\mathfrak{p})$.
Let $\mathfrak{c}_g(e)$ denote the centralizer of $e$ in $\mathfrak{g}$
and
$\mathfrak{p}^\perp$ denote the nilradical of $\mathfrak{p}$.
Also define elements $h \in \mathfrak{g}_0$ and $f \in \mathfrak{g}_{-1}$
so that $(e,h,f)$ is an $\mathfrak{sl}_2$-triple in $\mathfrak{g}$
(taking $h=f=0$ in the degenerate case $e=0$).
By \cite[Lemma 8.1(ii)]{BK}, we have that
\begin{equation}
\mathfrak{p} = \mathfrak{c}_{\mathfrak{g}}(e) \oplus [\mathfrak{p}^\perp,f].
\end{equation}
The projection $\mathfrak{p} \twoheadrightarrow
\mathfrak{c}_{\mathfrak{g}}(e)$ along this direct sum decomposition
induces a homomorphism
$S(\mathfrak{p}) \twoheadrightarrow S(\mathfrak{c}_{\mathfrak{g}}(e))$.
Now the precise statement
is that the restriction of this homomorphism to $\gr W(\pi)$ 
is an  isomorphism
$\gr W(\pi) \stackrel{\sim}{\rightarrow} S(\mathfrak{c}_{\mathfrak{g}}(e))$ 
of graded algebras; see \cite[Theorem 2.3]{Ly}.

\section{Invariants}\label{ssinvariants}
For $1 \leq i,j \leq n$, $0 \leq x \leq n$ and $r \geq 1$
define
\begin{equation}\label{thedef}
T_{i,j;x}^{(r)}
:=
\sum_{s = 1}^r (-1)^{r-s}
\sum_{\substack{i_1,\dots,i_s\\j_1,\dots,j_s}}
(-1)^{\#\{t=1,\dots,s-1\:|\:\row(j_t) \leq x\}}
e_{i_1,j_1} \cdots e_{i_s,j_s}
\end{equation}
where the second sum is over all $1 \leq i_1,\dots,i_s,j_1,\dots,j_s \leq N$
such that
\begin{enumerate}
\item[(a)] $\col(j_1)-\col(i_1)+\cdots+\col(j_s)-\col(i_s) +s= r$;
\item[(b)] $\col(i_t) \leq \col(j_t)$ for each $t=1,\dots,s$;
\item[(c)] if $\row(j_t) > x$ then
$\col(j_t) < \col(i_{t+1})$ for each
$t=1,\dots,s-1$;
\item[(d)]
if $\row(j_t) \leq x$ then $\col(j_t) \geq \col(i_{t+1})$
for each
$t=1,\dots,s-1$;
\item[(e)] $\row(i_1)=i$, $\row(j_s) = j$;
\item[(f)]
$\row(j_t)=\row(i_{t+1})$ for each $t=1,\dots,s-1$.
\end{enumerate}
Also set
\begin{equation}
T_{i,j;x}^{(0)} := 
\left\{
\begin{array}{ll}
1 &\hbox{if $i= j > x$,}\\
-1&\hbox{if $i=j \leq x$,}\\
0&\hbox{if $i \neq j$,}
\end{array}
\right.
\end{equation}
and introduce the generating function
\begin{equation}\label{earlier}
T_{i,j;x}(u) := \sum_{r \geq 0} T_{i,j;x}^{(r)} u^{-r}
\in U(\mathfrak p) [[u^{-1}]].
\end{equation}
These remarkable elements (or rather their images under the automorphism
$\eta$) were first introduced in \cite[(9.6)]{BK}.
As we will explain in the next section, certain of the elements
(\ref{thedef}) in fact generate the finite $W$-algebra
$W(\pi)$.

Here is a quite
different description of the elements
$T_{i,j;0}^{(r)}$ in the spirit of \cite[(12.6)]{BK}.
If either the $i$th or the $j$th row of the diagram is empty 
then we have simply that $T_{i,j;0}(u) = \delta_{i,j}$.
Otherwise, let $a \in \{1,\dots,l\}$ be minimal such that
$i > n - q_a$ and let $b \in \{2,\dots,l+1\}$ be maximal such that
$j > n - q_{b-1}$. 
Using the shorthand $\pi(r,c)$ for the entry
$(q_1+\cdots+q_{c}+r-n)$
in the $r$th row and the $c$th column of the diagram of $\pi$
(which makes sense only if $r > n - q_c$),
we have that
\begin{multline}\label{hattij}
T_{i,j;0}(u) =\\
u^{-S_{i,j}} \sum_{m=1}^{S_{i,j}} (-1)^{S_{i,j}-m}\!\!\!\!
\sum_{\substack{r_0,\dots,r_m \\ c_0,\dots,c_m}}
\prod_{t=1}^m
\big(e_{\pi(r_{t-1},c_{t-1}),\pi(r_t,c_t-1)}
+\delta_{r_{t-1},r_t} \delta_{c_{t-1},c_t-1}
u
\big)
\end{multline}
where the second summation is over all
rows $r_0,\dots,r_m$ and columns $c_0,\dots,c_m$ such that
 $a = c_0 < \cdots < c_m = b$,
 $r_0 = i$ and $r_m = j$, and
$\max(n - q_{c_t-1}, n - q_{c_t}) < r_t \leq n$ for each $t=1,\dots,m-1$.
This identity is proved by multiplying out the parentheses and comparing with
(\ref{thedef}).

\section{Finite $W$-algebras are quotients of shifted Yangians}\label{ssquotients}
Now we can formulate the main theorem from \cite{BK} precisely.
First, 
\cite[Theorem 10.1]{BK} asserts that the elements
\begin{align}
\{T_{i,i;i-1}^{(r)}\:&|\:i=1,\dots,n, r > s_{i,i}\},\label{ne1}\\
\{T_{i,i+1;i}^{(r)}\:&|\:i=1,\dots,n-1,r > s_{i,i+1}\},\label{ne2}\\
\{T_{i+1,i;i}^{(r)}\:&|\:i=1,\dots,n-1, r > s_{i+1,1}\}\label{ne3}
\end{align}
of $U(\mathfrak{p})$ from (\ref{thedef}) generate the subalgebra
$W(\pi)$.
Moreover, 
there is a unique surjective homomorphism 
\begin{equation}\label{thetadef}
\kappa:Y_n(\sigma) \twoheadrightarrow W(\pi)
\end{equation}
under which the generators (\ref{gen1})--(\ref{gen3}) of
$Y_n(\sigma)$ map to the corresponding generators (\ref{ne1})--(\ref{ne3})
of $W(\pi)$, i.e.
$$
\kappa(D_i^{(r)}) = T_{i,i;i-1}^{(r)},
\quad\kappa(E_i^{(r)}) = T_{i,i+1;i}^{(r)},
\quad\kappa(F_i^{(r)}) = T_{i+1,i;i}^{(r)}.
$$
The kernel of $\kappa$ is the two-sided ideal 
of $Y_n(\sigma)$ generated by 
$\{D_1^{(r)}\:|\:r > p_1\}$.
Finally, viewing $Y_n(\sigma)$ as a filtered algebra via the canonical
filtration and $W(\pi)$ as a filtered algebra via the Kazhdan filtration,
we have that
$\kappa(\F_d Y_n(\sigma)) =  \F_d W(\pi)$.

From now onwards we will abuse notation by using exactly the same notation
for the elements of $Y_n(\sigma)$ (or $Y_n(\sigma)[[u^{-1}]]$)
 introduced in chapter \ref{syangians}
as for their images in $W(\pi)$ (or $W(\pi)[[u^{-1}]]$)
under the map $\kappa$, relying on
context to decide which we mean. So in particular we will denote
the invariants 
(\ref{ne1})--(\ref{ne3}) from now on just by $D_i^{(r)},
E_i^{(r)}$ and $F_i^{(r)}$.
Thus, $W(\pi)$ is generated by these elements subject only to the
relations (\ref{r2})--(\ref{r13}) together with the one additional relation
\begin{align}\label{extra}
D_1^{(r)} &= 0&&\hbox{for $r > p_1$.}\\\intertext{More generally, given an admissible shape
$\nu = (\nu_1,\dots,\nu_m)$ for $\sigma$, $W(\pi)$ is generated
by the parabolic generators
(\ref{pgen1})--(\ref{pgen3}) subject only to the relations
from \cite[(3.3)--(3.14)]{BK} together with the one additional relation}
D_{1;i,j}^{(r)} &= 0&&\hbox{for 
$1 \leq i,j \leq \nu_1$ and $r > p_1$.}\label{extra2}
\end{align}
These parabolic generators
of $W(\pi)$ are also equal to certain
of the $T_{i,j;x}^{(r)}$'s; see \cite[Theorem 9.3]{BK} 
for the precise statement here.

We should also mention the special case that
the pyramid $\pi$ is an $n \times l$ rectangle, when
the nilpotent matrix $e$ consists of $n$ Jordan blocks all of the same 
size $l$ and the shift matrix $\sigma$ is the zero matrix.
In this case,
the relation (\ref{extra2}) implies 
 that $W(\pi)$ is the quotient
of the usual Yangian $Y_n$ from $\S$\ref{ssgrownup}
by the two-sided ideal generated by
$\{T_{i,j}^{(r)}\:|\:1 \leq i,j \leq n, r > l\}$. 
Hence in this case $W(\pi)$ is isomorphic to 
the {\em Yangian of level $l$} introduced by Cherednik \cite{Ch0,Ch}, as
was first noticed  by Ragoucy and Sorba \cite{RS}.

\section{More automorphisms}\label{ssmore}
Suppose that $\dot\pi$ is another pyramid with the same row lengths as $\pi$, 
and choose a shift matrix 
$\dot\sigma = (\dot s_{i,j})_{1 \leq i,j \leq n}$
corresponding to $\dot\pi$ as in $\S$\ref{sspyramids}.
Then, viewing $W(\pi)$ as a quotient of $Y_n(\sigma)$ and
$W(\dot\pi)$ as a quotient of $Y_n(\dot\sigma)$, the
automorphism $\iota$ from (\ref{iotadef}) factors
through the quotients to induce an isomorphism
\begin{equation}\label{iotadef2}
\iota:W(\pi)\rightarrow W(\dot\pi).
\end{equation}
Hence, the isomorphism type of the algebra $W(\pi)$
only actually depends on the conjugacy class of the nilpotent matrix
$e$, i.e. on the row lengths
$(p_1,\dots,p_n)$ of $\pi$, not on the pyramid $\pi$ itself.
We remark that there is a more conceptual explanation of this last statement;
see \cite{BGo}. 
Although we are not going to give any details here, this is
the reason we have 
modified the definition of $\iota$ in (\ref{iotadef}) compared to 
\cite{BK}: the modified
$\iota$ arises in an invariant way
that does not rely on the explicit generators and relations. 

In a similar fashion, the map $\tau$ from (\ref{taudef})
induces an anti-isomorphism
\begin{align}\label{indta}
\tau:&W(\pi) \rightarrow W(\pi^t),
\end{align}
where here $\pi^t$ denotes the {\em transpose pyramid}
$(q_l,\dots,q_1)$ obtained by reversing the order of the columns of $\pi$.
There is another way to define this map, as follows.
Let $w_\pi \in S_N$ denote the permutation which when applied to the entries
of the diagram $\pi$ numbered in the standard way down columns from left to
right gives the numbering down columns from right to left.
For example, if $\pi$ is as in (\ref{neeg}) then
$w_\pi = (1\,11)(2\,9\,3\,10\,4\,5\,6\,7\,8)$.
Let
\begin{align}\label{taup}
\tau:U(\mathfrak{g}) \rightarrow U(\mathfrak{g})
\end{align}
be the algebra antiautomorphism mapping
$x \in \mathfrak{g}$ to $w_\pi x^t w_\pi^{-1}$,
where $x^t$ is the usual
transpose matrix.
Letting $\mathfrak{p}'$ denote the parabolic subalgebra of $\mathfrak{g}$
associated to the pyramid $\pi^t$, 
the map $\tau$ sends $U(\mathfrak{p})$ to
$U(\mathfrak{p}')$.
Considering the form of the definition (\ref{thedef}) explicitly, one 
checks
that
$\tau(T_{i,j;x}^{(r)}) = T_{j,i;x}^{(r)}$
for all $1 \leq i,j \leq n, 0 \leq x \leq n$ and $r \geq 0$, where
the element $T_{i,j;x}^{(r)} \in U(\mathfrak{p})$ on the left hand side is defined
using $\pi$ and the element $T_{j,i;x}^{(r)} \in U(\mathfrak{p}')$
on the right hand side is defined using $\pi^t$.
Combining this with the results of $\S$\ref{ssquotients}, it follows that
$\tau$ maps the subalgebra $W(\pi)$ of $U(\mathfrak{p})$
to the subalgebra $W(\pi^t)$ of $U(\mathfrak{p}')$,
and its restriction to $W(\pi)$ 
coincides with  (\ref{indta}).

This discussion
has the following surprising consequence, for which we have been
unable to find a direct proof
(i.e. without using the explicit
generators).
Recalling (\ref{etadef}), 
let $\overline{\eta}:U(\mathfrak{p}) \rightarrow U(\mathfrak{p})$ be the
algebra automorphism defined by
\begin{equation}\label{baretadef}
\overline{\eta}(e_{i,j}) =  e_{i,j} +
\delta_{i,j}(n-q_{1} - q_2-\cdots - q_{\col(j)})
\end{equation}
for each $e_{i,j} \in \mathfrak{p}$.

\begin{Lemma}\label{miracle}
The subalgebra $W(\pi)$ of $U(\mathfrak{p})$ is equal to
$$
\{u \in U(\mathfrak{p})\:|\:
[\overline{\eta}(u), x] \in I_\chi U(\mathfrak{g})
\text{ for all }x \in \mathfrak{m}\}.
$$
\end{Lemma}

\begin{proof}
This follows by applying the antiautomorphism $\tau^{-1}$ to 
the definition (\ref{origdef}) of $W(\pi^t)$.
\end{proof}

There is one more useful automorphism of $W(\pi)$.
For a scalar $c \in \C$, let 
\begin{equation}
\eta_c:U(\mathfrak{g}) \rightarrow U(\mathfrak{g})
\end{equation}
be the algebra automorphism mapping
$e_{i,j} \mapsto e_{i,j} + \delta_{i,j} c$ for each 
$1 \leq i,j \leq N$.
It is obvious from the definitions in $\S$\ref{sswalgebras} that this leaves
the subalgebra $W(\pi)$ invariant, hence it restricts
to an algebra automorphism
\begin{equation}\label{etac}
\eta_c:W(\pi) \rightarrow W(\pi)
\end{equation}
The following lemma gives a description of
$\eta_c$ in terms of the generators of $W(\pi)$.

\begin{Lemma}\label{eta}
For any $c \in \C$, the following equations hold:
\begin{enumerate}
\item[\rm(i)]
$\eta_c(u^{p_i} D_{i}(u)) = (u+c)^{p_{i}}D_{i}(u+c)$
for $1 \leq i \leq n$;
\item[\rm(ii)]
$\eta_c(u^{s_{i,j}} E_{i,j}(u)) = (u+c)^{s_{i,j}}E_{i,j}(u+c)$
for $1 \leq i < j\leq n$;
\item[\rm(iii)]
$\eta_c(u^{s_{j,i}} F_{i,j}(u)) = (u+c)^{s_{j,i}}F_{i,j}(u+c)$
for $1 \leq i < j \leq n$;
\item[(\rm(iv)]
$\eta_c(u^{S_{i,j}} T_{i,j}(u)) = (u+c)^{S_{i,j}} T_{i,j}(u+c)$
for $1 \leq i,j \leq n$.
\end{enumerate}
\end{Lemma}

\begin{proof}
It is immediate from (\ref{hattij}) that
$$
\eta_c(u^{S_{i,j}} T_{i,j;0}(u))
= (u+c)^{S_{i,j}} T_{i,j;0}(u+c).
$$
We will deduce the lemma from this formula. 
To do so, let $\widehat T(u)$ denote the $n \times n$ matrix with
$ij$-entry $T_{i,j;0}(u)$. Consider the Gauss factorization
$\widehat T(u) = \widehat F(u) \widehat D(u) \widehat E(u)$
where 
$\widehat D(u)$ is a diagonal matrix with $ii$-entry
$\widehat D_i(u) \in U(\mathfrak{p})[[u^{-1}]]$,
$\widehat E(u)$ is an upper unitriangular matrix with $ij$-entry
$\widehat E_{i,j}(u) \in U(\mathfrak{p})[[u^{-1}]]$ and
$\widehat F(u)$ is a lower unitriangular matrix with $ji$-entry
$\widehat F_{i,j}(u) \in U(\mathfrak{p})[[u^{-1}]]$.
Thus,
$$
T_{i,j;0}(u) = \sum_{k=1}^{\min(i,j)}
\widehat F_{k,i}(u) \widehat D_k(u) \widehat E_{k,j}(u).
$$
Since $S_{i,j} = s_{i,k}+p_k+s_{k,j}$, it follows that
$$
\eta_c(T_{i,j;0}(u)) = \!\!\sum_{k=1}^{\min(i,j)}\!
(1+cu^{-1})^{s_{i,k}}
 \widehat F_{k,i}(u) 
(1+cu^{-1})^{p_k}
\widehat D_k(u) 
(1+cu^{-1})^{s_{k,j}}
\widehat E_{k,j}(u).
$$
From this equation we can read off immediately 
the Gauss factorization of the matrix
$\eta_c(\widehat T(u))$, hence 
the matrices $\eta_c(\widehat D(u)),
\eta_c(\widehat E(u))$ and
$\eta_c(\widehat F(u))$, to get that
\begin{align*}
\eta_c(u^{p_i} \widehat D_i(u)) &= (u+c)^{p_i} \widehat D_i(u+c),\\
\eta_c(u^{s_{i,j}} \widehat E_{i,j}(u)) &= (u+c)^{s_{i,j}} 
\widehat E_{i,j}(u+c),\\
\eta_c(u^{s_{j,i}} \widehat F_{i,j}(u)) &= (u+c)^{s_{j,i}} 
\widehat F_{i,j}(u+c).
\end{align*}
The first of these equations gives (i), since by \cite[Corollary 9.4]{BK}
we have that 
$\widehat D_i(u) = D_i(u)$ in $U(\mathfrak{p})[[u^{-1}]]$. 
Similarly, (ii) and (iii) for $j=i+1$
follow from the second and third equations,
looking just at the negative powers of $u$ and using
\cite[Corollary 9.4]{BK} again.
Then (ii) and (iii) for general $j$ follow using 
(\ref{eij})--(\ref{fij}).
Finally (iv) now follows from (i)--(iii) and 
the definition (\ref{tij}).
\end{proof}

\section{Miura transform}\label{ssmiura}
Recall from the definition 
that $W(\pi)$ is a subalgebra
of $U(\mathfrak{p})$, where $\mathfrak{p}$
is the parabolic subalgebra with Levi factor
$\mathfrak h =
\mathfrak{gl}_{q_1}\oplus\cdots\oplus \mathfrak{gl}_{q_l}$.
We will often identify $U(\mathfrak{h})$ with
$U(\mathfrak{gl}_{q_1})
\otimes \cdots \otimes U(\mathfrak{gl}_{q_{l}})$.
Let $\xi:U(\mathfrak{p}) \twoheadrightarrow U(\mathfrak{h})$
be the algebra homomorphism induced by the natural projection
$\mathfrak{p} \twoheadrightarrow \mathfrak{h}$.
We call the restriction
\begin{equation}\label{miura}
\xi:W(\pi) \rightarrow U(\mathfrak{h})
\end{equation}
of $\xi$ to $W(\pi)$ the {\em Miura transform}. By \cite[Theorem 11.4]{BK}
or \cite[Corollary 2.3.2]{Ly}, this restriction
is an injective algebra homomorphism,
allowing us to view $W(\pi)$ as a subalgebra of $U(\mathfrak{h})$.

Suppose that $l = l'+l''$ 
for non-negative integers $l',l''$,
and let
$\pi' := (q_1,\dots,q_{l'})$ and
$\pi'' := (q_{l'+1},\dots,q_l)$. We write 
$\pi = \pi' \otimes \pi''$
whenever a pyramid is cut into two in this way. 
Letting $\mathfrak{h}' := \mathfrak{gl}_{q_1}\oplus\cdots\oplus
\mathfrak{gl}_{q_{l'}}$ and
$\mathfrak{h}'' := 
\mathfrak{gl}_{q_{l'+1}}\oplus\cdots\oplus \mathfrak{gl}_{q_l}$,
the Miura transform allows us to 
view the algebras $W(\pi')$ and $W(\pi'')$ as subalgebras
of $U(\mathfrak{h}')$ and $U(\mathfrak{h}'')$, 
respectively. Moreover, identifying $\mathfrak{h}$ with
$\mathfrak{h}' \oplus \mathfrak{h}''$ hence
$U(\mathfrak{h})$ with $U(\mathfrak{h'})\otimes 
U(\mathfrak{h''})$, it follows from the definition
\cite[(11.2)]{BK} and injectivity of the Miura transforms
that the subalgebra $W(\pi)$ of $U(\mathfrak{h})$
is contained in the subalgebra
$W(\pi') \otimes W(\pi'')$ of $U(\mathfrak{h'})\otimes 
U(\mathfrak{h}'')$.
We denote the resulting inclusion homomorphism by
\begin{equation}\label{deltall}
\Delta_{l',l''}:W(\pi) \rightarrow W(\pi') \otimes W(\pi'').
\end{equation}
This is the {comultiplication}
from \cite[$\S$11]{BK} (modified slightly since 
we have shifted the definition of $W(\pi)$ by $\eta$). 
It is coassociative in 
an obvious sense; see \cite[Lemma 11.2]{BK}. 
The Miura transform $\xi$
for general $\pi$ may be recovered 
by iterating this comultiplication a total of $(l-1)$ 
times to split $\pi$ into its individual columns.

Let us explain the relationship between $\Delta_{l',l''}$ and
the comultiplication $\Delta$ from (\ref{dsub}).
Let $\dot\pi'$ be the right-justified pyramid
with the same row lengths as $\pi'$, and let $\dot\pi''$ be the left-justified
pyramid with the same row lengths as $\pi''$.
So $\dot\pi := \dot\pi' \otimes \dot\pi''$ is a pyramid with the same
row lengths as $\pi$. 
Read off a shift matrix $\dot\sigma = (\dot s_{i,j})_{1 \leq i,j \leq n}$ 
from the pyramid $\dot\pi$
by choosing the integer $k$ in (\ref{pyr3})
to be $l'$.
Finally define $\dot\sigma'$ 
resp. $\dot\sigma''$ to be the strictly lower
resp. upper triangular matrices with $\dot\sigma = \dot\sigma'+\dot\sigma''$.
Then
$W(\dot\pi)$ is naturally a quotient of the shifted Yangian $Y_n(\dot\sigma)$
and similarly $W(\dot\pi') \otimes W(\dot\pi'')$ is a quotient of
$Y_n(\dot\sigma') \otimes Y_n(\dot\sigma'')$.
Composing these quotient maps with the
isomorphisms
$$W(\dot\pi) \stackrel{\sim}{\rightarrow} W(\pi),
\qquad
W(\dot\pi') \otimes W(\dot\pi'') 
\stackrel{\sim}{\rightarrow} W(\pi') \otimes W(\pi'')
$$
arising from (\ref{iotadef2}), we obtain epimorphisms
$$
Y_n(\dot\sigma) \twoheadrightarrow W(\pi),\qquad
Y_n(\dot\sigma') \otimes Y_n(\dot\sigma'') \twoheadrightarrow W(\pi') \otimes 
W(\pi'').
$$
Now the second part of \cite[Theorem 11.9]{BK}
together with \cite[Remark 11.10]{BK} asserts that the following diagram 
commutes:
\begin{equation}\label{trick}
\begin{CD}
Y_n(\dot\sigma) &@>\Delta>>&Y_n(\dot\sigma')\otimes Y_n(\dot\sigma'')\\
@VVV&&@VVV\\
W(\pi) &@>\Delta_{l',l''}>>&W(\pi') \otimes W(\pi'').
\end{CD}
\end{equation}
Using this diagram, the results about $\Delta$ obtained
in $\S$\ref{ssgrownup} imply analogous statements for the maps
$\Delta_{l',l''}:W(\pi)\rightarrow W(\pi')\otimes W(\pi'')$ in general.
For example, (\ref{opp}) implies that
\begin{equation}\label{oppp}
\Delta_{l'',l'} \circ \tau
= P  \circ \tau \otimes \tau \circ \Delta_{l',l''},
\end{equation}
equality of maps from $W(\pi)$ to
$W((\pi'')^t) \otimes W((\pi')^t)$.
This can also be seen directly 
from the
alternative description of $\tau$ as the restriction of the map 
(\ref{taup}).

Note finally that the trivial $Y_n(\sigma)$-module 
from (\ref{triv})
factors through the
quotient map $\kappa$ to induce a one dimensional $W(\pi)$-module
on which all $D_i^{(r)}, E_i^{(r)}$ and $F_i^{(r)}$ act as zero. 
We call this the
{\em trivial $W(\pi)$-module}. It is clear from (\ref{thedef})
that it is simply the restriction of the trivial $U(\mathfrak{p})$-module
to the subalgebra $W(\pi)$.

\section{Vanishing of higher $T_{i,j}^{(r)}$'s}\label{svan}
We wish next to show that 
$T_{i,j}^{(r)}$ (viewed as an element of $W(\pi)$)
is zero 
whenever $r > S_{i,j}$. In order to prove this, we
derive a recursive formula
for $T_{i,j}^{(r)}$ as an element of $U(\mathfrak{p})$
which is of independent interest.

Recall the fixed choice of $k$ from (\ref{pyr1}).
Given $k \leq m \leq l$,
let $\pi_m$ denote the pyramid
$(q_1,\dots,q_m)$ of level $m$, i.e. the first $m$ columns of $\pi$.
Let $\sigma_m$ be the shift matrix for $\pi_m$
defined according to the formula (\ref{pyr3}),
using the same choice of $k$. Let $\mathfrak{g}_m$ denote the Lie algebra
$\mathfrak{gl}_{q_1+\cdots+q_m}$.
The usual embedding of
$\mathfrak{g}_m$ into the top left hand 
corner of $\mathfrak{g}$ 
induces an embedding $I_m: U(\mathfrak{g}_m) \hookrightarrow
U(\mathfrak{g})$ of universal enveloping algebras.
We need now to consider elements both of $W(\pi) \subseteq U(\mathfrak{g})$ 
and of $W(\pi_{m}) \subseteq U(\mathfrak{g}_{m})
\subseteq U(\mathfrak{g})$.
To avoid any confusion, we will always preceed the latter by the
embedding $I_m$. 
For instance, recalling the definitions from $\S$\ref{ssparabolic},
 the notation $I_{l-1}(\bar E_{a,b;i,j}^{(r)})$
 in the following lemma means the image of the element
$\bar E_{a,b;i,j}^{(r)}$ of $W(\pi_{l-1})$ under the
embedding $I_{l-1}$.
We always work relative to the minimal admissible shape 
$\nu = (\nu_1,\dots,\nu_m)$ for $\sigma$.

\begin{Lemma}\label{eorl}
Assume that $q_1 \geq q_l$ and $k \leq l-1$.
Then, for all meaningful $a,b,i,j$ and $r$, we have that
\begin{align*}
D_{a;i,j}^{(r)} &= 
\left\{
\begin{array}{ll}
I_{l-1}(D_{a;i,j}^{(r)})&\text{if $a<m$}\\
I_{l-1}(D_{m;i,j}^{(r)})
\\
\quad\quad
 + \displaystyle\sum_{h=1}^{\nu_m}
I_{l-1}(D_{m;i,h}^{(r-1)}) e_{q_1+\cdots+q_{l-1}+h,q_1+\cdots+q_{l-1}+j}\\
\quad\quad-
\left[I_{l-1} (D_{m;i,j}^{(r-1)}), e_{q_1+\cdots+q_{l-1}+j-q_l,q_1+\cdots+q_{l-1}+j}\right]&\text{if $a=m$,}
\end{array}\right.\\
E_{a,b;i,j}^{(r)}
&= 
\left\{
\begin{array}{ll}
I_{l-1} (E_{a,b;i,j}^{(r)})&\text{if $b < m$}\\
I_{l-1} (\bar E_{a,m;i,j}^{(r)})\\
\quad\quad+ 
\displaystyle\sum_{h=1}^{\nu_m}
I_{l-1}(E_{a,m;i,h}^{(r-1)}) e_{q_1+\cdots+q_{l-1}+h,q_1+\cdots+q_{l-1}+j}\\
\quad\quad-\left[I_{l-1} (E_{a,m;i,j}^{(r-1)}), e_{q_1+\cdots+q_{l-1}+j-q_l,q_1+\cdots+q_{l-1}+j}\right]&\text{if $b=m$,}
\end{array}\right.\\
F_{a,b;i,j}^{(r)} &= I_{l-1}(F_{a,b;i,j}^{(r)}).
\end{align*}
\end{Lemma}

\begin{proof}
The first equation involving $D_{a;i,j}^{(r)}$
and the second two equations in the case $b=a+1$
follow immediately from \cite[Lemma 10.4]{BK}.
The second two equations for $b > a+1$ may then be deduced in 
exactly the same way as \cite[Lemma 4.2]{BK}.
In the difficult case when $b=m$, 
one needs to use Lemma~\ref{barel} and also the observation that
\begin{equation*}
\left[I_{l-1}(E_{a,m-1;i,h}^{(r-s_{m-1,m}(\nu))}), 
e_{q_1+\cdots+q_{l-1}+j-q_l,q_1+\cdots+q_{l-1}+j}\right] = 0
\end{equation*}
for any $1 \leq h \leq \nu_m$ 
along the way. The latter fact is checked by considering the expansion of
$E_{a,m-1;i,h}^{(r-s_{m-1,m}(\nu))}$ using \cite[Theorem 9.3]{BK} and
Lemma~\ref{thee}.
\end{proof}

\begin{Lemma}\label{eorl2}
Assume that $q_1 \geq q_l$ and $k \leq l-1$.
Then, for all $1 \leq i,j \leq n$ and $r > 0$, we have that
\begin{multline*}
T_{i,j}^{(r)} = I_{l-1}(T_{i,j}^{(r)})
-
\sum_{\substack{1 \leq h \leq n-q_l \\ s_{h,j} \leq r}}
I_{l-1}(T_{i,h}^{(r-s_{h,j})})I_{l-1}( T_{h,j}^{(s_{h,j})})\\
+
\sum_{n-q_l< h \leq n}
I_{l-1}(T_{i,h}^{(r-1)})
e_{q_1+\cdots+q_l+h-n,q_1+\cdots+q_l+j-n}\\
-
\left[
I_{l-1}(T_{i,j}^{(r-1)}), e_{q_1+\cdots+q_{l-1}+j-n,
q_1+\cdots+q_l+j-n}
\right],
\end{multline*}
omitting the last three terms on the right hand side if $j \leq n - q_l$.
\end{Lemma}

\begin{proof}
Take $1 \leq a,b \leq m, 1 \leq i \leq \nu_a, 1 \leq j \leq \nu_b$
and $r > 0$.
By definition,
$$
T_{\nu_1+\cdots+\nu_{a-1}+i,\nu_1+\cdots+\nu_{b-1}+j}(u)
=
\sum_{c=1}^{\min(a,b)}
\sum_{s,t=1}^{\nu_c}
F_{c,a;i,s}(u) D_{c;s,t}(u) E_{c,b;t,j}(u).
$$
Now apply Lemma~\ref{eorl} to rewrite the terms on the right hand side
then simplify using the definition (\ref{bare}). 
\end{proof}

\begin{Theorem}\label{tvan}
The generators $T_{i,j}^{(r)}$ of $W(\pi)$ are zero
for all $1 \leq i,j \leq n$ and $r > S_{i,j}$.
\end{Theorem}

\begin{proof}
Proceed by induction on the level $l$. The base case $l=1$ is easy to
verify directly from the definitions.
For $l > 1$, we may assume by applying the antiautomorphism $\tau$ if necessary
that $q_1 \geq q_l$. Moreover we may assume that $k \leq l-1$.
Noting that $S_{i,j}-s_{h,j} = S_{i,h}$
for $i,h \leq j$, the induction hypothesis implies that 
all the terms on the right hand side of Lemma~\ref{eorl2} are zero
if $r > S_{i,j}$.
Hence $T_{i,j}^{(r)} = 0$.
\end{proof}

Finally we describe some PBW bases for the algebra $W(\pi)$.
Recalling the definition of the Kazhdan filtration
on $W(\pi)$ from $\S$\ref{sswalgebras}, \cite[Theorem 6.2]{BK} shows that
the associated graded algebra $\gr W(\pi)$ is free commutative on generators
\begin{align}\label{grb1}
\{\gr_r D_i^{(r)}\:&|\:1 \leq i \leq n, s_{i,i} < r \leq S_{i,i}\},\\
\{\gr_r E_{i,j}^{(r)}\:&|\:1 \leq i < j \leq n, s_{i,j} < r 
\leq S_{i,j}\},\label{grb2}\\
\{\gr_r F_{i,j}^{(r)}\:&|\:1 \leq i< j \leq n, s_{j,i} < r \leq 
S_{j,i}\}.\label{grb3}
\end{align}
Hence, as in \cite[Corollary 6.3]{BK}, the monomials in the elements
\begin{align}
\{D_i^{(r)}\:&|\:1 \leq i \leq n, s_{i,i} < r \leq S_{i,i}\},\label{b1}\\
\{E_{i,j}^{(r)}\:&|\:1 \leq i < j \leq n, s_{i,j} < r 
\leq S_{i,j}\},\label{b2}\\
\{F_{i,j}^{(r)}\:&|\:1 \leq i< j \leq n, s_{j,i} < r \leq 
S_{j,i}\}\label{b3}
\end{align}
taken in some fixed order give a basis for the algebra $W(\pi)$.

\begin{Lemma}\label{tb2}
The associated graded algebra $\gr W(\pi)$ is free commutative
on generators
$\{\gr_r T_{i,j}^{(r)}\:|\:1 \leq i,j \leq n, 
s_{i,j} < r \leq S_{i,j}\}$. Hence, the monomials
in the elements $\{T_{i,j}^{(r)}\:|\:1 \leq i,j \leq n, 
s_{i,j} < r \leq S_{i,j}\}$ taken in some 
fixed order form a basis for $W(\pi)$.
\end{Lemma}

\begin{proof}
Similar to the proof of Lemma~\ref{tbb}, but using
Theorem~\ref{tvan} too.
\end{proof}

\section{Harish-Chandra homomorphisms}\label{sscent}
Finally in this chapter we review the classical
description of the center $Z(U(\mathfrak{g}))$ of the universal enveloping
algebra of $\mathfrak{g} = \mathfrak{gl}_N$.
Recalling the notation (\ref{rdetdef})--(\ref{cdetdef}), 
define a monic polynomial
\begin{equation}\label{zNu}
Z_N(u) =\sum_{r=0}^N Z_N^{(r)}u^{N-r}
\in U(\mathfrak{g})[u]
\end{equation}
by setting
\begin{align}
Z_N(u) :=&
\rdet
\left(
\begin{array}{ccccc}
e_{1,1}+u-N+1&\cdots&e_{1,N-1}&e_{1,N}\\
\vdots&\ddots&\vdots&\vdots\\
e_{N-1,1}&\cdots&e_{N-1,N-1}+u-1&e_{N-1,N}\\
e_{N,1}&\cdots&e_{N,N-1}&e_{N,N}+u
\end{array}
\right)\label{znr}\\
=&\cdet\label{znc}
\left(
\begin{array}{ccccc}
e_{1,1}+u&e_{1,2}&
\cdots&e_{1,N}\\
e_{2,1}&e_{2,2}+u-1&
\cdots&e_{2,N}\\
\vdots&\vdots&\ddots&\vdots\\
e_{N,1}&e_{N,2}&\cdots&e_{N,N}+u-N+1
\end{array}
\right).
\end{align}
Then the coefficients
$Z_N^{(1)},\dots,Z_N^{(N)}$ of this polynomial
are algebraically independent and generate 
the center $Z(U(\mathfrak{g}))$.
For a proof, see
\cite[$\S$2.2]{CL} where this is deduced from the classical Capelli identity
or \cite[Remark 2.11]{MNO} where it is deduced from (\ref{cnr})--(\ref{cnc}).

So it is natural to parametrize the central
characters of $U(\mathfrak{g})$ by monic polynomials
in $\C[u]$ of degree $N$,
 the polynomial $f(u)$ corresponding to the central
character $Z(U(\mathfrak{g})) \rightarrow \C,
Z_N(u) \mapsto f(u)$.
Let $P$ denote the free abelian group
\begin{equation}\label{Pdef}
P = \bigoplus_{a \in \C} \Z \gamma_a.
\end{equation}
Given a monic $f(u) \in \C[u]$ of degree $N$,
we associate
the element 
\begin{equation}\label{pe1}
\theta = \sum_{a \in \C} c_a \gamma_a \in P
\end{equation}
whose coefficients $\{c_a\:|\:a \in \C\}$ are 
defined from the factorization
\begin{equation}\label{pe2}
f(u) = \prod_{a \in \C} (u+a)^{c_a}.
\end{equation}
This defines a bijection between the set of monic polynomials
of degree $N$ and the set of elements $\theta \in P$ whose coefficients
are non-negative integers summing to $N$.
We will from now on always use this latter set to parametrize central
characters.

Let us compute the images of
$Z_N^{(1)},\dots,Z_N^{(N)}$ under the Harish-Chandra
homomorphism. 
Let $\mathfrak{d}$ denote the standard Cartan
subalgebra of $\mathfrak{g}$ on basis
$e_{1,1},\dots,e_{N,N}$ and let $\delta_1,\dots,\delta_N$\label{del}
be the dual basis for $\mathfrak{d}^*$.
We often represent an element $\alpha \in \mathfrak{d}^*$
simply as an $N$-tuple $\alpha = (a_1,\dots,a_N)$
of elements of the field $\C$, defined from
$\alpha = \sum_{i=1}^N a_i \delta_i$.
Also let $\mathfrak{b}$ be the standard Borel
subalgebra consisting of upper triangular matrices.
We will parametrize highest weight 
modules already in ``$\rho$-shifted notation'':
for a weight $\alpha \in \mathfrak{d}^*$, let $M(\alpha)$
denote the {\em Verma module} of 
highest weight $(\alpha-\rho)$,
namely, the module
\begin{equation}\label{vermamod}
M(\alpha) := U(\mathfrak{g}) \otimes_{U(\mathfrak{b})}
\C_{\alpha-\rho}
\end{equation}
induced from the one dimensional $\mathfrak{b}$-module
of weight $(\alpha-\rho)$,
where $\rho$ here means the weight
$(0,-1,-2,\dots,1-N)$.
Thus, if $\alpha = (a_1,\dots,a_N)$,
the diagonal matrix $e_{i,i}$ acts on the highest weight space of $M(\alpha)$
by the scalar $(a_i + i-1)$.
Viewing the symmetric algebra
$S(\mathfrak{d})$ as an algebra of functions
on $\mathfrak{d}^*$, with the symmetric group $S_N$
acting by $w \cdot e_{i,i} := e_{wi,wi}$ as usual,
the Harish-Chandra homomorphism
\begin{equation}\label{hchom}
\Psi_{\!N}:
Z(U(\mathfrak{g})) \stackrel{\sim}{\rightarrow} 
S(\mathfrak{d})^{S_N}
\end{equation}
may be defined as the map sending
$z \in Z(U(\mathfrak{g}))$ to the 
unique element of $S(\mathfrak{d})$ with the property
that $z$ acts on $M(\alpha)$ by the scalar
$(\Psi_{\!N}(z))(\alpha)$ for each $\alpha \in \mathfrak{d}^*$.
Using (\ref{znc}) it is easy to see directly from this definition
that
\begin{equation}\label{hcim}
\Psi_{\!N}(Z_N(u)) = (u+e_{1,1}) (u+e_{2,2}) \cdots (u+e_{N,N}).
\end{equation}
The coefficients on the right hand side are the elementary
symmetric functions.
Define the {\em content} $\theta(\alpha)$
of the weight $\alpha = (a_1,\dots,a_N) \in \mathfrak{d}^*$ by setting
\begin{equation}\label{contdef0}
\theta(\alpha) := \gamma_{a_1}+\cdots + \gamma_{a_N} \in P.
\end{equation}
By (\ref{hcim}), the central character
of the Verma module $M(\alpha)$ 
is precisely the central character parametrized by 
$\theta(\alpha)$.

Now return to the setup of $\S$\ref{sswalgebras}.
Let 
$\psi:U(\mathfrak{g}) \rightarrow
U(\mathfrak{p})$
be the linear map defined as the composite
first of the projection
$U(\mathfrak{g}) \rightarrow U(\mathfrak{p})$
along the direct sum decomposition
$U(\mathfrak{g}) = U(\mathfrak{p}) \oplus U(\mathfrak{g}) I_\chi$
then the inverse $\eta^{-1}$ of the automorphism
$\eta$ from (\ref{etadef}).
The restriction of $\psi$
to $Z(U(\mathfrak{g}))$ gives a well-defined algebra homomorphism
\begin{equation}\label{psidef}
\psi:Z(U(\mathfrak{g})) \rightarrow
Z(W(\pi))
\end{equation} 
with image contained in the center of $W(\pi)$.
Applying this to the polynomial
$Z_N(u)$ we obtain elements
$\psi(Z_N^{(1)}),\dots,\psi(Z_N^{(N)})$ of
$Z(W(\pi))$.
The following lemma explains 
the relationship between these elements and the
elements
$C_n^{(1)}, C_n^{(2)},\dots$ of $Z(W(\pi))$ defined
by the formula (\ref{centelts}).

\begin{Lemma}\label{psieq}
$\psi(Z_N(u)) = u^{p_1} (u-1)^{p_2} \cdots (u-n+1)^{p_n}
C_n(u)$.
\end{Lemma}

\begin{proof}
A calculation using (\ref{znr})
shows that the image of $\psi(Z_N(u))$ under the Miura
transform $\xi:W(\pi) \rightarrow U(\mathfrak{gl}_{q_1})
\otimes\cdots\otimes U(\mathfrak{gl}_{q_l})$ from (\ref{miura})
is equal to $Z_{q_1}(u+q_1-n) \otimes \cdots \otimes Z_{q_l}(u+q_l-n)$.
So, since $\xi$ is injective, we have to check that
$\xi(u^{p_1} (u-1)^{p_2} \cdots (u-n+1)^{p_n} C_n(u))$ also equals
$Z_{q_1}(u+q_1-n) \otimes \cdots \otimes Z_{q_l}(u+q_l-n)$. 
By (\ref{centercom}) and (\ref{trick}), 
we have that $\xi(C_n(u)) = C_n(u) \otimes\cdots\otimes C_n(u)$
($l$ times). Therefore it just remains to observe 
in the special case that $\pi$ consists of a single column of
height $m \leq n$, i.e. when $W(\pi) = U(\mathfrak{gl}_m)$, that
$$
(u-n+m) \cdots (u-n+2) (u-n+1) C_n(u) = Z_m(u+m-n).
$$
This follows by a direct calculation from (\ref{cnc}),
exactly as in \cite[Remark 2.11]{MNO} (which is the
case $m=n$).
\end{proof}

We can also consider the Harish-Chandra homomorphism
\begin{equation}
\Psi_{\!q_1} \otimes\cdots\otimes \Psi_{\!q_l}
:Z(U(\mathfrak{h}))
\stackrel{\sim}{\rightarrow}
S(\mathfrak{d})^{S_{q_1} \times\cdots\times S_{q_l}}
\end{equation}
for $\mathfrak{h} = \mathfrak{gl}_{q_1}\oplus\cdots\oplus
\mathfrak{gl}_{q_l}$,
identifying
$U(\mathfrak{h})$ with
$U(\mathfrak{gl}_{q_1}) \otimes\cdots\otimes
U(\mathfrak{gl}_{q_l})$.
By (\ref{hcim}) and 
the explicit computation
of $\xi(\psi(Z_N(u)))$ made in the proof of Lemma~\ref{psieq},
the following diagram commutes:
\begin{equation}\label{hcfact}
\begin{CD}
Z(U(\mathfrak{g}))
& @>\sim>\Psi_{\!N}> &S(\mathfrak{d})&^{S_N\phantom{\cdots\times S_{q_l}ll}}\\
@V\xi\circ \psi VV&&@VVV\\
Z(U(\mathfrak{h}))&@>\sim>
\Psi_{\!q_1}\otimes\cdots\otimes \Psi_{\!q_l}>&S(\mathfrak{d})&^{S_{q_1}\times \cdots\times S_{q_l}}
\end{CD}
\end{equation}
where the right hand map is the inclusion
arising from the restriction of the automorphism
$S(\mathfrak{d}) \rightarrow S(\mathfrak{d}),
e_{i,i} \mapsto e_{i,i}+q_{\col(i)}-n$.
Hence the Harish-Chandra homomorphism $\Psi_{\!N}$
factors through the map $\psi$,
as has been observed in much greater generality than this by
Lynch \cite[Proposition 2.6]{Ly} and Premet \cite[6.2]{P}.
In particular this shows that $\psi$ is injective,
so the elements $\psi(Z_N^{(1)}),\dots,\psi(Z_N^{(N)})$
of $Z(W(\pi))$ are actually algebraically independent.

\chapter{Dual canonical bases}\label{sbases}

The appropriate setting for the combinatorics underlying the representation theory of the algebras $W(\pi)$ is provided by certain
dual canonical bases for representations of the Lie algebra
$\mathfrak{gl}_{\infty}$.
In this chapter we review these matters following \cite{qla} closely.
Throughout, $\pi$ denotes a fixed pyramid
$(q_1,\dots,q_l)$ with row lengths $(p_1,\dots,p_n)$, 
and $N = p_1+\cdots+p_n=q_1+\cdots+q_l$.

\section{Tableaux}\label{sstableaux}
By a {\em $\pi$-tableau}
we mean a filling of the boxes of the 
diagram of $\pi$ with arbitrary elements of the ground field $\C$.
Let $\Tab(\pi)$\label{tabpi} denote the set of all such $\pi$-tableaux.
If $\pi = \pi' \otimes \pi''$ for pyramids $\pi'$ and $\pi''$
and we are given a $\pi'$-tableau $A'$ and a $\pi''$-tableau $A''$, 
we write $A' \otimes A''$ for the $\pi$-tableau obtained by concatenating
$A'$ and $A''$. For example,
$$
A=\Diagram{1\cr 0&3&2\cr 4&3&1\cr}\,
=
\,\Diagram{1\cr 0&3\cr 4&3\cr}
\,\otimes\,
{\substack{\phantom{A^1_{3_{1}}}\\\Diagram{2\cr 1\cr}}}
\,=\,
\Diagram{1\cr 0\cr 4\cr}
\,\otimes\,
{\substack{\phantom{A^1_{3_{1}}}\\\Diagram{3\cr 3\cr}}}
\,\otimes\,
{\substack{\phantom{A^1_{3_{1}}}\\\Diagram{2\cr 1\cr}}}\,.
$$
We always number the rows of $A \in \Tab(\pi)$ by $1,\dots,n$ from top to bottom and the columns by $1,\dots,l$ from left to right,
like for the diagram of $\pi$.
We let $\gamma(A)$ denote the weight
$\alpha = (a_1,\dots,a_N) \in \C^N$ obtained from $A$ 
by {\em column reading} the entries of $A$ down columns starting with
the leftmost column.
For example, if $A$ is as above, then $\gamma(A) = (1,0,4,3,3,2,1)$.
Define the {\em content} $\theta(A)$\label{contdef1} of $A$
to be the content of the weight $\gamma(A)$
in the sense of (\ref{contdef0}),
an element of the free abelian group
$P = \bigoplus_{a \in \C} \Z \gamma_a$.

We say that two $\pi$-tableaux $A$ and $B$ are {\em row equivalent},
written $A \sim_{\ro} B$,
if one can be obtained from the other by permuting entries within rows.
The notion $\sim_{\co}$ of {\em column equivalence} is defined similarly.
Let $\Row(\pi)$ denote the set of all row equivalence classes
of $\pi$-tableaux. We refer to elements of $\Row(\pi)$ as
{\em row symmetrized $\pi$-tableaux}.\label{colpi}
Let $\Col(\pi)$ denote the set of all {\em column strict $\pi$-tableaux},
namely, the $\pi$-tableaux whose entries are strictly increasing up 
columns from bottom to top according to 
the partial order $\geq$ on $\C$ defined 
by $a \geq b$ if $(a-b) \in \N$. 
We stress the deliberate asymmetry of these definitions:
$\Col(\pi)$ is a subset of $\Tab(\pi)$ but $\Row(\pi)$ is a quotient.

Let us recall the usual definition of 
the {\em Bruhat ordering}\label{bruhord} on the set $\Row(\pi)$.
Given $\pi$-tableaux $A$ and $B$,
write $A \downarrow B$ if $B$ is obtained from $A$ by swapping
an entry $x$ in the $i$th row of $A$ with an entry $y$ in the $j$th
row of $A$, and moreover we have that $i < j$ and $x > y$.
For example,
$$
\Diagram{${2}$&${5}$\cr${7}$&${7}$\cr${3}$&${3}$&${5}$\cr}
\:
\downarrow
\:
\Diagram{${2}$&${3}$\cr${7}$&${7}$\cr${3}$&${5}$&${5}$\cr}
\:
\downarrow
\:
\Diagram{${2}$&${3}$\cr${7}$&${3}$\cr${7}$&${5}$&${5}$\cr}\:.
$$
Now for $A, B \in \Row(\pi)$, the notation
$A \geq B$ means that there exists $r \geq 1$ and
$\pi$-tableaux $A_1,\dots,A_r$ such that
\begin{equation}\label{bruhdef}
A \sim_\ro A_1 \downarrow \cdots \downarrow A_r \sim_\ro B.
\end{equation}
It is obvious that if $A \geq B$ then
$\theta(A) = \theta(B)$ (where the content $\theta(A)$
of a row-symmetrized $\pi$-tableau means the content of any representative
for $A$).

It just remains to introduce notions of
{\em dominant} and of {\em standard} $\pi$-tableaux.
The first of these is easy:
call an element $A \in \Row(\pi)$ {\em dominant} if it has a 
representative belonging to $\Col(\pi)$ and let $\Dom(\pi)$\label{dompi} denote the set 
of all such dominant row symmetrized $\pi$-tableaux.
The notion of a standard $\pi$-tableau is more subtle.
Suppose first that $\pi$ is left-justified,
when its diagram is a Young diagram in the usual sense.
In that case, a $\pi$-tableau $A \in \Col(\pi)$ with entries
$a_{i,1},\dots,a_{i,p_i}$ in its $i$th row read from left to right is
called {\em standard} if 
 $a_{i,j} \not > a_{i,k}$ for all
$1 \leq i \leq n$ and $1 \leq j < k \leq p_i$.
If $A$ has integer entries 
(rather than arbitrary elements of $\C$)
this is just saying that the entries of $A$ are 
strictly increasing up columns from bottom to top and
weakly increasing along rows from 
left to right, i.e. it is the usual notion of standard tableau.

\begin{Lemma}\label{lb}
Assume that $\pi$ is left-justified. Then any element $A \in \Dom(\pi)$
has a representative that is standard.
\end{Lemma}

\begin{proof}
By definition, we can choose a representative for $A$ that is column strict.
Let $a_{i,1},\dots,a_{i,p_i}$ be the entries on the $i$th row of this representative read from left to right, for each $i=1,\dots,n$. We need to show that we
can permute entries within rows 
so that it becomes standard. Proceed by induction on 
$$
\#\{(i,j,k)\:|\: 1 \leq i \leq n, 1 \leq j < k \leq p_i \hbox{ such that }a_{i,j} 
> a_{i,k}\}.
$$
If this number is zero then our tableau is already standard.
Otherwise we can pick
$1 \leq i \leq n$ and $1 \leq j< k \leq p_i$ such that
$a_{i,j} > a_{i,k}$, none of $a_{i,j+1},\dots,a_{i,k-1}$
lie in the same coset of $\C$ modulo $\Z$ as $a_{i,j}$, 
and either $i=n$ or 
$a_{i+1,j} \not > a_{i+1,k}$.
Then define $1 \leq h \leq i$ to be minimal 
so that $k \leq p_h$ and
$a_{r,j} > a_{r,k}$ for all $h \leq r \leq i$.
Thus our tableau contains the following entries:
$$
\begin{array}{ccc}
a_{h-1,j} & \leq  & a_{h-1,k}\\
a_{h,j} & > & a_{h,k}\\
a_{h+1,j} & > & a_{h+1,k}\\
\vdots & \vdots&\vdots \\
a_{i,j} & > & a_{i,k}\\
a_{i+1,j} & \leq &a_{i+1,k},
\end{array}
$$
where entries on the $(h-1)$th and/or $(i+1)$th 
rows should be omitted if they do not exist.
Now swap the entries
$a_{h,j} \leftrightarrow a_{h,k},
a_{h+1,j} \leftrightarrow a_{h+1,k}, \dots
a_{i,j} \leftrightarrow a_{i,k}$ and observe that 
the resulting tableau is still column strict.
Finally by the induction hypothesis we get that the new tableau is row equivalent to a standard tableau.
\end{proof}

To define what it means for $A \in \Col(\pi)$ to be standard
for more general pyramids $\pi$ we need to recall the
notion of row insertion; see e.g. \cite[$\S$1.1]{fulton}.
Suppose we are given a weight $(a_1,\dots,a_N) \in \C^N$.
We decide if it is admissible, and if so construct an element of
$\Row(\pi)$, according to the following algorithm.
Start from the diagram of $\pi$ with all boxes empty.
Insert $a_1$ into some box in the bottom ($n$th) row.
Then if $a_2\not< a_1$ insert $a_2$ into the bottom row too;
else if $a_2 < a_1$ replace the entry $a_1$ by $a_2$
and insert $a_1$ into the next row up instead.
Continue in this way: at the $i$th step the pyramid $\pi$ has $(i-1)$ boxes
filled in and we need to insert the entry $a_i$ into the bottom row.
If $a_i$ is $\not<$ all of the entries in this row, simply add
it to the row; else find the smallest entry $b$ in the row that is
strictly larger than $a_i$, replace this entry $b$ with $a_i$, then 
insert $b$ into the next row up in similar fashion.
If at any stage of this process one gets more than $p_i$ entries 
in the $i$th row for some $i$, 
the algorithm terminates and the weight $(a_1,\dots,a_N)$
is inadmissible; else, 
the weight $(a_1,\dots,a_N)$ is admissible and we have successfully computed
a tableau $A \in \Row(\pi)$.

Now, for any pyramid $\pi$, we say that
$A \in \Col(\pi)$ is {\em standard} if the weight
$\gamma(A)$ obtained from the column reading of $A$ is
admissible. 
Let $\Std(\pi)$\label{stdpi} denote the set of all such standard $\pi$-tableaux.
For $A \in \Std(\pi)$, we define the {\em rectification} 
$R(A) \in \Row(\pi)$ to be the row symmetrized $\pi$-tableau computed 
from the weight $\gamma(A)$
by the algorithm described in the previous paragraph.
In the special case that $\pi$ is left-justified, 
it is straightforward to check that 
the new definition of
standard tableau agrees with the one given before Lemma~\ref{lb}, and
moreover in this case the map $R$ is simply the map sending 
a tableau to its row equivalence class.
In general, it is clear from the algorithm that $R(A)$ belongs to $\Dom(\pi)$,
i.e. it has a representative that is column strict,
so rectification gives a map
\begin{equation}\label{rmap}
R:\Std(\pi) \rightarrow \Dom(\pi).
\end{equation}
Define an equivalence relation
$\parallel$ on $\Col(\pi)$ by declaring
that $A \parallel 
B$ if $B$ can be obtained from $A$ by shuffling columns of equal height in 
such a way that the relative position of all columns belonging to the 
same coset of $\C$ modulo $\Z$ remains the same.
Then the map $R:\Std(\pi) \rightarrow \Dom(\pi)$
is surjective, and $R(A) = R(B)$
if and only if $A \parallel B$, i.e. the fibres of $R$ are
precisely the $\parallel$-equivalence classes.
This follows in the left-justified case
using Lemma~\ref{lb}, and then in general
by a result of Lascoux and Sch\"utzenberger \cite{LS};
see \cite[$\S$A.5]{fulton} and \cite[$\S$2]{qla}.

We have now introduced all the sets
$\Tab(\pi), \Row(\pi), \Col(\pi),
\Dom(\pi)$ and
$\Std(\pi)$ of tableaux 
which will be needed later on to parametrize the various bases/modules
that we will meet.
We write
$\Tab_0(\pi),\Row_0(\pi), \Col_0(\pi), \Dom_0(\pi)$ and $\Std_0(\pi)$
for the subsets of
$\Tab(\pi),
\Row(\pi), \Col(\pi), \Dom(\pi)$ and $\Std(\pi)$ consisting just of the
tableaux all of whose entries are integers.
In fact, most of the problems that we will meet
are reduced in a straightforward fashion to
this special situation.
Finally, we define the {\em row reading}\label{roread}
$\rho(A)$
of $A \in \Row_0(\pi)$ to be the
weight $\alpha = (a_1,\dots,a_N) \in \Z^N$ obtained by
reading the entries in each row of $A$ in 
weakly increasing order, starting with the top row.
For example, if $A$ is the row equivalence class of the tableau
displayed in the first paragraph,
then $\rho(A) = (1,0,2,3,1,3,4)$.
 
\section{Dual canonical bases}\label{ssdcb}
Now let $\mathfrak{gl}_\infty$ denote the Lie algebra of
matrices with rows and columns labelled by $\Z$, all but finitely many
entries of which are zero.
It is generated by the usual Chevalley generators
$e_i, f_i$, i.e. the matrix units $e_{i,i+1}$
and $e_{i+1,i}$, together with the diagonal matrix
units $d_i = e_{i,i}$, for each $i \in \Z$.
The associated integral weight lattice $P_\infty$\label{pinf} 
is the free abelian
group with basis $\{\gamma_i\:|\:i \in \Z\}$ dual to
$\{d_i\:|\:i \in \Z\}$, and the
simple roots are $\gamma_i-\gamma_{i+1}$ for $i \in \Z$.
We will view $P_\infty$ as a subgroup of the group $P$ from (\ref{Pdef}).
Let $U_\Z$\label{uz} be the Kostant $\Z$-form for the 
universal enveloping algebra $U(\mathfrak{gl}_\infty)$,
generated by the divided powers $e_i^r / r!$, $f_i^r / r!$
and the elements $\binom{d_i}{r} = \frac{d_i(d_i-1)\cdots (d_i-r+1)}{r!}$ for all $i \in \Z, r \geq 0$.
Let $V_\Z$ be the {\em natural $U_\Z$-module}, that is, the $\Z$-submodule
of the natural $\mathfrak{gl}_\infty$-module 
generated by the standard basis vectors $v_i\:(i \in \Z)$.

Consider to start with the $U_\Z$-module arising as
the $N$th tensor power\label{ts}
$T^N(V_\Z)$ of $V_\Z$.
It is a free $\Z$-module with the monomial
basis $\{M_\alpha\:|\:\alpha \in \Z^N\}$
defined from $M_\alpha = v_{a_1}\otimes\cdots\otimes v_{a_N}$\label{malpha}
for $\alpha = (a_1,\dots,a_N) \in \Z^N$.
We also need the {\em dual canonical basis}
$\{L_\alpha\:|\:\alpha \in \Z^N\}$. The best way to define this is to first
quantize, then define $L_\alpha$ using a natural bar involution
on the $q$-tensor space, then specialize to $q=1$ at the end. 
We refer to \cite[$\S$4]{qla} for the details
of this construction (which is due to Lusztig \cite[ch.27]{Lubook}); 
the only significant difference is that 
in \cite{qla} the Lie algebra $\mathfrak{gl}_n$ is used in place of
the Lie algebra $\mathfrak{gl}_\infty$ here.
We just content ourselves with writing down an explicit formula
for the expansion of $M_\alpha$ as a linear combination of $L_\beta$'s
in terms of the usual Kazhdan-Lusztig polynomials $P_{x,y}(q)$ associated to
the symmetric group $S_N$ from \cite{KL} evaluated at $q=1$. 
To do this, let $S_N$ act on the right on the set 
$\Z^N$ by place permutation
in the natural way, and given any $\alpha \in \Z^N$ define
$d(\alpha) \in S_N$ to be the unique element of minimal length with the property
that $\alpha \cdot d(\alpha)^{-1}$ is a weakly increasing sequence.
Then, by \cite[$\S$4]{qla}, we have that
\begin{equation}\label{mpl}
M_\alpha = \sum_{\beta \in \Z^N} P_{d(\alpha) w_0, d(\beta) w_0}(1) L_\beta,
\end{equation}
writing $w_0$ for the longest element of $S_N$.

We also need to consider certain tensor products
of symmetric and exterior powers of $V_\Z$.
Let $S^N(V_\Z)$ denote the $N$th symmetric power of $V_\Z$,
defined as a quotient of $T^N(V_\Z)$ in the usual way.
Also let $\bw^N(V_\Z)$ denote the $N$th exterior power of $V_\Z$,
viewed unusually as the subspace of $T^N(V_\Z)$ consisting of
all skew-symmetric tensors. Recalling the fixed pyramid $\pi$,
let
\begin{align}\label{sdef}
S^\pi(V_\Z) := S^{p_1}(V_\Z) \otimes\cdots\otimes S^{p_n}(V_\Z),\\
\bw^\pi(V_\Z) := \bw^{q_1}(V_\Z) \otimes\cdots\otimes \bw^{q_l}(V_\Z).
\end{align}
Identifying $T^N(V_\Z) = T^{p_1}(V_\Z)\otimes\cdots\otimes
T^{p_n}(V_\Z) = T^{q_1}(V_\Z)\otimes\cdots\otimes T^{q_l}(V_\Z)$, we observe that
$S^\pi(V_\Z)$ is a quotient of
$T^N(V_\Z)$, while
$\bw^\pi(V_\Z)$ is a subspace.
Following \cite[$\S$5]{qla}, 
both of these free $\Z$-modules have two natural bases, 
a monomial basis and a dual canonical basis, parametrized
by the sets $\Row_0(\pi)$ and $\Col_0(\pi)$, respectively.

First we define these two bases for the space $S^\pi(V_\Z)$.
For $A \in \Row_0(\pi)$, define $M_A$\label{madef} to be the image of $M_{\rho(A)}$
and $L_A$ to be the image of $L_{\rho(A)}$
under the canonical quotient map
$T^N(V_\Z) \twoheadrightarrow S^\pi(V_\Z)$.
The monomial basis for $S^\pi(V_\Z)$ is then the set
$\{M_A\:|\:A \in \Row_0(\pi)\}$, and the dual canonical basis is
$\{L_A\:|\:A \in \Row_0(\pi)\}$.

Now we define the two bases for the space $\bw^\pi(V_\Z)$.
For $A \in \Col_0(\pi)$, let
\begin{equation}\label{nadef}
N_A := \sum_{B \sim_{\co} A}(-1)^{\ell(A,B)} M_{\gamma(B)},
\end{equation}
where  $\ell(A,B)$ denotes the minimal number of transpositions of adjacent
elements in the same column needed to get from $A$ to $B$. 
Also let $K_A$ denote the vector $L_{\gamma(A)} \in T^N(V_\Z)$.
Then both $N_A$ and $K_A$ belong to the subspace
$\bw^\pi(V_\Z)$ of $T^N(V_\Z)$; see \cite[$\S$5]{qla}.
Moreover, $\{N_A\:|\:A \in \Col_0(\pi)\}$ and
$\{K_A\:|\:A \in \Col_0(\pi)\}$ are bases
for $\bw^\pi(V_\Z)$, giving the monomial basis and the
dual canonical basis, respectively. 

The following formulae, derived in \cite[$\S$5]{qla}
as consequences of (\ref{mpl}), express the monomial bases
in terms of the dual canonical bases and certain 
Kazhdan-Lusztig polynomials:
\begin{align}\label{kl1}
M_A &= \sum_{B \in \Row_0(\pi)}
P_{d(\rho(A)) w_0, d(\rho(B)) w_0}(1) L_B,\\
N_A &= \sum_{B \in \Col_0(\pi)} 
\left(
\sum_{C \sim_\co A} (-1)^{\ell(A,C)} P_{d(\gamma(C)) w_0, d(\gamma(B))w_0}(1)
 \right) K_B,\label{nak}
\end{align}
for $A \in \Row_0(\pi)$ and $A \in \Col_0(\pi)$, respectively.

Note that $S^\pi(V_\Z)$ is a summand of the
commutative algebra 
$S(V_\Z)\otimes\cdots\otimes S(V_\Z)$, that is, the tensor product
of $n$ copies of the symmetric algebra $S(V_\Z)$.
In particular, if $\pi = \pi' \otimes \pi''$, the multiplication
in this algebra defines a $U_\Z$-module homomorphism
\begin{equation}\label{themapmu}
\mu:S^{\pi'}(V_\Z) \otimes S^{\pi''}(V_\Z) \rightarrow S^\pi(V_\Z).
\end{equation}
If we decompose $\pi$ into its individual columns
as $\pi = \pi_1 \otimes \cdots \otimes \pi_l$,
and then iterate the map (\ref{themapmu}) a total of
$(l-1)$ times, we get a multiplication map
$$
S^{\pi_1}(V_\Z) \otimes\cdots\otimes S^{\pi_l}(V_\Z)
\rightarrow S^\pi(V_\Z).
$$
Identifying 
$S^{\pi_1}(V_\Z) \otimes\cdots\otimes S^{\pi_l}(V_\Z)$
with $T^N(V_\Z)$ in the obvious fashion, this 
map gives us a surjective homomorphism
\begin{equation}\label{vdef}
\V:T^N(V_\Z) \twoheadrightarrow S^\pi(V_\Z)
\end{equation}
which is different from the canonical quotient map:
$\V$ maps
$M_{\gamma(A)}$ to $M_B$,
for $A \in \Tab_0(\pi)$ with row equivalence class $B$.
Define $P^\pi(V_\Z)$ to be the image of the subspace
$\bw^\pi(V_\Z)$ of $T^N(V_\Z)$ under this map $\V$.
Thus, the restriction of $\V$ defines a surjective
homomorphism
\begin{equation}\label{themapf}
\V:\bw^\pi(V_\Z) \twoheadrightarrow P^\pi(V_\Z).
\end{equation}
The $U_\Z$-module $P^\pi(V_\Z)$ is a well known 
$\Z$-form for the irreducible polynomial representation 
of $\mathfrak{gl}_\infty$ parametrized by the partition
$\lambda = (p_1,\dots,p_n)$.
For any $A \in \Col_0(\pi)$, define
\begin{equation}\label{vadef}
V_A := \V(N_A).
\end{equation}
By a classical result,
$P^\pi(V_\Z)$ is a free $\Z$-module
with {\em standard monomial basis} given by the
vectors $\{V_A\:|\:A \in \Std_0(\pi)\}$; see 
 \cite[Theorem 26]{qla} for a non-classical proof. 
Moreover, 
for $A \in \Col_0(\pi)$, we have that
\begin{equation}\label{fka}
\V(K_A) = \left\{
\begin{array}{ll}
L_{R(A)}&\text{if $A \in \Std_0(\pi)$,}\\
0&\text{otherwise,}
\end{array}
\right.
\end{equation}
recalling the rectification map $R$ from (\ref{rmap}).
The vectors $\{L_A\:|\:A \in \Dom_0(\pi)\}$ give
another basis for the submodule $P^\pi(V_\Z)$, which is the
{\em dual canonical basis} of Lusztig, or Kashiwara's {\em upper global 
crystal basis}. Finally, by (\ref{nak}) and (\ref{fka}), we have 
for any $A \in \Col_0(\pi)$ that
\begin{equation}\label{vae}
V_A = \sum_{B \in \Std_0(\pi)} 
\left(
\sum_{C \sim_\co A} (-1)^{\ell(A,C)} P_{d(\gamma(C)) w_0, d(\gamma(B))w_0}(1)
 \right) L_{R(B)}.
\end{equation}

\section{Crystals}\label{sscrystals}
In this section, we introduce the crystals underlying the
$U_\Z$-modules $T^N(V_\Z)$, 
$\bw^\pi(V_\Z)$, $S^\pi(V_\Z)$ and $P^\pi(V_\Z)$.
First, we define a crystal
$(\Z^N, \tilde e_i, \tilde f_i, \eps_i, \phi_i, \theta)$ 
in the sense of Kashiwara \cite{Ka}, as follows.
Take $\alpha = (a_1,\dots,a_N) \in \Z^N$ and $i \in \Z$.
The {\em $i$-signature} of $\alpha$ is the 
tuple $(\sigma_1,\dots,\sigma_N)$ defined from
\begin{equation}
\sigma_j = \left\{
\begin{array}{ll}
+&\hbox{if $a_j = i$,}\\
-&\hbox{if $a_j = i+1$,}\\
0&\hbox{otherwise.}
\end{array}\right.
\end{equation}
From this the {\em reduced $i$-signature} is computed
by successively replacing 
subsequences of the form $-+$ (possibly separated by $0$'s)
in the signature with $0 0$
until no $-$ appears to the left of 
a $+$.
Recall $\delta_j$ denotes the 
weight $(0,\dots,0,1,0,\dots,0) \in \C^N$
where $1$ appears in the $j$th place.
Define
\begin{align}\label{dc1}
\tilde e_i(\alpha) &:= \left\{
\begin{array}{ll}
\emptyset&\hbox{if there are no $-$'s in the reduced $i$-signature},\\
\alpha - \delta_j&\hbox{if the leftmost $-$ is in position $j$;}
\end{array}\right.\\
\tilde f_i(\alpha) &:= \left\{
\begin{array}{ll}
\emptyset&\hbox{if there are no $+$'s in the reduced $i$-signature},\\
\alpha + \delta_j&\hbox{if the rightmost $+$ is in position $j$;}
\end{array}\right.\\
\eps_i(\alpha) &= \hbox{the total number of $-$'s in the 
reduced $i$-signature};\\
\phi_i(\alpha) &= \hbox{the total number of $+$'s in the reduced 
$i$-signature}.\label{dc2}
\end{align}
Finally define
$\theta:\Z^N \rightarrow P_\infty$ to be the restriction of the map
(\ref{contdef0}).
This completes the definition of the crystal
$(\Z^N, \tilde e_i, \tilde f_i, \eps_i, \phi_i, \theta)$.
It is the $N$-fold tensor product of the
usual crystal associated to the natural module $V_\Z$ (but
for the opposite tensor product to the one used in \cite{Ka}).
This crystal carries information about the action of the Chevalley generators
of $U_\Z$ on the dual canonical basis $\{L_\alpha\:|\:\alpha \in \Z^N\}$ of
$T^N(V_\Z)$, thanks to the following result of 
Kashiwara \cite[Proposition 5.3.1]{Kaduke}:
for $\alpha \in \Z^N$, we have that 
\begin{align}\label{eil}
e_i L_\alpha &= 
\eps_i(\alpha)
L_{\tilde e_i (\alpha)} + 
\sum_{\substack{\beta \in \Z^N\\ \eps_i(\beta) < \eps_i(\alpha)-1}}
x_{\alpha,\beta}^i L_\beta\\
f_i L_\alpha &= 
\phi_i(\alpha)
L_{\tilde f_i (\alpha)} + 
\sum_{\substack{\beta \in \Z^N \\ \phi_i(\beta) < \phi_i(\alpha)-1}}
y_{\alpha,\beta}^i L_\beta\label{fil}
\end{align}
for $x_{\alpha,\beta}^i, y_{\alpha,\beta}^i \in \Z$.
The right hand side of (\ref{eil}) resp. (\ref{fil}) should be 
interpreted as zero if $\eps_i(\alpha) = 0$ resp. $\phi_i(\alpha) = 0$.

There are also crystals attached to the modules $S^\pi(V_\Z)$ and $\bw^\pi(V_\Z)$. To define them, identify $\Row_0(\pi)$ with a subset of $\Z^N$ by
row reading $\rho:\Row_0(\pi) \hookrightarrow \Z^N$, 
and identify $\Col_0(\pi)$ with a subset of $\Z^N$ by column reading
$\gamma:\Col_0(\pi) \hookrightarrow \Z^N$. In this way, 
both $\Row_0(\pi)$ and $\Col_0(\pi)$ become identified with 
subcrystals of the crystal
$(\Z^N, \tilde e_i, \tilde f_i, \eps_i, \phi_i, \theta)$.
This defines crystals
$(\Row_0(\pi), \tilde e_i, \tilde f_i, \eps_i, \phi_i, \theta)$
and 
$(\Col_0(\pi), \tilde e_i, \tilde f_i, \eps_i, \phi_i, \theta)$.
These crystals control the action of the Chevalley generators of
$U_\Z$ on the dual canonical bases
$\{L_A\:|\:A \in \Row_0(\pi)\}$ and
$\{K_A\:|\:A \in \Col_0(\pi)\}$, just like in (\ref{eil})--(\ref{fil}).
First, for $A \in \Row_0(\pi)$, we have that
\begin{align}\label{itf1}
e_i L_A &= \eps_i(A) L_{\tilde e_i(A)} + \sum_{\substack{B \in \Row_0(\pi) \\ \eps_i(B) < \eps_i(A) - 1}} x_{\rho(A), \rho(B)}^i L_B,\\
f_i L_A &= \phi_i(A) L_{\tilde f_i(A)} + \sum_{\substack{B \in \Row_0(\pi) \\ \phi_i(B) < \phi_i(A) - 1}} y_{\rho(A), \rho(B)}^i L_B.\label{itf2}
\end{align}
Second, for $A \in \Col_0(\pi)$, we have that
\begin{align}
e_i K_A &= \eps_i(A) K_{\tilde e_i(A)} + \sum_{\substack{B \in \Col_0(\pi) \\ \eps_i(B) < \eps_i(A) - 1}} x_{\gamma(A), \gamma(B)}^i K_B,\\
f_i K_A &= \phi_i(A) K_{\tilde f_i(A)} + \sum_{\substack{B \in \Col_0(\pi) \\ \phi_i(B) < \phi_i(A) - 1}} y_{\gamma(A), \gamma(B)}^i K_B.
\end{align}

Finally, there is a well known crystal attached to the polynomial
representation $P^\pi(V_\Z)$.
This has various different realizations, in terms of either 
the set $\Dom_0(\pi)$ 
or the set $\Std_0(\pi)$;
the realization as $\Std_0(\pi)$ when $\pi$ is
left-justified is the usual description from \cite{KN}. 
In the first case, we note that
$\Dom_0(\pi)$ is a subcrystal of the crystal
$(\Row_0(\pi), \tilde e_i, \tilde f_i, \eps_i, \phi_i, \theta)$, indeed it is the connected component of this crystal generated by the row equivalence class
of the
{\em ground-state tableau} $A_0$,\label{gst} 
that is, the tableau having 
all entries on row $i$ equal to $(1-i)$.
In the second case, as explained in \cite[$\S$2]{qla},
$\Std_0(\pi)$ is a subcrystal of the crystal
$(\Col_0(\pi), \tilde e_i, \tilde f_i, \eps_i, \phi_i, \theta)$, indeed again
it is the connected component of this crystal generated by the ground-state
tableau $A_0$.
In this way, we obtain two new crystals
$(\Dom_0(\pi), \tilde e_i, \tilde f_i, \eps_i, \phi_i, \theta)$ and
$(\Std_0(\pi), \tilde e_i, \tilde f_i, \eps_i, \phi_i, \theta)$.
The rectification map $R:\Std_0(\pi) \rightarrow \Dom_0(\pi)$
is the unique isomorphism between these crystals,
and it sends the ground-state tableau $A_0$ to its
row equivalence class.

\section{Consequences of the Kazhdan-Lusztig conjecture}\label{sskl}
In this section, we record a 
representation theoretic interpretation of the dual canonical basis of the
spaces $T^N(V_{\Z})$ and 
$\bigwedge^\pi(V_\Z)$, which is a well known reformulation 
of the Kazhdan-Lusztig conjecture \cite{BB,BrK}
in type $A$.
Later on in the article we will formulate analogous interpretations
for the dual canonical bases of the spaces $S^\pi(V_\Z)$ (conjecturally)
and $P^\pi(V_\Z)$.
Go back to the notation from $\S$\ref{sscent}, 
so $\mathfrak{g} = \mathfrak{gl}_N$,
$\mathfrak{d}$ is the standard Cartan subalgebra of diagonal matrices
and $\mathfrak{b}$ is the standard Borel subalgebra of upper triangular matrices.
Let $\mathcal O$\label{ocat} denote the \cite{BGG} category of all
finitely generated
$\mathfrak{g}$-modules which are locally finite over $\mathfrak{b}$
and semisimple over $\mathfrak{d}$.
The basic objects in $\mathcal O$ are the Verma modules $M(\alpha)$
and their unique irreducible quotients $L(\alpha)$
for $\alpha = (a_1,\dots,a_N) \in \C^N$, using the $\rho$-shifted notation
explained by (\ref{vermamod}).
Also recall
that we have parametrized the central characters
of $U(\mathfrak{g})$ by the set of elements $\theta$ of
$P = \bigoplus_{a \in \C} \Z \gamma_a$ whose coefficients are 
non-negative integers summing to $N$.

For $\theta \in P$, 
let $\mathcal O(\theta)$\label{ocatp} denote the full subcategory of $\mathcal O$
consisting of the objects all of whose composition factors
are of central character $\theta$, setting $\mathcal O(\theta) = 0$
by convention if the coefficients of $\theta$ are not
non-negative integers summing to $N$.
The category $\mathcal O$ has the following 
{\em ``block'' decomposition}:
\begin{equation}
\mathcal O = \bigoplus_{\theta \in P} \mathcal O(\theta).
\end{equation}
(For non-integral central characters
our ``blocks'' are not necessarily indecomposable.)
We will write $\pr_\theta:\mathcal O \rightarrow \mathcal O(\theta)$
for the natural projection functor. To be absolutely explicit,
if the coefficients of $\theta \in P$ 
are non-negative integers
summing to $N$ so $\theta$ corresponds
to the polynomial
$f(u) = u^N + f^{(1)} u^{N-1} + \cdots  + f^{(N)} \in \C[u]$
according to (\ref{pe1})--(\ref{pe2}), 
we have that
\begin{equation}\label{prtheta1}
\pr_\theta (M) = 
\left\{v \in M\:\bigg|\:
\begin{array}{ll}\hbox{for each 
$r=1,\dots,N$ there exists $p > 0$}\\
\hbox{such that
$(Z_N^{(r)} - f^{(r)})^p v = 0$}\end{array}\right\}.
\end{equation}
We have already observed in $\S$\ref{sscent} that the Verma module $M(\alpha)$ 
is of central character $\theta(\alpha)$. Hence, for
any $\theta \in P$, the 
modules $\{L(\alpha)\:|\:\alpha \in \C^N\text{ with }\theta(\alpha) = \theta\}$
form a complete set of pairwise non-isomorphic irreducibles
in the category $\mathcal O(\theta)$. 

Recall that the integral weight lattice 
$P_\infty$ of $\mathfrak{gl}_\infty$ is the subgroup
$\bigoplus_{i \in \Z} \Z \gamma_i$ of $P$.
Let us restrict our attention from now on to the full subcategory
\begin{equation}\label{vincent}
\mathcal O_0 = \bigoplus_{\theta \in P_\infty \subset P} \mathcal O(\theta)
\end{equation}
of $\mathcal O$
corresponding just to 
{\em integral} central characters.
The Grothendieck group $[\mathcal O_0]$ of this category has
the two natural bases
$\{[M(\alpha)]\:|\:\alpha \in \Z^N\}$ and
$\{[L(\alpha)]\:|\:\alpha \in \Z^N\}$.
Define a $\Z$-module isomorphism
\begin{equation}\label{canid}
j:T^N(V_\Z) \rightarrow [\mathcal O_0],\qquad
M_\alpha \mapsto [M(\alpha)].
\end{equation}
Note this isomorphism sends the $\theta$-weight space of
$T^N(V_{\Z})$ isomorphically onto 
the block component $[\mathcal O(\theta)]$
of $[\mathcal O_0]$, for each $\theta \in P_\infty$.
The Kazhdan-Lusztig conjecture \cite{KL}, 
proved in
\cite{BB, BrK},
can be formulated as follows
for the special case of the Lie algebra
$\mathfrak{gl}_N$.

\begin{Theorem}\label{klconj}
The map $j$ sends the dual canonical basis element
$L_{\alpha}$ of $T^N(V_{\Z})$ to the
class $[L(\alpha)]$ of the irreducible module $L(\alpha)$.
\end{Theorem}

\begin{proof}
In view of (\ref{mpl}), it suffices to show 
for $\alpha,\beta \in \Z^N$ 
that the composition
multiplicity of $L(\beta)$ in the Verma module $M(\alpha)$ is
given by the formula
$$
[M(\alpha):L(\beta)] = P_{d(\alpha)w_0, d(\beta) w_0}(1).
$$
This is well known consequence of the Kazhdan-Lusztig conjecture
combined with the translation principle for singular weights, 
or see \cite[Theorem 3.11.4]{BGS}.
\end{proof}

Using (\ref{canid}) we can view the action of $U_\Z$ on $T^N(V_\Z)$
instead as an
action on the Grothendieck group $[\mathcal O_0]$. The resulting
actions of the Chevalley generators
$e_i, f_i$ of $U_\Z$ on $[\mathcal O_0]$ are in fact induced 
by some exact functors 
$e_i, f_i:\mathcal O_0 \rightarrow \mathcal O_0$
on the category $\mathcal O_0$ itself. 
Like in \cite{BKtf}, 
these functors are certain {\em translation functors} arising 
from tensoring with the 
natural $\mathfrak{g}$-module or its dual 
then projecting onto certain blocks.
To be precise, let $V$ denote the natural $N$-dimensional $\mathfrak{g}$-module
of column vectors and let $V^*$ be its dual.
Then, for $i \in \Z$, we have that
\begin{align}\label{tf1}
e_i = \bigoplus_{\theta \in P_\infty}
\pr_{\theta + (\gamma_i - \gamma_{i+1})} \circ (? \otimes V^*) \circ
\pr_\theta,\\
f_i = \bigoplus_{\theta \in P_\infty}
\pr_{\theta - (\gamma_i - \gamma_{i+1})} \circ (? \otimes V) \circ
\pr_\theta.\label{tf2}
\end{align}
These exact functors are both left and right adjoint to each 
other in a canonical way (induced by the standard adjunctions between
$? \otimes V$ and $? \otimes V^*$). The next lemma is a well known consequence
of the tensor identity.

\begin{Lemma}\label{tid}
For $\alpha \in \C^N$, the module  
$M(\alpha) \otimes V$
has a filtration with
factors $M(\beta)$ for all weights $\beta \in \C^N$ obtained from 
$\alpha$ by adding $1$ to one of its entries.
Similarly, the module
$M(\alpha) \otimes V^*$
has a filtration with
factors $M(\beta)$ for all weights $\beta \in \C^N$ obtained from 
$\alpha$ by subtracting $1$ from one of its entries.
\end{Lemma}

Taking blocks and passing to the Grothendieck group, we deduce for $\alpha 
\in \Z^N$ and $i \in \Z$
that 
\begin{equation}\label{tide}
[e_i M(\alpha)] = \sum_{\beta} [M(\beta)]
\end{equation}
summing over all weights $\beta 
\in \Z^N$ obtained from $\alpha$ by replacing
an entry equal to $(i+1)$ by an $i$, and
\begin{equation}\label{tidf}
[f_i M(\alpha)] = \sum_{\beta} [M(\beta)]
\end{equation}
summing over all weights $\beta \in \Z^N$ obtained from $\alpha$ 
by replacing
an entry equal to $i$ by an $(i+1)$.
This verifies 
that the maps on the Grothendieck group $[\mathcal O_0]$ 
induced by the exact functors $e_i, f_i$
really do 
coincide with the action of the Chevalley generators of $U_\Z$
from (\ref{canid}).

Here is an alternative definition
of the functors $e_i$ and $f_i$, explained in detail in 
\cite[$\S$7.4]{CR}.
Let $\Omega = \sum_{i,j=1}^N e_{i,j} \otimes e_{j,i}
\in U(\mathfrak{g}) \otimes U(\mathfrak{g})$.
This element centralizes the image of
$U(\mathfrak{g})$ under the comultiplication
$\Delta:U(\mathfrak{g}) \rightarrow U(\mathfrak{g})\otimes
U(\mathfrak{g})$.
For any $M \in \mathcal O_0$,
$f_i M$ is precisely the generalized $i$-eigenspace
of the operator $\Omega$ acting on $M \otimes V$,
for any $M \in \mathcal O_0$.
Similarly, $e_i M$ is the generalized
$-(N+i)$-eigenspace of $\Omega$ acting on $M \otimes V^*$.

We need to recall a little more of the setup from \cite{CR}.
Define an endomorphism $x$ of the functor
$? \otimes V$ by letting $x_M:M \otimes V
\rightarrow M \otimes V$ be left multiplication by $\Omega$,
for all $\mathfrak{g}$-modules $M$.
Also define an endomorphism $s$ of the functor
$? \otimes V \otimes V$ by letting
$s_M: M \otimes V \otimes V \rightarrow
M \otimes V \otimes V$ be the permutation
$m \otimes v \otimes v' \mapsto m \otimes v' \otimes v$.
By \cite[Lemma 7.21]{CR}, we have that
\begin{equation}\label{cr1}
s_M \circ (x_M \otimes \id_{V})
= x_{M \otimes V} \circ s_M - \id_{M \otimes V \otimes V}
\end{equation}
for any $\mathfrak{g}$-module $M$, equality of maps
from $M \otimes V \otimes V$ to itself.
It follows that $x$ and $s$ restrict to well-defined
endomorphisms of the 
functors $f_i$ and $f_i^2$; we 
denote these restrictions
by $x$ and $s$ too.
Moreover, we have that
\begin{align}\label{cr2}
(s1_{f_i}) \circ (1_{f_i}s) \circ (s1_{f_i})&=
(1_{f_i} s) \circ (s 1_{f_i}) \circ (1_{f_i} s),\\
s^2 &= 1_{f_i^2},\\
s \circ (1_{f_i} x) &= (x 1_{f_i}) \circ s - 1_{f_i^2},\label{cr3}
\end{align}
equality of 
endomorphisms of $f_i^3, f_i^2$ and $f_i^2$, respectively.
In the language of 
\cite[$\S$5.2.1]{CR}, this shows that the category $\mathcal O_0$
equipped with the adjoint pair of functors $(f_i, e_i)$ 
and the endomorphisms $x \in \End(f_i)$ and $s \in \End(f_i^2)$
is an {\em $\mathfrak{sl}_2$-categorification} for each $i \in \Z$.
This has a number of important consequences, explored in detail
in \cite{CR}. 
We just record one more thing here, our proof of which 
also depends on Theorem~\ref{klconj}; see \cite{Ku}
for an independent proof.
Recall for the statement the definition of the
crystal
$(\Z^N, \tilde e_i, \tilde f_i, \eps_i, \varphi_i, \theta)$
from (\ref{dc1})--(\ref{dc2}).

\begin{Theorem}\label{kuj}
Let $\alpha \in \Z^N$ and $i \in \Z$.
\begin{enumerate}
\item[\rm(i)]
If $\eps_i(\alpha) = 0$ then $e_i L(\alpha) = 0$.
Otherwise, $e_i L(\alpha)$ is an indecomposable module with
irreducible socle and cosocle isomorphic to $L(\tilde e_i(\alpha))$.
\item[\rm(ii)]
If $\phi_i(\alpha) = 0$ then $f_i L(\alpha) = 0$.
Otherwise, $f_i L(\alpha)$ is an indecomposable module with
irreducible socle and cosocle isomorphic to $L(\tilde f_i(\alpha))$.
\end{enumerate}
\end{Theorem}

\begin{proof}
(i) For $\alpha \in \Z^N$, 
let $\eps'_i(\alpha)$ be the maximal integer $k \geq 0$
such that $(e_i)^k L(\alpha) \neq 0$.
If $\eps'_i(\alpha) > 0$, then \cite[Proposition 5.23]{CR}
shows that $e_i L(\alpha)$ is an indecomposable module with
irreducible socle and cosocle isomorphic to 
$L(\tilde e'_i(\alpha))$ for some $\tilde e'_i(\alpha) \in \Z^N$.
Moreover, using \cite[Lemma 4.3]{CR} 
too, $\eps'_i(\tilde e'_i(\alpha)) = \eps'_i(\alpha) - 1$
and all remaining composition factors of $e_i L(\alpha)$
not isomorphic to $L(\tilde e'_i(\alpha))$
are of the form $L(\beta)$ for $\beta \in \Z^N$
with $\eps'_i(\beta) < \eps'_i(\alpha)-1$.

Observe from (\ref{eil}) 
that $\eps_i(\alpha)$ is the maximal integer $k \geq 0$
such that $(e_i)^k L_{\alpha} \neq 0$, and assuming
$\eps_i(\alpha) > 0$ we know that
$e_i L_\alpha = \eps_i(\alpha) L_{\tilde e_i(\alpha)}$ 
plus a linear combination 
of $L_{\beta}$'s with $\eps_i(\beta) < \eps_i(\alpha)-1$.
Applying Theorem~\ref{klconj} and comparing with the preceeding paragraph, 
it follows immediately that $\eps_i(\alpha) = \eps'_i(\alpha)$,
in which case $\tilde e_i(\alpha) = \tilde e'_i(\alpha)$. 
This completes the proof.

(ii) Similar, or follows
from (i) using adjointness.
\end{proof}

It just remains to extend all of this to 
the parabolic case.
Continuing with the fixed pyramid $\pi = (q_1,\dots,q_l)$, recall from
(\ref{sswalgebras}) that $\mathfrak{h}$ denotes the standard Levi subalgebra
$\mathfrak{gl}_{q_1} \oplus\cdots\oplus \mathfrak{gl}_{q_l}$ 
of $\mathfrak{g}$ and $\mathfrak{p}$ is the corresponding
standard parabolic subalgebra of $\mathfrak{g}$.
Let $\mathcal O(\pi)$ denote the {\em parabolic category $\mathcal O$} consisting of all finitely generated $\mathfrak{g}$-modules that
are locally finite dimensional over $\mathfrak{p}$ and semisimple over 
$\mathfrak{h}$.
Note $\mathcal O(\pi)$ is a full subcategory
of the category $\mathcal O$.
To define the basic modules in 
$\mathcal O(\pi)$, let
$A \in \Col(\pi)$ be a column strict $\pi$-tableau and let
$\alpha=(a_1,\dots,a_N)\in\C^N$ denote the weight 
$\gamma(A)$ obtained from column reading
$A$ as in $\S$\ref{sstableaux}.
Let $V(\alpha)$ denote the
usual finite dimensional 
irreducible $\mathfrak{h}$-module
of highest weight $$
\alpha-\rho=(a_1,a_2+1,\dots,a_N+N-1).
$$ 
View $V(\alpha)$ as a $\mathfrak{p}$-module through the natural
projection $\mathfrak{p}\twoheadrightarrow \mathfrak{h}$, then form the
{\em parabolic Verma module}
\begin{equation}\label{parv}
N(A) := U(\mathfrak{g}) \otimes_{U(\mathfrak{p})} V(\alpha).
\end{equation}
The unique irreducible quotient of $N(A)$ is denoted $K(A)$;
by comparing highest weights we have that $K(A) \cong L(\alpha)$.
In this way, we obtain two natural bases
$\{[N(A)]\:|\:A \in \Col(\pi)\}$ and $\{[K(A)]\:|\:A \in \Col(\pi)\}$
for the Grothendieck group $[\mathcal O(\pi)]$ of $\mathcal O(\pi)$.
The vectors
$\{[N(A)]\:|\:A \in \Col_0(\pi)\}$ and $\{[K(A)]\:|\:A \in \Col_0(\pi)\}$
form bases for the Grothendieck group $[\mathcal O_0(\pi)]$ of the 
full subcategory
$\mathcal O_0(\pi) := \mathcal O(\pi) \cap \mathcal O_0$.
Moreover, the translation functors $e_i, f_i$ from (\ref{tf1})--(\ref{tf2}) 
send modules in $\mathcal O_0(\pi)$ to modules in $\mathcal O_0(\pi)$, hence 
the Grothendieck group $[\mathcal O_0(\pi)]$ is a $U_\Z$-submodule
of $[\mathcal O_0(\pi)]$.
Also recall the definition of the crystal structure
on $\Col_0(\pi)$ from $\S$\ref{sscrystals}.

\begin{Theorem}\label{parkl}
There is a unique $U_\Z$-module isomorphism $i:\bigwedge^\pi(V_\Z)
\rightarrow [\mathcal O_0(\pi)]$ 
such that
$i(N_A) = [N(A)]$ 
and
$i(K_A) = [K(A)]$ 
for each $A \in \Col_0(\pi)$.
Moreover, for $A \in \Col_0(\pi)$ and $i \in \Z$, the following
properties hold:
\begin{enumerate}
\item[\rm(i)]
If $\eps_i(A) = 0$ then $e_i K(A) = 0$.
Otherwise, $e_i K(A)$ is an indecomposable module with
irreducible socle and cosocle isomorphic to $K(\tilde e_i(A))$.
\item[\rm(ii)]
If $\phi_i(A) = 0$ then $f_i K(A) = 0$.
Otherwise, $f_i K(A)$ is an indecomposable module with
irreducible socle and cosocle isomorphic to $K(\tilde f_i(A))$.
\end{enumerate}
\end{Theorem}

\begin{proof}
Define a $\Z$-module isomorphism
$i:\bw^{\pi}(V_\Z) \rightarrow [\mathcal O_0(\pi)]$
by setting $i(N_A) := [N(A)]$ for each $A \in \Col_0(\pi)$.
We observe that the following diagram commutes:
\begin{equation}\label{cdcd}
\begin{CD}
\bw^{\pi}(V_\Z)&@>>>&T^N(V_\Z)\\
@ViVV&&@VVjV\\
[\mathcal O_{0}(\pi)]&@>>>&[\mathcal O_0]
\end{CD}
\end{equation}
where the horizontal maps are the natural inclusions. 
This is checked by computing the image  
either way round the diagram of $N_A$: one way round one uses 
the definitions (\ref{nadef}) 
and (\ref{canid}); the other way round uses the Weyl character formula to 
express $[V(\alpha)]$ 
as a linear combination of Verma modules over $\mathfrak h$, 
then exactness of the functor $U(\mathfrak{g}) \otimes_{U(\mathfrak{p})} ?$
to express $[N(A)]$ as a linear combination of Verma modules over
$\mathfrak{g}$.
Since we already know that all of the maps apart from $i$
are $U_\Z$-module homomorphisms, it then follows that $i$ is too.
To complete the proof of the first statement of the theorem,
it just remains to show that
$i(K_A) = [K(A)]$.
This follows by Theorem~\ref{klconj} because $K(A) \cong L(\gamma(A))$ and
$K_A = L_{\gamma(A)}$.
The remaining statements (i) and (ii) follow from 
Theorem~\ref{kuj}.
\end{proof}

\chapter{Highest weight theory}\label{shw}

In this chapter, we set up the usual machinery of highest weight theory
for the shifted Yangian $Y_n(\sigma)$, exploiting its 
triangular decomposition.
Fix throughout a shift matrix $\sigma = (s_{i,j})_{1 \leq i,j \leq n}$.

\section{Admissible modules}\label{rootsys}

Recall the definition of the Lie subalgebra
$\mathfrak{c}$ of $Y_n(\sigma)$ from $\S$\ref{ssgenerators}, and the root decomposition
(\ref{qgrading}).
Given a $\mathfrak{c}$-module $M$ and a weight $\alpha \in \mathfrak{c}^*$,
the {\em generalized $\alpha$-weight space} of $M$ is the subspace
\begin{equation}
M_\alpha := \left\{v \in M\:\bigg|\: \begin{array}{l}
\hbox{for each $i=1,\dots,n$ there exists $p > 0$}
\\\hbox{such that
$(D_i^{(1)} - \alpha(D_i^{(1)}))^p v = 0$}\end{array}\right\}.
\end{equation}
We say that $M$ is {\em admissible} if
\begin{enumerate}
\item[(a)] 
$M$ is the direct sum of its generalized weight spaces, i.e.
$M = \bigoplus_{\alpha \in \mathfrak{c}^*} M_\alpha$;
\item[(b)]
each $M_\alpha$ is finite dimensional;
\item[(c)]
the set of all $\alpha \in \mathfrak{c}^*$ such that
$M_\alpha$ is non-zero is contained in a finite union of sets
of the form $D(\beta)
:= \{\alpha \in \mathfrak{c}^*\:|\:\alpha \leq \beta\}$ for
$\beta \in \mathfrak{c}^*$.
\end{enumerate}
Given a $\mathfrak{c}$-module $M$ satisfying (a),
we define its {\em restricted dual}
\begin{equation}\label{rdual}
\overline{M} := \bigoplus_{\alpha \in \mathfrak c^*} (M_\alpha)^*
\end{equation}
to be the direct sum of the duals of its generalized weight spaces.

By an {\em admissible $Y_n(\sigma)$-module}, 
we mean a left $Y_n(\sigma)$-module which is admissible
when viewed as a $\mathfrak{c}$-module by restriction.
In that case, 
$\overline{M}$ is naturally a {\em right} $Y_n(\sigma)$-module
with action $(fx)(v) = f(xv)$ for
$f \in \overline{M}, v \in M$ and $x \in Y_n(\sigma)$.
Hence twisting with the inverse of the anti-isomorphism
$\tau:Y_n(\sigma)\rightarrow Y_n(\sigma^t)$ 
from (\ref{taudef}) we can view $\overline{M}$ instead as a
left $Y_n(\sigma^t)$-module, which we denote by $M^\tau$.
It is obvious that $M^\tau$ is also admissible.
Indeed, making the obvious definition on morphisms, 
$?^\tau$ can be viewed as a contravariant equivalence
between the categories of admissible $Y_n(\sigma)$- and $Y_n(\sigma^t)$-modules.

\section{Gelfand-Tsetlin characters}\label{sscharblock}

Next, let
$\mathscr P_n$ denote the set of all power series
$A(u) = A_1(u_1) A_2(u_{2}) \cdots A_n(u_n)$ 
in indeterminates $u_1,\dots,u_n$
such that
each $A_i(u)$ belongs to $1 + u^{-1}\C[[u^{-1}]]$.
Note that $\mathscr P_n$ is an abelian group under multiplication.
For $A(u) \in \mathscr P_n$, we always 
write $A_i(u)$ for the $i$th power series 
defined from the equation $A(u) = A_1(u_1)\cdots A_n(u_n)$
and $A_i^{(r)}$ for the $u^{-r}$-coefficient of $A_i(u)$. The
associated weight of $A(u) \in \mathscr P_n$ is defined by
\begin{equation}
\wt A(u) := A_1^{(1)} \eps_1 + A_2^{(1)}\eps_2+\cdots+A_n^{(1)} \eps_n
\in \mathfrak{c}^*.
\end{equation}
Now we form the {\em completed group algebra}
$\widehat\Z[\mathscr P_n]$. The elements of 
$\widehat \Z[\mathscr P_n]$ consist of formal sums
$S = \sum_{A(u) \in \mathscr P_n} m_{A(u)} [A(u)]$
for integers $m_{A(u)}$ with the property that
\begin{enumerate}
\item[(a)] the set $\{\wt A(u)\:|\:A(u) \in \operatorname{supp} S\}$
is contained in
a finite union of sets of the form $D(\beta)$ for $\beta \in \mathfrak{c}^*$;
\item[(b)] for each $\alpha \in \mathfrak{c}^*$ the
set $\{A(u) \in \operatorname{supp} S\:|\:\wt A(u) = \alpha\}$
is finite,
\end{enumerate}
where 
$\operatorname{supp} S$ denotes 
$\{A(u) \in \mathscr P_n\:|\:m_{A(u)} \neq 0\}$.
There is an obvious multiplication on $\widehat\Z[\mathscr P_n]$
extending the rule $[A(u)] [B(u)] = [A(u)B(u)]$.

Given an admissible $Y_n(\sigma)$-module $M$ 
and $A(u) \in \mathscr P_n$, the corresponding 
{\em Gelfand-Tsetlin subspace}
of $M$ is defined by
\begin{equation}\label{gws}
M_{A(u)} := \left\{v \in M\:\bigg|\:
\begin{array}{ll}\hbox{for each 
$i=1,\dots,n$ and $r > 0$ there exists}\\
\hbox{$p > 0$ such that
$(D_i^{(r)} - A_i^{(r)})^p v = 0$}\end{array}\right\}.
\end{equation}
Since the weight spaces of $M$ are finite dimensional and the operators
$D_i^{(r)}$ commute with each other, we have 
for each $\alpha \in \mathfrak{c}^*$ that
\begin{equation}
M_\alpha = \bigoplus_{\substack{A(u) \in \mathscr P_n \\ \wt A(u) = \alpha}}
M_{A(u)}.
\end{equation}
Hence, since $M$ is the direct sum of its generalized weight spaces, 
it is also the direct sum of its Gelfand-Tsetlin subspaces:
$M = \bigoplus_{A(u) \in \mathscr P_n} M_{A(u)}$.
Now we are ready to introduce a notion of {\em Gelfand-Tsetlin character}
of an admissible $Y_n(\sigma)$-module $M$, which is analogous to the 
characters of Knight \cite{Knight} for Yangians in general
and of Frenkel and Reshetikhin \cite{FR} in the setting of 
quantum affine algebras: set
\begin{equation}
\ch M := \sum_{A(u) \in \mathscr P_n} (\dim M_{A(u)}) [A(u)].
\end{equation}
By the definition of admissibility, $\ch M$ belongs to
the completed group algebra $\widehat \Z[\mathscr P_n]$.
For example, the Gelfand-Tsetlin character of the trivial
$Y_n(\sigma)$-module is $[1]$.

For the first lemma, recall the comultiplication
$\Delta:Y_n(\sigma) \rightarrow Y_n(\sigma') \otimes
Y_n(\sigma'')$ from (\ref{dsub}), where $\sigma'$ resp.
$\sigma''$ is the strictly lower resp. upper triangular matrix
such that $\sigma = \sigma' + \sigma''$. This allows us to view the
tensor product of a 
$Y_n(\sigma')$-module
$M'$ and a $Y_n(\sigma'')$-module $M''$ as a $Y_n(\sigma)$-module.
We will always denote this ``external'' tensor product 
by $M' \boxtimes M''$, to avoid confusion with
the usual ``internal'' tensor product of 
$\mathfrak{g}$-modules which we will also exploit later on.
We point out that $\Delta(D_i^{(1)}) = D_i^{(1)} \otimes 1 + 
1 \otimes D_i^{(1)}$, so the generalized $\alpha$-weight space
of $M \boxtimes N$ is equal to $\sum_{\beta \in \mathfrak{c}^*}
M_\beta \otimes M_{\alpha-\beta}$.

\begin{Lemma}\label{ten2}
Suppose that $M'$ is an admissible $Y_n(\sigma')$-module
and $M''$ is an admissible $Y_n(\sigma'')$-module.
Then $M' \boxtimes M''$ is an admissible
$Y_n(\sigma)$-module, and
$$
\ch (M' \boxtimes M'') = (\ch M') (\ch M'').
$$
\end{Lemma}

\begin{proof}
The fact that $M' \boxtimes M''$ is admissible is obvious.
To compute its character, 
order the set of weights of $M'$
as $\alpha_1, \alpha_2, \dots$
so that $\alpha_j > \alpha_k \Rightarrow j < k$. Let $M_{j}'$
denote $\sum_{1 \leq k \leq j} M_{\alpha_k}'$.
Then Theorem~\ref{utah}(i) implies that the subspace $M_{j}' \otimes M''$ 
of $M' \otimes M''$
is invariant under the action of all $D_i^{(r)}$. Moreover
in order to compute the 
Gelfand-Tsetlin character of $M' \boxtimes M''$, we can replace
it by
$$
\bigoplus_{j \geq 1} (M_j' \otimes M'') / (M_{j-1}' \otimes M'')
=
\bigoplus_{j \geq 1} (M_j' / M_{j-1}') \otimes M''
$$
with $D_i(u)$ acting as $D_i(u) \otimes D_i(u)$.
\end{proof}

The next lemma is concerned with the duality $?^\tau$
on admissible modules.

\begin{Lemma}\label{tau}
For an admissible $Y_n(\sigma)$-module $M$, we have that
$\ch (M^\tau) = \ch M$.
\end{Lemma}

\begin{proof}
$\tau(D_i^{(r)}) = D_i^{(r)}$.
\end{proof}

\section{Highest weight modules}\label{sshwm}

For $A(u) \in \mathscr P_n$, a vector $v$ in a $Y_n(\sigma)$-module $M$ is called a 
{\em highest weight vector} of {\em type} 
$A(u)$ if
\begin{enumerate}
\item[(a)] $E_i^{(r)} v = 0$ for all $i=1,\dots,n-1$ and
$r > s_{i,i+1}$;
\item[(b)] $D_i^{(r)} v = A_i^{(r)} v$ for all $i=1,\dots,n$
and $r > 0$.
\end{enumerate}
We call $M$ a {\em highest weight
module} of {\em type} $A(u)$ if it is generated by such a highest weight
vector. The following lemma gives an equivalent way to state these definitions
in terms of the elements $T_{i,j}^{(r)}$ from (\ref{tij}).

\begin{Lemma}\label{thw}
A vector $v$ in a $Y_n(\sigma)$-module is a
highest weight vector of type $A(u)$ if and only if
$T_{i,j}^{(r)} v = 0$ for all $1 \leq i < j \leq n$ and $r > s_{i,j}$, and 
$T_{i,i}^{(r)} v = A_i^{(r)} v$ for all $i=1,\dots,n$ and $r > 0$.
\end{Lemma}

\begin{proof}
By the definition (\ref{tij}), the left ideal of
$Y_n(\sigma)$ generated by $$
\{E_i^{(r)}\:|\:i=1,\dots,n-1,r > s_{i,i+1}\}
$$
coincides with the left ideal generated by
$$
\{T_{i,j}^{(r)}\:|\:1 \leq i < j \leq n, r > s_{i,j}\}.
$$
Moreover, $T_{i,i}^{(r)} \equiv D_i^{(r)}$ modulo this left ideal.
\end{proof}

In the next lemma, we 
write $\sigma = \sigma'+\sigma''$ where $\sigma'$ resp. $\sigma''$
is strictly lower resp. upper triangular.

\begin{Lemma}\label{ten1}
Suppose $v$ is a highest weight vector in a $Y_n(\sigma')$-module $M$
of type $A(u)$ and $w$ is a highest weight vector in a 
$Y_n(\sigma'')$-module $N$ of type $B(u)$. Then
$v \otimes w$ is a highest weight vector 
in the $Y_n(\sigma)$-module $M \boxtimes N$
of type $A(u)B(u)$.
\end{Lemma}

\begin{proof}
Apply Theorem~\ref{utah}.
\end{proof}

To construct the universal highest weight module of type $A(u)$,
let $\C_{A(u)}$ denote the one dimensional $Y_{(1^n)}$-module on
which $D_i^{(r)}$ acts as the scalar $A_i^{(r)}$.
Inflating through the epimorphism
$Y_{(1^n)}^{\sharp}(\sigma) \twoheadrightarrow
Y_{(1^n)}$ from (\ref{epis}), we can view
$\C_{A(u)}$ instead as a $Y_{(1^n)}^{\sharp}(\sigma)$-module.
Now form the induced module
\begin{equation}\label{mndef}
M(\sigma, A(u)) 
:= Y_n(\sigma) \otimes_{Y_{(1^n)}^{\sharp}(\sigma)} \C_{A(u)}.
\end{equation}
This is a highest weight module of type $A(u)$,
generated by the highest weight vector $v_+ := 1 \otimes 1$.
Clearly it is the universal such module, i.e. all other highest weight
modules of this type are quotients of $M(\sigma,A(u))$.
In the next theorem we record two natural bases for $M(\sigma,A(u))$.

\begin{Theorem}\label{vb1}
For any $A(u) \in \mathscr P_n$, 
the following sets of vectors give bases for the module 
$M(\sigma,A(u))$:
\begin{enumerate}
\item[\rm(i)] $\{x v_+\:|\:x \in X\}$,
where $X$ denotes the collection of all monomials in the elements
$\{F_{i,j}^{(r)}\:|\:1 \leq i< j \leq n, s_{j,i} < r\}$ taken
in some fixed order;
\item[\rm(ii)] $\{y v_+\:|\:y \in Y\}$,
where $Y$ denotes the collection of all monomials in the elements
$\{T_{j,i}^{(r)}\:|\:1 \leq i < j \leq n, s_{j,i} < r\}$ taken
in some fixed order.
\end{enumerate}
\end{Theorem}

\begin{proof}
Let $M := M(\sigma,A(u))$. 

(i) The isomorphism (\ref{tridec}) implies that
$Y_n(\sigma)$ is a free right
$Y_{(1^n)}^\sharp(\sigma)$-module 
with basis $X$.
Hence $M$ has basis
$\{x v_+\:|\:x \in X\}$.

(ii) Recall the definition of the canonical filtration
$\F_0 Y_n(\sigma) \subseteq \F_1 Y_n(\sigma) \subseteq \cdots$
of $Y_n(\sigma)$ from $\S$\ref{sspbw}.
In view of Lemma~\ref{tbb},
it may also be defined by declaring that all 
$T_{i,j}^{(r)}$ are of degree $r$.
Also introduce 
a filtration $\F_0 M \subseteq \F_1 M \subseteq \cdots$ of $M$
by setting $\F_d M := \F_d Y_n(\sigma) v_+$.
Let $X^{(d)}$ resp. $Y^{(d)}$ 
denote the set of all monomials
in the elements $X$ resp. $Y$ of total degree at most
$d$ in the canonical filtration.
Applying (i), one deduces at once that 
the set of all vectors of the form
$\left\{x v_+\:|\:x \in X^{(d)}\right\}$ form 
a basis for $\F_d M$. On the other hand using Lemmas~\ref{tbb}
and \ref{thw}, the vectors $\left\{y v_+\:|\:y \in Y^{(d)}\right\}$ 
span $\F_d M$. By dimension they must be linearly independent too.
Since $M = \bigcup_{d \geq 0} \F_d M$, this implies that the vectors
$\{y v_+\:|\:y \in Y\}$ give a basis for $M$ itself.
\end{proof}

This implies that 
the (generalized) $\wt A(u)$-weight space of $M(\sigma,A(u))$ is 
one dimensional, spanned by the vector $v_+$, while all other weights
are strictly lower in the dominance ordering.
Given this, the usual argument shows that $M(\sigma,A(u))$ has a unique
maximal submodule denoted $\rad M(\sigma,A(u))$. Set
\begin{equation}\label{lnsa}
L(\sigma,A(u)) := M(\sigma,A(u)) / \rad M(\sigma,A(u)).
\end{equation}
This is the unique (up to isomorphism) 
{\em irreducible highest weight module}
of type $A(u)$ for the algebra $Y_n(\sigma)$.
We also note that
\begin{equation}\label{oned}
\dim \End_{Y_n(\sigma)} (L(\sigma,A(u))) 
=
1
\end{equation}
for any $A(u) \in \mathscr P_n$.

\section{Classification of admissible irreducible representations}

A natural question arises at this point:
the module $M(\sigma,A(u))$ is certainly not admissible, 
since all of its generalized weight spaces 
other than the highest one are infinite dimensional, but the
irreducible quotient $L(\sigma,A(u))$ may well be.

\begin{Theorem}\label{tar}
For $A(u) \in \mathscr P_n$,
the irreducible $Y_n(\sigma)$-module $L(\sigma, A(u))$ is admissible if and only if
$A_i(u) / A_{i+1}(u)$ is a rational function for all
$i=1,\dots,n-1$.
\end{Theorem}

\begin{proof}
$(\Leftarrow)$.
Suppose that each $A_i(u) / A_{i+1}(u)$ is a rational function.
For $f(u) \in 1 + u^{-1}\C[[u^{-1}]]$, the twist of $L(\sigma, A(u))$
by the automorphism $\mu_{f}$ from (\ref{mufdef})
is isomorphic to
$L(\sigma, f(u_1\cdots u_n)A(u))$.
This allows us to reduce to the case that
each $A_i(u)$ is actually a polynomial in $u^{-1}$. 
Assuming this, we can find $l \geq s_{n,1}+s_{1,n}$ such that,
on setting $p_i := l - s_{n,i}-s_{i,n}$,
$u^{p_i}A_i(u)$ is a monic polynomial in 
$u$ of degree $p_i$ for each $i=1,\dots,n$. 
Let $\pi = (q_1,\dots,q_l)$ be the pyramid associated to the
shift matrix $\sigma$ and the level $l$.
For each $i=1,\dots,n$, 
factorize 
$u^{p_i} A_i(u)$ as $(u+a_{i,1}) \cdots (u + a_{i,p_i})$
for $a_{i,j} \in \C$, and write the numbers
$a_{i,1},\dots,a_{i,p_i}$ into the boxes on the $i$th row
of $\pi$ from left to right. 
For each $j=1,\dots,l$, 
let $b_{j,1},\dots,b_{j,q_j}$ denote the entries 
in the $j$th column
of the resulting $\pi$-tableau
read from top to bottom.
Let $M_j$ denote the usual
Verma module for the Lie algebra $\mathfrak{gl}_{q_j}$ of highest weight
$(b_{j,1},b_{j,2},\cdots,b_{j,q_j})$.
The tensor product $M_1 \boxtimes\cdots\boxtimes M_l$
is naturally a $W(\pi)$-module, hence a 
$Y_n(\sigma)$-module via the quotient map (\ref{thetadef}).
Applying Lemmas~\ref{ten2} and \ref{ten1}, 
it is an admissible $Y_n(\sigma)$-module
and it contains an obvious highest weight vector
of type $A(u)$.

$(\Rightarrow)$.
Assume to start with that the shift matrix
$\sigma$ is the zero matrix, i.e. $Y_n(\sigma)$ is just the usual Yangian
$Y_n$.
Suppose that $L(\sigma,A(u))$
is admissible for some $A(u) \in \mathscr{P}_n$.
In particular, for each $i=1,\dots,n-1$, the
$(\wt A(u)-\eps_i +\eps_{i+1})$-weight space of
$L(\sigma,A(u))$ is finite dimensional.
Given this an argument due originally to Tarasov \cite[Theorem 1]{Tarasov1},
see e.g. the proof of \cite[Proposition 3.5]{Molevfd}, 
shows that $A_i(u) / A_{i+1}(u)$ is a rational function
for each $i=1,\dots,n-1$.

Assume next that $\sigma$ is lower triangular, and 
consider the canonical 
embedding $Y_n(\sigma) \hookrightarrow Y_n$.
Given $A(u) \in \mathscr P_n$ such that
$L(\sigma,A(u))$ is admissible, 
the PBW theorem
implies that the induced module
$Y_n\otimes_{Y_n(\sigma)} L(\sigma,A(u))$
is also admissible and contains a non-zero 
highest weight vector of type $A(u)$. 
Hence by the preceeding paragraph
$A_i(u) / A_{i+1}(u)$ is a rational function 
for each $i=1,\dots,n-1$.

Finally suppose that $\sigma$ is arbitrary.
Recalling the isomorphism $\iota$ from (\ref{iotadef}), the
twist of a highest weight module by $\iota$ is again a highest weight
module of the same type, and the twist of an admissible module
is again admissible. So
the conclusion in general follows from the lower triangular case.
\end{proof}

In view of this result, let us define
\begin{equation}\label{Qn}
\mathscr{Q}_n := \left\{A(u) 
\in \mathscr P_n\:\bigg|\: 
\begin{array}{l}
A_i(u) / A_{i+1}(u)
\hbox{ is a rational function}\\\hbox{for each }i=1,\dots,n-1\end{array}
\right\}.
\end{equation}
Then Theorem~\ref{tar} implies that the
modules $\{L(\sigma,A(u))\:|\:A(u) \in \mathscr{Q}_n\}$
give a full set of pairwise non-isomorphic 
admissible irreducible $Y_n(\sigma)$-modules.

\begin{Remark}\label{tarrem}
The construction
explained in the proof of Theorem~\ref{tar} shows moreover
that every admissible irreducible $Y_n(\sigma)$-module
can be obtained 
from an admissible irreducible $W(\pi)$-module
via the homomorphism
\begin{equation}\label{infmap}
\kappa \circ \mu_f: Y_n(\sigma) \twoheadrightarrow W(\pi),
\end{equation}
for some pyramid $\pi$ associated to the shift matrix $\sigma$
and some $f(u) \in 1 + u^{-1}\C[[u^{-1}]]$.
\end{Remark}

\section{Composition multiplicities}

The final job in this chapter is to make precise
the sense in which Gelfand-Tsetlin characters characterize admissible modules.
We need to be a little careful here since admissible modules
need not possess a composition series.
Nevertheless, given admissible $Y_n(\sigma)$-modules $M$ and  $L$
with $L$ irreducible, we define the
{\em composition multiplicity} of $L$ in $M$ by
\begin{equation}
[M:L] := \sup \#\{i=1,\dots,r\:|\:M_i / M_{i-1} \cong L\}
\end{equation}
where the supremum is over all finite filtrations
$0 = M_0 \subset M_1\subset\cdots\subset M_r = M$.
By general principles, this multiplicity is additive on short exact sequences.
Now we repeat some standard arguments from \cite[ch. 9]{Kac}.

\begin{Lemma}\label{kac1}
Let $M$ be an admissible $Y_n(\sigma)$-module
and 
$\alpha \in \mathfrak{c}^*$ be a fixed weight. There is 
a
filtration $0 = M_0 \subset M_1 \subset\cdots\subset M_r = M$ and a subset
$I \subseteq \{1,\dots,r\}$
such that
\begin{enumerate}
\item[\rm(i)] for each $i \in I$, we have that 
$M_i / M_{i-1} \cong L(\sigma, A^{(i)}(u))$ for $A^{(i)}(u) \in \mathscr{Q}_{n}$ with $\wt A^{(i)}(u) \geq \alpha$;
\item[\rm(ii)] for each $i \notin I$, we have that
$(M_i / M_{i-1})_{\beta} = 0$ for all $\beta \geq \alpha$.
\end{enumerate}
In particular, given $A(u) \in \mathscr{Q}_n$ with 
$\wt A(u) \geq \alpha$, we have that
$$
[M:L(\sigma,A(u))] = \#\{i \in I\:|\:A^{(i)}(u) =A(u) \}.
$$
\end{Lemma}

\begin{proof}
Adapt the proof of \cite[Lemma 9.6]{Kac}.
\end{proof}

\begin{Corollary}\label{span}
For an admissible $Y_n(\sigma)$-module $M$, we have that
$$
\ch M = \sum_{A(u) \in \mathscr{Q}_{n}} 
[M:L(\sigma,A(u))] \ch L(\sigma,A(u)).
$$
\end{Corollary}

\begin{proof}
Argue using the lemma exactly as in \cite[Proposition 9.7]{Kac}.
\end{proof}

\begin{Theorem}\label{inj}
Let $M$ and $N$ be admissible $Y_n(\sigma)$-modules such that
 $\ch M = \ch N$. 
Then $M$ and $N$ have all the same composition multiplicities.
\end{Theorem}

\begin{proof}
This follows from Corollary~\ref{span} once we check that the
$\ch L(\sigma,A(u))$'s are linearly independent in an appropriate sense.
To be precise we need to show, given
$$
S = \sum_{A(u) \in \mathscr{Q}_n} m_{A(u)} \ch L(\sigma,A(u))
\in \widehat \Z[\mathscr P_n]
$$ 
for coefficients
$m_{A(u)} \in \Z$ satisfying
the conditions from $\S$\ref{sscharblock}(a)--(b), that 
$S = 0$ implies each $m_{A(u)} = 0$.
Suppose for a contradiction that $m_{A(u)} \neq 0$ for some $A(u)$. 
Amongst all such $A(u)$'s, pick
one with $\wt A(u)$ maximal in the dominance ordering.
But then, since $\ch L(\sigma,A(u))$ equals
$[A(u)]$ plus a (possibly infinite) linear combination of 
$[B(u)]$'s for $\wt(B(u)) < \wt(A(u))$, 
the coefficient of $[A(u)]$ in
$\sum_{A(u) \in \mathscr{Q}_{n}} m_{A(u)} \ch L(\sigma,A(u))$ 
is non-zero, which is the desired contradiction.
\end{proof}

\begin{Corollary}\label{taudual}
For $A(u) \in \mathscr Q_n$, we have that
$L(\sigma,A(u))^\tau \cong L(\sigma^t, A(u))$.
\end{Corollary}

\begin{proof}
Using (\ref{iotadef}), it is clear that
$L(\sigma,A(u))$ and $L(\sigma^t, A(u))$ have the same
formal characters. Hence by Lemma~\ref{tau}
so do $L(\sigma,A(u))^\tau$ and
$L(\sigma^t, A(u))$. 
\end{proof}

\chapter{Verma modules}\label{sverma}

Now we turn our attention to studying highest weight modules
over the algebras $W(\pi)$ themselves.
Fix throughout the chapter a pyramid $\pi = (q_1,\dots,q_l)$ 
of height $\leq n$, let
$(p_1,\dots,p_n)$ be the tuple of row lengths, and choose
a corresponding shift matrix
$\sigma = (s_{i,j})_{1 \leq i,j \leq n}$ as usual.
Notions of weights, highest weight vectors and so forth
are exactly as in the previous chapter,
viewing $W(\pi)$-modules as $Y_n(\sigma)$-modules
via the quotient map $\kappa:Y_n(\sigma) \twoheadrightarrow W(\pi)$
from (\ref{thetadef}).

\section{Parametrization of highest weights}\label{ssp}
Our first task is to understand the universal highest weight module
of type $A(u) \in \mathscr P_n$
for the algebra $W(\pi)$.
This module is obviously the unique largest
quotient of the $Y_n(\sigma)$-module $M(\sigma,A(u))$ from (\ref{mndef})
on which the kernel 
of the 
homomorphism $\kappa:Y_n(\sigma) \twoheadrightarrow W(\pi)$ from (\ref{thetadef})
acts as zero.
In other words, it is the $W(\pi)$-module
$W(\pi) \otimes_{Y_n(\sigma)} M(\sigma,A(u))$.
We will abuse notation and write simply $v_+$ 
instead of $1 \otimes v_+$
for the highest weight vector in $W(\pi) \otimes_{Y_n(\sigma)} 
M(\sigma,A(u))$.

\begin{Theorem}\label{vermathm}
For $A(u) \in \mathscr P_n$,
$W(\pi) \otimes_{Y_n(\sigma)} M(\sigma,A(u))$
is non-zero if and only if
$u^{p_i} A_i(u) \in \C[u]$ for each $i=1,\dots,n$.
In that case,
the following sets of vectors give bases for
$W(\pi) \otimes_{Y_n(\sigma)} M(\sigma,A(u))$:
\begin{enumerate}
\item[\rm(i)]
$\{x  v_+\:|\:x \in X\}$,
where $X$ denotes the collection of all monomials in the elements
$\{F_{i,j}^{(r)}\:|\:1 \leq i< j \leq n, s_{j,i} < r \leq S_{j,i}\}$ 
taken in some fixed order;
\item[\rm(ii)]
$\{y  v_+\:|\:y \in Y\}$,
where $Y$ denotes the collection of all monomials in the elements
$\{T_{j,i}^{(r)}\:|\:1 \leq i < j \leq n, s_{j,i} < r \leq S_{j,i}\}$ 
taken in some fixed order.
\end{enumerate}
\end{Theorem}

\begin{proof}
(ii) Let us work with the following
reformulation of the definition (\ref{mndef}):
the module $M(\sigma,A(u))$ is the quotient of
$Y_n(\sigma)$ by the left ideal $J$ generated by the elements
$$
\{E_i^{(r)}\:|\:i=1,\dots,n-1, r > s_{i,i+1}\}
\cup \{D_i^{(r)} - A_i^{(r)}\:|\:i=1,\dots,n, r > 0\}.
$$
Equivalently, by Lemma~\ref{thw}, $J$ is the left ideal of $Y_n(\sigma)$
generated by
the elements 
$$
P := 
\{T_{i,j}^{(r)}\:|\:1 \leq i < j \leq n, s_{i,j} < r\}
\cup \{T_{i,i}^{(r)} - A_i^{(r)}\:|\:1 \leq i \leq n, s_{i,i} < r\}.
$$
Also let 
$Q := \{T_{j,i}^{(r)}\:|\:1 \leq i < j \leq n,
s_{j,i} < r\}.$
Pick an ordering on $P \cup Q$ so that
the elements of $Q$ preceed the elements of
$P$. Obviously all ordered monomials in the elements
$P \cup Q$ containing at least one element of
$P$ belong to $J$. Hence by
Lemma~\ref{tbb} and Theorem~\ref{vb1}(ii),
the ordered monomials in the elements
$P \cup Q$ containing at least one element of
$P$ in fact form a basis for $J$.

Now it is clear that
$W(\pi) \otimes_{Y_n(\sigma)} M(\sigma,A(u))$
is the quotient of $W(\pi)$ by the image $\bar J$ of 
$J$ under the map $\kappa:Y_n(\sigma) \twoheadrightarrow W(\pi)$.
If $A_i^{(r)} \neq 0$ for some $1 \leq i \leq n$ and $r > p_i$,
i.e. $u^{p_i} A_i(u) \notin \C[u]$,
then the image of $T_{i,i}^{(r)} - A_i^{(r)}$ gives us a unit 
in $\bar J$ by Theorem~\ref{tvan}, hence 
$W(\pi) \otimes_{Y_n(\sigma)} M(\sigma,A(u)) = 0$ in this case.
On the other hand, if all $u^{p_i} A_i(u)$ belong to $\C[u]$,
we let
\begin{align*}
\bar P &:= 
\{T_{i,j}^{(r)}\:|\:1 \leq i < j \leq n, s_{i,j} < r \leq S_{i,j}\}\\
&\qquad\qquad\qquad\qquad\cup \{T_{i,i}^{(r)} - A_i^{(r)}\:|\:1 \leq i \leq n, s_{i,i} < r \leq S_{i,i}\},\\
\bar Q &:= \{T_{j,i}^{(r)}\:|\:1 \leq i < j \leq n,
s_{j,i} < r \leq S_{j,i}\}.
\end{align*}
Then Theorem~\ref{tvan} implies that 
$\bar J$ is spanned by all ordered monomials in the elements
$\bar P \cup \bar Q$ containing at 
least one element of
$\bar P$. By Lemma~\ref{tb2}, these monomials are
also linearly independent, hence form a basis for $\bar J$. 
It follows that the image of $Y$ gives a basis for $W(\pi) / \bar J$, 
proving (ii).

(i) This follows from (ii) by reversing the argument used to deduce (ii) 
from (i) in the proof of Theorem~\ref{vb1}.
\end{proof}

Now suppose that $v_+$ is a non-zero highest weight vector
in some $W(\pi)$-module $M$.
By Theorem~\ref{vermathm}, there exist
elements $(a_{i,j})_{1 \leq i \leq n, 1 \leq j \leq p_i}$ of $\C$
such that
\begin{align}\label{t1}
u^{p_1} D_1(u) v_+ &= (u+a_{1,1}) (u+a_{1,2}) \cdots (u+a_{1,p_1}) v_+,\\
(u-1)^{p_2} D_2(u-1) v_+ &= (u+a_{2,1}) (u+a_{2,2}) \cdots (u+a_{2,p_2}) 
v_+,\\\intertext{\vspace{-3mm}$$\vdots$$\vspace{-6mm}}
(u-n+1)^{p_n} D_n(u-n+1) v_+ &= (u+a_{n,1}) (u+a_{n,2}) \cdots (u+a_{n,p_n}) 
v_+.\label{t3}
\end{align}
In this way, the highest weight vector $v_+$
defines a row symmetrized $\pi$-tableau $A$ in the sense of $\S$\ref{sstableaux}, namely, 
the unique element of $\Row(\pi)$ 
with entries $a_{i,1},\dots,a_{i,p_i}$ on its $i$th row.
From now on, we will say simply that the highest weight vector $v_+$
is of {\em type} $A$
if these equations hold.

Conversely, 
suppose that we are given $A \in \Row(\pi)$ with entries
$a_{i,1},\dots,a_{i,p_i}$ on its $i$th row.
Define the corresponding {\em generalized Verma module} $M(A)$ 
to be the universal highest weight module
of type $A$, i.e.
\begin{equation}\label{genv}
M(A) := W(\pi) \otimes_{Y_n(\sigma)} M(\sigma,A(u))
\end{equation}
where  $A(u) = A_1(u_1) \cdots A_n(u_n)$ is defined
from
$$
(u-i+1)^{p_i}
A_i(u - i + 1) = (u+a_{i,1}) (u+a_{i,2}) \cdots (u+a_{i,p_i})
$$
for each $i=1,\dots,n$. 
Theorem~\ref{vermathm} then shows that the vector $v_+\in M(A)$
is a non-zero highest weight vector of type $A$.
Moreover, 
$M(A)$ is admissible and, as in $\S$\ref{sshwm},
it has a unique maximal submodule
denoted $\rad M(A)$. The quotient
\begin{equation}\label{lad}
L(A) := M(A) / \rad M(A) \cong
 W(\pi) \otimes_{Y_n(\sigma)} L(\sigma,A(u))
\end{equation}
is the unique (up to isomorphism) irreducible highest weight 
module of type $A$. The modules
$\{L(A)\:|\:A \in \Row(\pi)\}$
give a complete set of pairwise non-isomorphic irreducible admissible
representations of the algebra $W(\pi)$.

Let us describe in detail the situation when
the pyramid $\pi$ consists of a single column of height $m \leq n$.
In this case we have simply that 
$W(\pi) = U(\mathfrak{gl}_m)$ according to the definition (\ref{origdef}).
Let $A$ be a $\pi$-tableau with entries
$a_1,\dots,a_m \in \C$ read from top to bottom.
A highest weight vector for
$W(\pi)$ of type $A$ means a vector $v_+$ with the properties
\begin{enumerate}
\item[(a)] $e_{i,j} v_+ = 0$ for all $1\leq i < j \leq n$;
\item[(b)] $e_{i,i}v_+ = (a_i +n-m+ i-1)v_+$ for all $i=1,\dots,n$.
\end{enumerate}
Hence the module $M(A)$ here coincides with the Verma module
$M(\alpha)$ from (\ref{vermamod}) with
$\alpha = (a_1+n-m, \dots, a_m+n-m)$.

For another example, the trivial $W(\pi)$-module, which we defined earlier
to be the restriction of the trivial $U(\mathfrak{p})$-module, 
is isomorphic to the module
$L(A_0)$ where $A_0$ is the ground-state tableau 
from $\S$\ref{sscrystals}, i.e. the tableau having all entries on its $i$th row equal to $(1-i)$.

\section{Characters of Verma modules}\label{sschars}
By the {character} $\ch M$ of an admissible $W(\pi)$-module $M$,
we mean its Gelfand-Tsetlin character  
when viewed as a $Y_n(\sigma)$-module
via $\kappa:Y_n(\sigma) \twoheadrightarrow W(\pi)$.
Thus $\ch M$ is an element of the completed group algebra
$\widehat\Z[\mathscr P_n]$ from $\S$\ref{sscharblock}.

Given a decomposition $\pi = \pi' \otimes \pi''$ with
$\pi'$ of level $l'$ and $\pi''$ of level $l''$,
the comultiplication $\Delta_{l',l''}$ from (\ref{deltall}) 
allows us to view the tensor product of a
$W(\pi')$-module $M'$ and a $W(\pi'')$-module $M''$
as a $W(\pi)$-module, denoted $M' \boxtimes M''$.
Assuming $M'$ and $M''$ are both admissible,
Lemma~\ref{ten2} and (\ref{trick}) imply
that $M' \boxtimes M''$ is also admissible and
\begin{equation}\label{multy}
\ch (M' \boxtimes M'') = (\ch M') (\ch M'').
\end{equation}
Lemma~\ref{ten1} also carries over in an obvious way to this 
setting.

Introduce the following shorthand for some special elements of 
the completed group algebra $\widehat\Z[\mathscr P_n]$:
\begin{align}
x_{i,a} &:= [1 + (u_i+a+i-1)^{-1}],\label{xia}\\
y_{i,a} &:= [1 + (a+i-1)u_i^{-1}],\label{yia}
\end{align}
for $1 \leq i \leq n$ and $a \in \C$. We note that
\begin{equation}\label{xy}
y_{i,a} / y_{i,a-k} = x_{i,a-1} x_{i,a-2} \cdots x_{i,a-k}
\end{equation}
for any $k \in \N$. The following theorem implies in particular
that the character of any admissible $W(\pi)$-module actually belongs to the
completion of the subalgebra of $\widehat\Z[\mathscr P_n]$
generated just by the elements $\{y_{i,a}^{\pm 1}\:|\:i=1,\dots,n, a \in \C\}$.

\begin{Theorem}\label{nasty}
For $A\in \Row(\pi)$
with entries $a_{i,1},\dots,a_{i,p_i}$ on its $i$th row
for each $i=1,\dots,n$, we have that
$$
\ch M(A) =
\sum_c
\prod_{i=1}^n
\prod_{j=1}^{p_i}
\left\{
y_{i,a_{i,j}-(c_{i,j,i+1}+\cdots+c_{i,j,n})}
\prod_{k=i+1}^n
\frac{y_{k,a_{i,j}-(c_{i,j,k+1}+\cdots+c_{i,j,n})}}{y_{k,a_{i,j}-(c_{i,j,k}+\cdots+c_{i,j,n})}}
\right\}
$$
where the sum is over all tuples 
$c = (c_{i,j,k})_{1 \leq i < k \leq n, 1 \leq j \leq p_i}$
of natural numbers.
\end{Theorem}

The proof of this is more technical than conceptual,
so we postpone it to
$\S$\ref{nastypf}, preferring to illustrate its importance with some
applications first. 

\begin{Corollary}\label{tee}
Let $A_1,\dots,A_l$ be the columns of any representative of 
the row-symmetrized $\pi$-tableau 
$A \in \Row(\pi)$, so that 
$A \sim_{\ro} A_1 \otimes\cdots\otimes A_l$.
Then
$$
\ch M(A) = (\ch M(A_1))\times \cdots\times (\ch M(A_l))
=
\ch (M(A_1)\boxtimes\cdots\boxtimes M(A_l)).
$$
\end{Corollary}

\begin{proof}
This follows from the theorem on interchanging the first two products
on the right hand side.
\end{proof}

In order to derive the next corollary we need to explain 
an alternative way of managing the combinatorics in
Theorem~\ref{nasty}. Continue with
$A \in \Row(\pi)$ 
with entries $a_{i,1},\dots,a_{i,p_i}$ on its $i$th row
as in the statement of the theorem.
By a {\em tabloid} we mean an array
$t = (t_{i,j,a})_{1 \leq i \leq n, 1 \leq j \leq p_i, a < a_{i,j}}$ of integers from the set $\{1,\dots,n\}$
such that
\begin{enumerate}
\item[(a)] $\cdots \leq t_{i,j,a_{i,j}-3} \leq t_{i,j,a_{i,j}-2}\leq t_{i,j,a_{i,j}-1}$;
\item[(b)] $t_{i,j,a} = i$ for $a \ll a_{i,j}$;
\end{enumerate}
for each $1 \leq i \leq n$, $1 \leq j \leq p_i$.
Draw a diagram with rows parametrized by pairs
$(i,j)$ for $1 \leq i \leq n, 1 \leq j \leq p_i$ such that the
$(i,j)$th row consists of a strip of infinitely many boxes,
one in each of the columns parametrized by the numbers
$\dots, a_{i,j}-3, a_{i,j}-2, a_{i,j}-1$.
Then the tabloid $t$ can be recorded on the diagram 
by writing the number $t_{i,j,a}$ into the box in the
$a$th column of the
$(i,j)$th row. In this way tabloids can be thought of as fillings
of the boxes of the diagram by integers from the set $\{1,\dots,n\}$
so that the
entries on each row 
are weakly increasing and 
all but finitely many entries on row $(i,j)$ are equal to $i$.

Given a tabloid $t= (t_{i,j,a})_{1 \leq i \leq n, 1 \leq j \leq p_i, a < a_{i,j}}$, define
$c = (c_{i,j,k})_{1 \leq i < k \leq n, 1 \leq j \leq p_i}$ by
declaring that $c_{i,j,k} = \#\{a < a_{i,j}\:|\:t_{i,j,a} = k\}$, 
i.e. $c_{i,j,k}$ counts the number of
entries equal to $k$ appearing in the $(i,j)$th row of the tabloid $t$.
In this way we obtain a bijection $t \mapsto c$ from the set of all 
tabloids to the set of all tuples
of natural numbers as in the statement of Theorem~\ref{nasty}.
Moreover, for $t$ corresponding to $c$ via this bijection, the identity
(\ref{xy}) implies that
\begin{multline*}
\prod_{i=1}^n
\prod_{j=1}^{p_i}
\left\{
y_{i,a_{i,j}-(c_{i,j,i+1}+\cdots+c_{i,j,n})}
\prod_{k=i+1}^n
\frac{y_{k,a_{i,j}-(c_{i,j,k+1}+\cdots+c_{i,j,n})}}{y_{k,a_{i,j}-(c_{i,j,k}+\cdots+c_{i,j,n})}}
\right\}
\\= 
\prod_{i=1}^n \prod_{j=1}^{p_i} \prod_{a < a_i} x_{t_{i,j,a}, a},
\end{multline*}
where the
infinite product on the right hand side is interpreted using
the convention that
$x_{i,a-1} x_{i,a-2} \cdots = y_{i,a}$
for any $i=1,\dots,n$ and $a \in \C$.
Now we can restate Theorem~\ref{nasty}:
\begin{equation}\label{better}
\ch M(A) = 
\sum_t
\prod_{i=1}^n \prod_{j=1}^{p_i} \prod_{a < a_i}
x_{t_{i,j,a},a}
\end{equation}
where the first summation is over all tabloids
$t = (t_{i,j,a})_{1 \leq i \leq n, 1 \leq j \leq p_i, a < a_{i,j}}$.

\begin{Corollary}\label{ub}
For any $A \in \Row(\pi)$, all Gelfand-Tsetlin subspaces of
$M(A)$ are of dimension less than or equal
to $p_1! (p_1+p_2)! \cdots (p_1+p_2+\cdots+p_{n-1})!.$
\end{Corollary}

\begin{proof}
Two different tabloids $t$ and $t'$ contribute
the same monomial to the right hand side of (\ref{better})
if and only if
they have the same number of entries equal to $i$ appearing in
column $a$ for each $i=1,\dots,n$ and $a \in \C$.
So, given non-negative integers $k_{i,a}$ for
each $i=1,\dots,n$ and $a \in \C$, we need to show by (\ref{better}) 
that there are at
most $p_1! (p_1+p_2)! \cdots (p_1+\cdots + p_{n-1})!$ 
different tabloids with $k_{i,a}$ entries equal to $i$ in column $a$
for each $i=1,\dots,n$ and $a \in \C$.
Given such a tabloid, all entries in the rows parametrized by
$(n,1),\dots,(n,p_n)$ must equal to $n$, while in every other row there
are only finitely many entries equal to $n$ and 
all these entries must form a connected strip at the end of the row. 
So on removing all the boxes containing the entry $n$ we 
obtain a smaller diagram with rows indexed by pairs
$(i,j)$ for $i=1,\dots,n-1, j=1,\dots,p_i$. By induction there are at most
$p_1! (p_1+p_2)! \cdots (p_1+\cdots+p_{n-2})!$ admissible ways of filling the
boxes of this smaller diagram with $k_{i,a}$ entries equal to $i$ in column
$a$ for each $i=1,\dots,n-1$ and $a \in \C$. Therefore we just need to 
show that there are at most $(p_1+\cdots+p_{n-1})!$ admissible ways of 
inserting 
$k_{n,a}$ entries equal to $n$ into column $a$ for each $a \in \C$.
This follows from the following claim:

\vspace{2mm}
{\em \noindent
Suppose we are given $a_1,\dots,a_N \in \C$ and
non-negative integers $k_a$ for each $a \in \C$,
all but finitely many of which are zero. 
Draw a diagram with rows numbered $1,\dots,N$ such that
the $i$th row consists of an infinite strip of boxes, one in each of the
columns parametrized by $\dots,a_i-3,a_i-2,a_i-1$.
Then there are at most $N!$ different ways of deleting boxes from the 
ends of each row in such a way that a total of $k_a$ boxes are 
removed from column $a$ for each $a \in \C$.}

\vspace{2mm}
\noindent
This may be proved by  reducing first to the case that all $a_i$ belong to the same coset of $\C$ modulo $\Z$, then to the case that all $a_i$ are equal.
After these reductions 
it follows from the obvious fact that there are at most $N!$
different $N$-part compositions with prescribed transpose partition.
\end{proof}

\begin{Remark}
On analyzing the proof of the corollary more carefully,
one sees that this upper bound
$p_1! (p_1+p_2)! \cdots (p_1+\cdots+p_{n-1})!$ 
for the dimensions of the Gelfand-Tsetlin
subspaces of $M(A)$ is attained if and only if
all entries in the first $(n-1)$ rows 
of the tableau $A$ belong to the same coset of $\C$ modulo $\Z$.
At the other extreme, all Gelfand-Tsetlin subspaces of $M(A)$ are one 
dimensional if and only if all entries in the first $(n-1)$ rows
of the tableau $A$ belong to different cosets 
of $\C$ modulo $\Z$.
\end{Remark}

\section{The linkage principle}\label{sslink}
Our next application of Theorem~\ref{nasty} is to prove a 
``linkage principle'' showing that the row ordering from (\ref{bruhdef}) 
controls the types of composition factors that can occur in a 
generalized Verma module. 
In the special case that $\pi$ consists of a single column of height $n$, 
i.e. $W(\pi) = U(\mathfrak{gl}_n)$, this result is
\cite[Theorem A1]{BGG2}; even in this case the proof given here is quite
different.

\begin{Lemma}\label{alink}
Suppose $A \downarrow B$.
Then $\ch M(A) = \ch M(B) + (*)$ where $(*)$ is the character
of some admissible $W(\pi)$-module.
\end{Lemma}

\begin{proof}
In view of Corollary~\ref{tee}, it suffices prove this in the
special case that $\pi$ consists of a single column,
i.e. $W(\pi) = U(\mathfrak{gl}_m)$ for some $m$. 
But in that case it is well known
that $A \downarrow B$ implies that there is an embedding
$M(B) \hookrightarrow M(A)$;
see \cite{BGGold} or 
\cite[Lemma 7.6.13]{Dix}.
\end{proof}

\begin{Theorem}\label{bruhat}
Let $A, B \in \Row(\pi)$ with entries
$a_{i,1},\dots,a_{i,p_i}$ and
$b_{i,1},\dots,b_{i,p_i}$ on their $i$th rows, respectively.
The following are equivalent:
\begin{enumerate}
\item[\rm(i)] $A \geq B$;
\item[\rm(ii)] $[M(A):L(B)] \neq 0$;
\item[\rm(iii)] there exists a  tuple 
$c = (c_{i,j,k})_{1 \leq i < k \leq n, 1 \leq j \leq p_i}$
of natural numbers such that
\begin{equation*}
\prod_{i=1}^n
\prod_{j=1}^{p_i}
y_{i,b_{i,j}} =
\prod_{i=1}^n
\prod_{j=1}^{p_i}
\left\{
y_{i,a_{i,j}-(c_{i,j,i+1}+\cdots+c_{i,j,n})}
\prod_{k=i+1}^n
\frac{y_{k,a_{i,j}-(c_{i,j,k+1}+\cdots+c_{i,j,n})}}{y_{k,a_{i,j}-(c_{i,j,k}+\cdots+c_{i,j,n})}}
\right\}.
\end{equation*}
\end{enumerate}
\end{Theorem}

\begin{proof}
(i)$\Rightarrow$(ii).
If $A \geq B$ then Lemma~\ref{alink}
implies that $$
\ch M(A) = \ch M(B) + (*),
$$ 
where $(*)$ is the character
of some admissible $W(\pi)$-module.
Hence we get that $[M(A):L(B)] \geq [M(B):L(B)] = 1$.

(ii)$\Rightarrow$(iii).
Suppose that $[M(A):L(B)] \neq 0$.
The highest weight vector $v_+$ of $L(B)$
contributes 
$\prod_{i=1}^n \prod_{j=1}^{p_i} y_{i,b_{i,j}}$
to the formal character $\ch L(B)$.
Hence, by Corollary~\ref{span}, we see that 
$\ch M(A)$ also involves
$\prod_{i=1}^n \prod_{j=1}^{p_i} y_{i,b_{i,j}}$
with non-zero coefficient. In view of Theorem~\ref{nasty}
this implies (iii).

(iii)$\Rightarrow$(i).
Suppose that
\begin{equation*}
\prod_{i=1}^n
\prod_{j=1}^{p_i}
y_{i,b_{i,j}} =
\prod_{i=1}^n
\prod_{j=1}^{p_i}
\left\{
y_{i,a_{i,j}-(c_{i,j,i+1}+\cdots+c_{i,j,n})}
\prod_{k=i+1}^n
\frac{y_{k,a_{i,j}-(c_{i,j,k+1}+\cdots+c_{i,j,n})}}{y_{k,a_{i,j}-(c_{i,j,k}+\cdots+c_{i,j,n})}}
\right\}
\end{equation*}
for some tuple 
$c = (c_{i,j,k})_{1 \leq i < k \leq n, 1 \leq j \leq p_i}$
of natural numbers.
We show by induction on $\sum c_{i,j,k}$ that $A \geq B$.
If $\sum c_{i,j,k} = 0$ this is trivial since then $A \sim_{\ro} B$.
Otherwise, let $i_2$ be maximal such that
$c_{i,j,i_2} \neq 0$ for some $1 \leq i < i_2$
and $1 \leq j \leq p_i$. Considering the $y_{i_2,?}$'s on either side of
our equation gives that
$$
\prod_{j=1}^{p_{i_2}} y_{i_2,b_{i_2,j}}
=
\prod_{j=1}^{p_{i_2}} y_{i_2,a_{i_2,j}}
\times
\prod_{i=1}^{i_2-1} \prod_{j=1}^{p_i}
\frac{y_{i_2,a_{i,j}}}{y_{i_2,a_{i,j}-c_{i,j,i_2}}}.
$$
Hence there exist $1 \leq i_1 < i_2$,
$1 \leq j_1 \leq p_{i_1}$ and $1 \leq j_2 \leq p_{i_2}$
such that $a_{i_2,j_2} = a_{i_1,j_1} - c_{i_1,j_1,i_2} \neq a_{i_1,j_1}$.
Let $\bar A = (\bar a_{i,j})_{1 \leq i \leq n, 1\leq j \leq p_i}$ be the
$\pi$-tableau obtained from $A$ by swapping the entries
$a_{i_1,j_1}$ and $a_{i_2,j_2}$.
Define a new tuple 
$(\bar c_{i,j,k})_{1 \leq i < j \leq n, 1 \leq j \leq p_i}$
from
$$
\bar c_{i,j,k} = \left\{
\begin{array}{ll}
c_{i,j,k}&\hbox{if $(i,j,k) \neq (i_1,j_1,i_2)$},\\
0&\hbox{if $(i,j,k) = (i_1,j_1,i_2)$}.
\end{array}\right.
$$
Now using the maximality of the choice of $i_2$, one checks that
\begin{multline*}
\prod_{i=1}^n
\prod_{j=1}^{p_i}
\left\{
y_{i,\bar a_{i,j}-(\bar c_{i,j,i+1}+\cdots+\bar c_{i,j,n})}
\prod_{k=i+1}^n
\frac{y_{k,\bar a_{i,j}-(\bar c_{i,j,k+1}+\cdots+\bar c_{i,j,n})}}{y_{k,\bar a_{i,j}-
(\bar c_{i,j,k}+\cdots+\bar c_{i,j,n})}}
\right\}
=\\
\prod_{i=1}^n
\prod_{j=1}^{p_i}
\left\{
y_{i,a_{i,j}-(c_{i,j,i+1}+\cdots+c_{i,j,n})}
\prod_{k=i+1}^n
\frac{y_{k,a_{i,j}-(c_{i,j,k+1}+\cdots+c_{i,j,n})}}{y_{k,a_{i,j}-(c_{i,j,k}+\cdots+c_{i,j,n})}}
\right\}
=
\prod_{i=1}^n
\prod_{j=1}^{p_i}
y_{i,b_{i,j}}.
\end{multline*}
Since $\sum \bar c_{i,j,k} < \sum c_{i,j,k}$ we deduce by induction that
$\bar A \geq B$. Since $A \downarrow \bar A$ this completes the proof.
\end{proof}

\begin{Corollary}\label{gen}
For $A \in \Row(\pi)$
with entries $a_{i,1},\dots,a_{i,p_i}$ on its $i$th row, 
the following are equivalent:
\begin{enumerate}
\item[\rm(i)] $M(A)$ is irreducible;
\item[\rm(ii)]
$A$ is minimal with respect to the ordering $\geq$;
\item[\rm(iii)] $a_{i_1,j_1}\not > a_{i_2,j_2}$ for every
$1 \leq i_1 < i_2 \leq n$,
$1 \leq j_1 \leq p_{i_1}$ and $1 \leq j_2 \leq p_{i_2}$.
\end{enumerate}
Moreover, assuming (i)--(iii) hold, 
let $A_1,\dots,A_l$ be the columns of any representative of $A$
read from left to right,
so that $A \sim_{\ro} A_1\otimes\cdots\otimes A_l$.
Then we have that
$$
M(A) \cong M(A_1) \boxtimes\cdots\boxtimes M(A_l).
$$
\end{Corollary}

\begin{proof}
The equivalence of (i) and (ii) follows from Theorem~\ref{bruhat}.
The equivalence of (ii) and (iii) is clear from the definition 
of the Bruhat ordering.
The final statement follows from
Corollary~\ref{tee} and Theorem~\ref{inj}.
\end{proof}

\section{The center of $W(\pi)$}\label{sscenter}
Our final application of Theorem~\ref{nasty} is to 
prove that the center $Z(W(\pi))$ is a polynomial
algebra on generators $\psi(Z_N^{(1)}),\dots,\psi(Z_N^{(N)})$,
notation as in $\S$\ref{sscent}.
In the case that $\pi$ is an $n \times l$ rectangle,
when $W(\pi)$ is the Yangian of level $l$,
this result is due to Cherednik \cite{Ch0,Ch}; see 
also \cite[Corollary 4.1]{Molevcas}. 
For the first lemma, we point out that 
the usual Verma modules for the Lie algebra
$\mathfrak{h}$ 
are precisely the outer tensor product modules
$M(A_1) \boxtimes\cdots\boxtimes M(A_l)$
for $A \in \Tab(\pi)$ with columns $A_1,\dots,A_l$.
Moreover, 
if $\gamma(A) = (a_1,\dots,a_N)$ then
$M(A_1)\boxtimes\cdots\boxtimes M(A_l)$
is of usual highest weight
$(a_1+\row(1)-1,\dots,a_N+\row(N)-1) \in \mathfrak{d}^*$.
Recall also the definition of the Miura transform
$\xi$
from (\ref{miura}). 

\begin{Lemma}\label{semicent}
$\xi(Z(W(\pi))) \subseteq Z(U(\mathfrak{h}))$.
\end{Lemma}

\begin{proof}
Take $z \in Z(W(\pi))$ and $u \in U(\mathfrak{h})$.
We need to show that $[\xi(z),u] = 0$.
This follows by \cite[Theorem 8.4.4]{Dix} 
as soon as we check that $[\xi(z),u]$ annihilates 
$M(A_1)\boxtimes \cdots \boxtimes M(A_l)$
for generic $A \in \Tab(\pi)$ with columns $A_1,\dots,A_l$,
Corollary~\ref{gen} shows that 
$M(A_1)\boxtimes\cdots\boxtimes M(A_l)$
is generically irreducible when viewed as a
$W(\pi)$-module via $\xi$. Hence 
$\xi(z)$ acts on it as a scalar
by (\ref{oned}). So certainly $[\xi(z),u]$ acts as zero.
\end{proof}

\begin{Theorem}\label{fiish}
The map 
$\psi:Z(U(\mathfrak{gl}_N)) \rightarrow
Z(W(\pi))$ from (\ref{psidef}) is an isomorphism.
Hence, the elements $\psi(Z_N^{(1)}),
\dots, \psi(Z_N^{(N)})$ are algebraically
independent and generate the center $Z(W(\pi))$.
\end{Theorem}

\begin{proof}
In view of Lemma~\ref{semicent} and the
commutativity of the diagram (\ref{hcfact}), 
we just need to show that the image
of $z \in Z(W(\pi))$ under
$(\Psi_{\!q_1} \otimes\cdots\otimes \Psi_{\!q_l})\circ \xi$
is a symmetric polynomial in 
$e_{1,1}+q_{\col(1)}-n,\dots,e_{N,N}+q_{\col(N)}-n$. Equivalently,
by the definition of the Harish-Chandra homomorphism, 
we need to show, whenever $A, B$ are $\pi$-tableaux
with the same content,
that the element $z$ acts on the modules
$M(A_1) \boxtimes \cdots\boxtimes M(A_l)$ and
$M(B_1) \boxtimes \cdots \boxtimes M(B_l)$ by the same scalar,
where $A_i$ resp. $B_i$ denotes the $i$th column of $A$ resp. $B$.
If $B$ is obtained from $A$ by permuting entries within columns,
this is immediate from Lemma~\ref{semicent}.
If $B$ is obtained from $A$ by permuting enties within rows,
it follows from Theorem~\ref{inj} and
Corollary~\ref{tee}.
The general case follows from these two special situations.
\end{proof}

We remark that 
there is now a quite different proof of this theorem, valid 
for finite $W$-algebras associated to arbitrary 
finite dimensional semisimple Lie algebras, 
due to Ginzburg. For a sketch of the argument,
see the footnote to \cite[Question 5.1]{premet2}.

\begin{Corollary}\label{centcor}
The elements $C_n^{(1)},C_n^{(2)},\dots$
of $Y_n(\sigma)$ 
are algebraically independent and generate
the center $Z(Y_n(\sigma))$.
Moreover, 
$\kappa:Y_n(\sigma) \twoheadrightarrow W(\pi)$ 
maps $Z(Y_n(\sigma))$ surjectively onto $Z(W(\pi))$.
\end{Corollary}

\begin{proof}
This is immediate from the theorem on recalling that
$Y_n(\sigma)$ is a filtered inverse limit of $W(\pi)$'s
as explained in \cite[Remark~6.4]{BK}.
\end{proof}

We are grateful to one of the referees of \cite{BK} for pointing
out that we are already in a position to apply \cite{FO} to obtain
the following generalization of a theorem of Kostant from \cite{K1}.
In the case $W(\pi)$ is the Yangian of level $l$ this result
is \cite[Theorem 2]{FO}.

\begin{Theorem}
The algebra $W(\pi)$ is free as a module over its center.
\end{Theorem}

\begin{proof}
Recall that the associated graded algebra
$\gr W(\pi)$ is free commutative on generators
(\ref{grb1})--(\ref{grb3}), in particular $W(\pi)$ is a 
special filtered algebra in the sense of \cite{FO}.
Let $A$ be the quotient of $\gr W(\pi)$
by the ideal generated by the elements (\ref{grb2})--(\ref{grb3}).
Let $d_i^{(r)}$ resp. 
$c_n^{(r)}$ 
denote the image of $\gr_r D_i^{(r)}$ resp. $\gr_r C_n^{(r)}$ in $A$.
Thus, $A$ is the free polynomial algebra
$\C[d_i^{(r)}\:|\:i=1,\dots,n, r=1,\dots,p_i]$. Moreover
by Theorem~\ref{tvan} and (\ref{tij}) we have that
$d_i^{(r)} = 0$ for $r > p_i$. It follows from this and
(\ref{centelts}) that if we set
\begin{align*}
d_i(u) &= \sum_{r=0}^{p_i} d_i^{(r)} u^{p_i-r},\\
c_n(u) &= \sum_{r=0}^{N} c_n^{(r)} u^{N-r}
\end{align*}
then
$c_n(u) = d_1(u) d_2(u) \cdots d_n(u)$.
Now applying \cite[Theorem 1]{FO} as in the proof of
\cite[Theorem 2]{FO}, it suffices to show that
$c_n^{(1)},c_n^{(2)},\dots,c_n^{(N)}$ is a regular sequence
in $A$, i.e. that the image of
$c_n^{(r)}$ in $A / (A c_n^{(1)} + \cdots +A c_n^{(r-1)})$
is not invertible and not a zero divisor for each $r=1,\dots,N$.
For this, by \cite[Proposition 1(5)]{FO}, we just need to check that
the variety $Z = V(c_n^{(1)},\dots,c_n^{(N)})$ is equidimensional
of dimension $0$. Consider the morphism
$\varphi:\C^N \rightarrow \C^N$ mapping a point
$(x_i^{(r)})_{1 \leq i \leq n, 1 \leq r \leq p_i}$
to the coefficients of the following monic polynomial:
$$
\prod_{i=1}^n (u^{p_i} + d_i^{(1)} u^{p_i-1} + \cdots + d_i^{(p_i)}).
$$
Obviously $Z = \varphi^{-1}(0)$. Since $\C[u]$ is a unique factorization domain, $u^N = u^{p_1} \cdots u^{p_n}$ is the unique decomposition of $u^{N}$
as a product of monic polynomials of degrees $p_1,\dots,p_n$.
Hence $Z = \{0\}$.
\end{proof}

In view of Theorem~\ref{fiish}, the
center of $W(\pi)$ is canonically isomorphic to the
center of $U(\mathfrak{g})$.
So we can parametrize the central characters of $W(\pi)$
in exactly the same way as we did for $U(\mathfrak{g})$
in $\S$\ref{sscent},
by the set of $\theta \in P = \bigoplus_{a \in \C} \Z \gamma_a$ 
whose coefficients are non-negative integers
summing to $N$.
Given such an element
$\theta$,
define 
$f(u) = u^N + f^{(1)} u^{N-1}+\cdots+f^{(N)} \in \C[u]$
according to (\ref{pe1})--(\ref{pe2}).
Then, for an admissible $W(\pi)$-module $M$, 
define
\begin{equation}\label{prtheta2}
\pr_\theta(M) := 
\left\{v \in M\:\bigg|\:
\begin{array}{ll}\hbox{for each 
$r=1,\dots,N$ there exists $p > 0$}\\
\hbox{such that
$\psi(Z_N^{(r)} - f^{(r)})^p v = 0$}\end{array}\right\}.
\end{equation}
Equivalently, by (\ref{centelts}) and Lemma~\ref{psieq}, we have that
\begin{equation}\label{cw}
\pr_\theta(M) = \bigoplus_{A(u)} M_{A(u)}
\end{equation}
where the direct sum is over all $A(u) \in \mathscr P_n$
such that
$$
u^{p_1} (u-1)^{p_2} \cdots (u-n+1)^{p_n}
A_1(u) A_2(u-1) \cdots A_n(u-n+1) = f(u).
$$
Since the admissible $W(\pi)$-module 
$M$ is the direct sum of its Gelfand-Tsetlin subspaces, 
it follows that
\begin{equation}\label{blockdec}
M = \bigoplus_{\theta \in P} \pr_\theta(M),
\end{equation}
with the convention that $\pr_\theta(M)= 0$ if 
the coefficients of $\theta$ are not non-negative integers summing to $N$.
This is clearly a decomposition of $M$ as a $W(\pi)$-module.

\begin{Lemma}\label{cx}
All highest weight $W(\pi)$-modules of type $A \in \Row(\pi)$ 
are of central character $\theta(A)$.
\end{Lemma}

\begin{proof}
Suppose that the entries on the $i$th row of $A$ 
are $a_{i,1},\dots,a_{i,p_i}$.
By (\ref{centelts}), Lemma~\ref{psieq} and the definition
(\ref{t1})--(\ref{t3}), $\psi(Z_N(u))$ acts on any highest weight
module of type $A$
as the scalar $\prod_{i=1}^n \prod_{j=1}^{p_i} (u+a_{i,j})$.
\end{proof}

\section{Proof of Theorem~6.2}\label{nastypf}

Let $\bar\pi$ denote the pyramid 
obtained from  $\pi$
by removing the bottom row.
The tuple of row lengths corresponding to the pyramid $\bar\pi$ 
is $(p_1,\dots,p_{n-1})$ and the submatrix $\bar\sigma = 
(s_{i,j})_{1 \leq i,j \leq n-1}$ of the
shift matrix $\sigma = (s_{i,j})_{1 \leq i,j \leq n}$ 
chosen for $\pi$ gives a shift matrix for $\bar\pi$.
By the relations,
 there is a homomorphism
$W(\bar \pi) \rightarrow W(\pi)$ mapping the generators
$D_i^{(r)} \:(i=1,\dots,n-1, r > 0)$,
$E_i^{(r)} \:(i=1,\dots,n-2, r > s_{i,i+1})$
and
$F_i^{(r)} \:(i=1,\dots,n-2, r > s_{i+1,i})$
of $W(\bar\pi)$ to the elements with the same names
in $W(\pi)$. By the PBW theorem this map is in fact injective,
allowing us to view $W(\bar\pi)$ as a subalgebra of $W(\pi)$.
We will in fact prove the following branching theorem for generalized 
Verma modules.

\begin{Theorem}\label{nastier}
Let $A
\in \Row(\pi)$ with entries
$a_{i,1},\dots,a_{i,p_i}$ on its $i$th row for each $i=1,\dots,n$.
There is a filtration 
$0 = M_0 \subset M_1 \subset \cdots$ of $M(A)$
as a $W(\bar\pi)$-module with $\bigcup_{i \geq 0} M_i = M(A)$
and subquotients
isomorphic to the generalized Verma modules $M(B)$ for
$B  \in \Row(\bar\pi)$ such that $B$ has the entries
$(a_{i,1}-c_{i,1}),\dots,(a_{i,p_i}-c_{i,p_i})$ on its $i$th row
for each $i=1,\dots,n-1$, one for each
tuple $(c_{i,j})_{1 \leq i \leq n-1, 1 \leq j \leq p_i}$
of natural numbers.
\end{Theorem}

Let us first explain how to deduce Theorem~\ref{nasty} from this.
Proceed by induction on $n$, the case $n=1$ being trivial.
For the induction step, we have by Theorem~\ref{nastier} and the 
induction hypothesis that
the character of
$\res^{W(\pi)}_{W(\bar\pi)} M(A)$ equals
$$
\sum_c \prod_{i=1}^{n-1} \prod_{j=1}^{p_i}
\left\{y_{i,a_{i,j}-(c_{i,j,i+1}+\cdots+c_{i,j,n})}
\prod_{k=i+1}^{n-1}
\frac{y_{k,a_{i,j}-(c_{i,j,k+1}+\cdots+c_{i,j,n})}}{y_{k,a_{i,j}-(c_{i,j,k}+\cdots+c_{i,j,n})}}\right\},
$$
where the first sum is over all tuples
$c = (c_{i,j,k})_{1 \leq i < k \leq n-1, 1 \leq j \leq p_i}$ of
natural numbers.
But just like in the proof of Lemma~\ref{cx}, 
$$u^{p_1}(u-1)^{p_2} \cdots (u-n+1)^{p_n} D_1(u) D_2(u-1)
\cdots D_n(u-n+1)$$ 
acts on $M(A)$ as the scalar
$\prod_{i=1}^n \prod_{j=1}^{p_i} (u+a_{i,j})$.
Hence recalling (\ref{yia}), 
each monomial appearing in the expansion of $\ch M(A)$
must simplify to
$\prod_{i=1}^n \prod_{j=1}^{p_i} (u+a_{i,j})$
on replacing $y_{i,a}$ by $(u+a)$ everywhere.
In this way we can recover $\ch M(A)$ uniquely from the above
expression to complete the proof of  Theorem~\ref{nasty}.

To prove Theorem~\ref{nastier}, 
{\em we will assume from now on that the shift matrix $\sigma$ is 
 upper triangular}; the result in general
then follows easily by twisting with the 
isomorphism $\iota$ from (\ref{iotadef}).
Exploiting this assumption, the following lemma can be checked using the
formulae in \cite[$\S$5]{BKdrinfeld} and some elementary inductive arguments.

\begin{Lemma}\label{rels}
The following relations hold in $W(\pi)$.
\begin{enumerate}
\item[\rm(i)]
For all $i < j$,
$[F_{i,j}(u) D_i(u), F_{i,j}(v)D_i(v)] = 0$.
\item[\rm(ii)]
For all $i < j < k$,
$(u-v) [F_{j,k}(u), F_{i,j}(v)]$ equals
$$
\sum_{r \geq 0} (-1)^r 
\sum_{\substack{i < i_1 < \cdots < i_r \leq j \\ i_{r+1}=k}}
F_{i_r, i_{r+1}}(u) \cdots F_{i_1,i_2}(u) (F_{i,i_1}(v)-F_{i,i_1}(u)).
$$
\item[\rm(iii)]
For all $i < j$ and $k < i$ or $k > j$,
$[D_k(u), F_{i,j}(v)]=0$.
\item[\rm(iv)]
For all $i < j$,
$(u-v) [D_i(u), F_{i,j}(v)]=(F_{i,j}(u) - F_{i,j}(v)) D_i(u)$.
\item[\rm(v)]
For all $i < j$, $(u-v) [D_j(u), F_{i,j}(v)]$ equals
$$
\sum_{r \geq 0} (-1)^r \sum_{\substack{i < i_1 < \cdots < i_r < j \\
i_{r+1}=j}}
F_{i_r, i_{r+1}}(u) \cdots F_{i_1,i_2}(u) (F_{i,i_1}(v)-F_{i,i_1}(u))D_j(u).
$$
\item[\rm(vi)]
For all $i < j < k$,
$(u-v) [D_j(u), F_{i,k}(v)]$ equals
$$
\sum_{r \geq 0} (-1)^r \sum_{i < i_1 < \cdots < i_{r+1} = j}
F_{i_r, i_{r+1}}(u) \cdots F_{i_1,i_2}(u) (F_{i,i_1}(v)-F_{i,i_1}(u))
F_{j,k}(u)D_j(u).
$$
\item[\rm(vii)]
For all $i < j < k$,
$(u-v) [F_{j,k}(u) D_j(u), F_{i,k}(v)]$ equals
$$
\sum_{r \geq 0} (-1)^r
\sum_{\substack{i < i_1 < \cdots < i_r < j \\ i_{r+1}=k}}
F_{i_r,i_{r+1}}(u) \cdots F_{i_1,i_2}(u) (F_{i,i_1}(v) - F_{i,i_1}(u))
F_{j,k}(u)D_j(u).
$$
\end{enumerate}
\end{Lemma}

Recalling Theorem~\ref{tvan}, introduce  the shorthand
\begin{equation}
L_i(u) = \sum_{r=0}^{p_i} L_i^{(r)} u^{p_i-r}
:= u^{p_i} T_{n,i}(u) \in W(\pi)[u]
\end{equation}
for each $1 \leq i < n$.
Also for $h \geq 0$ set
\begin{equation}
L_{i,h}(u) := \frac{1}{h!} \frac{d^h}{du^h} L_i(u).
\end{equation}
We will apply the following simple observation repeatedly from now on:
given a vector $m$ of generalized weight $\alpha$ 
in a $W(\pi)$-module $M$ with the property that 
$\alpha + \eps_j - \eps_i$
 is not a weight of $M$ for any
$1 \leq j < i$, we have 
by (\ref{tij}) that
$L_i(u) m = u^{p_i} F_{i,n}(u) D_i(u) m$. 

\begin{Lemma}\label{hard}
Suppose we are given $1 \leq i < n$ and a vector
$m$ of generalized weight $\alpha$ in a $W(\pi)$-module $M$
such that 
\begin{enumerate}
\item[\rm(i)]
$\alpha - d (\eps_i - \eps_n) + \eps_j - \eps_i$ is not a weight of $M$
for any $1 \leq j < i$ and $d \geq 0$;
\item[\rm(ii)]
$u^{p_i} D_i(u) m \equiv (u+a_1) \cdots (u+a_{p_i}) m
\pmod{M'[u]}$ for some $a_1,\dots,a_{p_i} \in \C$
and some subspace $M'$ of $M$.
\end{enumerate}
For $j = 1,\dots,p_i$, 
define
$m_j := L_{i,h(j)}(-a_j) m$
where $$
h(j) = \#\{k=1,\dots,j-1\:|\:a_k = a_j\}.
$$
Then we have that
\begin{multline*}
u^{p_i} D_i(u) m_j \equiv (u+a_1) \cdots (u+a_{j-1}) (u+a_j-1)
(u+a_{j+1}) \cdots (u+a_{p_i})m_j
\\- \sum_{\substack{k=1,\dots,j-1\\ a_k = a_j}}
\frac{(u+a_1) \cdots (u+a_{p_i})}{(u+a_j)^{h(j)-h(k)+1}} m_k
\pmod{\sum_{r=1}^{p_i} L_i^{(r)} M'[u]}.
\end{multline*}
Moreover, the subspace of $M$ spanned by the vectors
$m_1,\dots,m_{p_i}$ coincides with the subspace spanned by the vectors
$L_i^{(1)} m, \dots, L_i^{(p_i)} m$.
\end{Lemma}

\begin{proof}
By Lemma~\ref{rels}(iv) and the assumptions (i)--(ii), we have that
\begin{multline*}
(u-v) [u^{p_i} D_i(u), L_i(v)] m
\equiv
(v+a_1) \cdots (v+a_{p_i}) L_i(u) m\\
-
(u+a_1) \cdots (u+a_{p_i}) L_i(v) m
\pmod{\sum_{r=1}^{p_i} L_i^{(r)} M'[u,v]}.
\end{multline*}
Hence,
\begin{multline*}
u^{p_i} D_i(u) L_i(v) m
\equiv
(u+a_1) \cdots (u+a_{p_i}) L_i(v) m
-
\frac{(u+a_1) \cdots (u+a_{p_i}) L_i(v)}{u-v} m
\\+ \frac{
(v+a_1) \cdots (v+a_{p_i}) L_i(u)}{u-v} m.
\end{multline*}
Apply the operator
$\frac{1}{h(j)!} \frac{{d}^{h(j)}}{{d}v^{h(j)}}$ 
to both sides
using the
Leibniz rule then set $v := -a_j$ to deduce that
\begin{multline*}
u^{p_i} D_i(u) L_{i,h(j)}(-a_j) m
\equiv
(u+a_1)\cdots
(u+a_{p_i}) L_{i,h(j)}(-a_j) m
\\-
\sum_{k=0}^{h(j)}
\frac{(u+a_1)\cdots(u+a_{p_i})}{(u+a_j)^{h(j)-k+1}}
L_{i,k}(-a_j) m.
\end{multline*}
The left hand side equals $u^{p_i} D_i(u) m_j$ by definition.
The right hand side simplifies to give
$$
(u+a_1)\cdots
(u+a_j-1)\cdots
(u+a_{p_i}) m_j
-
\sum_{k=0}^{h(j)-1}
\frac{(u+a_1)\cdots(u+a_{p_i})}{(u+a_j)^{h(j)-k+1}}
L_{i,k}(-a_j) m
$$
which is exactly what we need to prove the first part of the lemma.

For the second part, we observe that the transition matrix between
the vectors $L_{i}(u_1)m,\cdots, L_i(u_{p_i})m$
and 
$L_i^{(1)}m,\cdots,L_i^{(p_i)}m$
is a Vandermonde matrix with determinant 
$\prod_{1 \leq j < k \leq p_i} (u_j-u_k)$.
Apply $\frac{1}{h(j)!} \frac{d^{h(j)}}{du_j^{h(j)}}$
for $j=1,\dots,p_i$ to deduce that the determinant of the transition
matrix between
$L_{i,h(1)}(u_1) m,\cdots, L_{i,h({p_i})}(u_{p_i})m$
and
$L_i^{(1)}m,\cdots,L_i^{(p_i)}m$ is
$$
\frac{1}{h(1)! h(2)! \cdots h(p_i)!}
\frac{d^{h(1)}}{du_1^{h(1)}}
\cdots \frac{d^{h({p_i})}}{du_{p_i}^{h({p_i})}}
\prod_{1 \leq j < k \leq  p_i} (u_j - u_k).
$$
Evaluate this expression at 
$u_j=-a_j$ for each $j=1,\dots,p_i$ to get
$$
(-1)^{h(1)+\cdots+h(p_i)}
\prod_{\substack{1 \leq j < k \leq p_i \\ a_j \neq a_k}} (a_k - a_j)
\neq 0.
$$
Hence the transition matrix 
between the vectors 
$m_1,\dots,m_{p_i}$
and
$L_{i}^{(1)}m,\cdots,L_{i}^{(p_i)}m$ is invertible, so they span the
same space.
\end{proof}

\begin{Lemma}\label{reallyhard}
Under the same assumptions as Lemma~\ref{hard},
let $C_d$ denote the set of all $p_i$-tuples
$c = (c_1,\dots,c_{p_i})$ of natural numbers summing to $d$.
Put a total order on $C_d$ so that
$c' < c$ if $c'$ is lexicographically greater than $c$.
For $c \in C_d$ let
$$
m_c := \prod_{j=1}^{p_i} \prod_{k=1}^{c_j}
L_{i,h_c(j,k)}(-a_j+k-1) m
$$
where $h_c(j,k) = \#\{l=1,\dots,j-1\:|\:a_l - c_l = a_j - k+1\}$.
Then
$$
u^{p_i} D_i(u) m_c \equiv (u+a_1-c_1) \cdots (u+a_{p_i}-c_{p_i})
m_c \pmod{M_c'[u]}
$$
where $M_c'$ is the subspace of $M$ spanned by all the vectors $m_{c'}$
for $c' < c$ and $L_i^{(r_1)} \cdots L_i^{(r_d)} M'$
for $1 \leq r_1,\dots,r_d \leq p_i$.
Moreover,
the vectors $\{m_c\:|\:c \in C_{d}\}$
span the same subspace of $M$ as the vectors
$L_{i}^{(r_1)} \cdots L_{i}^{(r_d)} m$
for all $1 \leq r_1,\dots,r_d \leq p_i$.
\end{Lemma}

\begin{proof}
Note first that the definition of the vectors $m_c$ does not depend
on the order taken in the products, thanks to Lemma~\ref{rels}(i).
Now 
proceed by induction on $d$, the case $d=1$ being precisely the result
of the previous lemma. 
For $d > 1$, 
define vectors $m_1,\dots,m_{p_i}$ 
according to the preceeding lemma.
For $r=1,\dots,p_i$, let $M_r'$ be the subspace spanned by
$m_1,\dots,m_{r-1}$ and $L_{i}^{(s)} M'$ for all $s=1,\dots,p_i$.
Then the preceeding lemma shows that
$$
u^{p_i}D_{i}(u) m_r \equiv (u+a_1) \cdots (u+a_r-1) \cdots (u+a_{p_i})
m_r \pmod{M_r'[u]}
$$
and that $m_1,\dots,m_{p_i}$ span the same space as the vectors
$L_{i}^{(1)} m,\dots,L_{i}^{(p_i)} m$.

For $c \in C_{d-1}$ and $r=1,\dots,p_i$, let
$$
m_{r,c} :=  \prod_{j=1}^{p_i}\prod_{k=1}^{c_j} 
L_{i,h_{r,c}(j,k)}(-a_j+\delta_{j,r}+k-1)m_r
$$
where $h_{r,c}(j,k) := 
\#\{l=1,\dots,j-1\:|\:a_l -\delta_{r,l}- c_l = a_j -\delta_{r,j}- k+1\}$.
Let $M_{r,c}'$ be the subspace of $M$ spanned by 
all $m_{r,c'}$ for $c' < c$ together with
$L_{i}^{(r_1)} \cdots L_{i}^{(r_{d-1})}M_r'$ for
all $1 \leq r_1,\dots,r_{d-1} \leq p_i$.
Then by the induction hypothesis,
$$
u^{p_i} D_{i}(u) m_{r,c} \equiv (u+a_1 - c_1) \cdots
(u+a_r - 1 - c_r) \cdots (u+a_{p_i}-c_{p_i}) m_{r,c} \pmod{M_{r,c}'[u]}.
$$
Moreover the vectors $\{m_{r,c}\:|\:c \in C_{d-1}\}$
span the same subspace of $M$ as the vectors
$L_{i}^{(r_1)} \cdots L_{i}^{(r_{d-1})} m_r$ for all
$1 \leq r_1,\dots,r_{d-1} \leq p_i$.
Now observe that if $c \in C_{d-1}$ 
satisfies $c_1=\dots=c_{r-1} = 0$, then
$m_{r,c} = m_{c+\delta_r}$
where $c+\delta_r \in C_{d}$ 
is the tuple $(c_1,\dots,c_{r-1},c_r+1,c_{r+1},\dots,c_{p_i})$;
otherwise, $m_{r,c}$ lies in the subspace spanned by the
$m_{s,c'}$ for $s < r, c' \in C_{d-1}$.
The lemma follows.
\end{proof}

At last we can complete the proof of Theorem~\ref{nastier}.
Let $C$ denote the set of all 
tuples $c = (c_{i,j})_{1 \leq i \leq n-1, 1 \leq j \leq p_i}$
of natural numbers. 
Writing $|c|_i$ for $\sum_{j=1}^{p_i} c_{i,j}$ and
$|c|$ for $|c|_1 + |c|_2+\cdots+|c|_{n-1}$,
we put a total order on $C$ so that $c' \leq c$ if any
of the following hold:
\begin{enumerate}
\item[(a)] $|c'| < |c|$;
\item[(b)]
$|c'| = |c|$ but
$|c'|_{n-1} = |c|_{n-1}$,$|c'|_{n-2} = |c|_{n-2}$,\dots,
$|c'|_{i+1} = |c|_{i+1}$ and
$|c'|_{i} > |c|_{i}$
for some $i \in \{1,\dots,n-1\}$;
\item[(c)]
$|c'|_i = |c|_i$ but
the tuple $(c_{i,1}',\dots,c_{i,p_i}')$ is lexicographically greater
than or equal to the tuple
$(c_{i,1},\dots,c_{i,p_i})$ 
for every $i=1,\dots,n-1$.
\end{enumerate}
Now let $M := M(A)$ for short.
For each $c \in C$, define a vector $m_c \in M$ by
$$
m_c := 
\prod_{i=1}^{n-1}
\left\{
\prod_{j=1}^{p_i}
\prod_{k=1}^{c_{i,j}}
L_{i,h_c(i,j,k)}(-a_{i,j}+k-1)
\right\} v_+
$$
where $h_c(i,j,k) = \#\{l=1,\dots,j-1\:|\:a_{i,l}-c_{i,l} = a_{i,j} - k+1\}$
and the first product is taken in order of increasing $i$ from left to right.
The second part of Lemma~\ref{reallyhard} and Theorem~\ref{vermathm}(ii) imply that the vectors
$\{y m_{c}\:|\:y \in Y, c \in C\}$
form a basis for $M$, where $Y$ 
here denotes the set of all monomials in the elements $\{T_{j,i}^{(r)}\:|\:1 \leq i < j \leq n-1, r = 1,\dots,p_i\}$.
For each $c \in C$, let $M_c$ resp. $M_c'$ denote the subspace of $M$ spanned by $\{y m_{c'}\:|\:y \in Y, c' \leq c\}$
resp. $\{y m_{c'}\:|\:y \in Y, c' < c\}$.
Clearly $M = \bigcup_{c \in C} M_c$.
Now we complete the proof of Theorem~\ref{nastier} 
by showing 
that each $M_c$ is actually a $W(\bar\pi)$-submodule of $M$
with $M_c / M_c' \cong M(B)$ for
$B \in \Row(\bar\pi)$ such that $B$ has entries
$(a_{i,1}-c_{i,1}),\dots,(a_{i,p_i}-c_{i,p_i})$ on its $i$th row for each $i=1,\dots,n-1$.

Proceeding by induction on the total ordering on $C$,
the induction hypothesis allows us to assume that
$M_c'$ is a $W(\bar\pi)$-submodule of $M$.
Then the vectors
$$
\{y m_{c'} + M_c'\:|\:y \in Y, c' \geq c\}
$$
form a basis for the $W(\bar\pi)$-module $M / M_c'$.
Hence the vector $\overline m_c :=
m_c + M_c'$ is a vector of maximal weight in $M / M_c'$, 
so it is annihilated by all $E_i^{(r)}$ for $i=1,\dots,n-2$ and $r > s_{i,i+1}$.
Moreover, using Lemma~\ref{rels}(iii),(vi) and (vii), 
Lemma~\ref{reallyhard} and the PBW theorem for $Y_{(1^{n})}^\flat(\sigma)$, 
one checks that
$$
u^{p_i} D_i(u) \overline m_c =
(u+a_{i,1}-c_{i,1}) \cdots (u+a_{i,p_i} - c_{i,p_i})
\overline m_c.
$$
Hence, $\overline m_c \in M / M_c'$ is a highest weight vector of type 
$B$ as claimed.
Now it follows easily using Theorem~\ref{vermathm}(ii) and the universal
property of generalized Verma modules that
$M_c$ is a $W(\bar\pi)$-submodule of $M$ and $M_c / M_c' \cong M(B)$.

\chapter{Standard modules}\label{sstandard}

In this chapter, we begin by classifying the 
finite dimensional
irreducible
representations of $W(\pi)$ and of $Y_n(\sigma)$, following
the argument
in the case of the Yangian $Y_n$ itself due to Tarasov
\cite{Tarasov2} and Drinfeld \cite{D3}.
Then we define and study another family of finite dimensional
$W(\pi)$-modules which we call standard modules.

\section{Two rows}\label{sstr}
In this section we assume that $n=2$ and let $\pi$
be any pyramid with just two rows
of lengths $p_1 \leq p_2$.
We will represent the $\pi$-tableau with entries $a_1,\dots,a_{p_1}$ on its first row and $b_1,\dots,b_{p_2}$ on its second row by
$\substack{a_1 \cdots a_{p_1} \\ b_1 \cdots b_{p_2}}$.
The first lemma is well known; see e.g. \cite{CP2}.
We reproduce here the detailed argument following 
\cite[Proposition 3.6]{Molevfd}
since we need to slightly weaken the hypotheses later on.

\begin{Lemma}\label{hw1}
Assume $p_1=p_2=l$ and 
$a_1,\dots,a_l,b_1,\dots,b_l,a,b \in \C$.
\begin{enumerate}
\item[\rm(i)] If $a_i > b$ implies that $a_i \geq a > b$
for each $i=1,\dots,l$, then all highest weight vectors
in $L(\substack{a_1 \cdots a_l \\ b_1 \cdots b_l}) 
\boxtimes L(\substack{a \\ b})$ are scalar multiples of $v_+ \otimes v_+$.
\item[\rm(ii)]
If $a > b_i$ implies that $a > b \geq b_i$
for each $i=1,\dots,l$, then all highest weight vectors
in $L(\substack{a \\ b})\boxtimes 
L(\substack{a_1 \cdots a_l \\ b_1 \cdots b_l})$ 
are scalar multiples of $v_+ \otimes v_+$.
\end{enumerate}
\end{Lemma}

\begin{proof}
(i) Abbreviate $e := e_{1,2}$, $d_2 := e_{2,2}$ and
$f := e_{2,1}$ in the Lie algebra $\mathfrak{gl}_2$.
Let $f^{(r)}$ denote $f^r / r!$.
Recall that
the irreducible $\mathfrak{gl}_2$-module
$L(\substack{a \\ b})$ of highest weight
$(a,b+1)$ has basis
$v_+, f v_+, f^{(2)} v_+, \dots$ if $a \not > b$ or
$v_+, f v_+, \dots, f^{(a-b-1)} v_+$ if $a > b$.
Also $e f^{(r+1)} v_+ = (a-b-r-1) f^{(r)} v_+$.

Suppose that 
$L(\substack{a_1 \cdots a_l \\ b_1 \cdots b_l}) 
\boxtimes L(\substack{a \\ b})$ contains a highest weight
vector $v$ that is not a scalar multiple of $v_+ \otimes v_+$.
We can write
$$
v = \sum_{i=0}^k m_i \otimes f^{(k-i)} v_+
$$
for vectors $m_0 \neq 0, m_1,\dots,m_k$ and 
$k \geq 0$ with $k < a-b$ in case $a > b$.
The element $T_{1,2}^{(r+1)}$ 
acts on the tensor product
as $T_{1,2}^{(r+1)}\otimes 1 
+ T_{1,2}^{(r)} \otimes d_2 + T_{1,1}^{(r)} \otimes e
\in W(\pi) \otimes U(\mathfrak{gl}_2)$.

Apply $T_{1,2}^{(r+1)}$
to the vector $v$ and compute the $? \otimes 
y^{(k)} v_+$-coefficient to deduce that
$$
T_{1,2}^{(r+1)} m_0 + (b+k+1) T_{1,2}^{(r)} m_0 = 0
$$
for all $r \geq 0$. It follows that $T_{1,2}^{(r)} m_0 = 0$ for all $r > 0$,
hence $m_0$ is a scalar multiple of the canonical highest weight vector $v_+$ of
$L(\substack{a_1 \cdots a_l \\ b_1 \cdots b_l})$.
Moreover we must in fact have that $k \geq 1$ since $v$ is not a multiple
of $v_+ \otimes v_+$.

Next compute the $? \otimes f^{(k-1)} v_+$-coefficient of 
$T_{1,2}^{(r+1)}v$ to get that
$$
T_{1,2}^{(r+1)} m_1 + (b+k) T_{1,2}^{(r)} m_1 + (a-b-k) T_{1,1}^{(r)} m_0 = 0.
$$
Multiply by $(-(b+k))^{l-r}$ and sum over $r=0,1,\dots,l$ to deduce that
$$
T_{1,2}^{(l+1)} m_1 + (a-b-k) \sum_{r=0}^l (-(b+k))^{l-r} T_{1,1}^{(r)} m_0 = 0.
$$
But $T_{1,2}^{(l+1)}  = 0$ in $W(\pi)$ by 
a trivial special case of Theorem~\ref{tvan}.
Moreover, by the definition (\ref{t1}), we have that 
$\sum_{r=0}^l u^{l-r} T_{1,1}^{(r)} m_0 =
(u+a_1) \cdots (u+a_l) m_0$. So we have shown that
$$
(a-b-k) (a_1 - b-k) (a_2 -b-k) \cdots (a_l - b-k) = 0.
$$
Since $k \geq 1$ and $k < a-b$ in case $a > b$, we have that
$(a-b-k) \neq 0$. Hence we must have that $a_i = b+k$ for some
$i=1,\dots,l$, i.e. $a_i > b$ and either $a \not > b$
or $a_i < a$. This is a contradiction.

(ii) Similar.
\end{proof}

\begin{Corollary}\label{hw2}
Assume $p_1=p_2=l$ and 
$a_1,\dots,a_l,b_1,\dots,b_l,a,b \in \C$.
\begin{enumerate}
\item[\rm(i)] If $b < a_i$ implies that $b < a \leq a_i$
for each $i=1,\dots,l$, then 
$L(\substack{a \\ b}) \boxtimes 
L(\substack{a_1 \cdots a_l \\ b_1 \cdots b_l})$ 
is a highest weight module generated by the 
highest weight vector $v_+ \otimes  v_+$.
\item[\rm(ii)]
If $b_i < a$ implies that $b_i \leq b < a$
for each $i=1,\dots,l$, then 
$L(\substack{a_1 \cdots a_l \\ b_1 \cdots b_l})\boxtimes
L(\substack{a \\ b})$ is a highest weight module
generated by the highest weight vector $v_+ \otimes  v_+$.
\end{enumerate}
\end{Corollary}

\begin{proof}
(i) By Lemma~\ref{hw1}(i),
$L(\substack{a_1 \cdots a_l \\ b_1 \cdots b_l})
\boxtimes L(\substack{a \\ b})$ has simple socle generated by the highest weight vector
$v_+ \otimes v_+$. Now apply the duality $?^\tau$ 
using 
Corollary~\ref{taudual} and (\ref{oppp}) to deduce that
$$
(L(\substack{a_1 \cdots a_l \\ b_1 \cdots b_l})
\boxtimes L(\substack{a \\ b}))^\tau \cong
L(\substack{a \\ b}) \boxtimes 
L(\substack{a_1 \cdots a_l \\ b_1 \cdots b_l})
$$
has a unique maximal submodule and that the highest weight vector
$v_+ \otimes v_+$ does not belong to this submodule.
Hence it is a highest weight module generated by the vector
$v_+ \otimes v_+$.

(ii) Similar.
\end{proof}

\begin{Remark}\label{r1}
The module
$L(\substack{a_1 \cdots a_l \\ b_1 \cdots b_l})$ in the statement of
Corollary~\ref{hw2} can in fact be replaced by any 
non-zero quotient of the generalized Verma module
$M(\substack{a_1 \cdots a_l \\ b_1 \cdots b_l})$.
This follows because the only property of
$L(\substack{a_1 \cdots a_l \\ b_1 \cdots b_l})$
needed for the proof of Lemma~\ref{hw1} 
is that all its highest weight vectors 
are scalar multiples of $v_+$; any non-zero submodule of 
$M(\substack{a_1 \cdots a_l \\ b_1 \cdots b_l})^\tau$ also has 
this property.
\end{Remark}

\begin{Lemma}\label{hw3}
Assume $p_1\leq p_2$ and 
 $a_1,\dots,a_{p_1}, b_1,\dots,b_{p_2}, b \in \C$.
\begin{enumerate}
\item[\rm(i)] If $a_i > b$ implies that $a_i > b_i \geq b$
for each $i=1,\dots,p_1$ then
all highest weight vectors in 
$L(\substack{a_1\cdots a_{p_1} \\ b_1 \cdots b_{p_2}})
\boxtimes L(\substack{\phantom{a}\\b})$
are scalar multiples of $v_+ \otimes v_+$.
\item[\rm(ii)]
All highest weight vectors in the module
$L(\substack{\phantom{a}\\b})
\boxtimes L(\substack{a_1\cdots a_{p_1} \\ b_1 \cdots b_{p_2}})$
are scalar multiples of $v_+ \otimes v_+$.
\end{enumerate}
\end{Lemma}

\begin{proof}
Let $\sigma = (s_{i,j})_{1 \leq i,j \leq 2}$ be a shift matrix
corresponding to the pyramid $\pi$.
Also note (since $n=2$) that 
$L(\substack{\phantom{a}\\b})$ is the 
one dimensional $\mathfrak{gl}_1$-module with basis $v_+$
such that $e_{1,1} v_+ = (b+1)v_+$.

(i) 
Suppose that $m \otimes v_+$ is a non-zero highest weight vector in 
$L(\substack{a_1\cdots a_{p_1} \\ b_1 \cdots b_{p_2}})
\boxtimes L(\substack{\phantom{a}\\b})$.
So we have that $E_1^{(r+1)} (m \otimes v_+) = 0$ for all
$r > s_{1,2}$ and 
$$
u^{p_1} D_1(u) (m \otimes v_+) = (u+c_1)(u+c_2) \cdots (u+c_{p_1}) (m \otimes v_+)
$$
for some scalars $c_1,\dots,c_{p_1} \in \C$.

Applying the Miura transform to Lemma~\ref{eorl}
(or see \cite[Lemma 11.3]{BK} and 
\cite[Theorem 4.1(i)]{BK}),
we have that $\Delta_{p_2,1}(E_1^{(r+1)}) = E_1^{(r+1)}  \otimes 1
 + E_1^{(r)} \otimes e_{1,1}$ for all $r > s_{1,2}$.
Hence
$E_1^{(r+1)} m + (b+1) E_1^{(r)} m = 0$
for all $r > s_{1,2}$. On setting
$m' := E_1^{(s_{1,2}+1)} m$, we deduce that
$E_1^{(s_{1,2}+r+1)} m = (-(b+1))^r m'$ for all $r \geq 0$, i.e.
$$
E_1(u) m = (1-(b+1)u^{-1} + (b+1)^2 u^{-2} - \cdots) u^{-s_{1,2}-1} m'
= \frac{u^{-s_{1,2}-1}}{1+(b+1)u^{-1}} m'.
$$
If $m' = 0$ then we have that $E_1^{(r)} m = 0$ for all
$r > s_{1,2}$, hence $m$ is a scalar multiple of $v_+$ as required.
So assume from now on that $m' \neq 0$ and aim for a contradiction.

Since $\Delta_{p_2,1}(D_1^{(r)}) = D_1^{(r)} \otimes 1$
for all $r > 0$ we have that
\begin{equation*}
D_1(u) m = (1+c_1u^{-1})(1+c_2 u^{-1})\cdots (1+c_{p_1}u^{-1}) m.
\end{equation*}
The last two equations and the identity
$[D_1(u), E_1^{(s_{1,2}+1)}] = u^{s_{1,2}} D_1(u) E_1(u)$ 
in $W(\pi)[[u^{-1}]]$ show that 
\begin{equation*}
D_1(u)  m' = 
\frac{(1+c_1 u^{-1}) \cdots (1+c_{p_1}u^{-1}) (1+(b+1)u^{-1})}{1+bu^{-1}}m'.
\end{equation*}
Since $D_1^{(r)} = 0$ for $r > p_1$ it follows from this that
$b = c_i$ for some $1 \leq i \leq p_1$.
Without loss of generality we may as well assume that
$b = c_1$. 
Then we have shown that
\begin{equation*}
D_1(u)  m' = 
(1+(c_1+1) u^{-1})(1+c_2 u^{-1}) 
\cdots (1+c_{p_1}u^{-1})m'.
\end{equation*}

Now we claim that if we have any non-zero vector
in $L(\substack{a_1 \cdots a_{p_1} \\
b_1 \cdots b_{p_2}})$ on which $D_1(u)$ acts as the scalar
$(1+d_1 u^{-1}) \cdots (1+d_{p_1} u^{-1})$ then there exists
a permutation $w \in S_{p_1}$ such that
$a_i \geq d_{wi}$ and moreover if $a_i > b_i$ then $d_{wi} > b_i$,
for each $i=1,\dots,p_1$.
To prove this, we may replace the module
$L(\substack{a_1 \cdots a_{p_1} \\
b_1 \cdots b_{p_2}})$ with the tensor product
$L(\substack{a_1 \\ b_1}) \boxtimes\cdots\boxtimes
L(\substack{a_{p_1} \\ b_{p_1}}) \boxtimes 
L(\substack{\phantom{a}\\b_{p_1+1}}) \boxtimes\cdots\boxtimes
L(\substack{\phantom{a}\\b_{p_2}})$, since that contains
$L(\substack{a_1 \cdots a_{p_1} \\
b_1 \cdots b_{p_2}})$ (possibly twisted by the isomorphism $\iota$) 
as a subquotient.
Now the claim follows from Lemma~\ref{ten2} and the
familiar fact that
if we have a non-zero vector in the irreducible 
$\mathfrak{gl}_2$-module
$L(\substack{a \\ b})$ on which $D_1(u)$ acts as the scalar
$(1+d u^{-1})$ then $a \geq d$ and moreover if $a > b$ then 
$d > b$.

Applying the claim to
the non-zero vectors $m$ and $m'$ of
$L(\substack{a_1 \cdots a_{p_1} \\
b_1 \cdots b_{p_2}})$, we deduce (after reordering
if necessary) that there exists a permutation $w \in S_{p_1}$ such that
\begin{enumerate}
\item[(a)]
$a_1 \geq c_1 + 1$ and moreover if $a_1 > b_1$ then $c_1+1 > b_1$;
$a_2 \geq c_2$ and moreover if $a_2 > b_2$ then $c_2 > b_{2}$; \:\dots\:;
$a_{p_1} \geq c_{p_1}$ and moreover if $a_{p_1} > b_{p_1}$
then $c_{p_1} > b_{p_1}$;
\item[(b)]
$a_1 \geq c_{w1}$ and moreover if $a_1 > b_1$ then $c_{w1} > b_1$;
$a_2 \geq c_{w2}$ and moreover if $a_2 > b_2$ then $c_{w2} > b_{2}$;
\:\dots\:; $a_{p_1} \geq c_{wp_1}$ and moreover if $a_{p_1} > b_{p_1}$
then $c_{wp_1} > b_{p_1}$.
\end{enumerate}
From this we can derive the required contradiction, as follows.
Suppose that we know that
$c_i > b$ for some $i$. 
Then $a_i \geq c_i > b$, hence by the hypothesis from the statement of the lemma $a_i \geq c_{wi} > b_i \geq b$. Hence $c_{wi} > b$.
Now we do know that $c_1 = b$. Hence
$a_1 \geq c_1+1 > b$, so $a_1 \geq c_{w1} > b_1 \geq b$.
Hence $c_{w1} > b$. Combining this with the preceeding observation we deduce
that $c_{w^k 1} > b$ for all $k \geq 1$, hence in particular $c_1 > b$.

(ii) We have that $\Delta_{1,p_2}(E_1^{(r)}) = 
1 \otimes E_1^{(r)}$ for all $r > s_{1,2}$. 
So 
if $v_+ \otimes m$ is a highest weight vector in
$L(\substack{\phantom{a}\\b})
\boxtimes 
L(\substack{a_1\cdots a_{p_1} \\ b_1 \cdots b_{p_2}})$
then $E_1^{(r)} m = 0$ for all $r > s_{1,2}$.
Hence $m$ is a scalar multiple of $v_+$ as required.
\end{proof}

\begin{Corollary}\label{hw4}
Assume $p_1\leq p_2$ and 
 $a_1,\dots,a_{p_1}, b_1,\dots,b_{p_2}, b \in \C$.
\begin{enumerate}
\item[\rm(i)] If $b < a_i$ implies that $b \leq b_i < a_i$
for each $i=1,\dots,p_1$ then the module 
$L(\substack{\phantom{a}\\b}) \boxtimes 
L(\substack{a_1\cdots a_{p_1} \\ b_1 \cdots b_{p_2}})$
is a highest weight module generated by the highest weight vector 
$v_+ \otimes v_+$.
\item[\rm(ii)]
The module
$L(\substack{a_1\cdots a_{p_1} \\ b_1 \cdots b_{p_2}})
\boxtimes L(\substack{\phantom{a}\\b})$
is a highest weight module generated by the highest weight vector
$v_+ \otimes v_+$.
\end{enumerate}
\end{Corollary}

\begin{proof}
Argue using the duality $?^\tau$ exactly as in the proof of Corollary~\ref{hw2}.
\end{proof}

\begin{Remark}\label{rem2}
As in Remark~\ref{r1}, the module
$L(\substack{a_1 \cdots a_{p_1} \\ b_1 \cdots b_{p_2}})$ in the statement of
Corollary~\ref{hw4}(ii) can be replaced by any 
non-zero quotient of the generalized Verma module
$M(\substack{a_1 \cdots a_{p_1} \\ b_1 \cdots b_{p_2}})$.
We cannot quite say the same thing for 
Corollary~\ref{hw4}(i), but by the proof we can at least replace
$L(\substack{a_1 \cdots a_{p_1} \\ b_1 \cdots b_{p_2}})$ by any 
non-zero quotient $M$ of the generalized Verma module
$M(\substack{a_1 \cdots a_{p_1} \\ b_1 \cdots b_{p_2}})$
with the property that all of its Gelfand-Tsetlin weights, i.e.
the $A(u) \in \mathscr P_2$ such that $M_{A(u)} \neq 0$, are
also Gelfand-Tsetlin weights of the module
$L(\substack{a_1\\ b_1})
\boxtimes\cdots\boxtimes
L(\substack{a_{p_1}\\ b_{p_1}})
\boxtimes
L(\substack{\phantom{a_1}\\ b_{p_1+1}})
\boxtimes\cdots\boxtimes
L(\substack{\phantom{a_1}\\ b_{p_2}}).$
\end{Remark}

Now we can prove the main theorem of the section. This is new only
if $p_1 \neq p_2$.

\begin{Theorem}\label{sl2}
Assume $p_1 \leq p_2$ and  
$a_1,\dots,a_{p_1}, b_1,\dots,b_{p_2} \in \C$ 
are scalars
such that the following property holds for each $i=1,\dots,p_1$: 
\begin{quote}
If the set 
$\{a_j - b_k \:|\: i \leq j \leq p_1, i \leq k \leq p_2
\text{ such that }a_j > b_k\}$ 
is non-empty then $(a_i - b_i)$ 
is its smallest element.
\end{quote}
Then the irreducible $W(\pi)$-module
$L(\substack{a_1 \cdots a_{p_1} \\ b_1 \cdots b_{p_2}})$ is isomorphic to the
tensor product of the modules
$$
L(\substack{a_1 \\ b_1}),\dots,
L(\substack{a_{p_1} \\ b_{p_1}}),
L(\substack{\phantom{a}\\ b_{p_1+1}}),\dots,
L(\substack{\phantom{a}\\ b_{p_2}})
$$
taken in any order that matches the shape of the pyramid $\pi$.
\end{Theorem}

\begin{proof}
Assume to start with that the pyramid $\pi$ is left-justified.
First we show for $p_1 > 0$ that
$$
L(\substack{a_1 \cdots a_{p_1} \\ b_1 \cdots b_{p_1}})
\cong
L(\substack{a_1  \\ b_1})
\boxtimes
L(\substack{a_2 \cdots a_{p_1} \\ b_2 \cdots b_{p_1}}).
$$
Since $a_1 > b_i$ implies that $a_1 > b_1 \geq b_i$ for all
$i=2,\dots,p_1$, Lemma~\ref{hw1}(ii) implies that $v_+\otimes v_+$
is the unique (up to scalars) 
highest weight vector in the module on the right hand side.
Since $b_1 < a_i$ implies $b_1 < a_1 \leq a_i$, Corollary~\ref{hw2}(i)
shows that this vector generates
the whole module. Hence it is irreducible, so isomorphic to
$L(\substack{a_1 \cdots a_{p_1} \\ b_1 \cdots b_{p_1}})$ by
Lemma~\ref{ten1}.
Next we show for $p_2 > p_1$ that
$$
L(\substack{a_1 \cdots a_{p_1} \\ b_1 \cdots b_{p_2}})
\cong
L(\substack{a_1 \cdots a_{p_1}\phantom{21} \\ b_1 \cdots b_{p_2-1}})
\boxtimes L(\substack{\phantom{a}\\ b_{p_2}}).
$$
Since $a_i > b_{p_2}$ implies $a_i > b_i \geq b_{p_2}$
Lemma \ref{hw3}(i) implies that $v_+ \otimes v_+$ is the unique
(up to scalars) highest weight vector in the module on the right hand side.
But by Corollary~\ref{hw4}(ii) this vector generates the whole module, hence it
is irreducible.
Using these two facts,
it follows by induction on $p_2$ that
$$
L(\substack{a_1 \cdots a_{p_1} \\ b_1 \cdots b_{p_2}})
\cong 
L(\substack{a_1 \\ b_1}) \boxtimes\cdots\boxtimes
L(\substack{a_{p_1} \\ b_{p_1}}) \boxtimes
L(\substack{\phantom{a}\\ b_{p_1+1}})\boxtimes\cdots\boxtimes
L(\substack{\phantom{a}\\ b_{p_2}}).
$$
This proves the theorem for one particular ordering of the tensor product
and for one particular choice of the pyramid $\pi$ with row lengths 
$(p_1,p_2)$.
The theorem for all other orderings and pyramids
follows from this by character considerations. 
\end{proof}

Suppose finally that 
we are
given an arbitrary two row tableau
$A$ with entries $a_1,\dots,a_{p_1}$ on row one
and $b_1,\dots,b_{p_2}$ on row two. We can always reindex the entries in the rows
so that the hypothesis of Theorem~\ref{sl2} is satisfied:
first reindex to
ensure if possible that $a_1-b_1$ is the minimal positive integer difference
amongst all the dfferences $a_i - b_j$, then inductively reindex the
remaining entries $a_2,\dots,a_{p_1}, b_2,\dots,b_{p_2}$.
Hence Theorem~\ref{sl2} shows that {\em every} irreducible 
admissible $W(\pi)$-module can be realized as a tensor product of
irreducible $\mathfrak{gl}_2$- and $\mathfrak{gl}_1$-modules. 
This remarkable observation was first made 
by Tarasov \cite{Tarasov2} in the case $p_1=p_2$.

\begin{Corollary}\label{stupid}
If the irreducible module
$L(\substack{a_1 \cdots a_{p_1} \\ b_1 \cdots b_{p_2}})$ is finite dimensional
for scalars $a_1,\dots,a_{p_1}, b_1,\dots,b_{p_2} \in \C$
then there exists a permutation $w \in S_{p_2}$ such that
$a_{1} > b_{w1}, a_{2} > b_{w2},\dots, a_{p_1} > b_{w p_1}$.
\end{Corollary}

\begin{proof}
Reindexing if necessary, we may assume that the hypothesis of Theorem~\ref{sl2} is satisfied. Then by the theorem we must have that
$L(\substack{a_i \\ b_i})$ is finite dimensional for each
$i=1,\dots,p_1$, i.e $a_i > b_i$ for each such $i$.
\end{proof}

\section{Classification of finite dimensional irreducible representations}

Now assume that $\pi = (q_1,\dots,q_l)$ is an arbitrary pyramid 
with row lengths $(p_1,\dots,p_n)$.
Let $\sigma = (s_{i,j})_{1 \leq i,j \leq n}$ be a shift matrix corresponding 
to $\pi$, so that
$W(\pi)$ is canonically a quotient of the shifted Yangian $Y_n(\sigma)$.
Recall the definitions of the sets 
$\Row(\pi)$ of row symmetrized $\pi$-tableaux, 
$\Col(\pi)$ of column strict $\pi$-tableaux 
and $\Dom(\pi)$ of dominant row symmetrized $\pi$-tableaux
from $\S$\ref{sstableaux}.

\begin{Theorem}\label{fd}
For $A \in \Row(\pi)$, the irreducible $W(\pi)$-module
$L(A)$ is finite dimensional if and only if $A$ is dominant,
i.e. it has a representative belonging to $\Col(\pi)$.
\end{Theorem}

\begin{proof}
Suppose first that $L(A)$ is finite dimensional.
For each $i=1,\dots,n-1$, let $\sigma_i$ denote the
$2 \times 2$ submatrix 
$$
\left( \begin{array}{ll}s_{i,i} 
&s_{i,i+1}\\ s_{i+1,i} & s_{i+1,i+1} \end{array}
\right)
$$ 
of the matrix $\sigma$.
Also let $a_{i,1},\dots,a_{i,p_i}$ be the entries in the $i$th row of $A$
for each $i=1,\dots,n$.
The map $\psi_{i-1}$ from (\ref{psisdef}) obviously 
induces an embedding of the 
shifted Yangian $Y_2(\sigma_i)$ into $Y_n(\sigma)$.
The highest weight vector $v_+ \in L(A)$ is also a 
highest weight vector in the restriction of $L(A)$ 
to $Y_2(\sigma_i)$ using this embedding.
Hence by Corollary~\ref{stupid} 
there exists $w \in S_{p_{i+1}}$ such that
$$
a_{i,1} > a_{i+1,w1},\:\: a_{i,2} > a_{i+1,w2},\:\dots\:,\:a_{i,p_i} > a_{i+1,w p_i},
$$
for each $i=1,\dots,n-1$.
Hence $A$ has a representative belonging to $\Col(\pi)$.

Conversely, suppose that $A$ has a representative belonging to $\Col(\pi)$.
Let $A_1,\dots,A_l$ be the columns of this representative, 
so that $A \sim_{\ro} A_1\otimes\cdots\otimes A_l$. 
Since $A_i$ is column strict, 
the irreducible module
$L(A_i)$ is finite dimensional. 
By Lemma~\ref{ten1} the tensor product
$L(A_1) \boxtimes \cdots \boxtimes L(A_l)$ is then 
a finite dimensional $W(\pi)$-module
containing a highest weight vector of type $A$. Hence $L(A)$ is 
finite dimensional.
\end{proof}

Hence, the modules
$\{L(A)\:|\:A \in \Dom(\pi)\}$
give a full set of pairwise non-isomorphic finite dimensional irreducible 
$W(\pi)$-modules. As a corollary, we have the following result
classifying the finite dimensional irreducible representations of
the shifted Yangians $Y_n(\sigma)$ themselves.
Since every finite dimensional $Y_n(\sigma)$-module
is admissible, it is enough for this to 
determine which of the irreducible modules $L(\sigma,A(u))$
from (\ref{lnsa}) is finite dimensional.

\begin{Corollary}\label{yfdc}
For $A(u) \in \mathscr P_n$, 
the irreducible $Y_n(\sigma)$-module $L(\sigma,A(u))$ is
finite dimensional if and only if
there exist (necessarily unique)
monic polynomials $P_1(u),\dots,P_{n-1}(u)$, $Q_1(u),\dots,Q_{n-1}(u) \in \C[u]$ such that
$(P_i(u), Q_i(u)) = 1$,
$Q_i(u)$ is of degree $d_i := s_{i,i+1}+s_{i+1,i}$, and
$$
\frac{A_i(u)}{A_{i+1}(u)}
=
\frac{P_i(u)}{P_i(u-1)} \times \frac{u^{d_i}}{Q_i(u)}
$$
for each $i=1,\dots,n-1$.
\end{Corollary}

\begin{proof}
Recall from Remark~\ref{tarrem}
that every admissible irreducible $Y_n(\sigma)$-module may be obtained
by inflating an admissible irreducible $W(\pi)$-module
through the map (\ref{infmap}), for some pyramid $\pi$ with shift matrix
$\sigma$ and some $f(u) \in 1 + u^{-1} \C[[u^{-1}]]$.
Given this and Theorem~\ref{fd}, we see that
$L(\sigma,A(u))$ is finite dimensional if and only if
there exist $l \geq s_{n,1}+s_{1,n}$, $f(u) \in 1+u^{-1}\C[[u^{-1}]]$ and
scalars $a_{i,j}\in \C$ for $1 \leq i \leq n, 1 \leq j \leq p_i := l - s_{n,i}-s_{i,n}$
such that
\begin{enumerate}
\item[(a)] $A_i(u) = f(u) (1+a_{i,1}u^{-1}) \cdots (1+a_{i,p_i}u^{-1})$
for each $i=1,\dots,n$;
\item[(b)] $a_{i,j} \geq a_{i+1,j}$ for each $i=1,\dots,n-1$
and $j=1,\dots,p_i$.
\end{enumerate}
Following the proof of \cite[Theorem 2.8]{Molevfd}, these conditions are
equivalent to the existence of monic polynomials
$P_1(u),\dots,P_{n-1}(u),Q_1(u),\dots,Q_{n-1}(u) \in \C[u]$
such that $Q_i(u)$ is of degree $d_i$ and
$$
\frac{A_i(u)}{A_{i+1}(u)}
=
\frac{P_i(u)}{P_i(u-1)} \times \frac{u^{d_i}}{Q_i(u)}
$$
for each $i=1,\dots,n-1$.
Finally to get uniqueness of the $P_i(u)$'s and $Q_i(u)$'s we
have to insist in addition that $(P_i(u), Q_i(u)) = 1$.
\end{proof}

\begin{Remark}
From Corollary~\ref{yfdc} and (\ref{syn}), it also follows that 
the isomorphism classes of irreducible
$SY_n(\sigma)$-modules are parametrized in the same fashion by
monic polynomials $P_1(u),\dots,P_{n-1}(u)$, 
$Q_1(u),\dots,Q_{n-1}(u) \in \C[u]$ such that
$Q_i(u)$ is of degree $d_i$ and $(P_i(u), Q_i(u)) = 1$
for each $i=1,\dots,n-1$.
In the case $\sigma$ is the zero matrix, each $Q_i(u)$ is 
of course just equal to $1$, so we recover the classification
from \cite{D3}
of finite dimensional irreducible representations
of the Yangian of $\mathfrak{sl}_n$ by their {\em Drinfeld polynomials} $P_1(u),\dots,P_{n-1}(u)$; see also \cite[$\S$2]{Molevfd} once more. 
\end{Remark}

\section{Tensor products}\label{sstp}
Continuing with the notation from the previous section, 
we set $m := q_l$ for short.
For $A \in \Col(\pi)$ with columns $A_1,\dots,A_l$ from left to right,
let
\begin{equation}\label{VAdef}
V(A) := L(A_1) \boxtimes \cdots \boxtimes L(A_l).
\end{equation}
We will refer to the modules $\{V(A)\:|\:A \in \Col(\pi)\}$
as {\em standard modules}.
As we observed already in the proof of Theorem~\ref{fd}, each $V(A)$ 
is a finite dimensional $W(\pi)$-module, and the vector
$v_+\otimes\cdots\otimes v_+ \in V(A)$ is a highest weight vector
of type equal to the row equivalence class of $A$.
We wish to give a sufficient condition for $V(A)$
to be a highest weight module generated by this highest weight vector,
following an argument due to Chari \cite{Chari} in the context of quantum affine algebras. The key step is provided by the following lemma; 
in its statement we work with the usual action of the symmetric group $S_m$
on finite dimensional irreducible $\mathfrak{gl}_m$-modules, 
and $s_1,\dots,s_{m-1} \in S_m$ denote the basic transpositions.

\begin{Lemma}\label{techy}
Suppose that we are given a $\pi$-tableau $A$ with columns $A_1,\dots,A_l$ from left to right,
together with $1 \leq t < m =q_l$ and
$w \in S_m$ such that $$
t \geq w^{-1} t < w^{-1} (t+1).
$$
Letting $a_1,\dots,a_p$ 
resp.
$c_1,\dots,c_q,b_1,\dots,b_p$ 
denote 
the entries in the $(n-m+t)$th 
resp. the $(n-m+t+1)$th row of $A$
read from left to right,
assume that
\begin{enumerate}
\item[\rm(i)] $a_i > b_i$ for each $i=1,\dots,p$;
\item[\rm(ii)] $a_i \not > a_j$ for each $1 \leq i < j \leq p$;
\item[\rm(iii)] 
either $c_i \not < a_j$ or $c_i \leq b_j$
for each $i=1,\dots,q$ and $j=1,\dots,p$;
\item[\rm(iv)] none of the elements $c_1,\dots,c_q$ lie in the same
coset of $\C$ modulo $\Z$ as $a_p$;
\item[\rm(v)] $A_l$ is column strict.
\end{enumerate}
Then the vector
$v_+\otimes\cdots\otimes v_+ \otimes s_t w v_+$
is an element of the $W(\pi)$-submodule of $L(A_1)\boxtimes\cdots\boxtimes
L(A_{l-1})\boxtimes L(A_l)$
generated by the vector $v_+\otimes\cdots\otimes v_+ \otimes w v_+$.
\end{Lemma}

Since this is technical, let us postpone the proof until the end of the section and explain the applications. For the first one,
recall from $\S$\ref{sstableaux}
the definition of the set $\Std(\pi)$ of standard $\pi$-tableaux
in the case that $\pi$ is left-justified.

\begin{Theorem}\label{charithm}
Assume that the pyramid $\pi$ is left-justified and let
$A \in \Std(\pi)$.
Then the $W(\pi)$-module $V(A)$
is a highest weight module generated by the highest weight vector
$v_+\otimes\cdots\otimes v_+$.
\end{Theorem}

\begin{proof}
Let $A_1,\dots,A_l$ denote the columns of $A$ from left to right, and set
$M := L(A_1) \boxtimes\cdots\boxtimes L(A_{l-1}), L := L(A_l)$ for short.
By induction on $l$, $M$ is a highest weight module generated
by the vector $v_+\otimes\cdots\otimes v_+$.
Fix the 
reduced expression
$w_0 = s_{i_h} \cdots s_{i_1}$
for the longest element of
the symmetric group
$S_m$ where
$$
(i_1,\dots,i_h) = (m-1;m-2,m-1;\dots;2,\dots,m-1;1,\dots,m-1).
$$
For $r=0,\dots,h$ let $v_r := s_{i_r} \cdots s_{i_1} v_+
\in L$. Note by the choice of reduced expression that
$i_{r+1} \geq s_{i_1} \cdots s_{i_r} (i_{r+1}) < s_{i_1} \cdots s_{i_r} (i_{r+1}+1)$. So, taking $w = s_{i_r} \cdots s_{i_1}$ and $t = i_{r+1}$ for
some $r=0,\dots,h-1$, the hypotheses of Lemma~\ref{techy}
are satisfied.
Hence the lemma implies that the vector
$v_+\otimes\cdots\otimes v_+ \otimes v_{r+1}$ lies in the
$W(\pi)$-submodule of $M \boxtimes L$ generated by the vector
$v_+\otimes\cdots\otimes v_+ \otimes v_r$.
This is true for all $r=0,\dots,h-1$, and $v_h = w_0 v_+$. So this
shows that the vector
$v_+\otimes\cdots\otimes v_+ \otimes w_0 v_+$ lies in the
$W(\pi)$-submodule of $M \boxtimes L$ generated by the highest weight vector
$v_+\otimes\cdots\otimes v_+ \otimes v_+$.

Now to complete the proof we show that $M \boxtimes L$
is generated as a $W(\pi)$-module by the vector
$v_+\otimes\cdots\otimes v_+ \otimes w_0 v_+$.
Let $M_d$ denote the span of all generalized weight spaces of $M$
of weight $\lambda - (\eps_{j_1}-\eps_{j_1+1})
 -\cdots-(\eps_{j_d}-\eps_{j_d+1})$
for $1 \leq j_1,\dots,j_d < n$, where
$\lambda \in \mathfrak{c}^*$ is the weight of the highest weight vector
$v_+\otimes\cdots\otimes v_+$ of $M$.
We will prove by induction on $d \geq 0$ that
$M_d \otimes L$
is contained in the $W(\pi)$-submodule of
$M \boxtimes L$
generated by the vector
$(v_+\otimes \cdots\otimes v_+) \otimes w_0 v_+$.
Note to start with for any vector $y \in L$
and $1 \leq i < m$ that 
$$
E_{n-m+i}^{(1)}((v_+\otimes\cdots\otimes v_+) \otimes y)
=
(v_+\otimes\cdots\otimes v_+) \otimes (e_{i,i+1} y).
$$
Since $L$ is
generated as a $\mathfrak{gl}_m$-module
by the lowest weight vector $w_0 v_+$
this is enough to verify the base case.
Now for the induction step we know already that
$M$ is a highest weight module, hence it suffices to show that 
every vector of the form $(F_i^{(r)} x) \otimes y$ 
for $1 \leq i < n, r > 0, x \in M_{d-1}$ and $y \in L$ 
lies in the $W(\pi)$-submodule of $M \boxtimes L$ 
generated by $M_{d-1}\otimes L$.
But for this we have that
$$
F_i^{(r)} (x \otimes y) \equiv (F_i^{(r)} x) \otimes y
\pmod{M_{d-1}\otimes L}
$$
by Theorem~\ref{utah}(iii).
\end{proof}

For the second application, we return
to an arbitrary pyramid $\pi = (q_1,\dots,q_l)$.
The following theorem reduces 
the problem of computing the characters of all finite dimensional
irreducible $W(\pi)$-modules to that of computing the characters 
just of the modules $L(A)$ where all
entries of $A$ lie in the same coset of $\C$ modulo $\Z$.
Twisting moreover with the automorphism
$\eta_c$ from (\ref{etac}) using Lemma~\ref{eta} one can reduce
further to the case that all entries of $A$ actually lie in 
$\Z$ itself, i.e. $A \in \Dom_0(\pi)$.

\begin{Theorem}\label{irret}
Suppose that $\pi = \pi' \otimes \pi''$ for pyramids $\pi'$ and $\pi''$,
and we are given $A' \in \Dom(\pi')$
and $A'' \in \Dom(\pi'')$ such that
no entry of $A'$ lies in the same
coset of $\C$ modulo $\Z$ as an entry of $A''$.
Then the $W(\pi)$-module
$L(A') \boxtimes L(A'')$
is irreducible 
with highest weight vector $v_+ \otimes v_+$.
\end{Theorem}

\begin{proof}
By character considerations, we may assume for the proof
that the pyramid $\pi'$ is right-justified of level $l'$
and the pyramid $\pi''$
is left-justified of level $l''$.
Pick a standard $\pi''$-tableau representing $A''$ and let 
$A_{l'+1},A_{l'+2},\dots,A_{l}$
be its columns read from left to right.
We claim that
$L(A') \boxtimes L(A_{l'+1}) \boxtimes \cdots \boxtimes L(A_{l})$ is
a highest weight module generated by the highest weight vector
$v_+ \otimes v_+\otimes\cdots\otimes v_+$.
The theorem follows from this claim as follows.
By Theorem~\ref{charithm}, $L(A'')$ is a quotient of 
$L(A_{l'+1})\boxtimes\cdots\boxtimes L(A_l)$.
Hence we get from the claim 
that $L(A') \boxtimes L(A'')$ is a highest weight module
generated by the highest weight vector $v_+ \otimes v_+$.
Similarly so is $L(A'')^\tau \boxtimes L(A')^\tau$, 
hence $v_+ \otimes v_+$ is actually
the unique (up to scalars) highest weight vector in $L(A') \boxtimes L(A'')$. 
Thus $L(A') \boxtimes L(A'')$ is irreducible.

To prove the claim, 
fix the same reduced expression
$w_0 = s_{i_h} \cdots s_{i_1}$ for the longest element of 
$S_m$ as in the proof of Theorem~\ref{charithm}. 
Let $v_r := s_{i_r} \cdots s_{i_1} v_+ \in L(A_l)$.
We are actually going to show that $v_{r+1}$ lies in the $W(\pi)$-submodule
of $L(A') \boxtimes L(A_{l'+1}) \boxtimes \cdots \boxtimes L(A_{l})$
generated by the vector
$v_+ \otimes v_+\otimes\cdots\otimes v_+ \otimes v_r$ for each
$r=0,\dots,h-1$.
Given this, it follows that 
$v_+\otimes v_+\otimes\cdots\otimes v_+ \otimes w_0 v_+$
lies in the $W(\pi)$-submodule generated by the highest weight 
vector.
Since we already know by induction that
$L(A') \boxtimes L(A_{l'+1}) \boxtimes \cdots \boxtimes L(A_{l-1})$
is highest weight, the argument can then be completed
in the same way as in last paragraph of the proof of Theorem~\ref{charithm}.

So finally fix a choice of $r = 0,\dots,h-1$.
Let $w := s_{i_r} \cdots s_{i_1}$ and $t := i_{r+1}$.
Pick a representative for $A'$ so that,
letting $a_1,\dots,a_p$ resp. $c_1,\dots,c_q,b_1,\dots,b_p$ denote the
entries in its $(n-m+t)$th resp. $(n-m+t+1)$th row read from left to right,
we have that
\begin{enumerate}
\item[\rm(a)] $a_i > b_i$ for each $i=1,\dots,p$;
\item[\rm(b)] $a_i \not > a_j$ for each $1 \leq i < j \leq p$;
\item[\rm(c)] 
either $c_i \not < a_j$ or $c_i \leq b_j$
for each $i=1,\dots,q$ and $j=1,\dots,p$.
\end{enumerate}
To see that this is possible, it is easy 
to arrange things so that
(a) and (b) are satisfied. If $p > 0$ we rearrange the
$(n-m+i+1)$th row so that $a_p - b_p$ is the smallest
positive integer in the set
$\{a_p-b_1,\dots,a_p-b_p,a_p-c_1,\dots,a_p-c_q\}$.
The condition (c) is then automatic for $j=p$, and the
remaining entries 
$c_1,\dots,c_q,b_1,\dots,b_{p-1}$
can then be rearranged 
inductively to get (c) in general.
Let $A_1,\dots,A_{l'}$ denote the columns of this representative
from left to right.
It then follows by Lemma~\ref{techy} that
$v_+\otimes \cdots \otimes v_+ \otimes v_{r+1}$
lies in the $W(\pi)$-submodule of
$L(A_1)\boxtimes\cdots\boxtimes L(A_{l-1})\boxtimes 
L(A_{l})$ generated by the vector $v_+\otimes\cdots\otimes v_+\otimes v_r$.
Since $L(A')$ is a quotient of the submodule of
$L(A_1)\boxtimes\cdots\boxtimes L(A_{l'})$ generated by the highest weight
vector $v_+\otimes\cdots\otimes v_+$, this completes the proof.
\end{proof}

We still need to explain the proof of Lemma~\ref{techy}.
Let the notation be as in the statement of the lemma
and abbreviate $n-m+t$ by $i$.
Let $\pi'$ be the pyramid consisting just of the
$i$th and $(i+1)$th rows of $\pi$.
The $2 \times 2$ submatrix
$\sigma'$ 
consisting just of the $i$th and $(i+1)$th rows and columns of $\sigma$
gives a choice of shift matrix for $\pi'$.
As in the proof of Theorem~\ref{fd}, the map $\psi_{i-1}$
from (\ref{psisdef})
induces an embedding
$\phi:Y_2(\sigma') \hookrightarrow Y_n(\sigma)$.
For $j=1,\dots,l$, let
$$
q_j' := \left\{
\begin{array}{ll}
2&\hbox{if $n-q_j < i$,}\\
1&\hbox{if $n-q_j = i$,}\\
0&\hbox{if $n-q_j > i$.}
\end{array}\right.
$$
So, numbering the columns of the pyramid $\pi'$
by $1,\dots,l$ in the same way as in the pyramid $\pi$,
its columns are of heights $q_1',q_2',\dots,q_l'$ from left to right
(including possibly some empty columns at the left hand edge).
Recall the quotient map
$\kappa
:Y_n(\sigma) \twoheadrightarrow W(\pi)$
and the
Miura transform
$\xi:W(\pi) \hookrightarrow
U(\mathfrak{gl}_{q_1}) \otimes\cdots\otimes U(\mathfrak{gl}_{q_l})$
from  (\ref{thetadef}) and
(\ref{miura}).
Similarly we have the quotient map
$\kappa':Y_2(\sigma') \twoheadrightarrow W(\pi')$
and the Miura transform
$\xi':W(\pi') \hookrightarrow
U(\mathfrak{gl}_{q_1'}) \otimes\cdots\otimes U(\mathfrak{gl}_{q_l'})$.
For each $j=1,\dots,l$, define an algebra embedding
$\phi_j:U(\mathfrak{gl}_{q_j'}) \hookrightarrow U(\mathfrak{gl}_{q_j})$
so that if $q_j' = 2$ then
\begin{align*}
e_{1,1} &\mapsto e_{q_j-n+i,q_j-n+i},
&e_{1,2} &\mapsto e_{q_j-n+i,q_j-n+i+1},\\
e_{2,1} &\mapsto e_{q_j-n+i+1,q_j-n+i},
&e_{2,2} &\mapsto e_{q_j-n+i+1,q_j-n+i+1},
\end{align*}
and if $q_j'=1$ then
$e_{1,1} \mapsto e_{1,1}$. 
We have now defined all the maps in the following diagram:
\begin{equation}
\begin{CD}
Y_2(\sigma') & @>\kappa'>>& W(\pi') &@>\xi' >>& U(\mathfrak{gl}_{q_1'})
\otimes\cdots\otimes U(\mathfrak{gl}_{q_l'})\\
@V\phi VV&&&&&&@VV\phi_1\otimes\cdots\otimes \phi_l V\\
Y_n(\sigma) & @>\kappa>>& W(\pi) &@>\xi> >& U(\mathfrak{gl}_{q_1})
\otimes\cdots\otimes U(\mathfrak{gl}_{q_l})
\end{CD}
\end{equation}
This diagram definitely does {\em not} commute.
So the two actions
of $Y_2(\sigma')$ on the 
$U(\mathfrak{gl}_{q_1})
\otimes\cdots\otimes U(\mathfrak{gl}_{q_l})$-module
$L(A_1) \boxtimes\cdots\boxtimes L(A_l)$ defined 
using the homomorphism
$\xi\circ\kappa\circ\phi$ or using the homomorphism
$\phi_1\otimes\cdots\otimes \phi_l \circ \xi'\circ\kappa'$
are in general different.
In the proof of 
the following lemma we will see that in fact the two actions coincide on
special vectors.

\begin{Lemma}\label{techy2}
The following subspaces
of $L(A_1)\boxtimes\cdots\boxtimes L(A_{l-1})\boxtimes L(A_l)$ are equal:
\begin{align*}
&(\xi\circ\kappa\circ\phi)(Y_2(\sigma')) 
(v_+\otimes\cdots\otimes v_+ \otimes wv_+),\\
&(\phi_1\otimes\cdots\otimes \phi_l \circ\xi'\circ\kappa')(Y_2(\sigma')) 
(v_+\otimes\cdots\otimes v_+ \otimes wv_+).
\end{align*}
\end{Lemma}

\begin{proof}
For $j=1,\dots,l-1$, let $v_j$
be an element of $L(A_j)$ whose weight is equal to the 
weight of the highest weight vector $v_+$ of $L(A_j)$ minus
some multiple of the $i$th simple root $\eps_i-\eps_{i+1} \in \mathfrak{c}^*$.
Also let $v_l$ be any element of $L(A_l)$.
We claim for 
any element $x$ of $Y_2(\sigma')$ that
$$
(\xi\circ\kappa\circ\phi)(x) 
(v_1\otimes\cdots\otimes v_l)
=(\phi_1\otimes\cdots\otimes \phi_l \circ\xi'\circ\kappa')(x) 
(v_1\otimes\cdots\otimes v_l).
$$
Clearly the lemma follows from this claim. The advantage of the claim is that
it suffices to prove it for $x$ running over a set of generators for the algebra
$Y_2(\sigma')$,
since the vector on the right hand side of the equation 
can obviously be expressed as 
a linear combination of vectors of the form
$v_1'\otimes\cdots\otimes v_l'$ where again the weight of
$v_j'$ is equal to the weight of $v_+$ minus some multiple of $\eps_i-\eps_{i+1}$
for each $j=1,\dots,l-1$, 

So now we proceed to prove the claim just for 
$x = D_1^{(r)}, D_2^{(r)}, E_1^{(r)}$
and $F_1^{(r)}$ and all meaningful $r$.
For each of these choices for $x$, explicit formulae for
$\kappa'(x) \in W(\pi')$ and $\kappa\circ\phi(x) 
\in W(\pi)$ are given by (\ref{thedef}).
On applying the Miura transforms one obtains explicit formulae for
$(\xi\circ\kappa\circ\phi)(x)$ and
$(\phi_1\otimes\cdots\otimes \phi_l \circ\xi'\circ\kappa')(x)$
as elements of $U(\mathfrak{gl}_{q_1}) 
\otimes\cdots\otimes U(\mathfrak{gl}_{q_l})$.
By considering these formulae directly, one observes finally that
$(\xi\circ\kappa\circ\phi)(x)-
(\phi_1\otimes\cdots\otimes \phi_l \circ\xi'\circ\kappa')(x)$ 
is a linear
combination of terms of the form $x_1\otimes\cdots\otimes x_l$
such that some $x_j$ ($j=1,\dots,l-1$) 
annihilates $v_j$ by weight considerations, which proves the claim.
Let us explain this last step in detail just in 
the case $x = D_2^{(r)}$, all the other
cases being entirely similar.
In this case, we have that 
$$
(\xi\circ\kappa\circ\phi-
\phi_1\otimes\cdots\otimes \phi_l \circ\xi'\circ\kappa')(x)
=
\sum_{\substack{i_1,\dots,i_r\\ j_1,\dots,j_r}}
(-1)^{\#\{s=1,\dots,r-1\:|\:\row(j_s) \leq i\}}
e_{i_1,j_1} \cdots e_{i_r,j_r},
$$
where we are identifying $U(\mathfrak{gl}_{q_1})
\otimes\cdots\otimes U(\mathfrak{gl}_{q_l})$ with $U(\mathfrak{h})$
as usual,
and the sum is over $1 \leq i_1,\dots,i_r,j_1,\dots,j_r \leq n$
with
\begin{enumerate}
\item[(a)] $\row(i_1) = \row(j_r) = i+1$;
\item[(b)] $\col(i_s) = \col(j_s)$ for all $s =1,\dots,r$;
\item[(c)] $\row(j_s) = \row(i_{s+1})$ for all $s=1,\dots,r-1$;
\item[(d)] if $\row(j_s) \geq i+1$ then $\col(j_s) < \col(i_{s+1})$ for all $s=1,\dots,r-1$;
\item[(e)] if $\row(j_s) \leq i$ then $\col(j_s) \geq \col(i_{s+1})$ for all $s=1,\dots,r-1$;
\item[(f)] $\row(j_s) \notin \{i,i+1\}$ for at least one $s = 1,\dots,r-1$.
\end{enumerate}
Take such a monomial
$e_{i_1,j_1} \cdots e_{i_r,j_r} \in
U(\mathfrak{gl}_{q_1})
\otimes\cdots\otimes U(\mathfrak{gl}_{q_l})$.
Let $c$ be minimal such that there exists
$j_s$ with $\col(j_s) = c$ and $\row(j_s) \notin \{i,i+1\}$, then
take the maximal such $s$.
Consider the component of
$e_{i_1,j_1} \cdots e_{i_r,j_r}$
in the $c$th tensor position $U(\mathfrak{gl}_{q_c})$.
If $\row(j_s) > i+1$, then by the choices of $c$ and $s$, this component
is of the form $u e_{i_s,j_s} u'$ where $\row(i_s) \leq i+1 < \row(j_s)$
and the weight of $u'$ is some multiple
of $\eps_i-\eps_{i+1}$.
Similarly if $\row(j_s) < i$ then  this
component is of the form $u e_{i_{s+1},j_{s+1}} u'$
where $\row(j_{s+1}) \geq i > \row(i_{s+1})$ and the weight of $u'$
is some multiple of $\eps_i-\eps_{i+1}$.
In either case, this component annihilates the vector
$v_c \in L(A_c)$ by weight considerations.
\end{proof}

Now let $L_j$ be the irreducible $U(\mathfrak{gl}_{q_j'})$-submodule
of $L(A_j)$ generated by the highest weight vector $v_+$
for each $j=1,\dots,l-1$, embedding
$U(\mathfrak{gl}_{q_j'})$ into $U(\mathfrak{gl}_{q_j})$
via $\phi_j$. Similarly, let $L_l$ be the
$U(\mathfrak{gl}_{q_l'})$-submodule of
$L(A_l)$ generated by the vector $w v_+$.
Recall by the hypotheses in Lemma~\ref{techy} that the tableau
 $A_l$ is column strict and 
$t \geq w^{-1}(t) < w^{-1}(t+1)$. It follows that 
the vector $w v_+ \in L_l$ is a highest weight vector
for the action of $U(\mathfrak{gl}_{q_l'})$ with $e_{1,1}$
acting as $(a+i-1)$ and $e_{2,2}$ acting as $(b+i)$,
for some $b < a \geq a_p$. In particular $L_l$ is also
irreducible.
So in our usual notation the $W(\pi')$-module 
$L_1 \boxtimes \cdots \boxtimes L_l$ is isomorphic to the tensor product
$$
L(\substack{\phantom{a_1}\\c_1+i-1}) \boxtimes \cdots \boxtimes L(\substack{\phantom{a_1}\\c_q+i-1}) 
\boxtimes L(\substack{a_1+i-1\\b_1+i-1})
\boxtimes\cdots\boxtimes L(\substack{a_{p-1}+i-1\\b_{p-1}+i-1}) 
\boxtimes L(\substack{a+i-1\\b+i-1})
$$
for some $b < a \geq a_p$.
Using the remaining hypotheses (i)--(iv) from Lemma~\ref{techy},
we apply Corollaries~\ref{hw2}(i) and \ref{hw4}(i), or rather the slightly
stronger versions of these corollaries described in Remarks~\ref{r1}
and \ref{rem2}, repeatedly to this tensor product working from right to left
to deduce that
$L_1\boxtimes\cdots\boxtimes L_l$ is actually a highest weight $W(\pi')$-module
generated by the highest weight vector
$v_+\otimes\cdots\otimes v_+ \otimes w v_+$.
Hence in particular, since $s_t w v_+ \in L_l$, we get that
$$
v_+\otimes\cdots\otimes v_+ \otimes s_t w v_+
\in 
W(\pi')(v_+\otimes \cdots \otimes w v_+).
$$
In view of Lemma~\ref{techy2}, this completes the proof of
Lemma~\ref{techy}.

\section{Characters of standard modules}\label{sscsm}
We wish to 
explain how to compute the Gelfand-Tsetlin characters of the standard modules
$\{V(A)\:|\:A \in \Col(\pi)\}$ from (\ref{VAdef}).
In view of (\ref{multy}) it suffices just to consider the special
case that $\pi$ consists of a single column of height $m \leq n$,
when $W(\pi) = U(\mathfrak{gl}_m)$.
Take $A \in \Col(\pi)$
with entries $a_1 > \dots > a_m$ read from top to bottom.
Choose an arbitrary scalar $c \in\C$
so that $a_m+m-1 \geq c$.
Then 
$$
(b_1,\dots,b_m) := 
(a_1-c,a_2+1-c,\dots,a_m+m-1-c)
$$ 
is a partition.
Draw its Young diagram in the usual English way and define the
{\em residue} of the box in the $i$th row and $j$th column to
be $(j-i)$.
For example, if $(b_1,b_2,b_3) = (5,3,2)$
then the Young diagram with boxes labelled by their residues
is as follows
$$
\Diagram{${0}$&${1}$&${2}$&${3}$&${4}$\cr -${1}$&${0}$&${1}$\cr -${2}$&-${1}$\cr}
$$
Given a filling $t$ of the boxes of this diagram with the integers
$\{1,\dots,m\}$ we associate the monomial
\begin{equation}
x(t) := 
\prod_{i=1}^m
\prod_{j=1}^{b_i}
x_{n-m+t_{i,j}, c+j-i}\in \widehat{\Z}[\mathscr P_n]
\end{equation}
where $t_{i,j}$ denotes the entry of $t$ in the $i$th row and $j$th column
and $x_{i,a}$ and $y_{i,a}$ are as in (\ref{xia})--(\ref{yia}).
Then we have that
\begin{equation}\label{st}
\ch V(A) = 
y_{n-m+1,c} y_{n-m+2,c-1} \cdots y_{n,c-m+1}\times 
\sum_t x(t)
\end{equation}
summing over all fillings $t$ of the boxes of the diagram
with integers $\{1,\dots,m\}$ such that the entries are weakly
increasing along rows from left to right and strictly increasing
down columns from top to bottom.
The proof of this formula is based like the proof of
Theorem~\ref{vermathm} on branching $V(A)$ from
$\mathfrak{gl}_m$ to $\mathfrak{gl}_{m-1}$.
This time however the restriction is completely understood by the classical
branching theorem for finite dimensional representations of $\mathfrak{gl}_m$, so 
everything is easy.
The closest reference that we could find in the literature is
\cite[Lemma 2.1]{NT}; see also \cite{GT, Ch0} and \cite[Lemma 4.7]{FM}
(the last of these references greatly influenced our choice of notation here).

For example, suppose that $m=n$ and that the entries of $A$
are $1,-1,-2,\dots,1-n$
from top to bottom. Then $V(A)$ is the $n$-dimensional
{natural representation} of $\mathfrak{gl}_n$.
Taking $c=0$, the possible fillings of the Young diagram
$\diagram{\cr}$ are
$\diagram{$\scriptstyle{1}$\cr}$,
$\diagram{$\scriptstyle{2}$\cr}$,\dots,
$\diagram{$\scriptstyle{n}$\cr}$.
Hence
\begin{equation}
\ch V(A) = x_{1,0}+x_{2,0}+\cdots+x_{n,0}.
\end{equation}
Let us make a few further comments, still assuming that $m=n$.
By (\ref{xy}), we have that
\begin{equation}
y_{i,a} = \left\{
\begin{array}{ll}
x_{i,-i+1} x_{i,-i+2}\cdots x_{i,a-1}&\hbox{if $a > 1-i$,}\\
1&\hbox{if $a = 1-i$,}\\
x_{i,a}^{-1} x_{i,a+1}^{-1}\cdots x_{i,-i}^{-1}
&\hbox{if $a < 1-i$}
\end{array}\right.
\end{equation}
for each $i=1,\dots,n$.
Hence if the scalar $c$ in (\ref{st}) is an integer, i.e. if the representation
$V(A)$ is a {\em rational representation} of $\mathfrak{gl}_n$,
then $\ch V(A)$ belongs to the
subalgebra $\Z[x_{i,a}^{\pm 1}\:|\:i=1,\dots,n, a \in \Z]$
of $\widehat{\Z}[\mathscr P_n]$.
Moreover, the character of a rational representation of 
$\mathfrak{gl}_n$ in the usual sense can be deduced from
its Gelfand-Tsetlin character by applying the algebra homomorphism
\begin{equation}
\Z[x_{i,a}^{\pm 1}\:|\:i=1,\dots,n, a \in \Z]
\rightarrow
\Z[x_{i}^{\pm 1}\:|\:i=1,\dots,n],
\quad
x_{i,a} \mapsto x_i.
\end{equation}
Finally, if one can choose the scalar $c$ in (\ref{st}) to be $0$, i.e.
if the representation $V(A)$
is actually a {\em polynomial representation} of $\mathfrak{gl}_n$,
then the formula (\ref{st}) is especially simple since
the leading monomial
$y_{1,c} y_{2,c-1} \cdots y_{n,c-n+1}$ is equal to $1$.
So the Gelfand-Tsetlin character of any polynomial representation of 
$\mathfrak{gl}_n$ belongs to the 
subalgebra $\Z[x_{i,a}\:|\:i=1,\dots,n, a \in \Z]$ of
$\widehat{\Z}[\mathscr P_n]$.

\section{Grothendieck groups}\label{ssgg}
Let us at long last introduce some categories of 
$W(\pi)$-modules. First, let $\mathcal M(\pi)$\label{mcat} denote the category
of all finitely generated, 
admissible $W(\pi)$-modules.
Obviously $\mathcal M(\pi)$ is an abelian category closed under taking
finite direct sums. Note that the duality $?^\tau$ 
defines a contravariant equivalence
$\mathcal M(\pi) \rightarrow \mathcal M(\pi^t)$.
Also, for any other pyramid $\dot\pi$ with the same row lengths as $\pi$,
the isomorphism $\iota$ from (\ref{iotadef2}) induces an
isomorphism $\mathcal M(\pi) \rightarrow \mathcal M(\dot\pi)$.

\begin{Lemma}\label{cl}
Every module in the category $\mathcal M(\pi)$ has a composition series.
\end{Lemma}

\begin{proof}
Copying the standard 
proof that modules in the usual category $\mathcal O$
have composition series, it suffices to prove the lemma for
the generalized Verma
module $M(A)$, $A \in \Row(\pi)$.
In that case it follows because all the weight spaces of $M(A)$ are 
finite dimensional, and moreover there are only finitely 
many irreducibles $L(B)$
with the same central character as $M(A)$ by Lemma~\ref{cx}.
\end{proof}

Hence, the Grothendieck group $[\mathcal M(\pi)]$ of the category
$\mathcal M(\pi)$
is the free abelian group with basis $\{[L(A)]\:|\:A \in \Row(\pi)\}$.
By Theorem~\ref{bruhat}, we have that
$[M(A)] = [L(A)]$ plus an $\N$-linear combination
of $[L(B)]$'s for $B < A$. It follows that
the generalized Verma modules $\{[M(A)]\:|\:A \in \Row(\pi)\}$ also form 
a basis for $[\mathcal M(\pi)]$.
By Theorem~\ref{inj},  the character map $\ch$ 
defines an injective map
\begin{equation}
\ch:[\mathcal M(\pi)] \hookrightarrow \widehat{\Z}[\mathscr P_n].
\end{equation}
Now suppose $\pi = \pi' \otimes \pi''$ for pyramids $\pi'$
and $\pi''$.
We claim that the tensor product $\boxtimes$ induces a multiplication
\begin{equation}\label{newmu}
\mu:[\mathcal M(\pi')] \otimes [\mathcal M(\pi'')]
\rightarrow [\mathcal M(\pi)].
\end{equation}
To see that this makes sense, 
we need to check that the tensor product
$M' \boxtimes M''$ of $M' \in \mathcal M(\pi')$ and
$M'' \in \mathcal M(\pi'')$ belongs to $\mathcal M(\pi)$.
In view of Lemma~\ref{cl}, it suffices to check this
for generalized Verma modules. So take
 $A' \in \Row(\pi')$ and $A'' \in \Row(\pi'')$.
Then, by Corollary~\ref{tee}, we have that
$$
\ch (M(A')\boxtimes M(A'')) = \ch M(A)
$$ where
$A \sim_{\ro} A' \otimes A''$. In view of Theorem~\ref{inj} and Lemma~\ref{cl}, this shows that $M(A') \boxtimes M(A'')$ has a composition series
with factors belonging to $\mathcal M(\pi)$, hence it belongs to
$\mathcal M(\pi)$ itself.
Moreover,
\begin{equation}
\mu([M(A')] \otimes [M(A'')]) = [M(A)].
\end{equation}
Recalling the decomposition (\ref{blockdec}), 
the category $\mathcal M(\pi)$ has the following {\em ``block'' decomposition}
\begin{equation}\label{blockdecc}
\mathcal M(\pi) = \bigoplus_{\theta \in P} \mathcal M(\pi,\theta)
\end{equation}
where $\mathcal M(\pi,\theta)$ is the full subcategory
of $\mathcal M(\pi)$ consisting of objects all of whose composition
factors are of central character
$\theta$; by convention, we set $\mathcal M(\pi,\theta) = 0$ 
if the coefficients of $\theta$ are not non-negative integers
summing to $N$.
Like in (\ref{vincent}), 
we now restrict our attention just to modules with integral central characters:
let 
\begin{equation}\label{blockdeccc}
\mathcal M_0(\pi) := \bigoplus_{\theta \in P_\infty \subset P}
\mathcal M(\pi,\theta).
\end{equation}
The Grothendieck group $[\mathcal M_0(\pi)]$ has the two natural bases
$\{[M(A)]\:|\:A \in \Row_0(\pi)\}$ and
$\{[L(A)]\:|\:A \in \Row_0(\pi)\}$.

Next recall the definition of the $U_\Z$-module 
$S^\pi(V_\Z)$ 
from (\ref{sdef}).
This is also a free abelian group, with two natural bases
$\{M_A\:|\:A \in \Row_0(\pi)\}$
and $\{L_A\:|\:A \in \Row_0(\pi)\}$.
Define an isomorphism of abelian groups
\begin{equation}\label{fountain}
k:S^\pi(V_\Z) \rightarrow [\mathcal M_0(\pi)],\qquad
M_A \mapsto [M(A)]
\end{equation}
for each $A \in \Row_0(\pi)$.
Under this isomorphism, the $\theta$-weight space of
$S^\pi(V_\Z)$ corresponds to the block component
$[\mathcal M(\pi,\theta)]$ of $[\mathcal M_0(\pi)]$,
for each $\theta \in P_\infty$.
Moreover, the isomorphism is compatible with the multiplications
$\mu$ arising from (\ref{themapmu}) and (\ref{newmu})
in the sense that for every decomposition $\pi = \pi' \otimes \pi''$
the following diagram commutes:
\begin{equation}\label{mad}
\begin{CD}
S^{\pi'}(V_\Z) \otimes S^{\pi''}(V_\Z) &@>\mu>>& S^\pi(V_\Z)\\
@Vk \otimes k VV&&@VV k V\\
[\mathcal M_0(\pi')] \otimes [\mathcal M_0(\pi'')] &@>>\mu>& [\mathcal M_0(\pi)]
\end{CD}
\end{equation} 
Now we can formulate the following conjecture, which may be viewed
as a more precise formulation in type $A$ of \cite{VD}.
Note this conjecture is true if $\pi$ consists of a single column; 
see Theorem~\ref{klconj}.
It is also true if $\pi$ has just two rows,
by comparing Theorem~\ref{sl2} and \cite[Theorem 20]{qla}.

\begin{Conjecture}\label{mainconj}
For each $A \in \Row_0(\pi)$, 
the map $k:S^\pi(V_\Z) \stackrel{\sim}{\rightarrow} [\mathcal M_0(\pi)]$
maps the dual canonical basis element $L_A$
to the class $[L(A)]$ of the irreducible module $L(A)$.
In other words, for every $A,B \in \Row_0(\pi)$, we have that
$$
[M(A):L(B)] = P_{d(\rho(A)) w_0, d(\rho(B)) w_0}(1),
$$
notation as in (\ref{kl1}).
\end{Conjecture}

Let us turn our attention to finite dimensional $W(\pi)$-modules.
Let $\mathcal F(\pi)$\label{fcat} denote the category of all finite dimensional 
$W(\pi)$-modules, a full subcategory of the category $\mathcal M(\pi)$.
Let $\mathcal F_0(\pi) = \mathcal F(\pi) \cap \mathcal M_0(\pi)$.
Like in (\ref{blockdecc})--(\ref{blockdeccc}), we have the block decompositions
\begin{align}
\mathcal F(\pi) &= \bigoplus_{\theta \in P} \mathcal F(\pi,\theta),\\
\mathcal F_0(\pi) &= \bigoplus_{\theta \in P_\infty \subset P}
\mathcal F(\pi,\theta).
\end{align}
By Theorem~\ref{fd}, the Grothendieck group $[\mathcal F(\pi)]$ has basis
$\{[L(A)]\:|\:A \in \Dom(\pi)\}$ coming from the simple modules.
Hence $[\mathcal F_0(\pi)]$ has basis
$\{[L(A)]\:|\:A \in \Dom_0(\pi)\}$;
we refer to these $L(A) \in \mathcal F_0(\pi)$
as the {\em rational} irreducible representations of $W(\pi)$.

Recall the subspace $P^\pi(V_\Z)$ of $S^\pi(V_\Z)$
from $\S$\ref{ssdcb}. 
Comparing (\ref{vadef}) and (\ref{VAdef}) and using (\ref{mad}), 
it follows that the map
$k:S^\pi(V_\Z) \rightarrow [\mathcal M_0(\pi)]$
maps $V_A$ to
$[V(A)]$.
Hence there is a well-defined map
$j:P^\pi(V_\Z) \rightarrow [\mathcal F_0(\pi)]$ such that
$V_A\mapsto [V(A)]$ for each $A \in \Col_0(\pi)$.
Moreover, the following diagram commutes:
\begin{equation}\label{we}
\begin{CD}
P^\pi(V_\Z) &@>>> & S^\pi(V_\Z)\\
@VjVV&&@VVkV\\
[\mathcal F_0(\pi)]&@>>>&[\mathcal M_0(\pi)]
\end{CD}
\end{equation}
where the horizontal maps are the natural inclusions.

\begin{Lemma}
The map $j:P^\pi(V_\Z) \rightarrow [\mathcal F_0(\pi)],
V_A \mapsto [V(A)]$ is an isomorphism of
abelian groups.
\end{Lemma}

\begin{proof}
Arguing with the isomorphism $\iota$,
it suffices to prove this in the special case that $\pi$ is left-justified.
In this case, recall from (\ref{rmap}) that 
$R(A)$ denotes the row equivalence class of
$A \in \Std_0(\pi)$.
By Theorem~\ref{charithm}, for each $A \in \Std_0(\pi)$ 
the standard module
$V(A)$ is a quotient of $M(R(A))$, hence
we have that $V(A) = L(R(A))$ plus an $\N$-linear combination of
$L(B)$'s for $B < A$. It follows that 
$\{[V(A)]\:|\:A \in \Std_0(\pi)\}$ is a basis
for $[\mathcal F_0(\pi)]$.
Since the map $j:P^\pi(V_\Z) \rightarrow [\mathcal F_0(\pi)]$
maps the basis $\{V_A\:|\:A \in \Std_0(\pi)\}$ of $P^\pi(V_\Z)$
onto this basis of $[\mathcal F_0(\pi)]$, 
it follows that
$j$ is indeed an isomorphism. 
\end{proof}

This lemma implies that $\{[V(A)]\:|\:A \in \Std_0(\pi)\}$
is a basis for the Grothendieck group
$[\mathcal F_0(\pi)]$. Hence, the Gelfand-Tsetlin character
of any module in  $\mathcal F_0(\pi)$ belongs to the subalgebra
$\Z[x_{i,a}^{\pm 1}\:|\:i=1,\dots,n, a \in \Z]$
of $\widehat{Z}[\mathscr P_n]$, since we know already that this is true for 
the standard modules.
In the next lemma we extend this ``standard basis'' from $[\mathcal F_0(\pi)]$
to all of the Grothendieck group $[\mathcal F(\pi)]$.
Recall for the statement the definition of the relation $\parallel$
on $\Std(\pi)$ from the paragraph after (\ref{rmap}).

\begin{Lemma}
For $A, B \in \Std(\pi)$ we have that
$[V(A)] = [V(B)]$ if and only if $A \parallel B$.
The elements of the set
 $\{[V(A)]\:|\:A \in \Std(\pi)\}$
form a basis for $[\mathcal F(\pi)]$.
In particular, the elements $\{[V(A)]\:|\:A \in \Std_0(\pi)\}$
form a basis for $[\mathcal F_0(\pi)]$.
\end{Lemma}

\begin{proof}
If $A \parallel B$, it is obvious 
from (\ref{multy}) that $[V(A)] = [V(B)]$.
We have already proved the last statement 
about $\mathcal F_0(\pi)$, and
we know that
$[V(A)] \neq [V(B)]$
for distinct $A, B \in \Std_0(\pi)$ 
The remaining parts of the lemma
are consequences of these two statements and
Theorem~\ref{irret}.
\end{proof}

\chapter{Character formulae}\label{sfdr}

Throughout the chapter, we 
fix a pyramid $\pi = (q_1,\dots,q_l)$ 
with associated shift matrix $\sigma = (s_{i,j})_{1 \leq i,j \leq n}$
as usual.
Conjecture~\ref{mainconj} 
immediately implies that the isomorphism 
$j:P^\pi(V_\Z) \rightarrow [\mathcal F_0(\pi)]$
from (\ref{we}) maps
the dual canonical basis 
of the polynomial representation
$P^\pi(V_\Z)$ to the basis of the Grothendieck group 
$[\mathcal F_0(\pi)]$ arising from irreducible modules.
In this chapter, we will give an independent proof of this
statement.
Hence we can in principle compute the Gelfand-Tsetlin characters of all finite dimensional
irreducible $W(\pi)$-modules.

\section{Skryabin's theorem}\label{ssskryabin}

We begin by recalling the relationship between
the algebra $W(\pi)$ and the representation theory of
$\mathfrak{g}$.
Let $\C_\chi$ denote the one dimensional $\mathfrak{m}$-module
defined by the character $\chi$.
Also recall the definitions
(\ref{etadef})--(\ref{origdef}).
Introduce the {\em generalized Gelfand-Graev representation}
\begin{equation}\label{gggr}
Q_\chi := U(\mathfrak{g}) / U(\mathfrak{g}) I_\chi \cong
U(\mathfrak{g}) \otimes_{U(\mathfrak{m})} \C_\chi.
\end{equation}
We write $1_\chi$ for the coset of $1 \in U(\mathfrak{g})$ in $Q_\chi$.
Often we work with the {\em dot action}
of $u \in U(\mathfrak{p})$ on $Q_\chi$ defined by
$u \cdot u'1_\chi := \eta(u) u'1_\chi$ for all $u' \in U(\mathfrak{g})$.
By the definition of $W(\pi)$, right multiplication by $\eta(w)$
leaves $U(\mathfrak{g}) I_\chi$ invariant for each $w \in W(\pi)$.
Hence, there is a well-defined
right $W(\pi)$-module structure on $Q_\chi$ 
such that $(u \cdot 1_\chi) w = uw \cdot 1_\chi$
for $u \in U(\mathfrak{p})$ and $w \in W(\pi)$.
This makes the $\mathfrak{g}$-module $Q_\chi$ into a
$(U(\mathfrak{g}), W(\pi))$-bimodule.
As explained in the introduction of \cite{BK}, 
the associated representation $W(\pi) \rightarrow 
\End_{U(\mathfrak{g})}(Q_\chi)^{\op}$ is actually an isomorphism.

Let ${\mathcal W}(\pi)$ 
denote the category of
{\em generalized Whittaker modules} of type $\pi$, that is, 
the category  of all $\mathfrak{g}$-modules on which
$(x-\chi(x))$ acts locally nilpotently for all $x\in \mathfrak{m}$. 
For any $\mathfrak{g}$-module $M$, let
\begin{equation}\label{whdef}
\Wh(M):=\{v\in M\:|\: xv=\chi(x)v\ \text{for all $x\in
\mathfrak{m}$}\}.
\end{equation}
Given $w\in W(\pi)$ and $v \in \Wh(M)$,
the vector $w \cdot v := \eta(w) v$ again belongs to $\Wh(M)$,
so $\Wh(M)$ is a left $W(\pi)$-module via the dot action.
In this way, we obtain a functor 
$\Wh$ from
${\mathcal W}(\pi)$ to the
category of all left $W(\pi)$-modules.
In the other direction, $Q_\chi \otimes_{W(\pi)} ?$ is a functor
from $W(\pi)\Mod$ to ${\mathcal W}(\pi)$.
The functor $\Wh$ is isomorphic in an obvious way to the functor
$\hom_{U(\mathfrak{g})}(Q_\chi, ?)$,
so adjointness of tensor and hom gives rise to
a canonical adjunction between
the functors $Q_\chi \otimes_{W(\pi)} ?$
and $\Wh$.
The unit and the counit of this canonical adjunction
are defined by
$M \mapsto \Wh(Q_\chi \otimes_{W(\pi)} M),
v \mapsto 1_\chi \otimes v$ for $M \in W(\pi)\Mod$ and $v \in M$, and by
$Q_\chi \otimes_{W(\pi)} \Wh(M)
\rightarrow M, u1_\chi \otimes v \mapsto uv$ for $M \in \mathcal W(\pi)$, 
$u \in U(\mathfrak{g})$ and $v \in \Wh(M)$, respectively.
{\em Skryabin's theorem} \cite{Skry} asserts that these maps are 
actually isomorphisms,
so that the functors $\Wh$ and
$Q_\chi\otimes_{W(\pi)} ?$ are quasi-inverse equivalences
between the categories $\mathcal W(\pi)$ and
$W(\pi)\Mod$. 

Skryabin also proved that $Q_\chi$ is a free right $W(\pi)$-module
and explained how to write down an explicit basis,
as we briefly recall.
Let $b_1,\dots,b_h$ be a homogeneous basis for $\mathfrak{m}$
such that each
$b_i$ of degree $-d_i$.
The elements $[b_1,e], \dots, [b_h,e]$ are again linearly independent,
and $[b_i,e]$ is of degree $(1-d_i)$.
Hence there exist elements $a_1,\dots,a_h \in \mathfrak{p}$
such that each
 $a_i$ is of degree $(d_i -1)$
and
\begin{equation}
([a_i, b_j],e) = (a_i, [b_j,e]) = \delta_{i,j}.
\end{equation}
Now it follows from \cite{Skry} that
the elements
$\{a_1^{i_1}\cdots a_h^{i_h}\cdot 1_\chi\:|\:
i_1,\dots,i_h \geq 0\}$ form a basis for
$Q_\chi$ as a free right $W(\pi)$-module.
Hence, there is a unique right $W(\pi)$-module homomorphism
$\p: Q_\chi \twoheadrightarrow W(\pi)$
defined by
\begin{equation}\label{tiger}
\p(a_1^{i_1} \cdots a_h^{i_h}\cdot 1_\chi)
= \delta_{i_1,0} \cdots \delta_{i_h,0}
\end{equation}
for all $i_1,\dots,i_h \geq 0$.
In particular, $\p(1_\chi) = 1$.

\section{Tensor identities}\label{ssti}
Throughout the section, we
let $V$ be a 
finite dimensional $\mathfrak{g}$-module
with fixed basis $v_1,\dots,v_r$.
Define the coefficient functions
$c_{i,j} \in U(\mathfrak{g})^*$ from the equation
\begin{equation}\label{cffuncs}
u v_j = \sum_{i=1}^r c_{i,j}(u)v_i
\end{equation} 
for all $u \in U(\mathfrak{g})$.
Given any $M \in \mathcal W(\pi)$, it is clear that
$M \otimes V$ (the usual tensor product of $\mathfrak{g}$-modules)
also belongs to the category $\mathcal W(\pi)$.
Thus $? \otimes V$ gives an exact functor 
from $\mathcal W(\pi)$ to $\mathcal W(\pi)$.
Using Skryabin's equivalence of categories, we can 
transport this functor directly to the category
$W(\pi)\Mod$: for a $W(\pi)$-module $M$, let
\begin{equation}\label{circledast}
M \circledast V := \Wh((Q_\chi \otimes_{W(\pi)} M) \otimes V).
\end{equation}
This defines an exact functor 
$?\circledast V:W(\pi)\Mod\rightarrow W(\pi)\Mod$.
The following lemma is a reformulation of \cite[Theorem 4.2]{Ly}.
For the statement, fix a right
$W(\pi)$-module homomorphism
$\p:Q_\chi \twoheadrightarrow W(\pi)$ with $\p(1_\chi) = 1$;
such maps exist by (\ref{tiger}).

\begin{Theorem}\label{el}
For any left $W(\pi)$-module $M$ and any $V$ as above,
the restriction of the map
$(Q_\chi \otimes_{W(\pi)} M) \otimes
V \rightarrow M \otimes V,
(u1_\chi \otimes m) \otimes v \mapsto \p(u1_\chi) m \otimes v$ 
defines a natural vector space isomorphism
$$
\chi_{M,V}:M \circledast V \stackrel{\sim}{\longrightarrow} M \otimes V.
$$
The inverse isomorphism maps $m \otimes v_j$ to 
$\sum_{i=1}^r (x_{i,j}\cdot 1_\chi \otimes m) \otimes v_i$,
where $(x_{i,j})_{1 \leq i,j \leq r}$ is the (necessarily invertible) matrix
with entries in $U(\mathfrak{p})$ determined uniquely by the properties
\begin{enumerate}
\item[\rm(i)] 
$\p(x_{i,j}\cdot 1_\chi) = \delta_{i,j}$;
\item[\rm(ii)]
$\displaystyle [x,\eta(x_{i,j})] 
+ \sum_{s=1}^r c_{i,s}(x) \eta(x_{s,j}) \in U(\mathfrak{g}) I_\chi$
for all $x \in \mathfrak{m}$.
\end{enumerate}
Finally, if $(x_{i,j}')_{1 \leq i, j \leq r}$ is another
matrix with entries in $U(\mathfrak{p})$ satisfying (ii) (primed),
then 
$x_{i,j}' = \sum_{k=1}^r x_{i,k} w_{k,j}$
where $(w_{i,j})_{1 \leq i, j \leq r}$ is the matrix
with entries in $W(\pi)$ defined from the equation
$\p(x_{i,j}' \cdot 1_\chi) = w_{i,j}$.
\end{Theorem}

\begin{proof}
For any vector space $M$, let
$E_M$ denote the space of all
linear maps $f:U(\mathfrak{m})\rightarrow M$
which annihilate $(I_\chi)^p$ for $p \gg 0$,
viewed as an $\mathfrak{m}$-module
with action defined by $(xf)(u) = f(ux)$ for $x \in \mathfrak{m}$,
$u \in U(\mathfrak{m})$.
In the case $M = \C$, we denote $E_M$ simply by $E$.
Skryabin proved the following fact in the course
of \cite{Skry}:
for $M \in W(\pi)\Mod$ there is a natural
$\mathfrak{m}$-module isomorphism
$$
\phi_M:Q_\chi \otimes_{W(\pi)} M \rightarrow 
E_M
$$
defined by $\phi_M(u'1_\chi \otimes m)(u) = \p(uu'1_\chi)m$ 
for $u \in U(\mathfrak{m})$,
$u' \in U(\mathfrak{g})$ and $m \in M$.
Using the fact that $\mathfrak{m}$ is nilpotent and $V$ is finite dimensional,
one checks that evaluation at $1$ defines a natural isomorphism
$$
\xi_V:\Wh(E \otimes V) \stackrel{\sim}{\longrightarrow} V,
\qquad
\sum_{i=1}^r f_i \otimes v_i \mapsto \sum_{i=1}^r f_i(1) v_i.
$$
Finally, there is an obvious
 natural isomorphism
$\psi_M:M \otimes E \rightarrow E_M$ mapping
$m \otimes f \in M \otimes E$ to the function
$u \mapsto f(u)m$.
Combining these things, we obtain
the following natural isomorphisms:
$$
\begin{CD}
M  \circledast V
=
\Wh((Q_\chi \otimes_{W(\pi)} M) \otimes V)
&@>\varphi_M \otimes \id_{V}>>&
\Wh(E_M \otimes V)\\
&@>\psi_M^{-1} \otimes \id_{V}>>
&M\otimes\Wh(E \otimes V)
&@>\id_M \otimes \xi_{V}>>
&M \otimes V.
\end{CD}
$$
Let $\chi_{M,V}:M \circledast V \rightarrow M \otimes V$
denote the composite isomorphism.

Assume in this paragraph that $M = W(\pi)$, the regular $W(\pi)$-module.
In this case, the inverse image of $1 \otimes v_j$
under the isomorphism $\chi_{M,V}$
can be written as
$\sum_{i=1}^r (x_{i,j}\cdot 1_\chi \otimes 1)
\otimes v_i$ for unique elements $x_{i,j} \in U(\mathfrak{p})$.
Now compute to see that 
$$
\xi_V^{-1}(v_j) = \sum_{i=1}^r f_{i,j} \otimes v_i \in \Wh(E \otimes V)
$$
for elements $f_{i,j} \in E$
with $f_{i,j}(u) = c_{i,j}(u^*)$ for $u \in U(\mathfrak{m})$.
Here, $*:U(\mathfrak{m}) 
\rightarrow U(\mathfrak{m})$
is the antiautomorphism with $x^* = \chi(x) - x$ for each $x \in
\mathfrak{m}$.
So,
$$
(\psi_{M} \otimes \id_V) \circ (\id_{M} \otimes \xi_V^{-1})
(1 \otimes v_j) = \sum_{i=1}^r \hat f_{i,j} \otimes v_i,
$$
where $\hat f_{i,j} \in E_M$ satisfies $\hat f_{i,j}(u) = c_{i,j}(u^*) 1$.
On the other hand,
$$
(\varphi_M \otimes \id_V)\bigg(\sum_{i=1}^r (x_{i,j} \cdot 1_\chi \otimes 1) \otimes v_i\bigg) = \sum_{i=1}^r g_{i,j} \otimes v_i
$$
where $g_{i,j}(u) = \p(u \eta(x_{i,j}) 1_\chi)$.
So each $x_{i,j}$ is determined by the property that 
\begin{equation}\label{theprop}
\p(u \eta(x_{i,j})1_\chi) = c_{i,j}(u^*)
\end{equation}
for all $u \in U(\mathfrak{m})$.
Taking $u =1$ in (\ref{theprop}), 
we see that $\p(x_{i,j}\cdot 1_\chi) = \delta_{i,j}$,
as in property (i). Moreover, $\sum_{i=1}^r (x_{i,j}\cdot 1_\chi \otimes 1)
\otimes v_i$ is a Whittaker vector, which is 
equivalent to property (ii). 
Conversely, one checks that 
properties (i) and (ii) imply (\ref{theprop}),
hence they also determine the $x_{i,j}$'s uniquely.

Now return to general $M$. Property (ii) implies
that $\sum_{i=1}^r (x_{i,j}\cdot 1_\chi \otimes m) \otimes v_i$
belongs to $M \circledast V$ for any $m \in M$.
By functoriality, the image of this element under the isomorphism
$\chi_{M,V}$ constructed in the first paragraph of the proof
must equal $m \otimes v_j$. By property (i)
this is also its image under 
the restriction of the map
$(Q_\chi \otimes_{W(\pi)} M) \otimes V \rightarrow M \otimes V,
(u1_\chi \otimes m) \otimes v \mapsto \p(u1_\chi)m \otimes v$.
This shows that the isomorphism $\chi_{M, V}$ constructed in the proof
coincides with the map $\chi_{M, V}$ from the statement of the theorem.

To see that the matrix 
$(x_{i,j})_{1 \leq i,j \leq r}$ 
is invertible,
we may assume without loss of generality that the basis $v_1,\dots,v_r$
has the property that $x v_i 
\in \C v_{1}+\cdots+\C v_{i-1}$
for each $i=1,\dots,r$ and $x \in \mathfrak{m}$, i.e.
$c_{i,j}(x) = 0$ for $i \geq j$.
But then, if one replaces $x_{i,j}$ by $\delta_{i,j}$ for all $i \geq j$, 
the new elements still satisfy (\ref{theprop}).
Hence by uniqueness we must already have that $x_{i,j} = \delta_{i,j}$
for $i \geq j$, i.e. the matrix $(x_{i,j})_{1 \leq i,j \leq r}$
is unitriangular, so it is invertible.

Finally suppose $(x_{i,j}')_{1 \leq i,j \leq r}$ is another 
matrix satisfying (ii) (primed).
Taking $M = W(\pi)$ once more, 
$\sum_{i=1}^r (x_{i,j}' \cdot 1_\chi \otimes 1) \otimes v_i$ belongs to
$\Wh((Q_\chi \otimes_{W(\pi)} M) \otimes V)$.
Hence, by what we have proved already,
there exist elements $w_{i,j} \in W(\pi)$ such that
$$
\sum_{i=1}^r (x_{i,j}' \cdot 1_\chi \otimes 1) \otimes v_i
=
\sum_{i,k=1}^r (x_{i,k} \cdot 1_\chi \otimes w_{k,j}) \otimes v_i.
$$
Equating coefficients gives that $x_{i,j}' = 
\sum_{k=1}^r x_{i,k} w_{k,j}$.
With a final application of the
right $W(\pi)$-module homomorphism $\p$ using (i), we get
that $w_{i,j} = \p(x_{i,j}' \cdot 1_\chi)$, which completes the proof.
\end{proof}

Now we can prove the following important ``tensor identity''.

\begin{Corollary}\label{tensorid}
For any $\mathfrak{p}$-module $M$ and any $V$ as above,
the restriction of the 
map
$(Q_\chi \otimes_{W(\pi)} M) \otimes V \rightarrow 
M \otimes V$ sending $(u \cdot 1_\chi \otimes m) \otimes v
\mapsto um \otimes v$ for each $u \in U(\mathfrak{p}), m \in M, v \in V$ 
defines a natural 
isomorphism
$$
\mu_{M,V}: M \circledast V \stackrel{\sim}{\longrightarrow} M \otimes V
$$
of $W(\pi)$-modules.
Here, we are viewing the $U(\mathfrak{p})$-modules
$M$ and $M \otimes V$ on the left and right hand sides as
$W(\pi)$-modules by restriction.
The inverse map sends
$m \otimes v_k$ to
$\sum_{i,j=1}^r
(x_{i,j} \cdot 1_\chi \otimes y_{j,k} m) \otimes v_i$,
where $(x_{i,j})_{1 \leq i,j \leq r}$ is the matrix
defined in Theorem~\ref{el} and
$(y_{i,j})_{1 \leq i,j \leq r}$ is the inverse matrix.
\end{Corollary}

\begin{proof}
Letting $U(\mathfrak{p})$ act on $Q_\chi \otimes_{W(\pi)} M$
via the dot action, the given map $(Q_\chi \otimes_{W(\pi)} M) \otimes V
\rightarrow M \otimes V$ is a $\mathfrak{p}$-module
homomorphism. Hence its restriction $\mu_{M,V}$
is a $W(\pi)$-module homomorphism.
To prove that $\mu_{M,V}$ is an isomorphism, note by Theorem~\ref{el}
that there is a well-defined map
$$
M \otimes V \rightarrow M \circledast V, \quad m \otimes v_k \mapsto
\sum_{i,j=1}^r
(x_{i,j} \cdot 1_\chi \otimes y_{j,k} m) \otimes v_i.
$$
This is a two-sided inverse to $\chi_{M,V}$.
\end{proof}

Let us make some comments about associativity of 
$\circledast$.
Suppose that we are 
given another 
finite dimensional $\mathfrak{g}$-module $V'$.
For any $W(\pi)$-module $M$, Skryabin's equivalence gives an isomorphism
\begin{align*}
(Q_\chi \otimes_{W(\pi)} \Wh((Q_\chi \otimes_{W(\pi)} M) \otimes V)) \otimes V'
&\stackrel{\sim}{\longrightarrow}
((Q_\chi \otimes_{W(\pi)} M) \otimes V) \otimes V',\\
u' 1_\chi \otimes x \otimes v' &\mapsto u'x \otimes v'
\end{align*}
for $u' \in U(\mathfrak{g}), x \in 
\Wh((Q_\chi \otimes_{W(\pi)} M) \otimes V)$
and $v' \in V'$.
So, 
in view of 
the natural associativity isomorphism at the level of $\mathfrak{g}$-modules,
we conclude that the restriction of the linear map
\begin{align*}
(Q_\chi \otimes_{W(\pi)} ((Q_\chi \otimes_{W(\pi)} M) \otimes V))\otimes V'
&\rightarrow
(Q_\chi \otimes_{W(\pi)} M) \otimes V \otimes V',\\
(u' 1_\chi \otimes ((u 1_\chi \otimes m) \otimes v)) \otimes v'
&\mapsto 
(u'((u 1_\chi \otimes m) \otimes v)) \otimes v'
\end{align*}
defines a natural isomorphism
\begin{equation}\label{aVV}
a_{M, V, V'}:(M \circledast V) \circledast V'
\rightarrow M \circledast (V \otimes V')
\end{equation}
of $W(\pi)$-modules.
If $M$ is actually a $\mathfrak{p}$-module, 
it is straightforward to check that the following diagram
commutes:
\begin{equation}\label{Second}
\begin{CD}
(M \circledast V) \circledast V'&@>a_{M, V, V'}>> &M \circledast (V \otimes V') 
\\
@V\mu_{M,V} \circledast \id_{V'}VV&&@VV\mu_{M, V \otimes V'}V\\
(M \otimes V) \circledast V'&@>>\mu_{M \otimes V, V'}> &M \otimes V \otimes V'. 
\end{CD}
\end{equation}
Also, given a third finite dimensional module $V'''$,
the following diagram commutes:
\begin{equation}\label{third}
\begin{CD}
((M \circledast V) \circledast V') \circledast V''
&@>a_{M \circledast V, V', V''}>>& 
(M \circledast V) \circledast (V' \otimes V'')\\
@Va_{M, V, V'} \circledast \id_{V''}VV&&@VVa_{M, V, V' \otimes V''}V\\
(M \circledast (V \otimes V')) \circledast V''
&@>>a_{M, V \otimes V', V''}>& M \circledast (V \otimes V' \otimes V'').
\end{CD}
\end{equation}
Writing $\C$ for the trivial $\mathfrak{g}$-module,
there is for each $W(\pi)$-module $M$ a natural isomorphism
$i_M: M \circledast \C \rightarrow M$ mapping
$(1_\chi \otimes m) \otimes 1 \mapsto m$ for each $m \in M$.
There are also some commutative 
triangles
arising from the compatibility of $i$ with $a$ and with $\mu$, but they are quite 
obvious so we omit them.

Finally, we translate the canonical adjunction between
the functors $? \otimes V$ and $? \otimes V^*$
into an adjunction between $? \circledast V$ and $? \circledast V^*$, where
$V^*$ here denotes the usual dual $\mathfrak{g}$-module.
Let $\overline{v}_1,\dots,\overline{v}_r$ 
be the basis for $V^*$ dual to the basis
$v_1,\dots,v_r$ for $V$.
Then the unit of the canonical adjunction is the map
$\iota:\operatorname{Id} \rightarrow
(? \circledast V) \circledast V^*$ 
defined on a $W(\pi)$-module $M$  to be the composite
\begin{equation}\label{adj1}
\begin{CD}
M&@>i_M^{-1}>\phantom{\displaystyle T}>&M \circledast \C &@>>>&M \circledast (V \otimes V^*)
&@>a_{M, V, V^*}^{-1}>>& (M \circledast V) \circledast V^*,
\end{CD}
\end{equation}
where the second map is $(1_\chi \otimes m)  \otimes 1
\mapsto \sum_{i=1}^r (1_\chi \otimes m) \otimes v_i \otimes \overline{v}_i$.
The counit of the canonical adjunction is the map
$\eps:(? \circledast V^*) \circledast V \rightarrow \operatorname{Id}$
defined on a $W(\pi)$-module $M$ to be the composite
\begin{equation}\label{adj2}
\begin{CD}
(M \circledast V^*) \circledast V &@>a_{M, V^*, V}>\phantom{t}>
&M \circledast (V^* \otimes V)
&@>>>& M \circledast \C &@>i_M >> &M,
\end{CD}
\end{equation}
where the second map is the 
restriction of 
$(u 1_\chi \otimes m) \otimes f \otimes v
\mapsto (u1_\chi \otimes m) \otimes f(v)$.

\section{Translation functors}\label{sstf}
In this section we extend the definition of the
translation functors $e_i, f_i$ from $\S$\ref{sskl}
to the category $\mathcal M_0(\pi)$ from $\S$\ref{ssgg}.
Throughout the section, we let $V$ denote the natural $N$-dimensional
$\mathfrak{g}$-module of column vectors, with standard basis
$v_1,\dots,v_N$.
We first define an endomorphism
\begin{equation}\label{xendo}
x:?\circledast V\rightarrow ?\circledast V
\end{equation}
of the functor $?\circledast V:  
W(\pi)\Mod \rightarrow W(\pi)\Mod$.
On a $W(\pi)$-module $M$,
$x_M$ is the endomorphism of
$M \circledast V = \Wh((Q_\chi \otimes_{W(\pi)} M) \otimes V)$
defined by left multiplication by
$\Omega = \sum_{i,j=1}^N e_{i,j} \otimes e_{j,i}
\in U(\mathfrak{g}) \otimes U(\mathfrak{g})$.
Here, we are treating the $\mathfrak{g}$-module
$Q_\chi \otimes_{W(\pi)} M$ as the first tensor position and $V$
as the second,
so $\Omega((u 1_\chi \otimes m) \otimes v)$ means 
$\sum_{i,j=1}^N (e_{i,j} u 1_\chi \otimes m) \otimes e_{j,i} v$.
Next, we define an endomorphism
\begin{equation}\label{sendo}
s: (? \circledast V) \circledast V \rightarrow (? \circledast V) \circledast V
\end{equation}
of the functor $(? \circledast V) \circledast V:
W(\pi)\Mod \rightarrow W(\pi)\Mod$.
Recalling (\ref{aVV}), we take
$s_M: (M \circledast V) \circledast V
\rightarrow (M \circledast V) \circledast V$ 
to be the composite $a_{M, V, V}^{-1} \circ 
\widehat s_M \circ a_{M, V, V}$, where
$\widehat s_M$ is the endomorphism of 
$M \circledast (V \otimes V)=
\Wh((Q_\chi \otimes_{W(\pi)} M) \otimes V \otimes V)$
defined by left multiplication by 
$\Omega^{[2,3]}$, i.e. $\Omega$ acting on the second and third tensor positions
so $\Omega^{[2,3]} ((u 1_\chi \otimes m)\otimes v \otimes v')$
means $\sum_{i,j=1}^N
(u1_\chi \otimes m) \otimes e_{i,j} v \otimes e_{j,i} v'$
(which equals 
$(u1_\chi \otimes m) \otimes v' \otimes v$).
Actually these definitions are just the natural translations 
through Skryabin's equivalence of categories of the
endomorphisms $x$ and $s$ from $\S$\ref{sskl}
of the functors $? \otimes V$ and $? \otimes V \otimes V$.

More generally, suppose that we are given $d \geq 1$,
and introduce the following endomorphisms of the 
$d$th power 
$(? \circledast V)^d$:
for $1 \leq i \leq d$ and $1 \leq j < d$, let
\begin{equation}\label{xxii}
x_i := (1_{? \circledast V})^{d-i} x (1_{? \circledast V})^{i-1},\qquad
s_j := (1_{? \circledast V})^{d-j-1} s (1_{? \circledast V})^{j-1}.
\end{equation}
There is an easier description of these endomorphisms.
To formulate this, we exploit the natural isomorphism 
\begin{equation}
a_d:(? \circledast V)^d
\stackrel{\sim}{\longrightarrow}? \circledast V^{\otimes d}
\end{equation} 
obtained by iterating
the associativity isomorphism from
(\ref{aVV}).
For $1 \leq i \leq d$ and $1 \leq j < d$, let
$\widehat x_i$ and $\widehat s_j$ denote the endomorphisms
of the functor $? \circledast V^{\otimes d}$ defined by left multiplication
by the elements $\sum_{h=1}^i \Omega^{[h,i+1]}$ and 
$\Omega^{[j+1,j+2]}$, respectively, notation as above.
Then we have that
\begin{equation}\label{xisj}
x_i = a_d^{-1} \circ \widehat x_i \circ a_d,
\qquad
s_j = a_d^{-1} \circ \widehat s_j \circ a_d.
\end{equation}
Using this alternate description, the following 
identities are straightforward to check:
\begin{align}
x_i x_j &= x_j x_i,\label{dhr1}\\
s_i x_i &= x_{i+1} s_i - 1,\\
s_i x_j &= x_j s_i \qquad\qquad\text{if }j \neq i,i+1,\\
s_i^2 &= 1,\label{dhr3}\\
s_i s_j &= s_j s_i \qquad\qquad\text{\,if }|i-j| > 1,\\
s_i s_{i+1} s_i &= s_{i+1} s_i s_{i+1}.\label{dhr5}
\end{align}
These 
are the defining relations of the degenerate affine Hecke algebra
$H_d$.

Let us next bring the adjoint
functor $? \circledast V^*$ into the picture,
where $V^*$ is the dual $\mathfrak{g}$-module.

\begin{Lemma}\label{tid2}
The functors $? \circledast V$ and $? \circledast V^*$
map 
objects in $\mathcal M(\pi)$ to objects in $\mathcal M(\pi)$.
Moreover,
for $A \in \Row(\pi)$, we have that
\begin{enumerate}
\item[\rm(i)] $\ch  (M(A) \circledast V) =
\sum_{i=1}^N \ch M(B_i)$ where $B_i$ is 
the row equivalence class of the tableau obtained by adding $1$ to the 
$i$th entry of a fixed representative for $A$;
\item[\rm(ii)] $\ch  (M(A) \circledast V^*) =
\sum_{i=1}^N \ch M(B_i)$ where $B_i$ is 
the row equivalence class of the tableau obtained by subtracting $1$ from the 
$i$th entry of a fixed representative for $A$.
\end{enumerate}
\end{Lemma}

\begin{proof}
We just prove the statements about $? \circledast V$,
since the ones for $? \circledast V^*$ are similar.
Recall from Corollary~\ref{tee} and Theorem~\ref{inj}
that $M(A)$ has all the same composition factors
as $M(A_1) \boxtimes\cdots\boxtimes M(A_l)$, 
where $A_1 \otimes\cdots\otimes A_l$ is some representative for $A$.
To prove that $? \circledast V$ sends objects in 
$\mathcal M(\pi)$ to objects in $\mathcal M(\pi)$,
it suffices by exactness of
the functor to check that
$(M(A_1) \boxtimes\cdots\boxtimes M(A_l)) \circledast V$
belongs to $\mathcal M(\pi)$. 
Since $M(A_1) \boxtimes\cdots\boxtimes M(A_l)$
is the restriction of a $\mathfrak{p}$-module $M$,
Corollary~\ref{tensorid} implies that
$M \circledast V \cong M \otimes V$ as $W(\pi)$-modules.
Now observe that 
$V$ has a filtration as a 
$\mathfrak{p}$-module with factors
$V_1,\dots,V_{l}$ being the natural modules of the
components $\mathfrak{gl}_{q_1},\dots,\mathfrak{ql}_{q_l}$
of $\mathfrak{h}$, respectively.
Hence $M \otimes V$ has a filtration with factors
$M \otimes V_{i}$. Now apply Lemma~\ref{tid} to each of these
factors in turn, to deduce that 
$M \circledast V$ has a filtration with factors isomorphic to
$$
M(A_1)\boxtimes \cdots\boxtimes M(A_{i-1})
\boxtimes M(B_i) \boxtimes M(A_{i+1}) \boxtimes\cdots\boxtimes M(A_l),
$$
one for each $i=1,\dots,l$ and each $B_i$  obtained from the 
column tableau $A_i$ by adding $1$ to one of its entries.
Hence it belongs to $\mathcal M(\pi)$.
Taking Gelfand-Tsetlin characters gives (i) as well.
\end{proof}

For $\theta \in P_\infty$, let
$\pr_\theta:\mathcal M_0(\pi) \rightarrow
\mathcal M(\pi,\theta)$ be the projection functor
along the decomposition (\ref{blockdeccc}).
Explicitly, for a module $M \in \mathcal M_0(\pi)$, we have that
$\pr_\theta(M)$ is the summand of $M$ defined by (\ref{prtheta2}),
or $\pr_\theta(M) = 0$ if the coefficients of $\theta$ are not non-negative
integers summing to $N$.
In view of Lemma~\ref{tid2}, it makes sense to define
exact functors $e_i, f_i:\mathcal M_0(\pi)
\rightarrow \mathcal M_0(\pi)$ by setting
\begin{align}\label{tf3}
e_i &:= \bigoplus_{\theta \in P_\infty}
\pr_{\theta + (\eps_i - \eps_{i+1})} \circ (? \circledast V^*) \circ
\pr_\theta,\\
f_i &:= \bigoplus_{\theta \in P_\infty}
\pr_{\theta - (\eps_i - \eps_{i+1})} \circ (? \circledast V) \circ
\pr_\theta.\label{tf4}
\end{align}
Note $e_i$ is right adjoint to $f_i$, indeed,
the canonical adjunction 
between
$? \circledast V$ and $? \circledast V^*$ 
from (\ref{adj1})--(\ref{adj2})
induces a canonical adjunction between $f_i$ and $e_i$.
Similarly, $e_i$ is also left adjoint to $f_i$.
Moreover, applying Lemma~\ref{tid2} and taking blocks, we see
for $A \in \Row_0(\pi)$ and $i \in \Z$ that
\begin{equation}
[e_i M(A)] = \sum_{B} [M(B)]
\end{equation}
summing over all $B$ obtained from $A$ by replacing
an entry equal to $(i+1)$ by an $i$, and
\begin{equation}
[f_i M(A)] = \sum_{B} [M(B)]
\end{equation}
summing over all $B$ 
obtained from $A$ by replacing
an entry equal to $i$ by an $(i+1)$; cf. (\ref{tide})--(\ref{tidf}).
Hence if we identify the Grothendieck group 
$[\mathcal M_0(\pi)]$ with the $U_\Z$-module
$S^\pi(V_\Z)$ via the isomorphism
(\ref{fountain}), the maps
on the Grothendieck group
induced by the exact functors $e_i, f_i$
coincide with the action of 
$e_i, f_i\in U_\Z$.
Moreover, for any $M \in \mathcal M_0(\pi)$, we have that
\begin{equation}
M \circledast V = \bigoplus_{i \in \Z} f_i M,\qquad\qquad
M \circledast V^* = \bigoplus_{i \in \Z} e_i M.
\end{equation}

\begin{Lemma}
For $M \in \mathcal M_0(\pi)$,
$f_i M$ coincides  with the generalized $i$-eigenspace of
$x_M\in \End_{W(\pi)}(M \circledast V)$.
\end{Lemma}

\begin{proof}
It suffices to check this on a generalized Verma module $M(A)$
for $A \in \Row_0(\pi)$.
Say the entries of $A$ in some order are $a_1,\dots,a_N$ and let
$B$ be obtained from $A$ by replacing the entry
$a_t$ by $a_t+1$, for some $1 \leq t \leq N$.
Recall the elements
\begin{align*}
Z_N^{(1)} &= \sum_{i=1}^N (e_{i,i}-N+i),\\
Z_N^{(2)} &= \sum_{i < j} \left(
(e_{i,i}-N+i)(e_{j,j}-N+j) - e_{i,j} e_{j,i}
\right) 
\end{align*}
of $Z(U(\mathfrak{g}))$ from (\ref{zNu}).
For any $\mathfrak{g}$-module $M$, 
the operator $\Omega$ acts on $M \otimes V$
in the same way as $Z_N^{(2)} \otimes 1
+ Z_N^{(1)} \otimes 1-\Delta(Z_N^{(2)})$.
Also by Lemma~\ref{cx}, $\psi(Z_N^{(1)})$ acts on $M(A)$ as
$\sum_{r=1}^N a_r$ and 
$\psi(Z_N^{(2)})$ acts as $\sum_{r < s} a_r a_s$.
It follows that $x_{M(A)}$ stabilizes any $W(\pi)$-submodule of
$M(A) \circledast V$, and it acts on any irreducible subquotient 
having the same central character as $L(B)$
by scalar multiplication by
$$
a_t = \sum_{r < s} a_r a_s + \sum_{r=1}^N a_r-\sum_{r < s} (a_r + \delta_{r,t}) (a_s + \delta_{s,t}).
$$
Since $M(A) \circledast V = \bigoplus_{i \in \Z} f_i M(A)$
and all irreducible subquotients of $f_i M(A)$ have the same central
character as $L(B)$ for some $B$ obtained from $A$ by replacing an
entry $i$ by an $(i+1)$,
this identifies $f_i M(A)$ as the generalized $i$-eigenspace
of $x_{M(A)}$.
\end{proof}

As in \cite[$\S$7.4]{CR}, this lemma together with the relations
(\ref{dhr1})--(\ref{dhr5}) imply that the endomorphisms
$x$ and $s$ 
restrict to well-defined endomorphisms also denoted $x$ and $s$
of the functors $f_i$ and $f_i^2$, respectively. 
Moreover, the identities (\ref{cr2})--(\ref{cr3})
also hold in this setting.
This means that the
category $\mathcal M_0(\pi)$ equipped with the adjoint pair
of functors $(f_i, e_i)$ and the endomorphisms
$x \in \End(f_i)$ and $s \in \End(f_i^2)$ 
is an $\mathfrak{sl}_2$-categorification
in the sense of \cite{CR}, for all $i \in \Z$. So we can appeal
to all the general results developed in \cite{CR} in our study
of the category $\mathcal M_0(\pi)$. 

\begin{Theorem}\label{crlem}
Let $A \in \Row_0(\pi)$ and $i \in \Z$.
\begin{enumerate}
\item[\rm(i)] Define $\eps'_i(A)$ to be the maximal integer
$k \geq 0$ such that $(e_i)^k L(A) \neq 0$.
Assuming $\eps_i'(A) > 0$,
$e_i L(A)$ has irreducible socle and cosocle
isomorphic to $L(\tilde e'_i (A))$ for some
$\tilde e'_i (A) \in \Row_0(\pi)$ with $\eps'_i(\tilde e'_i (A)) = 
\eps'_i(A) - 1$.
The multiplicity of
$L(\tilde e'_i (A))$ as a composition factor of
$e_i L(A)$ is equal to $\eps_i'(A)$, and all other composition
factors are of the form $L(B)$ for $B \in \Row_0(\pi)$ with
$\eps'_i(B) < \eps'_i(A) - 1$.
\item[\rm(ii)] Define $\phi'_i(A)$ to be the maximal integer
$k \geq 0$ such that $(f_i)^k L(A) \neq 0$.
Assuming $\phi_i'(A) > 0$,
$f_i L(A)$ has irreducible socle and cosocle
isomorphic to $L(\tilde f'_i (A))$ for some
$\tilde f'_i (A) \in \Row_0(\pi)$ with $\phi'_i(\tilde f'_i (A)) = 
\phi'_i(A) - 1$.
The multiplicity of
$L(\tilde f'_i (A))$ as a composition factor of
$f_i L(A)$ is equal to $\phi_i'(A)$, and all other composition
factors are of the form $L(B)$ for $B \in \Row_0(\pi)$ with
$\phi'_i(B) < \phi'_i(A) - 1$.
\end{enumerate}
\end{Theorem}

\begin{proof}
This follows from \cite[Lemma 4.3]{CR} and \cite[Proposition 5.23]{CR},
as in the first paragraph of the proof of Theorem~\ref{kuj}.
\end{proof}

\begin{Remark}
This theorem gives a representation theoretic definition of
a crystal structure
$(\Row_0(\pi), \tilde e'_i, \tilde f'_i, \eps'_i, \phi'_i, 
\theta)$
on the set $\Row_0(\pi)$. In $\S$\ref{sscrystals},
we gave a combinatorial definition of another crystal
structure 
$(\Row_0(\pi), \tilde e_i, \tilde f_i, \eps_i, \phi_i, \theta)$
on the same underlying set. If Conjecture~\ref{mainconj} is true, then it 
follows by (\ref{itf1})--(\ref{itf2})
(as in the proof of Theorem~\ref{kuj}) that these two crystal
structures are in fact {\em equal}, that is,
$\eps'_i(A) = \eps_i(A), \phi'_i(A) = \phi_i(A), \tilde e'_i(A) = \tilde e_i(A)$
and $\tilde f'_i(A) = \tilde f_i(A)$ for all $A \in \Row_0(\pi)$.
\end{Remark}

\begin{Remark}
Even without Conjecture~\ref{mainconj}, one can show using
\cite[Lemma 4.3]{CR} and \cite[Theorem 5.37]{BeK}
that the two crystals
$(\Row_0(\pi), \tilde e'_i, \tilde f'_i, \eps'_i, \phi'_i, \theta)$
and
$(\Row_0(\pi), \tilde e_i, \tilde f_i, \eps_i, \phi_i, \theta)$
are at least {\em isomorphic}. 
However, 
there is an identification problem here: without invoking
Conjecture~\ref{mainconj}
we do not know how to prove
that the identity map
on the underlying set $\Row_0(\pi)$
is an isomorphism between the two crystals.
An analogous identification problem arises in a number of other
situations; compare for example \cite{BKclifford} and \cite{BKsergeev}.
\end{Remark}

\section{Translation commutes with duality}\label{ssduality}
There is a right module analogue of Skryabin's theorem.
To formulate it quickly, recall Lemma~\ref{miracle} and
the automorphism $\overline{\eta}$ from (\ref{baretadef}).
Let
\begin{equation}
\overline{Q}_\chi := U(\mathfrak{g}) / I_\chi U(\mathfrak{g}).
\end{equation}
We write $\overline{1}_\chi$ for the coset of $1$ in $\overline{Q}_\chi$,
and define the dot action of $u \in U(\mathfrak{p})$
on $\overline{Q}_\chi$ by
$\overline{1}_\chi u' \cdot u := \overline{1}_\chi u' \overline{\eta}(u)$.
Make the right $U(\mathfrak{g})$-module
$\overline{Q}_\chi$ into a $(W(\pi), U(\mathfrak{g}))$-bimodule
so that $w \overline{1}_\chi\cdot u = \overline{1}_\chi \cdot wu$ for 
each $w \in W(\pi)$ and $u \in U(\mathfrak{p})$.
Let $\overline{\mathcal{W}}(\pi)$ denote the category of all right
$U(\mathfrak{g})$-modules on which $(x-\chi(x))$ acts locally nilpotently
for all $x \in \mathfrak{m}$. 
For $M \in\overline{\mathcal{W}}(\pi)$, let
\begin{equation}\label{bwh}
\overline{\Wh}(M) := \{v \in M\:|\:v x = \chi(x) v\text{ for all }x \in \mathfrak{m}\},
\end{equation}
naturally a right $W(\pi)$-module with dot
action $v \cdot w := v \overline{\eta}(w)$ for $v \in \overline{\Wh}(M)$
and $w \in W(\pi)$. This defines an equivalence of categories
$\overline{\Wh}:\overline{\mathcal{W}}(\pi) \rightarrow \text{mod-}W(\pi)$
with quasi-inverse $? \otimes_{W(\pi)} \overline{Q}_\chi:
\text{mod-}W(\pi)\rightarrow \overline{\mathcal{W}}(\pi)$.
The quickest way to see this is to use the antiautomorphism
$\tau$ from (\ref{taup}) to identify
the category $\overline{\mathcal{W}}(\pi)$
with $\mathcal{W}(\pi^t)$ and
the category
$\text{mod-}W(\pi)$ with
$W(\pi^t)\text{-mod}$.
When that is done, the functor
$\overline{\Wh}:\overline{\mathcal{W}}(\pi) \rightarrow \text{mod-}W(\pi)$
becomes identified with Skryabin's original 
equivalence of categories
$\Wh:\mathcal{W}(\pi^t) \rightarrow W(\pi^t)\text{-mod}$
from $\S$\ref{ssskryabin}.

Given a finite dimensional $\mathfrak{g}$-module $V$ as in $\S$\ref{ssti},
there is also a right module analogue $? \circledast \overline{V}$
of the functor $? \circledast V$.
Here, $\overline{V}$ denotes the dual vector space $V^*$
viewed as a right
 $U(\mathfrak{g})$-module  via $(fx)(v) = f(xv)$ for
$f \in \overline{V}, v \in V$ and $x \in \mathfrak{g}$. 
Then,
by definition, $? \circledast \overline{V}
:\text{mod-}W(\pi) \rightarrow \text{mod-}W(\pi)$ is the functor
defined on objects by
\begin{equation}
M \circledast \overline{V} := 
\overline{\Wh}((M \otimes_{W(\pi)} \overline{Q}_\chi) \otimes \overline{V}).
\end{equation}
Moreover, given another finite dimensional $\mathfrak{g}$-module
$V'$ and any right $W(\pi)$-module $M$, 
there is an associativity isomorphism
\begin{equation}\label{aVV2}
a_{M, \overline{V}, \overline{V}'}:
(M \circledast \overline{V}) \circledast \overline{V}'
\stackrel{\sim}{\longrightarrow}
M \circledast (\overline{V} \otimes \overline{V}')
\end{equation}
defined in an analogous way to (\ref{aVV}).
Another way to think about the functor
$? \circledast \overline{V}$ is to first identify
right $W(\pi)$-modules with left $W(\pi^t)$-modules using
the antiautomorphism $\tau$, then
$? \circledast \overline{V}:
\text{mod-}W(\pi) \rightarrow \text{mod-}W(\pi)$ is 
naturally isomorphic to the functor
$? \circledast V: W(\pi^t)\text{-mod} \rightarrow W(\pi^t)\text{-mod}$
defined as in $\S$\ref{ssskryabin}.
For an admissible left $W(\pi)$-module $M$, 
recall the restricted dual $\overline{M}$ from (\ref{rdual}).
Assuming $\pi$ is left-justified, we are going to prove
that $\circledast$ commutes with duality in the sense
that $\overline{M} \circledast \overline{V} \cong
\overline{M \circledast V}$; equivalently,
$M^\tau \circledast V \cong (M \circledast V)^\tau$.
Although not proved here, this is
true even without the assumption that $\pi$ is left-justified; 
see Remark~\ref{gencase}.

For the proof, 
we say that a (necessarily finite dimensional) $\mathfrak{g}$-module
$V$ is {\em dualizable}
if there is a basis $v_1,\dots,v_r$ 
for $V$ and a pair of mutually inverse matrices 
$(x_{i,j})_{1 \leq i,j \leq r}$
and $(y_{i,j})_{ 1 \leq i,j \leq r}$ with entries in $U(\mathfrak{p})$
such that
\begin{enumerate}
\item[(a)]
$\displaystyle [x, \eta(x_{i,j})] + \sum_{s=1}^r c_{i,s}(x) \eta(x_{s,j})
\in U(\mathfrak{g})I_\chi$
for all $1 \leq i,j \leq r$ and 
$x \in \mathfrak{m}$;
\item[(b)]
$\displaystyle [\overline{\eta}(y_{i,j}),x] + \sum_{s=1}^r \overline{\eta}(y_{i,s})c_{s,j}(x)\in I_\chi
U(\mathfrak{g})$
for all $1 \leq i,j \leq r$ and 
$x \in \mathfrak{m}$.
\end{enumerate}
Here, $c_{i,j} \in U(\mathfrak{g})^*$ is the coefficient function defined
by (\ref{cffuncs}).

\begin{Lemma}\label{all}
Suppose that $V$ is dualizable.
Let $v_1,\dots,v_r$ be any basis for $V$
and $(x_{i,j})_{1 \leq i,j \leq r}$ be any invertible matrix
with entries in $U(\mathfrak{p})$ satisfying property (a) above.
Then the inverse matrix $(y_{i,j})_{1 \leq i,j \leq r}$
satisfies property (b) above.
\end{Lemma}

\begin{proof}
Since $V$ is dualizable there exists some basis
$v_1',\dots,v_r'$ for $V$ and some pair of
mutually inverse matrices
$(x_{i,j}')_{1 \leq i,j \leq r}$
and $(y_{i,j}')_{1 \leq i,j \leq r}$ satisfying properties (a) and (b)
(primed).
Conjugating by an invertible scalar matrix
if necessary, we can assume that
$v_1'=v_1,\dots,v_r' = v_r$.
The last part of Theorem~\ref{el} implies that there is an invertible matrix
$(w_{i,j})_{1 \leq i,j \leq r}$ with entries in $W(\pi)$ such that
$x_{i,j} = \sum_{k=1}^r x_{i,k}' w_{k,j}$.
Let $({v}_{i,j})_{1 \leq i,j \leq r}$ be the inverse matrix.
Then $y_{i,j} = \sum_{k=1}^r {v}_{i,k} y_{k,j}'$.
Using Lemma~\ref{miracle} together with property (b) for $y_{k,j}'$, 
we get for $x \in \mathfrak{m}$ that
\begin{align*}
[\overline{\eta}(y_{i,j}), x] + \sum_{s=1}^r \overline{\eta}(y_{i,s})
c_{s,j}(x)
&=
\sum_{k=1}^r \left(
[\overline{\eta}({v}_{i,k}y_{k,j}'), x] + \sum_{s=1}^r \overline{\eta}(v _{i,k}y_{k,s}')
c_{s,j}(x)\right)\\
&\equiv
\sum_{k=1}^r
\overline{\eta}({v}_{i,k})
\left([\overline{\eta}(y_{k,j}'), x] + \sum_{s=1}^r \overline{\eta}(y_{k,s}')
c_{s,j}(x)\right)\\
&\equiv 0 \pmod{I_\chi U(\mathfrak{g})}.
\end{align*}
Hence $(y_{i,j})_{1 \leq i,j \leq r}$ satisfies property (b).
\end{proof}

\begin{Lemma}\label{barchi}
For any right $W(\pi)$-module $M$ and any dualizable 
$\mathfrak{g}$-module $V$, there is a natural
vector space isomorphism
$$
{\chi}_{M, \overline{V}}: M \circledast \overline{V} 
\rightarrow M \otimes \overline{V}
$$
determined uniquely by the following property. Let $v_1,\dots,v_r$ 
be any basis for $V$, let
$\overline{v}_1,\dots,\overline{v}_r$ be the dual basis for $\overline{V}$,
let $(x_{i,j})_{1 \leq i,j \leq r}$
be the matrix defined in Theorem~\ref{el}
and let
$(y_{i,j})_{1 \leq i,j \leq r}$ be the inverse matrix.
Then 
${\chi}_{M, \overline{V}}$ maps
$\sum_{j=1}^r
(m \otimes \overline{1}_\chi \cdot y_{i,j}) \otimes \overline{v}_j$ to
$m \otimes \overline{v}_i$,
for each $m \in M$ and $1 \leq i \leq r$.
\end{Lemma}

\begin{proof}
By Lemma~\ref{all}, the elements $y_{i,j}$ satisfy property (b).
Therefore, twisting the
conclusion of Theorem~\ref{el} (with $\pi$ replaced by $\pi^t$) by the antiautomorphism $\tau$, we deduce that
the map $\sum_{j=1}^r
(m \otimes \overline{1}_\chi \cdot y_{i,j}) \otimes \overline{v}_j
\mapsto m \otimes \overline{v}_i$ is a vector space isomorphism
$\chi_{M, \overline{V}}:M \circledast \overline{V} \rightarrow
M \otimes \overline{V}$.
Moreover, this definition 
is independent of the initial choice of basis.
It remains to check naturality. Clearly it is natural in $M$.
To see that it is natural in $V$, let $V'$ be another dualizable
$\mathfrak{g}$-module and $f:V \rightarrow V'$ be a $\mathfrak{g}$-module
homomorphism. 
Let $f^*:\overline{V}' \rightarrow \overline{V}$
be the dual map.
We need to show that the following diagram commutes:
$$
\begin{CD}
M \circledast \overline{V}' &@>\chi_{M, \overline{V}'}>>& M \otimes \overline{V}'\phantom{\,.}\\
@V\operatorname{id}_M 
\circledast f^*VV&&@VV\operatorname{id}_M \otimes f^*V\\
M \circledast \overline{V} &@>>\chi_{M, \overline{V}}>& M \otimes \overline{V} \,.
\end{CD}
$$
Pick a basis $v_1',\dots,v_s'$ for $V'$
and let $\overline{v}_1',\dots,\overline{v}_s'$ be the dual basis for
$\overline{V}'$.
Say $f(v_j) = \sum_{i=1}^s a_{i,j} v_i'$, so
$f^*(\overline{v}_i') = \sum_{j=1}^r a_{i,j} \overline{v}_j$.
Let $(x_{i,j}')_{1 \leq i,j \leq s}$ be the matrix defined by 
applying Theorem~\ref{el}
to the chosen basis of $V'$, and let $(y_{i,j}')_{1 \leq i,j \leq s}$
be the inverse matrix.
By the naturality in Theorem~\ref{el}, we have that $\sum_{k=1}^s x_{i,k}' a_{k,j}
= \sum_{k=1}^r a_{i,k} x_{k,j}$.
Hence, $\sum_{k=1}^r a_{i,k} y_{k,j} = \sum_{k=1}^s
y_{i,k}' a_{k,j}$.
This is exactly what is needed.
\end{proof}

\begin{Theorem}\label{maindual}
For any admissible left $W(\pi)$-module $M$ and any 
dualizable $\mathfrak{g}$-module $V$,
there is a natural isomorphism
$\omega_{M,V}: \overline{M} \circledast \overline{V}
\rightarrow \overline{M \circledast V}$
of right $W(\pi)$-modules such that the following diagram commutes
$$
\begin{CD}
\overline{M} \circledast \overline{V} &@>\omega_{M,V}>>& \overline{M \circledast V} \phantom{\:.}\\
@V{\chi}_{\overline{M},\overline{V}} VV&&@AA\chi^*_{M,V}A\\
\overline{M} \otimes \overline{V}&@>>>& \overline{M \otimes V}\:.
\end{CD}
$$
Here, the left hand map is the isomorphism from
Lemma~\ref{barchi}, the right hand map is the dual of the
isomorphism from Theorem~\ref{el}, and the bottom map 
sends $f \otimes g$ to the function
$m \otimes v \mapsto f(m)g(v)$.
Moreover, given another dualizable $\mathfrak{g}$-module $V'$,
the following diagram commutes:
$$\begin{CD}
(\overline{M} \circledast \overline{V}) \circledast \overline{V}'
&@>\omega_{M, V} \circledast \id_{\overline{V}'}>>
&
\overline{M \circledast V} \circledast \overline{V}'
&@>\omega_{M \circledast V, V'}>>&\overline{(M \circledast V) \circledast V'}\\
@V a_{\overline{M}, \overline{V}, \overline{V}'}VV&&&&&&@AA a^*_{M, V, V'} A\\
\overline{M} \circledast (\overline{V} \otimes \overline{V}')
&@=&\overline{M} \circledast \overline{V \otimes V}'&@>>\omega_{M, V \otimes V'}>&
\overline{M \circledast (V \otimes V')}
\end{CD}
$$
where $a_{\overline{M}, \overline{V}, \overline{V}'}$ is the map from (\ref{aVV2}) and $a^*_{M, V, V'}$
is the dual of the map from (\ref{aVV}).
\end{Theorem}

\begin{proof}
Define 
$\omega_{M,V}: \overline{M} \circledast \overline{V}
\rightarrow \overline{M \circledast V}$ so that the given diagram commutes.
The resulting map 
is natural in both $M$ and $V$, since the other three maps
in the diagram are. 
We need to check that it is a right
$W(\pi)$-module homomorphism.
Fix a basis $v_1,\dots,v_r$ for $V$ and define a matrix
$(x_{i,j})_{1 \leq i,j \leq r}$ as in Theorem~\ref{el}.
Let  $(y_{i,j})_{1 \leq i,j \leq r}$ be the inverse matrix.
Let 
\begin{equation}\label{ddef}
\delta:U(\mathfrak{p}) \rightarrow U(\mathfrak{p})\otimes\End_{\C}(V)
\end{equation} 
be the 
composite $(\id_{U(\mathfrak{p})} \otimes \rho) \circ \Delta$
where $\Delta:U(\mathfrak{p}) \rightarrow U(\mathfrak{p}) \otimes U(\mathfrak{p})$ is the comultiplication and $\rho: U(\mathfrak{p}) \rightarrow
\End_{\C}(V)$ is the representation of $\mathfrak{p}$ on $V$.
Take $w\in W(\pi)$ and let
$\delta(w) = \sum_{i,j} w'_{i,j} \otimes e_{i,j}$.
For any left $W(\pi)$-module $M$ and any $m \in M$, we have that
$$
w \cdot \left(\sum_{i=1}^r (x_{i,j} \cdot 1_\chi \otimes m) \otimes v_i\right)
=
\sum_{i,k=1}^r (w'_{i,k} x_{k,j} \cdot 1_\chi \otimes m) \otimes v_i
\in M \circledast V.
$$
In the special case $M = W(\pi)$ and $m=1$, 
this must equal
$\sum_{i,k=1}^r x_{i,k} \cdot 1_\chi \otimes w_{k,j} \otimes v_i$
for elements $w_{k,j} \in W(\pi)$, with
$\sum_{k=1}^r x_{i,k}w_{k,j} = \sum_{h,k=1}^r w'_{i,k} x_{k,j}$.
Hence in the general case too, we have that
$$
w \cdot \left(\sum_{i=1}^r (x_{i,j} \cdot 1_\chi \otimes m) \otimes v_i\right)
=
\sum_{i,k=1}^r (x_{i,k} \cdot 1_\chi \otimes w_{k,j}m) \otimes v_i.
$$
Using this formula we can now 
lift the dot action of $W(\pi)$ on $M \circledast V$
directly to the vector space $M \otimes V$ 
via the isomorphism $\chi_{M, V}$, to make $M \otimes V$ itself
into a left $W(\pi)$-module with action defined by
\begin{equation}\label{lact}
w (m \otimes v_j) = \sum_{i=1}^r w_{i,j} m \otimes v_i
\end{equation}
where the elements $w_{i,j} \in W(\pi)$
are defined 
from
$\delta(w) = \sum_{h,k=1}^r x_{i,h} w_{h,j} y_{h,j} \otimes e_{i,j}$.
Instead, let $\overline{v}_1,\dots,\overline{v}_r$ be the dual basis
for $\overline{V}$. By a similar argument to the above,
we lift the
dot action of $W(\pi)$ on $\overline{M} \circledast \overline{V}$
to the vector space $\overline{M} \otimes \overline{V}$
via the isomorphism $\chi_{\overline{M}, \overline{V}}$.
This makes $\overline{M} \otimes \overline{V}$ into a right
$W(\pi)$-module with action defined by
\begin{equation}\label{ract}
(f \otimes \overline{v}_i) w = \sum_{j=1}^r f w_{i,j} \otimes \overline{v}_j.
\end{equation}
Under these identifications, the statement that $\omega_{M, V}$
is a module homomorphism is 
equivalent to saying that the
natural map $\overline{M} \otimes \overline{V}
\rightarrow \overline{M \otimes V}$ is a module homomorphism,
which is easily checked given (\ref{lact})--(\ref{ract}).

The commutativity of the final diagram is checked by
a direct calculation which we leave as an exercise;
the matrices (\ref{hot1})--(\ref{hot2}) from the proof of Lemma~\ref{d2}
below are useful in doing this.
\end{proof}

We do not yet have {\em any} examples of dualizable
$\mathfrak{g}$-modules.

\begin{Lemma}\label{d2}
Finite direct sums, direct summands, tensor products and duals 
of dualizable
modules are dualizable.
\end{Lemma}

\begin{proof}
It is obvious that direct sums of dualizable modules are dualizable.

Consider direct summands.
Let $V$ be dualizable and suppose that $V = V' \oplus V''$
as a $\mathfrak{g}$-module.
Let $v_1,\dots,v_s$ be a basis for $V'$ and $v_{s+1},\dots,v_r$
be a basis for $V''$.
Let $(x_{i,j})_{1 \leq i,j \leq r}$ be the matrix obtained by
applying Theorem~\ref{el}
 to the basis $v_1,\dots,v_r$ for $V$.
By Lemma~\ref{all} and the assumption that $V$ is dualizable, 
the inverse matrix
$(y_{i,j})_{1 \leq i,j \leq r}$ satisfies property (b) above.
Note also that $c_{i,j}= c_{j,i} = 0$ if $1 \leq i \leq s < j \leq r$.
Using this and the uniqueness in Theorem~\ref{el}, we deduce 
that $x_{i,j} =x_{j,i}=0$ 
if $1 \leq i \leq s < j\leq r$ too.
Hence, 
$(y_{i,j})_{1 \leq i,j \leq s}$
is the inverse of the matrix
 $(x_{i,j})_{1 \leq i,j \leq s}$.
Since the matrices $(x_{i,j})_{1 \leq i,j \leq r}$ and 
$(y_{i,j})_{1 \leq i,j \leq r}$ satisfy properties (a) and (b),
the submatrices
$(x_{i,j})_{1 \leq i,j \leq s}$ and 
$(y_{i,j})_{1 \leq i,j \leq s}$ do to.
Hence $V'$ is dualizable.

Next consider tensor products.
Let $V$ and $V'$ be dualizable, with bases
$v_1,\dots,v_r$ and $v_1',\dots,v_s'$, respectively.
Let $\overline{v}_1,\dots, \overline{v}_r$
and $\overline{v}_1',\dots, \overline{v}_s'$ be the dual bases
for $\overline{V}$ and $\overline{V}'$, respectively.
Write $e_{i,j}$ for the $ij$-matrix unit in
$\End_{\C}(V) = \End_{\C}(\overline{V})^{\op}$
and $e'_{p,q}$ for the $pq$-matrix unit in
$\End_{\C}(V') = \End_{\C}(\overline{V}')^{\op}$.
Let $x = 
\sum_{i,j=1}^r x_{i,j} \otimes e_{i,j} \in U(\mathfrak{p}) \otimes \End_{\C}(V)$ be the matrix obtained by applying Theorem~\ref{el}
to the given basis for $V$ and let
$y = \sum_{i,j=1}^r y_{i,j} \otimes e_{i,j}$ be the inverse matrix.
Similarly, define $x' = \sum_{p,q=1}^s x_{p,q}' \otimes e_{p,q}'$
by applying Theorem~\ref{el} to the given basis for $V'$ and let
$y' = \sum_{p,q=1}^s y_{p,q}' \otimes e_{p,q}'$ be the inverse.
Let 
$\delta:U(\mathfrak{p}) \rightarrow U(\mathfrak{p})\otimes\End_{\C}(V)$
be the map (\ref{ddef}) from the proof of Theorem~\ref{maindual}.
Consider
 the following elements of
$U(\mathfrak{p}) \otimes \End_{\C}(V) \otimes \End_{\C}(V')$:
\begin{align}\label{hot1}
\sum_{i,j=1}^r \sum_{p,q=1}^s
x_{i,p;j,q} \otimes e_{i,j} \otimes e_{p,q}'
:=((\delta \otimes \id_{\End_\C(V')})(x')) (x \otimes 1),\\
\sum_{i,j=1}^r \sum_{p,q=1}^s
y_{i,p;j,q} \otimes e_{i,j} \otimes e_{p,q}'
:=
(y \otimes 1)((\delta \otimes \id_{\End_\C(V')})(y')).\label{hot2}
\end{align}
Clearly these are mutual inverses.
Now let $M$ be any left $W(\pi)$-module.
Recall the isomorphisms $\chi_{M,V}$
and $\chi_{M \circledast V, V'}$ from
Theorem~\ref{el} and the associativity isomorphism 
$a_{M, V,V'}: (M \circledast V) \circledast V' \rightarrow M \circledast (V \otimes V')$
from (\ref{aVV}).
The image of
$m \otimes v_j \otimes v_q'$
under the composite map
$$
a_{M, V, V'} \circ
\chi_{M \circledast V, V'}^{-1}  \circ (\chi_{M, V}^{-1} \otimes \id_{V'})
:
M \otimes V \otimes V' \rightarrow M \circledast (V \otimes V')
$$
is equal to $\sum_{i=1}^r \sum_{p=1}^s 
x_{i,p;j,q} \cdot 1_\chi \otimes m \otimes v_i \otimes v_p'$.
As in the proof of Theorem~\ref{el}, the fact that this is
a Whittaker vector implies that 
the matrix
$(x_{i,p;j,q})_{1 \leq i,j \leq r, 1 \leq p,q \leq s}$
satisfies property (a) 
with respect to the basis 
$\{v_i \otimes v_p'\:|\:i=1,\dots,r, p = 1,\dots,s\}$
for $V \otimes V'$.
Instead, take any right $W(\pi)$-module $M$.
Recalling (\ref{aVV2}) and the isomorphisms $\chi_{M,\overline{V}}$
and $\chi_{M \circledast{\overline{V}}, \overline{V}'}$ from
Lemma~\ref{barchi},
the image of
$m \otimes \overline{v}_i \otimes \overline{v}_p'$
under the map
$$
a_{M,\overline{V}, \overline{V}'} \circ
\chi_{M \circledast \overline{V}, 
\overline{V}'}^{-1}  \circ (\chi_{M, \overline{V}}^{-1} \otimes \id_{\overline{V}'})
:
M \otimes \overline{V} \otimes \overline{V}' \rightarrow M \circledast (\overline{V} \otimes \overline{V}')
$$
is equal to $\sum_{j=1}^r \sum_{q=1}^s 
m \otimes \overline{1}_\chi\cdot y_{i,p;j,q} \otimes \overline{v}_j \otimes \overline{v}_q'$.
The fact that this is a Whittaker vector implies that
the matrix
$(y_{i,p;j,q})_{1 \leq i,j \leq r, 1 \leq p,q \leq s}$
satisfies property (b). 
Hence $V \otimes V'$ is dualizable.

Finally we consider the dual $\mathfrak{g}$-module $V^*$,
assuming that $V$ is dualizable of dimension $r$.
Note that
$V^* \cong D \otimes \bigwedge^{r-1} (V)$
where $D$ is a one-dimensional representation.
Since $V^{\otimes (r-1)}$ is dualizable by the preceeding paragraph,
and $\bigwedge^{r-1} (V)$ is a summand of it, 
it follows that $\bigwedge^{r-1}(V)$ is dualizable.
It is obvious that any one dimensional representation is dualizable.
Hence $V^*$ is too.
\end{proof}

\begin{Lemma}\label{d3}
If $\pi$ is left-justified,
the natural $\mathfrak{g}$-module $V$ is dualizable.
\end{Lemma}

\begin{proof}
Let $v_1,\dots,v_N$ be the standard basis for $V$.
We are going to write down mutually inverse matrices
$(x_{i,j})_{1 \leq i,j \leq N}$ and 
$(y_{i,j})_{1 \leq i,j \leq N}$ and verify that they
satisfy properties (a) and (b) by brute force.
Since $\pi$ is left-justified,
we can take $k=0$ in (\ref{pyr1}).
All other notation throughout the proof is as in $\S$\ref{svan}.

For $1 \leq i,j \leq N$, define
\begin{equation}
x_{i, j} := 
(-1)^{\col(i) - \col(j)}
I_{\col(i)-1} (T_{\row(j),\row(i)}^{(\col(i)-\col(j))}),
\end{equation}
interpreted as $\delta_{i,j}$ if $\col(i) \leq \col(j)$.
We claim for all $1 \leq i,j,r,s \leq N$
with $\col(s) = \col(r)-1$ that
\begin{itemize}
\item[(i)] $[e_{r,s},\eta(x_{i,j})] +\delta_{i,r}\eta(x_{s,j}) \in
U(\mathfrak{g})I_\chi$;
\item[(ii)]
 $[e_{r,s},\overline{\eta}(x_{i,j})] +\delta_{i,r}\overline{\eta}(x_{s,j})
\in I_\chi U(\mathfrak{g})$.
\end{itemize}
We just explain the argument to check (i), since (ii) is entirely similar
given Lemma~\ref{miracle}.
We may as well assume that $\col(i) > \col(j)$, since it is trivial otherwise.
If $\col(i) < \col(r)$ then $[e_{r,s}, \eta(x_{i,j})] = 0$
obviously, while if $\col(i) > \col(r)$
then $[e_{r,s},\eta(x_{i,j})] \in U(\mathfrak{g}) I_\chi$, 
as $x_{i,j}$ belongs to $W(\pi_{\col(i)-1})$.
So assume that $\col(j) < \col(i) = \col(r)$.
In that case, we expand $\eta(x_{i,j})$
using Lemma~\ref{eorl2} (with $l= \col(i)-1$) then commute with $e_{r,s}$.
Almost all the resulting
terms are zero.
The only term that possibly contributes comes from the
third term on the right hand side of Lemma~\ref{eorl2}
when $q_1+\cdots+q_{\col(i)-1} + h-n=s$,
from which we get exactly
$-\delta_{i,r}\eta(x_{s,j})$ modulo $U(\mathfrak{g}) I_\chi$,
as required.

Since $\mathfrak{m}$ is generated by the 
elements $e_{r,s}$ with $\col(s) = \col(r)-1$,
formula (i) is all that is needed to verify that the matrix
$(x_{i,j})_{1 \leq i,j \leq N}$ satisfies property (a).
The inverse matrix
$(y_{i,j})_{1 \leq i,j \leq N}$ is given explicitly by
\begin{equation}\label{yij}
y_{i,j} = \sum_{t=0}^{\col(i)-\col(j)}
\sum_{i_0,\dots,i_t}
(-1)^t
x_{i_0,i_1} x_{i_1,i_2} \cdots x_{i_{t-1},i_t},
\end{equation}
where the summation is over all $1 \leq i_0,\dots,i_t \leq N$ such that
$i_0 = i, i_t = j$ and
$\col(i_0) > \cdots > \col(i_t)$.
It just remains to check that this matrix satisfies property (b).
For this, it is enough to show that
$[\overline\eta(y_{i,j}), e_{r,s}] + \overline{\eta}(y_{i,r})
\delta_{s,j}
 \in I_\chi U(\mathfrak{g})$
when $\col(s)=\col(r)-1$.
Using formula (ii), we get for $i_0,\dots,i_t$ as in (\ref{yij}) that
$$
[\overline{\eta}(x_{i_0,i_1} x_{i_1,i_2} \cdots x_{i_{t-1},i_t}),e_{r,s}] 
\equiv
\overline{\eta}(x_{i_0,i_1} \cdots x_{i_{h-1},r}x_{s, i_{h+1}} \cdots x_{i_{t-1},i_t})
$$
modulo
$I_\chi U(\mathfrak{g})$
if $i_h = r$ for some $h=0,\dots,t-1$, and it is congruent to $0$ otherwise.
Now a calculation using this and (\ref{yij})
completes the proof.
\end{proof}

\begin{Theorem}\label{leftcase}
If $\pi$ is left-justified,
 every finite dimensional $\mathfrak{g}$-module
is dualizable.
\end{Theorem}

\begin{proof}
It is easy to see that any finite dimensional $\mathfrak{g}$-module
on which $\mathfrak{g}' = \mathfrak{sl}_N$ acts trivially
is dualizable. 
Every finite dimensional $\mathfrak{g}$-module
is a direct sum of summands of tensor products of such modules
and copies of the natural module.
Hence every finite dimensional $\mathfrak{g}$-module
is dualizable by Lemmas~\ref{d2} and \ref{d3}.
\end{proof}

\begin{Remark}\label{gencase}
In fact Theorem~\ref{leftcase}
is true for an arbitrary pyramid.
The only way we have found to see this is by reducing the general
case to the left-justified case treated above. 
In order to do this, the key point is that
the functor arising from 
twisting with $\iota$ 
commutes with the bifunctor
$? \circledast ?$.
This can be proved by an argument in the spirit of \cite{GG,BGo}, 
using the invariant definition of $\iota$ mentioned briefly
in $\S$\ref{ssmore}.
\end{Remark}

Finally, we return to the natural $\mathfrak{g}$-module $V$ and check
one more technical fact which will allow us to descend
from the functor $? \circledast V$ to the translation functors
$e_i,f_i$ as defined $\S$\ref{sstf}.
To formulate this, 
we need the endomorphism $x$ of the functor
$? \circledast \overline{V}$
that is the right module analogue of (\ref{xendo}).
So, for a right $W(\pi)$-module $M$,
${x}_M: M \circledast \overline{V}
\rightarrow M \circledast \overline{V}$
is the map induced by right multiplication by
$\Omega = \sum_{i,j=1}^N e_{i,j} \otimes e_{j,i}$.

\begin{Lemma}\label{nasty2}
Assume that the natural $\mathfrak{g}$-module $V$ is dualizable
(which is true e.g. if $\pi$ is left-justified).
Then, for any admissible left 
$W(\pi)$-module $M$, the following diagram commutes:
$$
\begin{CD}
\overline{M} \circledast \overline{V} &@>\omega_{M, V}>>& \overline{M \circledast V}\phantom{,}\\
@Vx_{\overline{M}}
VV&&@VV x^*_M V\\
\overline{M} \circledast \overline{V} &@>>\omega_{M,V}>&\overline{M \circledast V},
\end{CD}
$$
where $\omega_{M, V}$ is as in Theorem~\ref{maindual} and 
$x_M^*$ denotes 
the dual map to $x_M$.
\end{Lemma}

\begin{proof}
Letting $(x_{i,j})_{1 \leq i,j \leq N}$ be the matrix from
Theorem~\ref{el} and $(y_{i,j})_{1 \leq i,j \leq N}$ be its inverse
as usual, 
we have for any $m \in M$ that 
\begin{align*}
\Omega\left(\sum_{i=1}^N (x_{i,j}\cdot 1_\chi \otimes m) \otimes v_i\right)
&=
\sum_{i,k=1}^N (e_{i,k} \eta(x_{i,j}) 1_\chi \otimes m )\otimes v_k\\
&=
\sum_{\col(i) \leq \col(k)} (e_{i,k} \eta(x_{i,j}) 1_\chi \otimes m) \otimes v_k\\
&\qquad+ \sum_{\col(i) > \col(k)} ((\eta(x_{i,j}) e_{i,k} - \eta(x_{k,j})) 1_\chi
\otimes m) \otimes v_k.
\end{align*}
Considering the special case $M = W(\pi)$ first then using naturality,
this must equal $\sum_{i,k=1}^N x_{k,i} \cdot 1_\chi \otimes w_{j,i} m \otimes v_k$
for some elements $w_{j,i} \in W(\pi)$.
Equating coefficients, we get that
\begin{multline}
w_{j,i} = \sum_{\col(h) \leq \col(k)} (-1)^{\col(k)-\col(h)} y_{i,k} e_{h,k} x_{h,j}
+
\sum_{k=1}^N (q_{\col(k)}-n) y_{i,k} x_{k,j}\\
+
\sum_{\substack{\row(h)=\row(k) \\\col(h) = \col(k)-1}}y_{i,h} x_{k,j}.
\end{multline}
Now we can lift 
the endomorphism of $M \circledast V$
to an endomorphism of the vector space $M \otimes V$
through the isomorphism $\chi_{M, V}$.
We obtain the endomorphism of $M \otimes V$
defined simply by left multiplication by
$\sum_{i,j=1}^N w_{j,i} \otimes e_{i,j} \in W(\pi) \otimes \End_\C(V)$.
With an entirely similar calculation, we lift the endomorphism
of $\overline{M} \circledast \overline{V}$
to an endomorphism of the vector space $\overline{M} \otimes \overline{V}$
through the isomorphism $\chi_{\overline{M}, \overline{V}}$.
We obtain the endomorphism of
$\overline{M} \otimes \overline{V}$ defined by right multiplication by
the same element
$\sum_{i,j=1}^N w_{j,i} \otimes e_{i,j} \in W(\pi) \otimes \End_\C(\overline{V})^{\op}$.
Using these descriptions, 
the desired commutativity of the diagram is now easy to check.
\end{proof}

\section{Whittaker functor}
\label{ssmain}
Recall that $\mathfrak{c}$ denotes the Lie subalgebra of $W(\pi)$
spanned by $D_1^{(1)},\dots,D_n^{(1)}$.
We point out that as elements of $U(\mathfrak{g})$, we have simply that
\begin{equation}
D_i^{(1)} = 
\sum_{\substack{1 \leq j \leq N \\ \row(j) = i}} e_{j,j}
\end{equation}
for each $i=1,\dots,n$.
Hence, $\mathfrak{c}$ is a subalgebra of the 
standard Cartan subalgebra $\mathfrak{d}$ of $\mathfrak{g}$,
indeed, $\mathfrak{c}$ is the centralizer of
$e$ in $\mathfrak{d}$.

Let $M$ be a $\mathfrak{g}$-module 
which is the direct sum of its 
generalized $\mathfrak{c}$-weight spaces,
i.e. $M = \bigoplus_{\alpha \in \mathfrak{c}^*} M_\alpha$. 
We do not assume that each $M_\alpha$ is finite dimensional.
Set
\begin{equation}\label{wdef}
\V(M) := \overline{\overline{\Wh}(\overline{M})},
\end{equation}
where $\overline{M}$ denotes the restricted
dual $\bigoplus_{\alpha \in \mathfrak{c}^*} M_\alpha^*$
as in (\ref{rdual})
viewed as a right $U(\mathfrak{g})$-module
with action $(fu)(v) := f(uv)$,
$\overline{\Wh}(\overline{M})$ denotes the 
right $W(\pi)$-module
obtained by applying the functor $\overline{\Wh}$
from (\ref{bwh}), and
finally $\overline{\overline{\Wh}(\overline{M})}$ denotes the
left $W(\pi)$-module obtained by taking the restricted dual 
one more time.
There is an obvious definition on morphisms, making
$\V$ into a (covariant) right exact functor.

For the first lemma, recall the automorphism $\overline{\eta}:U(\mathfrak{p})
\rightarrow U(\mathfrak{p})$ from (\ref{baretadef}).

\begin{Lemma}\label{wh2}
Let $M$ be a $\mathfrak{p}$-module 
such that $M = \bigoplus_{\alpha \in \mathfrak{c}^*} M_\alpha$
and each $M_\alpha$ is finite dimensional.
Then there is a natural $W(\pi)$-module isomorphism 
between $\V(U(\mathfrak{g}) \otimes_{U(\mathfrak{p})}
 M)$
and the $W(\pi)$-module equal to $M$ as a vector space with 
action defined by
$u \circ v := \overline{\eta}(u) v$ 
for $u \in W(\pi), v \in M$.
\end{Lemma}

\begin{proof}
Let $I := U(\mathfrak{g}) \otimes_{U(\mathfrak{p})} M$. Note that
$I = U(\mathfrak{m}) \otimes M$ as a left $U(\mathfrak{m})$-module.
So for $\alpha \in \mathfrak{c}^*$, the generalized $\alpha$-weight space
of $I$ is $$
I_\alpha = \bigoplus_{\beta \in \mathfrak{c}^*} U(\mathfrak{m})_{\beta} \otimes M_{\alpha-\beta},
$$ 
where $U(\mathfrak{m})_\beta$ is the 
$\beta$-weight
space of $U(\mathfrak{m})$ with respect to the adjoint action of $\mathfrak{c}$.
By definition, $$
\overline{\Wh}(\overline{I})
=
\{f \in \hom_{\mathfrak{m}}(I, \C_\chi)\:|\:
f(I_\alpha) = 0 \text{ for all but finitely many $\alpha \in \mathfrak{c}^*$}\}.
$$
The restriction of the obvious isomorphism
$\hom_{\mathfrak{m}}(I, \C_\chi)
\stackrel{\sim}{\rightarrow} M^*$
to the subspace
$\overline{\Wh}(\overline{I})$
gives an injective linear map
$\varphi:\overline{\Wh}(\overline{I}) \hookrightarrow \overline{M}$.
We claim that $\varphi$ is also surjective, hence an isomorphism.
To see this,
take any $f \in \overline{M}$.
Its inverse image in $\hom_{\mathfrak{m}}(I, \C_\chi)$ is the map 
$\hat f$ sending
$u \otimes m \mapsto \chi(u) f(m)$ 
for each $u \in U(\mathfrak{m})$ and 
$m \in M$.
Since $\chi(u) = 0$ if $u \notin U(\mathfrak{m})_0$, we get that
$\hat f$ vanishes on $I_\alpha$ for all but finitely many
$\alpha$.
Hence $\hat f \in \overline{\Wh}(\overline{I})$, proving the claim.
The dual map  to $\varphi$ now gives a natural vector space
isomorphism
$\overline\varphi:M \stackrel{\sim}{\rightarrow} \V(I)$.
It just remains to check that the $W(\pi)$-module structure on
$\V(I)$ corresponds under this isomorphism to the circle action
of $W(\pi)$ on $M$
\end{proof}

Recall the category $\mathcal{O}$ of $\mathfrak{g}$-modules
from $\S$\ref{sskl}, and the Verma modules $M(\alpha)$
for each $\alpha \in \C^N$ from (\ref{vermamod}).

\begin{Lemma}
Take any weight $\alpha \in \C^N$. Let $A$ be the $\pi$-tableau
with $\gamma(A) = \alpha$.
Let $A_i$ denote the $i$th column of $A$.
Then
$$
\V(M(\alpha)) \cong M(A_1) \boxtimes \cdots \boxtimes M(A_l).
$$
Moreover, $\V$ maps
short exact sequences of $\mathfrak{g}$-modules
with finite Verma flags to short exact sequences
of $W(\pi)$-modules.
Hence, $\V$ maps objects in $\mathcal O$ to objects in $\mathcal M(\pi)$.
\end{Lemma}

\begin{proof}
The first statement follows by Lemma~\ref{wh2}
since $M(\alpha) \cong U(\mathfrak{g}) \otimes_{U(\mathfrak{p})} M$
where $M$ is the $\mathfrak{p}$-module
whose pull-back through the automorphism
$\overline{\eta}$ 
is isomorphic to $M(A_1)\boxtimes\cdots\boxtimes M(A_l)$.
The second statement follows because all
short exact sequences of $\mathfrak{g}$-modules
with finite Verma flags are split when viewed as short exact sequences of $\mathfrak{m}$-modules.
For the final statement, take any $M \in \mathcal O$ and 
let $P \twoheadrightarrow M$ be its projective cover in $\mathcal O$.
Since $\V$ is right exact, it suffices to show that $\V(P)$ 
belongs to $\mathcal M(\pi)$.
This follows because $P$ has a finite Verma flag,
so $\V(P)$ has a finite filtration with factors
of the form $M(A_1) \boxtimes \cdots \boxtimes M(A_l)$
for $A \in \Tab(\pi)$.
We have already observed several times 
that the latter modules belong to $\mathcal M(\pi)$
thanks to Corollary~\ref{tee}.
\end{proof}

In view of the lemma, the functor $\V$ restricts to a well-defined
right exact functor
\begin{equation}\label{Co}
\V:\mathcal O \rightarrow \mathcal M(\pi).
\end{equation}
Moreover, $\V$ preserves central characters, so it
also sends the subcategory $\mathcal O_0$ of $\mathcal O$
consisting of all modules with integral central character
to the subcategory $\mathcal M_0(\pi)$ of $\mathcal M(\pi)$.
For the next lemma, recall from Remark~\ref{gencase}
that every finite dimensional
$\mathfrak{g}$-module is dualizable (though we have only proved
that here if $\pi$ is left-justified).

\begin{Lemma}\label{tfc}
For any $M \in \mathcal O$ and any  dualizable 
$V$, there is a natural isomorphism
$$
\nu_{M, V}: \V(M \otimes V) \stackrel{\sim}{\longrightarrow} \V(M) \circledast V
$$
of $W(\pi)$-modules.
Moreover, given another dualizable module $V'$, the following
diagram commutes:
$$
\begin{CD}
\V(M \otimes V \otimes V')
&@>\nu_{M \otimes V, V'} >>&\V(M \otimes V) \circledast V'\phantom{.}\\
@V \nu_{M, V \otimes V'} VV &&@VV \nu_{M, V} \circledast \id_{V'} V\\
\V(M) \circledast (V \otimes V') &@<< \phantom{n}a_{\V(M), V, V'}\phantom{n} < &(\V(M) \circledast V) \circledast V'.
\end{CD}
$$
Finally, letting $V^*$ denote the dual $\mathfrak{g}$-module
(which is dualizable by Lemma~\ref{d2}) the following diagram commutes:
$$
\begin{CD}
\V(M \otimes V^*)&@>\phantom{aaaa}\nu_{M, V^*}\phantom{aaaa}>>&\V(M) \circledast V^*\\
@V \iota_{\V(M \otimes V^*)} VV&&@AA \V(\eps_M) \circledast \id_{V^*} A\\
(\V(M \otimes V^*) \circledast V) \circledast V^*
&@>> \nu^{-1}_{M \otimes V^*, V} \circledast id_{V^*}>&\V(M \otimes V^* \otimes V) \circledast V^*
\end{CD}
$$
where $\iota$ is the unit of the adjunction between $? \circledast V$ and $\circledast V^*$
from (\ref{adj1}), 
and $\eps$
is the counit of the canonical adjunction between
$? \otimes V$ and $? \otimes V^*$.
\end{Lemma}

\begin{proof}
Take a module $M \in \mathcal O$ and a dualizable $\mathfrak{g}$-module $V$.
Set $N := \V(M)
= \overline{\overline{\Wh}(\overline{M})}$.
Theorem~\ref{maindual} gives us a natural isomorphism
$$
\omega_{N,V}:\overline{N} \circledast \overline{V}
\stackrel{\sim}{\rightarrow} \overline{N \circledast V}.
$$
By definition, $\overline{N} \circledast \overline{V}
=
\overline{\Wh}((\overline{\Wh}(\overline{M}) \otimes_{W(\pi)}
\overline{Q}_\chi) \otimes \overline{V})$, so from the canonical isomorphism
$\overline{\Wh}(\overline{M}) \otimes_{W(\pi)} \overline{Q}_\chi
\stackrel{\sim}{\rightarrow} \overline{M}$ we get induced an
isomorphism $\overline{N} \circledast \overline{V}
\stackrel{\sim}{\rightarrow} 
\overline{\Wh}(\overline{M} \otimes \overline{V})$.
Finally, there is an obvious isomorphism
$\overline{\Wh}(\overline{M} \otimes \overline{V})
\stackrel{\sim}{\rightarrow} \overline{\Wh}(\overline{M \otimes V})$.
Composing these maps, we have constructed a natural isomorphism
$$
\overline{\V(M) \circledast V}
=
\overline{N \circledast V}
\stackrel{\omega_{N,V}^{-1}}{\longrightarrow}
\overline{N} \circledast \overline{V}
\stackrel{\sim}{\longrightarrow}
\overline{\Wh}(\overline{M} \otimes \overline{V})
\stackrel{\sim}{\longrightarrow}
\overline{\Wh}(\overline{M \otimes V}) =
\overline{\V(M \otimes V)}.
$$
Let $\nu_{M,V}: \V(M \otimes V) \rightarrow \V(M) \circledast V$
be the dual map.
This is a natural isomorphism of $W(\pi)$-modules.

Now we consider the commutativity of the two diagrams.
The first one is checked using
the commutative diagram from Theorem~\ref{maindual}.
For the second, consider 
{\small
$$
\begin{CD}
\V(M \otimes V^*) \circledast V
&@<\nu_{M \otimes V^*, V}<<&\V(M \otimes V^* \otimes V)&@> \V(\id_M \otimes e)>>& \V(M \otimes \C)&@>\V(i_M)>> &\V(M)\phantom{.}\\
@VV\nu_{M, V^*} \circledast \id_{V}V&&@VV\nu_{M, V^* \otimes V}V&&@VV\nu_{M, \C}V&&@|\\
(\V(M) \circledast V^*) \circledast V&@>>\!\!\!a_{\V(M), V^*, V}\!\!\!>&\V(M) \circledast (V^* \otimes V) &@>>\id_{\V(M)} \circledast e >&\V(M) \circledast \C &@>>i_{\V(M)}> &\V(M).
\end{CD}
$$
}
Here, $e: V^*\otimes V \rightarrow \C$ is evaluation
$f \otimes v \mapsto f(v)$,
 $i_M: M \otimes \C \rightarrow M$ is the multiplication
$m \otimes c \mapsto cm$, and $i_{\V(M)}$ is as in (\ref{adj2}).
This diagram commutes: the left hand square commutes thanks to the first
 diagram
just checked, the middle square commutes
by naturality of $\nu$, and the right hand square is easy.
The composite $\V(M \otimes V^* \otimes V) \rightarrow \V(M)$ along the top 
of the diagram is precisely the map $\V(\eps_M)$, while the composite
$\eps_{\V(M)}:(\V(M) \circledast V^*) \circledast V \rightarrow \V(M)$
along the bottom is the definition of counit of the adjunction
from (\ref{adj2}).
Hence, we have shown that the following diagram commutes:
$$
\begin{CD}
\V(M\otimes V^* \otimes V)&@>\nu_{M \otimes V^*, V}>>&\V(M \otimes V^*) \circledast V\phantom{.}\\
@V\V(\eps_M) VV&&@VV\nu_{M, V^*}\circledast \id_{V} V\\
\V(M)&@<<\phantom{H}\eps_{\V(M)}\phantom{H}<&(\V(M) \circledast V^*)\circledast V.
\end{CD}
$$
This implies the commutativity of the second diagram;
see \cite[Lemma 5.3]{CR}.
\end{proof}

Recall the translation functors
$e_i, f_i:\mathcal O_0 \rightarrow \mathcal O_0$ from $\S$\ref{sskl},
and their counterparts on the category
$\mathcal{M}_0(\pi)$ from $\S$\ref{sstf}.

\begin{Lemma}\label{worst}
Assume that the natural $\mathfrak{g}$-module $V$ is dualizable
(which is true e.g. if $\pi$ is left-justified).
Then the functor $\V:\mathcal O_0 \rightarrow \mathcal M_0(\pi)$ 
commutes with the translation functors $f_i, e_i$ for all $i \in \Z$,
i.e. there are natural isomorphisms
$\nu^+:\V \circ f_i \stackrel{\sim}{\rightarrow} 
f_i \circ \V$ and $\nu^-:\V \circ e_i \stackrel{\sim}{\rightarrow} 
e_i \circ \V$.
In fact, $(\V, \nu^+, \nu^-)$ is a morphism of 
$\mathfrak{sl}_2$-categorifications in the sense of \cite[5.2.1]{CR}.
\end{Lemma}

\begin{proof}
Recall the endomorphism $x$ 
of the functor
$? \otimes V$
and the endomorphism $s$ of the functor
$(? \otimes V) \circ (? \otimes V)$
from $\S$\ref{sskl}, and the analogous endomorphisms
of $? \circledast V$ and $(? \circledast V) \circ (? \circledast V)$
from $\S$\ref{ssduality}.
We claim that the following diagrams commute
for all $M \in \mathcal O$:
\begin{equation}\label{eq1}
\begin{CD}
\V(M \otimes V) &@>\nu_{M,V}>>& \V(M) \circledast V\phantom{,}\\
@V\V(x_M)VV&&@VVx_{\V(M)}V\\
\V(M \otimes V)&@>>\nu_{M,V}>&\V(M) \circledast V,
\end{CD}
\end{equation}
\begin{equation}\label{eq2}
\begin{CD}
\V(M \otimes V \otimes V)&@>\nu_{M \otimes V,V}>>&\V(M\otimes V) \circledast V
&@>\nu_{M,V} \circledast \id_V>>&(\V(M) \circledast V) \circledast V\phantom{.}\\
@V\V(s_M)VV&&&&&&@VVs_{\V(M)}V\\
\V(M \otimes V \otimes V)&@>>\nu_{M \otimes V,V}>&\V(M\otimes V) \circledast V
&@>>\nu_{M,V} \circledast \id_V>&(\V(M) \circledast V) \circledast V.
\end{CD}
\end{equation}
The commutativity of the first of these is checked using the definition of
$\nu_{M,V}$ from the proof of Lemma~\ref{tfc}, together with Lemma~\ref{nasty2}.
The commutativity of the second diagram follows
immediately from the naturality of
the isomorphism $\nu_{M, V \otimes V}$, the commutative diagram from
Lemma~\ref{tfc} and the definitions of the maps $s_M$ and $s_{\V(M)}$.

Now let us prove the lemma.
Recalling the definitions (\ref{tf3})--(\ref{tf4}), 
the isomorphisms
$\nu_{M, V}: \V(M \otimes V) \rightarrow \V(M) \circledast V$
and $\nu_{M, V^*}: \V(M \otimes V^*) \rightarrow \V(M) \circledast V^*$
restrict to give natural isomorphisms
$\nu^+_M: \V(f_i M) \rightarrow f_i \V(M)$
and $\nu^-_M: \V(e_i M) \rightarrow e_i \V(M)$
for each $M \in \mathcal{O}_0$.
This defines the natural isomorphisms $\nu^{\pm}$.
The fact that the triple $(\V, \nu^+, \nu^-)$ is a
morphism of $\mathfrak{sl}_2$-categorifications follows from 
(\ref{eq1})--(\ref{eq2})
together with the final commutative diagram
from Lemma~\ref{tfc}.
\end{proof}

In the generality of (\ref{Co}), the right exact functor $\V$ 
is usually not exact.
However, by a result of Lynch \cite[Lemma 4.6]{Ly} (which
Lynch attributes originally to N. Wallach)
$\V$ is exact on short exact sequences of 
$\mathfrak{g}$-modules that are finitely
generated over $\mathfrak{m}$; see the next lemma.
For this reason, we are going to restrict our attention 
from now on to the
parabolic category $\mathcal{O}(\pi)$ from $\S$\ref{sskl}
and the category $\mathcal{F}(\pi)$ of finite dimensional
$W(\pi)$-modules from $\S$\ref{ssgg}.

\begin{Lemma}\label{exactness}
The restriction of the functor $\V$ to
$\mathcal{O}(\pi)$ defines an {exact} functor
$$
\V:\mathcal{O}(\pi) \rightarrow \mathcal{F}(\pi).
$$
Moreover, $\V$ maps the
parabolic Verma module $N(A)$ from (\ref{parv})
to the standard module $V(A)$ from (\ref{VAdef}),
for any $A \in \Col(\pi)$.
\end{Lemma}

\begin{proof}
The second statement is immediate from Lemma~\ref{wh2}.
For the first statement, take any $M \in \mathcal O(\pi)$.
Note to start with that $M$ is finitely generated as a 
$U(\mathfrak{m})$-module.
This follows because the parabolic Verma modules are finitely
generated as $U(\mathfrak{m})$-modules. 
By definition,
$$
\overline{\Wh}(\overline{M}) = 
\{f \in \hom_{\mathfrak{m}}(M, \C_\chi)\:|\:
f(M_\alpha) = 0\text{ for all but finitely many }\alpha \in \mathfrak{c}^*\}.
$$
It is already clear from this that
$\V(M)$ is finite dimensional, i.e. it lies in $\mathcal F(\pi)$,
because
$\hom_{\mathfrak{m}}(M, \C_\chi)$ certainly is by the finite generation.
We claim that in fact 
$$
\overline{\Wh}(\overline{M}) = \hom_{\mathfrak{m}}(M, \C_\chi).
$$
To see this, it suffices to show that any
$f \in \hom_{\mathfrak{m}}(M, \C_\chi)$ vanishes on $M_\alpha$
for all but finitely many $\alpha \in \mathfrak{c}^*$.
Pick weights $\alpha_1,\dots,\alpha_r \in \mathfrak{c}^*$
such that $M$ is generated as a $U(\mathfrak{m})$-module
by $M_{\alpha_1}\oplus\cdots\oplus M_{\alpha_r}$.
For any $\alpha \in \mathfrak{c}^*$, the weight space $M_\alpha$
is spanned by terms of the form $u_i m_i$ for 
$u_i \in U(\mathfrak{m})_{\alpha - \alpha_i}$ and $m_i \in M_{\alpha_i}$.
But $f(u_i m_i) = \chi(u_i) f(m_i)$
and $\chi(u_i) = 0$ unless $\alpha = \alpha_i$,
so we deduce that $f(M_\alpha) = 0$ 
unless $\alpha \in \{\alpha_1,\dots,\alpha_r\}$.

To complete the proof of the lemma, we now show that
$\hom_{\mathfrak{m}}(?, \C_\chi)$ is an exact functor
on the category of $\mathfrak{g}$-modules that are finitely generated
over $\mathfrak{m}$.
Let $E$ denote the space of linear maps
$f:U(\mathfrak{m}) \rightarrow \C$ 
which annihilate $(I_\chi)^p$ for $p \gg 0$,
viewed as an $\mathfrak{m}$-module
by $(xf)(u) = f(ux)$ for $f \in E,
x \in \mathfrak{m}$ and $u \in U(\mathfrak{m})$.
By \cite[Assertion 2]{Skry}, $E$
is an injective $\mathfrak{m}$-module,
so the functor
$\hom_{\mathfrak{m}}(?, E)$ is exact.
For any $\mathfrak{g}$-module $M$,
$\hom_{\mathfrak{m}}(M, E)$
is naturally a right $U(\mathfrak{g})$-module with action
$(fu)(v) = f(uv)$ for $f \in \hom_{\mathfrak{m}}(M, E),
u \in U(\mathfrak{g})$ and $v \in M$.
Moreover, if $M$ is finitely generated as a 
$U(\mathfrak{m})$-module,
then $\hom_{\mathfrak{m}}(M, E)$ belongs to the category
$\overline{\mathcal{W}}(\pi)$.
It remains to observe 
that the $\mathfrak{m}$-module
$\Wh(E)$ can be identified with $\C_\chi$, so that 
for any $\mathfrak{g}$-module $M$
$$
\hom_{\mathfrak{m}}(M, \C_\chi)
=
\hom_{\mathfrak{m}}(M, \Wh(E))
=
\overline{\Wh}(\hom_{\mathfrak{m}}(M, E)).
$$
Hence, 
on the category of $\mathfrak{g}$-modules that are finitely generated
over $\mathfrak{m}$,
the functor $\hom_{\mathfrak{m}}(?, \C_\chi)$
factors as the composite of 
the exact functor
$\hom_{\mathfrak{m}}(?, E)$
and Skryabin's equivalence of categories
$\overline{\Wh}:\overline{\mathcal{W}}(\pi)
\rightarrow W(\pi)\text{-mod}$, so
it is exact.
\end{proof}

For a while now, we will restrict our attention to integral 
central characters.
By Lemma~\ref{exactness},
the restriction of $\V$ to $\mathcal{O}_0(\pi)$
gives an exact functor
\begin{equation}
\V: 
\mathcal O_0(\pi) \rightarrow \mathcal F_0(\pi).
\end{equation}
Also let $\I:\mathcal F_0(\pi) \rightarrow \mathcal M_0(\pi)$
be the natural inclusion functor.
We use the same notation
$\V$ and $\I$
for the induced maps at the level of Grothendieck groups.
Recall also the 
isomorphism
$i:\bigwedge^\pi(V_\Z) \rightarrow [\mathcal O_0(\pi)],
N_A \mapsto [N(A)]$
from the proof of Theorem~\ref{parkl},
and the isomorphisms $j:P^\pi(V_\Z) \rightarrow [\mathcal F_0(\pi)],
V_A \mapsto [V(A)]$ 
and $k:S^\pi(V_\Z) \rightarrow [\mathcal M_0(\pi)],
M_A \mapsto [M(A)]$ from (\ref{we}).
We observe that the following diagram 
commutes:
\begin{equation}\label{lastsq}
\begin{CD}
\bigwedge^\pi(V_\Z) &@>\V>>&P^\pi(V_\Z) &@>\I>> & S^\pi(V_\Z)\\
@ViVV&&@VVjV&&@VVkV\\
[\mathcal O_0(\pi)]&@>>\V>&[\mathcal F_0(\pi)]&@>>\I>&[\mathcal M_0(\pi)],
\end{CD}
\end{equation}
where the top $\V$ is the map from (\ref{themapf}) and the 
top $\I$ is the natural inclusion.
To see this, 
we already checked in (\ref{we}) that the right hand square commutes,
and the fact that $\V(N(A)) \cong V(A)$ from Lemma~\ref{exactness}
is exactly what is needed to check that the left hand
square does too.
Now we are ready to invoke Theorem~\ref{parkl}, or rather,
to invoke the Kazhdan-Lusztig conjecture, since Theorem~\ref{parkl}
was a direct consequence of it.
For the statement of the following theorem, recall the definition of
 the bijection
$R:\Std_0(\pi) \rightarrow \Dom_0(\pi)$ from (\ref{rmap});
in the case that $\pi$ is left-justified the
rectification $R(A)$ of a standard $\pi$-tableau $A$ simply means its
row equivalence class.

\begin{Theorem}\label{main}
For $A \in \Col_0(\pi)$, we have that
\begin{equation*}
\V(K(A)) \cong \left\{
\begin{array}{ll}
L(R(A))&\text{if $A$ is standard,}\\
0&\text{otherwise.}
\end{array}\right.
\end{equation*}
\end{Theorem}

\begin{proof}
Note that it suffices to prove the theorem in the special case
that $\pi$ is left-justified. 
Indeed, in view of Theorem~\ref{parkl},
the properties of the homomorphism
$\V:\bigwedge^\pi(V_\Z) \rightarrow P^\pi(V_\Z)$ 
and the commutativity of the 
left hand square in (\ref{lastsq}), the theorem follows if we can
show that $j(L_A) = [L(A)]$ for every $A \in \Dom_0(\pi)$. 
This last statement is independent of the particular choice of $\pi$,
thanks to the existence of the isomorphism $\iota$.
So assume from now on that $\pi$ is left-justified.

Using Theorem~\ref{parkl} and 
the commutativity of the left hand square in (\ref{lastsq}) again,
we know already for $A \in \Col_0(\pi)$
that $\V(K(A)) \neq 0$ if and only if  $A \in \Std_0(\pi)$.
Let $A_0 \in \Col_0(\pi)$ be the ground-state tableau, with all entries on 
row $i$ equal to $(1-i)$.
Since the crystal $(\Std_0(\pi), \tilde e_i, \tilde f_i, 
\eps_i, \varphi_i, \theta)$ 
is connected, it makes sense to define the {\em height} 
of $A \in \Std_0(\pi)$ to be the minimal number of 
applications of the operators
$\tilde f_i$, $\tilde e_i$ ($i \in \Z$) needed to map
$A$ to $A_0$. We proceed to prove that
$\V(K(A)) \cong L(R(A))$ for $A \in \Std_0(\pi)$
by induction on height.
For the base case, observe that
no other elements of $\Col_0(\pi)$ have the same content
as $A_0$, hence $N(A_0) = K(A_0)$.
Similarly, $V(R(A_0)) = L(R(A_0))$.
Hence by 
Lemma~\ref{exactness},
we have that
$\V(K(A_0)) \cong L(R(A_0))$.

Now for the induction step, take $B \in \Std_0(\pi)$
of height $> 0$. 
We can write $B$ as either $\tilde f_i(A)$ or
as $\tilde e_i(A)$,
where $A \in \Std_0(\pi)$ is of strictly
smaller height. We will assume that the first case holds,
i.e. that $B = \tilde f_i(A)$, since the argument in the
second case is entirely similar. By the induction hypothesis,
we know already that $\V(K(A)) \cong L(R(A))$.
We need to show that $\V(K(B)) \cong L(R(B))$.

Note by 
Lemma~\ref{exactness}
that $\V(N(B)) \cong V(B)$, and 
by exactness of the functor $\V$, we know that $\V(K(B))$ is
a non-zero quotient of $V(B)$.
Since $B \in \Std_0(\pi)$, Theorem~\ref{charithm} shows that 
$V(B)$ is a highest weight module
of type $R(B)$.
Hence $\V(K(B))$ is a highest weight module of type $R(B)$ too.
Also $K(B)$ is both a quotient and a submodule of $f_i K(A)$
by Theorem~\ref{parkl}.
Hence by Lemma~\ref{worst}, 
$\V(K(B))$ is both a quotient and a submodule of
$\V(f_i K(A)) \cong f_iL(R(A))$.
In particular, $L(R(B))$ is a quotient of $f_i L(R(A))$
and $\V(K(B))$ is a non-zero submodule of it.

Finally, we know by Theorem~\ref{crlem} that the socle and cosocle of
$f_i L(R(A))$ are irreducible and isomorphic to each other.
Since we know already that $L(R(B))$ is a quotient of 
$f_i L(R(A))$, it follows that 
the socle of $f_i L(R(A))$ is isomorphic to $L(R(B))$.
Since $\V(K(B))$ embeds into $f_i L(R(A))$, this means
that $\V(K(B))$ has irreducible socle isomorphic to $L(R(B))$
too. But $\V(K(B))$ is a highest weight module of type $R(B)$.
These too statements together imply that
$\V(K(B))$ is indeed irreducible.
\end{proof}

\begin{Corollary}\label{c1}
The isomorphism $j:P^\pi(V_\Z) \rightarrow
[\mathcal F_0(\pi)]$ maps 
$L_A$ to $[L(A)]$ for each $A \in \Dom_0(\pi)$.
Hence, for $A \in \Col_0(\pi)$ and $B \in \Std_0(\pi)$
$$
[V(A):L(R(B))] = 
\sum_{C \sim_\co A} (-1)^{\ell(A,C)} P_{d(\gamma(C))w_0, d(\gamma(B)) w_0}(1),
$$
notation as in (\ref{vae}).
\end{Corollary}

\begin{proof}
The first statement follows from the theorem and the
commutativity of the diagram (\ref{lastsq}).
The second two statement then follows by (\ref{vae}).
\end{proof}

\begin{Corollary}
For $A \in \Dom_0(\pi)$ and $i \in \Z$, the following
properties hold.
\begin{enumerate}
\item[\rm(i)]
If $\eps_i(A) = 0$ then $e_i L(A) = 0$.
Otherwise, $e_i L(A)$ is an indecomposable module with
irreducible socle and cosocle isomorphic to $L(\tilde e_i(A))$.
\item[\rm(ii)]
If $\phi_i(A) = 0$ then $f_i L(A) = 0$.
Otherwise, $f_i L(A)$ is an indecomposable module with
irreducible socle and cosocle isomorphic to $L(\tilde f_i(A))$.
\end{enumerate}
\end{Corollary}

\begin{proof}
Argue using (\ref{itf1})--(\ref{itf2}), 
Theorem~\ref{crlem} and Corollary~\ref{c1}, like in the
proof of Theorem~\ref{kuj}.
\end{proof}

Since the Gelfand-Tsetlin characters of standard modules are known, 
one can now in principle
compute the characters
of the finite dimensional irreducible $W(\pi)$-modules
with integral central character, by inverting the 
unitriangular square submatrix
$([V(A):L(R(B))])_{A, B \in \Std_0(\pi)}$ of the
decomposition matrix from Corollary~\ref{c1}.
Using Theorem~\ref{irret} too, one can deduce from this
the characters
of arbitrary finite dimensional irreducible $W(\pi)$-modules.
All the other combinatorial results just formulated can also 
be extended to arbitrary central characters in similar fashion.
We just record here the 
extension of Theorem~\ref{main} itself to arbitrary
central characters.

\begin{Corollary}
For $A \in \Col(\pi)$, we have that
\begin{equation*}
\V(K(A)) \cong \left\{
\begin{array}{ll}
L(R(A))&\text{if $A$ is standard,}\\
0&\text{otherwise.}
\end{array}\right.
\end{equation*}
Thus, the functor $\V:\mathcal O(\pi)
\rightarrow \mathcal F(\pi)$ sends irreducible modules to irreducible
modules or to zero.
Every finite dimensional irreducible $W(\pi)$-module arises in this way.
\end{Corollary}

\begin{proof}
This is a consequence of Corollary~\ref{c1}, 
Lemma~\ref{exactness}, Theorem~\ref{irret} 
and \cite[Proposition 5.12]{BG}.
\end{proof}

The final result gives a criterion for the irreducibility
of the standard module $V(A)$, in the spirit of \cite{LNT}.
Note as a special case of this corollary, we recover 
the main result of \cite{Molevirr} concerning Yangians.
Following \cite[Lemma 3.8]{LZ}, 
we say that two sets $A = \{a_1,\dots,a_r\}$ and $B = \{b_1,\dots,b_s\}$
of numbers from $\C$ are {\em separated} if
\begin{enumerate}
\item[(a)] $r < s$ and there do not exist 
$a,c \in A - B$ and $b \in B - A$ such that $a < b < c$;
\item[(b)] $r = s$ and there do not exist
$a,c \in A  - B$ and $b,d \in B - A$ such that $a < b < c < d$
or $a > b > c > d$;
\item[(c)] $r > s$ and there do not exist 
$c \in A - B$ and $b,d \in B - A$ such that $b > c > d$.
\end{enumerate}
Say that a $\pi$-tableau $A \in \Col(\pi)$
is {\em separated} if the sets $A_i$ and $A_j$ of entries in the $i$th and
$j$th columns of $A$, respectively,  are separated for each $1 \leq i < j \leq l$.

\begin{Theorem}
For $A \in \Col(\pi)$, the standard module $V(A)$ is irreducible if and only if
$A$ is separated, in which case it is isomorphic to $L(B)$
where $B \in \Dom(\pi)$ is the row equivalence class of $A$.
\end{Theorem}

\begin{proof}
Using Theorem~\ref{irret}, the proof reduces to the special case
that $A$ belongs to $\Col_0(\pi)$.
In that case, we apply \cite[Theorem 1.1]{LZ} and
the main result of Leclerc, Nazarov and Thibon \cite[Theorem 31]{LNT};
see also \cite{Caldero}.
These references imply that
$V_A$ is equal to $L_B$ for some $B \in \Dom_0(\pi)$
if and only if $A$ is separated. Actually, the references cited
only prove the $q$-analog of this statement, but it 
follows at $q=1$ too by the positivity of the structure constants
from \cite[Remark 24]{qla};
see the argument from the proof of \cite[Proposition 15]{LNT}.
By Theorem~\ref{main}, this 
shows that $V(A)$ is irreducible if and only if $A$ is separated.
Finally, when this happens, we must have that
$V(A) \cong L(B)$ where $B$ is the row equivalence class
of $A$, since $V(A)$ always contains a highest weight vector of that type.
\end{proof}

\backmatter

\chapter*{Notation}

\begin{tabular}{llr}
$\parallel$&Equivalence relation on $\Col(\pi)$&\pageref{rmap}\\
$\leq$&Bruhat order on $\Row(\pi)$&\pageref{bruhord}\\
$\bigwedge^\pi(V_\Z)$&$\bigwedge^{q_1}(V_\Z)\otimes\cdots\otimes
\bigwedge^{q_l}(V_\Z)$&\pageref{sdef}\\
$\mathfrak{b}$&Upper triangular matrices in $\mathfrak{g}$&\pageref{vermamod}\\
$\mathfrak{c}$&Lie algebra spanned by $D_1^{(1)},\dots,D_n^{(1)}$&\pageref{c}\\
$\mathfrak{d}$&Diagonal matrices in $\mathfrak{g}$&\pageref{pe2}\\
$\mathfrak{g} = \bigoplus_{j \in \Z} \mathfrak{g}_j$&
Lie algebra $\mathfrak{gl}_N$; $e_{i,j}$ has degree $(\col(j)-\col(i))$&\pageref{sswalgebras}\\
$\mathfrak{h}$&Levi subalgebra
$\mathfrak{g}_0 \cong \mathfrak{gl}_{q_1} \oplus\cdots\oplus \mathfrak{gl}_{q_l}$&\pageref{edef}\\
$\mathfrak{m}$&Nilpotent subalgebra $\bigoplus_{j < 0} \mathfrak{g}_j$&\pageref{edef}\\
$\mathfrak{p}$&Parabolic subalgebra $\bigoplus_{j \geq 0} \mathfrak{g}_j$&\pageref{edef}\\
$\gamma(A)$&Column reading of a tableau&\pageref{sstableaux}\\
$\gamma_a, \gamma_i$ & Standard bases for $P, P_\infty$&\pageref{Pdef},\pageref{pinf}\\
$\Delta, \Delta_{l',l''}$&Comultiplications&\pageref{com2},\pageref{deltall}\\
$\delta_i$ & Standard basis for $\mathfrak{d}^*$&\pageref{del}\\
$\eps_i$ & Standard basis for $\mathfrak{c}^*$& \pageref{eps}\\
$\eta$&$e_{i,j} \mapsto e_{i,j}
+ \delta_{i,j}(n-q_{\col(j)}-q_{\col(j)+1}-\cdots - q_l)$\quad\qquad&\pageref{etadef}\\
$\overline{\eta}$&
$e_{i,j} \mapsto e_{i,j}
+ \delta_{i,j}(n-q_1-q_2-\cdots - q_{\col(j)})$&\pageref{baretadef}\\
$\theta(\alpha),\theta(A)$&Content of a weight, tableau&\pageref{contdef0},\pageref{contdef1}\\
$\iota,\tau,\mu_f, \eta_c$&Automorphisms of $Y_n(\sigma) / W(\pi)$&\pageref{taudef},\pageref{taup}\\
$\kappa$&Canonical surjection of $Y_n(\sigma)$ onto $W(\pi)$&\pageref{thetadef}\\
$\mu$&Multiplication on Grothendieck group&\pageref{themapmu},\pageref{newmu}\\
$\mu_{M, V}$&Isomorphism from tensor identity&\pageref{tensorid}\\
$\nu_{M, V}$&Isomorphism between $\V(M \otimes V)$, $\V(M) \circledast V$&\pageref{tfc}\\
$\xi$&Miura transform
embedding $W(\pi)$ into $U(\mathfrak{h})$&\pageref{miura}\\
$\pi = (q_1,\dots,q_l)$&Pyramid of level $l$ with row lengths $(p_1,\dots,p_n)$&\pageref{sspyramids}\\
$\rho(A)$&Row reading of a tableau&\pageref{roread}\\
$\sigma= (s_{i,j})_{1 \leq i,j \leq n}$&Shift matrix &\pageref{shift}\\
$\chi_{M, V}$&Isomorphism between $M \circledast V$ and $M \otimes V$&\pageref{el}\\
$\Psi_N$&Harish-Chandra homomorphism&\pageref{hchom}\\
$\psi$&Isomorphism between $Z(U(\mathfrak{g}))$ and $Z(W(\pi))$&\pageref{psidef}\\
$\omega_{M, V}$&Isomorphism between $\overline{M} \circledast \overline{V}$,
$\overline{M \circledast V}$&\pageref{maindual}\\
$a_{M, V, V'}$&Associativity isomorphism&\pageref{aVV}\\
$\col(i), \row(i)$&Column, row numbers of $i$ in $\pi$&\pageref{rl}\\
$e_i, f_i$&Generators of $U_\Z$, translation functors&\pageref{ssdcb},\pageref{tf1}\\
$\tilde e_i, \tilde f_i$&Crystal operators&\pageref{sscrystals}\\
$s_M$&Endomorphisms of $(M \otimes V) \otimes V$, $(M \circledast V) \circledast V$&\pageref{cr1},\pageref{sendo}\\
$x_M$&Endomorphisms of $M \otimes V$, $M \circledast V$&\pageref{cr1},\pageref{xendo}\\
$x_{i,a}, y_{i,a}$&Elements of $\widehat{Z}[\mathscr P_n]$&\pageref{xia}\\
\end{tabular}

\pagebreak

\begin{tabular}{llr}
$A_0$&Ground-state tableau&\pageref{gst}\\
$\Col(\pi)$&Column strict $\pi$-tableaux&\pageref{colpi}\\
$C_n^{(r)}$&Central elements of $Y_n(\sigma) / W(\pi)$&\pageref{centref}\\
$\Dom(\pi)$&Dominant row symmetrized $\pi$-tableaux&\pageref{dompi}\\
$D_i^{(r)}, E_i^{(r)}, F_i^{(r)}$&Generators of $Y_n(\sigma) / W(\pi)$&\pageref{def}\\
$E_{i,j}^{(r)}, F_{i,j}^{(r)}$&Higher root elements&\pageref{defu}\\
$\mathcal F(\pi)$&Finite dimensional representations of $W(\pi)$&\pageref{fcat}\\
$\mathbb{F}$&Algebraically closed field of characteristic $0$&\pageref{syangians}\\
$K_{(1^n)}^{\sharp/\flat}(\sigma)$
& Ideals of positive/negative Borel subalgebras&\pageref{kernels}\\
$K_A$&Dual canonical basis element of $\bw^\pi(V_\Z)$&\pageref{nadef}\\
$K(A)$&Irreducible highest weight module in $\mathcal O(\pi)$&\pageref{parv}\\
$L_\alpha, L_A$&Dual canonical bases of $T^N(V_\Z)$, $S^\pi(V_\Z)$&\pageref{mpl}\\
$L(A)$&Irreducible highest weight module in $\mathcal M(\pi)$&\pageref{lad}\\ 
$M(\alpha)$&Verma module of highest weight $(\alpha-\rho)$&\pageref{vermamod}\\
$M_\alpha, M_A$&Monomial bases of $T^N(V_\Z)$, $S^\pi(V_\Z)$&\pageref{malpha}\\
$M(A)$&Generalized Verma module&\pageref{genv}\\
$\mathcal M(\pi)$&Analogue of category $\mathcal O$ for $W(\pi)$&\pageref{mcat}\\
$N_A$&Monomial basis element of $\bw^\pi(V_\Z)$&\pageref{nadef}\\
$N(A)$&Parabolic Verma module&\pageref{parv}\\
$\mathcal O, \mathcal O(\pi)$&Category $\mathcal O$, parabolic category $\mathcal O$&\pageref{ocatp}\\
$P^\pi(V_\Z)$&Polynomial representation of $U_\Z$&\pageref{themapf}\\
$P, P_\infty$&Free $\Z$-modules on bases $\{\gamma_a\:|\:a \in \C\}$,
$\{\gamma_i\:|\:i \in \Z\}$&\pageref{ssdcb}\\
$P_n, Q_n$&Weight lattice, root lattice in $\mathfrak{c}^*$&\pageref{pn}\\
$\mathscr P_n, \mathscr Q_n$&Gelfand-Tsetlin weights, rational weights&\pageref{sscharblock},\pageref{Qn}\\
$Q_\chi$&Generalized Gelfand-Graev representation&\pageref{gggr}\\
$R(A)$&Rectification of a standard $\pi$-tableau&\pageref{rmap}\\
$\Row(\pi)$&Row symmetrized $\pi$-tableaux&\pageref{colpi}\\
$\Std(\pi)$&Standard $\pi$-tableaux&\pageref{stdpi}\\
$S_{i,j}$&$S_{i,j} = s_{i,j} + p_{\min(i,j)}$&\pageref{Sijdef}\\
$SY_n(\sigma)$&Special shifted Yangian&\pageref{synind}\\
$S^\pi(V_\Z)$&$S^{p_1}(V_\Z)\otimes\cdots\otimes S^{p_n}(V_\Z)$&\pageref{sdef}\\
$T^N(V_\Z)$&Tensor space&\pageref{ts}\\
$T_{i,j}^{(r)}$&Alternate generators of $Y_n(\sigma) / W(\pi)$&\pageref{tijind}\\
$T_{i,j;x}^{(r)}$&Invariants in $U(\mathfrak{p})$&\pageref{thedef}\\
$U_\Z$&Kostant $\Z$-form for $\mathfrak{gl}_\infty$&\pageref{uz}\\
$V_\Z$&$\Z$-form for natural $U_\Z$-module&\pageref{uz}\\
$\V$&Map $\bw^\pi(V_\Z) \rightarrow S^{\pi}(V_\Z)$, Whittaker functor&\pageref{themapf},\pageref{wdef}\\
$V$&Natural $\mathfrak{g}$-module of column vectors&\pageref{tf1}\\
$V_A$&Standard monomial basis element of $P^\pi(V_\Z)$&\pageref{vadef}\\
$V(A)$&Standard module&\pageref{VAdef}\\
$\Wh, \overline{\Wh}$&Left, right Whittaker vectors&\pageref{whdef},\pageref{bwh}\\
$W(\pi)$&Finite $W$-algebra&\pageref{origdef}\\
$Y_n(\sigma)$&Shifted Yangian &\pageref{yangian}\\
$Y_{(1^n)}^{\sharp/\flat}(\sigma)$&Positive/negative Borel
subalgebras
of $Y_n(\sigma)$&\pageref{triangular}\\
$Z_N^{(r)}$&Generators of $Z(U(\mathfrak{gl}_N))$&\pageref{zNu}\\
\end{tabular}

\end{document}